\documentclass[reqno,11pt]{amsbook}
\usepackage{amsmath,amssymb,amsfonts,
}

\usepackage{diagrams}

\diagramstyle[centredisplay,dpi=600,nohug]
\newif\ifhavegoodboldmath\havegoodboldmathtrue
\newarrow{Tto}----{->}
\newarrow{Tinc}{\ifhavegoodboldmath bold\fi hook}---{->}
\newarrow{Tincdash}{\ifhavegoodboldmath bold\fi hook}{dash}{}{dash}{->}
\newarrow{Tonto}----{->>}
\newarrow{Tdash}{}{dash}{}{dash}{->}

\theoremstyle{plain}
\newtheorem{thm}{Theorem}[section]
\newtheorem{lem}[thm]{Lemma}
\newtheorem{cor}[thm]{Corollary}
\newtheorem{propo}[thm]{Proposition}

\theoremstyle{definition}
\newtheorem{defn}[thm]{definition}
\newtheorem{ex}[thm]{Example}
\newtheorem{parraf}[thm]{}

\theoremstyle{remark}
\newtheorem*{rem}{Remark}
\newtheorem*{notaci}{Notations}
\newtheorem*{notacion}{Notation}
\newtheorem{cuesab}{Question}

\numberwithin{equation}{thm}

\newcommand{\HA}{\widehat{A}}
\newcommand{\HB}{\widehat{B}}
\newcommand{\hf}{\widehat{f}}
\newcommand{\hd}{\widehat{d}}
\newcommand{\HM}{\widehat{M}}
\newcommand{\hm}{\widehat{\mathfrak{m}}}

\newcommand{\BA}{\mathbb A}
\newcommand{\BD}{\mathbb D}
\newcommand{\NN}{\mathbb N}

\newcommand{\ZZ}{\mathbb Z}

\newcommand{\fm}{\mathfrak m}
\newcommand{\fn}{\mathfrak n}
\newcommand{\fp}{\mathfrak p}
\newcommand{\fq}{\mathfrak q}
\newcommand{\fr}{\mathfrak r}

\newcommand{\FJ}{\mathfrak J}
\newcommand{\FS}{\mathfrak S}
\newcommand{\FT}{\mathfrak T}
\newcommand{\FU}{\mathfrak U}
\newcommand{\FV}{\mathfrak V}

\newcommand{\FX}{\mathfrak X}
\newcommand{\FY}{\mathfrak Y}
\newcommand{\FZ}{\mathfrak Z}

\newcommand{\sch}{\mathsf {Sch}}
\newcommand{\sfa}{\mathsf {FS}}
\newcommand{\sfn}{\mathsf {NFS}}
\newcommand{\scha}{\mathsf {Sch}_{\mathsf {af}}}
\newcommand{\sfaa}{\sfa_{\mathsf {af}}}
\newcommand{\sfna}{\sfn_{\mathsf {af}}}

\newcommand{\CC}{\mathcal C}

\newcommand{\CF}{\mathcal F}
\newcommand{\CG}{\mathcal G}
\newcommand{\CH}{\mathcal H}
\newcommand{\CI}{\mathcal I}
\newcommand{\CJ}{\mathcal J}
\newcommand{\CK}{\mathcal K}
\newcommand{\CL}{\mathcal L}

\newcommand{\CO}{\mathcal O}
\newcommand{\CP}{\mathcal P}
\newcommand{\CQ}{\mathcal Q}

\newcommand{\dirlim}[1]{\begin{array}[t]{c} \mathrm{lim}\\[-7.5 pt]
 {\longrightarrow} \\[-7.5 pt] {\scriptstyle {#1}} \end{array}}
\newcommand{\invlim}[1]{\begin{array}[t]{c} \mathrm{lim}\\[-7.5 pt]
 {\longleftarrow} \\[-7.5 pt] {\scriptstyle {#1}} \end{array}}

\newcommand{\lto}{\longrightarrow}
\newcommand{\xto}{\xrightarrow}

\newcommand{\epi}{\twoheadrightarrow}
\newcommand{\mono}{\rightarrowtail}
\newcommand{\inc}{\hookrightarrow}
\newcommand{\iso}{\tilde{\to}}

\newcommand{\imp}{\Rightarrow}
\newcommand{\dimp}{\Leftrightarrow}

\newcommand{\tr}{\triangle}
\newcommand{\tc}{\widehat{\otimes}}
\newcommand{\om}{\widehat{\Omega}}

\DeclareMathOperator{\aut}{Aut}
\DeclareMathOperator{\Jac}{Jac}
\DeclareMathOperator{\C}{C}
\DeclareMathOperator{\B}{B}
\DeclareMathOperator{\Z}{Z}
\DeclareMathOperator{\h}{H}
\DeclareMathOperator{\rg}{rg}
\DeclareMathOperator{\htt}{ht}
\DeclareMathOperator{\Ima}{Im}
\DeclareMathOperator{\spec}{Spec}
\DeclareMathOperator{\Max}{Max}
\DeclareMathOperator{\spf}{Spf}
\DeclareMathOperator{\supp}{Supp}

\DeclareMathOperator{\Img}{Im}

\DeclareMathOperator{\Ker}{Ker}
\DeclareMathOperator{\Hom}{Hom}
\DeclareMathOperator{\Der}{Der}
\DeclareMathOperator{\Homcont}{Homcont}
\DeclareMathOperator{\Dercont}{Dercont}
\DeclareMathOperator{\ga}{\Gamma}

\DeclareMathOperator{\D}{D}
\DeclareMathOperator{\fD}{\mathfrak{D}}
\DeclareMathOperator{\dimtop}{dimtop}
\DeclareMathOperator{\tor}{Tor}
\DeclareMathOperator{\shom}{\CH\mathit{om}}

\DeclareMathOperator{\ext}{Ext}
\DeclareMathOperator{\rad}{J}

\DeclareMathOperator{\coh}{\mathsf{Coh}}
\DeclareMathOperator{\Com}{\mathsf{Comp}}
\DeclareMathOperator{\com}{\text{-}\mathsf{comp}}
\DeclareMathOperator{\modu}{\text{-}\mathsf{mod}}
\DeclareMathOperator{\Modu}{\mathsf{Mod}}
\DeclareMathOperator{\alg}{\text{-}\mathsf{Alg}}

\newcommand{\ie}{{\it i.e.} }

\begin{document}

\frontmatter
\title{Infinitesimal local study of formal schemes}

\author[L. Alonso]{Leovigildo Alonso Tarr\'{\i}o}
\address{Departamento de \'Alxebra\\
Facultade de Matem\'a\-ticas\\
Universidade de Santiago de Compostela\\
E-15782  Santiago de Compostela, SPAIN}
\email{leoalonso@usc.es}
\urladdr{http://web.usc.es/\~{}lalonso/}
\curraddr{Department of Mathematics\\
Purdue University\\
West Lafayette IN 47907, USA}

\author[A. Jerem\'{\i}as]{Ana Jerem\'{\i}as L\'opez}
\address{Departamento de \'Alxebra\\
Facultade de Matem\'a\-ticas\\
Universidade de Santiago de Compostela\\
E-15782  Santiago de Compostela, SPAIN}
\email{jeremias@usc.es}
\curraddr{Department of Mathematics\\
Purdue University\\
West Lafayette IN 47907, USA}

\author[M. P\'erez]{Marta P\'erez Rodr\'{\i}guez}
\address{Departamento de Matem\'a\-ticas\\
Escola Superior de En\-xe\-\~ne\-r\'{\i}a Inform\'atica\\
Campus de Ourense, Univ. de Vigo\\
E-32004 Ou\-ren\-se, Spain}
\email{martapr@uvigo.es}

\thanks{This work was partially supported by Spain's MCyT and E.U.'s
FEDER research project BFM2001-3241, supplemented by Xunta de Galicia grant PGDIT 01PX120701PR.}


\hyphenation{pseu-do}

\begin{abstract}
We make a systematic study of the infinitesimal lifting conditions of a pseudo finite type map of noetherian formal schemes. We recover the usual general properties in this context, and, more importantly, we uncover some new phenomena. We define a completion map of formal schemes as the one that arises canonically by performing the completion of a noetherian formal scheme along a subscheme, following the well-known pattern of ordinary schemes. These maps are \'etale in the sense of this work (but not adic). They allow us to give a local description of smooth morphisms. These morphisms can be factored locally as a completion map followed from a smooth \emph{adic} morphism. The latter kind of morphisms can be factored locally as an \'etale \emph{adic} morphism followed by a (formal) affine space. We also characterize \'etale \emph{adic} morphisms giving an equivalence of categories between the category of \'etale \emph{adic} formal schemes over a noetherian formal scheme $(\FX, {\CO}_{\mkern-2mu \FX})$ and the category of \'etale schemes over the ordinary scheme $(\FX, {\CO}_{\mkern-2mu \FX}/\CI)$, with $\CI$ an Ideal of definition. These results characterize completely the local behavior of formally smooth pseudo finite type (\ie \emph{smooth}) maps of noetherian formal schemes. In the final chapter, we study the basic deformation theory of (non necessarily adic) smooth morphisms.
\end{abstract}


\maketitle

\tableofcontents

\chapter*{Introduction}

In this work we study the infinitesimal lifting properties in the setting of locally noetherian formal schemes. Though this kind of study had not been undertaken before in a general way, it has been treated as needed in several scattered 
references. On the other hand, the novelty in our approach lies in the systematic consideration of a finiteness condition introduced by the first two authors, jointly with J. Lipman in an Oberwolfach seminar in 1996, see \cite{AJL1}. This condition was also considered, independently, by Yekutieli, see \cite{y}. In particular, this leads to the consideration of non-adic morphisms between formal schemes, something rather unusual in the literature. But let us start at the beginning and put our research in perspective.

\section*{Motivation and context of this work}

One of the characteristic features of Grothendieck's theory of schemes is that there may be nilpotent elements in the structural sheaf. From the function-theoretic point of view these looks unnatural, because nilpotents are non-zero functions vanishing at every point. However their consideration encodes useful infinitesimal information. For instance, different scheme structures on a point express tangent information of the embedding of the point in an ambient variety.

From this point of view it is only natural, given a closed subset of a variety, or more generally, of a scheme, to define a structure that encompasses all possible subscheme structures of the closed subset. These structures are called formal schemes and were introduced by Grothendieck in his Bourbaki talk ``G\'eom\'etrie formelle et g\'eom\'etrie alg\'ebrique'' \cite{gfga}. 

However, the origins of this notion can be traced back ---as lots of other modern constructions in algebraic geometry--- to Zariski. Specifically, the idea is present in his classical memoir ``Theory and applications of holomorphic functions on algebraic varieties over arbitrary ground fields'' \cite{z}. In this paper, Zariski establishes a theory of holomorphic functions along a subvariety, completely satisfactory in the affine case. His main application was the connectedness theorem. His inspiration was to look for an analogue of holomorphic functions along a subvariety in analytic geometry. In the global case, the lack of the language of sheaves was the technical problem that did not allow him to progress further.

Grothendieck's definition uses indeed the language of sheaves and is the natural extension to Zariski's ideas in the context of the theory of schemes. One of Grothendieck's achievements is the theorem of comparison between the cohomology of a proper map of schemes and the cohomology of its completion. This result gives the theorem on formal functions from where it follows Stein factorization and the connectedness theorem for schemes. Grothendieck developed some further applications as the Lefschetz-type theorems for Picard groups and fundamental groups in his 1962 seminar \cite{sga2}. 

Another question that arises is the algebraization problem. Given a morphism of formal schemes $f\colon\FX\to\FY$, is it the completion of a morphism of usual schemes $g \colon X \to Y$, \ie such that $f = \widehat{g}$? Under certain hypothesis the question is solved in \cite{gfga} for a proper map $f$. A modern exposition can be found in \cite{it}. Further development in the problem of algebraization has been done by Hironaka, Matsumura and Faltings.

After the initial period, formal schemes have always been present in the backstage of algebraic geometry but they were rarely studied in a systematic way. However, it has become more and more clear that the wide applicability of formal schemes in several areas of mathematics require such study. Let us cite a few of this applications. The construction of De Rham cohomology fo a scheme $X$ of zero characteristic embeddable in a smooth scheme $P$,  studied by Hartshorne \cite{ha2} (and, independently, by Deligne), is defined as the hypercohomology of the formal completion of the De Rham complex of the completion of $P$ along $X$. Formal schemes also play a key role in $p$-adic cohomologies (crystalline, rigid \dots).

Formal schemes are algebraic models of rigid analytic spaces. The rough idea goes back to Grothendieck but was developed by Raynaud and has been studied by Raynaud himself in collaboration with Bosch and L\"utkebohmert. Further development is due to Berthelot and de Jong.

Applications keep appearing in different branches of mathematics. We mention in this vein that Strickland \cite{st} has pointed out the importance of formal schemes in the context of (stable) homotopy theory.

We see then that the connections and applications of formal schemes are wide but after the important work exposed in the \textit{El\'{e}ments of  G\'{e}om\'{e}trie Alg\'{e}brique} (specifically, \cite{EGA1} and \cite{EGA31}) the general study of formal schemes has been very scarce. In particular a hypothesis that is almost always present in reference works is that the morphisms are adic, \ie that the topology of the sheaf of rings of the initial scheme is induced by the topology of the base formal scheme. This hypothesis is clearly natural because it guarantees that its fibers are usual schemes, therefore an adic morphism of formal schemes is, in the terminology of Grothendieck's school, a relative scheme over a base that is a formal scheme.

But there are important examples of maps of formal schemes that do not correspond to this situation. The first example that comes to mind is the natural map $\spf(A[[X]]) \to \spf(A)$ for an adic ring $A$. This morphism has a finiteness property that had not been made explicit until \cite{AJL1} (and independently, in \cite{y}). This property is called \emph{pseudo finite type}. The fact that pseudo finite type morphisms need not be adic allows fibers that are not usual schemes, and the structure of these maps is, in principle, subtler than the structure of adic maps.

Pseudo finite type maps arise in several significant contexts. To compute residues in a point of a variety over a perfect field we make the completion of the local ring at the given point. The geometric interpretation of this construction is to make a formal completion of the variety at the point and consider it a formal scheme over the base field. The resulting formal scheme is not adic but it is indeed of pseudo finite type. In the context of the extension of Grothendieck duality to formal schemes (see \cite{AJL1}) one generalizes the construction of relative residues and arrives to the following situation. We have a commutative diagram
\begin{diagram}[height=2.3em,w=2.3em,labelstyle=\scriptstyle]
Z 	         & \rTto^{\text{via }f}  & W \\
\dTinc^{i_1} &                       & \dTinc^{i_2} \\
X 	         & \rTto^f               & Y \\
\end{diagram}
where the horizontal map $f$ is a finite type map of formal schemes and the vertical maps, $i_1, i_2$ are closed immersions. We obtain a completion morphism $\widehat{f} \colon X_{/Z} \to Y_{/W}$ that is of pseudo finite type but it is adic only if $f^{-1}(W) = Z$, an important case but that does not always hold in practice.
Recently Lipman, Nayak and Sastry in \cite{LNS} made a canonical construction of Cousin complexes and its associated pseudofunctors for pseudo finite type maps of noetherian formal schemes.

One can interpret Hartshorne's approach to De Rham cohomology as associating a non adic formal scheme of pseudo finite type over the base field to a singular algebraic variety of characteristic zero whose cohomology behaves in a reasonable way. This non adic formal scheme turns out to be smooth in the sense explained below. 

In a different connection, Berthelot has pointed out the usefulness of non adic formal schemes of pseudo finite type over the formal spectrum of a discrete valuation ring as models of non quasi-compact rigid spaces.

A specially interesting aspect of pseudo finite type maps is their behavior with respect to the infinitesimal lifting properties. As in the case of usual schemes the theory contains more precise results when one adds a finiteness condition. In the present work we develop the infinitesimal lifting properties of pseudo finite type maps of formal schemes. We study the infinitesimal lifting conditions (formal smoothness, formal unramifiedness and formal \'etale condition) for formal schemes following the path (and using the results) of \cite{EGA44}. We call a morphism between locally noetherian formal schemes smooth, unramified or \'etale if it is of pseudo finite type and formally smooth, formally unramified or formally \'etale, respectively. Our definition of smooth maps agrees essentially with the one given in \cite[\S 2.4]{LNS}. Note also that condition (ii) in \cite[Definition 2.1]{y} corresponds to a smooth map in which the base is a ordinary noetherian scheme, so smooth formal embeddings are examples of smooth maps of formal schemes. We obtain characterizations analogous to the usual scheme case but, as a special feature of this work we define in general the completion of a formal scheme along a closed formal subscheme (see \ref{defcom}). A morphism arising this way is a completion map and can be characterized intrinsically as a certain kind of \'etale map (certainly non adic), see below. These kind of maps provide us the right tool to get a local and general description of \'etale and smooth morphisms. It is noteworthy to point out that a smooth map of formal schemes can have a non smooth usual subscheme as underlying reduced subscheme, as it is the case in Hartshorne's construction.

In this work, we also treat a basic version of deformation theory for smooth morphisms. We expect to be able to extend it to the non-smooth and non-separated cases in the future.
Also, in subsequent work we will apply our results to the decomposition  of the De Rham complex of a smooth formal scheme over a characteristic $p > 0$ field. We expect that other applications to the cohomology of singular schemes will follow.

\section*{Brief description of contents}

Though there is a description of every chapter and every section at their beginnings, it could be useful to give here an overview of the main concepts and results of this work.

In the first chapter we gather several results needed later, some of them generalizing well-known results for usual schemes that will be of use in the following chapters to study the infinitesimal lifting properties.

The first section recalls basic definitions from \cite[Chapter 10]{EGA1}. We discuss basic examples as formal affine spaces and formal disks that will be the source of examples and counterexamples used throughout the text. The section closes with the definition of the topological and algebraic dimension (or simply dimension) of a formal scheme. This will allow to control the rank of differentials for smooth morphisms and shows a difference between formal schemes and usual schemes. On usual schemes both dimensions agree, but on formal schemes the algebraic dimension is a property of the structural sheaf. For instance, let $K$ be a field, the underlying space of $\spf(K[[X]])$ is a point so its topological dimension is 0, but is algebraic dimension is 1 and this reflects an important aspect of the structure of this formal scheme.

In the next section we recall the notion of adic map and state in a systematic way its general properties. We define also the fibers of a map of formal schemes. Note that if a morphism is adic, then its fibers are usual schemes. We define relative dimensions at a point, both algebraic and topological, as the corresponding dimensions of the fibers. We introduce the class of pseudo closed embeddings as those maps of formal schemes that are a limit of a directed system of maps of ordinary schemes such that all of the maps in the system are closed embeddings. This class includes the closed embeddings defined in \cite[p. 442]{EGA1}. The section closes studying separated maps and open embeddings that, as already remarked in \emph{loc.~cit}.\/, behave in a similar fashion to the case of ordinary schemes.

The third section of the first chapter is devoted to a systematic exposition of the finiteness condition on morphisms of formal schemes that will be used in the sequel. We recall and state the main properties of finite, pseudo-finite, finite type and pseudo finite type morphisms. For the study of \'etale maps in chapter 3, one needs a reasonable non-adic analogous of quasi-finite maps. These are the \emph{quasi-coverings} (see definition \ref{defncuasireves}). We study its main properties and show that they are different from the pseudo finite type maps with finite fibers, that we call here \emph{pseudo quasi-finite} maps.

In the last section we discuss some basic properties of flat maps of formal schemes, including an interpretation of Grothendieck's local criterion of flatness in this context. We also define the completion of a formal scheme along a formal subscheme, a notion that has not appeared before but that it is modeled on the completion of a usual scheme along a closed subscheme. These maps will play an essential role later. We also discuss the invariance of certain properties of morphisms under completion.

The second chapter deals with the main tools and definitions of this work. In it, we define a suitable notion of the module of differentials. We recall Grothendieck construction for topological algebras and globalize it to get an appropriate relative cotangent sheaf for a pseudo finite type map of formal schemes. It turns out that the ``correct'' cotangent sheaf for the infinitesimal study of formal schemes is the completion of the sheaf of differentials defined in terms only of the structure of ringed spaces without topology on the structural sheaf. We study some general properties including a key one: its coherence. In the next section we see that the Fundamental Exact Sequences carry over this context.

In the fourth section of this chapter, we define the infinitesimal lifting conditions for formal schemes. We discuss the basic properties and prove, in down to earth terms, the fact that they are local on the base and \emph{on the source} (Proposition \ref{condinflocal}). Next, it is proved that a completion is an \'etale morphism, an statement that has no significant counterpart in the case of ordinary schemes and that will be a ``building block'' for the local characterization of smooth maps of formal schemes. It is also proven in this section that the infinitesimal conditions carry over to a completion of a morphism of formal schemes along suitable closed subschemes.

The closing section of this chapter relates the module of differentials studied in the first sections to the infinitesimal lifting conditions. It is proved in this context that the vanishing of the module of differentials is equivalent to the map being unramified, and that a smooth morphism is flat and its module of differentials is locally free. The section ends with Zariski's Jacobian criterion for formal schemes.

Chapter three contains the main results of this work, especially those that show the difference in infinitesimal behavior of formal schemes with respect to the behavior of usual schemes. In the first part, formed by the first three sections there is a local study of the infinitesimal lifting properties and we relate unramified, smooth and \'etale morphisms of formal schemes to the same properties of ordinary schemes. One can express a morphism of formal schemes $f \colon \FX \to \FY$ as a direct limit 
\[f = \dirlim  {n \in \NN} f_{n}	\]
of morphisms between the underlying ordinary schemes $f_n \colon X_n \to Y_n$.
Then, we can pose an important question: Do the infinitesimal  properties of $f_n$ imply those of $f$? We show that $f_n$ is unramified \emph{for any} $n$ if, and only if, $f$ is unramified. We show also that the fibers of an unramified map are usual schemes and characterize pseudo closed embeddings as unramified maps such that $f_0$ is unramified (with the previous notation). Given a morphism of noetherian formal schemes it is proved that if the underlying maps of formal schemes are smooth, then the original morphism of formal schemes is. The reciprocal result is false for smooth, and therefore, for \'etale morphisms. We characterize smooth morphisms as those maps that are flat and whose fibers are smooth formal schemes over a field. We also note that a smooth morphism factors locally as an \'etale map and a (formal) affine space and obtain a matricial Jacobian criterion for smoothness. All of this is applied to \'etale morphisms, and we obtain, for instance, that to be an \'etale map is equivalent to be flat and unramified and also to be smooth and a quasi-covering.

The next section gives three results that explain the structure of \'etale maps. First, we prove that an open embedding is the same as an adic \'etale radical morphism (Theorem \ref{caractencab}). This generalizes a well-known property of usual schemes, see \cite[Th\'eor\`eme (17.9.1)]{EGA44}. We characterize (Theorem \ref{caracmorfcompl}) completion morphisms as those \'etale maps whose underlying map of schemes is a closed embedding, which is also equivalent to be a flat pseudo closed embedding. This shows completion morphisms as a new class of \'etale morphisms of formal schemes. We also prove that, given a formal scheme $\FY$ there is an equivalence between the category of \'etale \emph{adic} formal schemes over $\FY$ and \'etale schemes over $Y_0$ where 
\[\FY = \dirlim  {n \in \NN} Y_{n}.\] 
The case in which $\FY$ is smooth over a base an ordinary scheme has appeared in \cite[Proposition 2.4]{y}.
This result (and the cited previous special case) is a reinterpretation in this context of \cite[Th\'eor\`eme (18.1.2)]{EGA44}, which, in fact, is used in the proof.

With all these results we prove the local structure theorem of smooth morphisms on formal schemes that says that if $f:\FX \to \FY$ be a morphism smooth at $x \in \FX$, then there exists an open subset $\FU \subset \FX$ with $x \in \FU$ such that  $f|_{\FU}$ factors as
\[
\FU \xto{\kappa} \FX' \xto{f'} \FY
\]
where $\kappa$ is a completion morphism and $f'$ is an \emph{adic} smooth  morphism. This is Theorem \ref{tppall}. In a sense, together with the previous result, it explains the local structure of a smooth map. An \'etale map can be locally factored as a completion followed by an adic \'etale map which are classified by the fundamental group of the underlying scheme of the base formal scheme. This gives a way for understanding the local properties of these morphisms.

The work closes with an study of basic deformation theory of smooth maps of formal schemes. The results are quite similar to those that are standard for ordinary schemes. There is an obstruction that lives in a cohomology group of order one for extending a map over an infinitesimal neighborhood over a smooth formal scheme. Also, we consider the following situation: a smooth morphism $f_{0}:\FX_{0} \to \FY_{0}$ and a closed immersion $\FY_{0} \inc \FY$ defined by a square zero Ideal. In this case: if there exists a smooth $\FY$-formal scheme $\FX$  such that  $\FX \times_{\FY} \FY_{0} = \FX_{0}$, it is unique when $\ext^{1}(\om^{1}_{\FX_{0}/\FY_{0}}, f_{0}^{*}\CI)=0$ (Proposition \ref{deform4}). Furthermore, we show that there exist such a smooth $\FY$-formal scheme $\FX$ if a certain element in $\ext^{2}(\om^{1}_{\FX_{0}/\FY_{0}}, f_{0}^{*}\CI)$ vanishes, see Proposition \ref{obstrext2}. In particular, if $f_{0}$ is a map of affine schemes, such a $\FX$ always exists.

\section*{A note on terminology }

Whenever a formal scheme morphism $f \colon \FX \to \FY$ is expressed as a direct limit of morphism of their underlying ordinary schemes $f_{n}:\FX_n \to \FY_n$ (with $n \in \NN$) the question about how a property $\CP$ that holds for all the maps in the system $\{f_n\}_{n \in \NN}$ is reflected on the map $f$ arises. It has been customary to call ``pseudo-$\CP$'' the properties of locally noetherian formal schemes such that $\CP$ holds for all $f_n$ but that, in general, do not hold for the corresponding morphisms $f$ of formal schemes unless $f$ is adic . This is the case for morphisms of pseudo finite type or pseudo closed embeddings.

In previous stages of this work we used the term \emph{pseudo-smooth} for what it is called now smooth. By the preceding discussion, the old terminology would not be adequate because, as remarked previously, the underlying maps of ordinary schemes of a smooth map of formal schemes may not be smooth. Also, in this way, our terminology agrees with the one used in \cite{LNS}. Note that the notion is new even though the terminology is classic.

\section*{Acknowledgements}
We have benefited form conversations on these topics and also on terminology with Joe Lipman, Suresh Nayak and Pramath Sastry. We also want to thank Luis Narv\'aez, Adolfo Quir\'os and Norbert Schappacher for encouragement and further comments. The authors thank the Mathematics department of Purdue University for hospitality and support. 

The diagrams were typeset with Paul Taylor's \texttt{diagrams.tex}.

\vspace{20pt}

\mainmatter

\chapter{Formal schemes}\label{cap1}
\setcounter{equation}{0}
In this chapter we  establish the concepts and notations about  formal schemes that will be useful for the reading of  this memory. In  Section \ref{sec11} we recall basic facts of the theory of formal schemes illustrated with  examples. Sections \ref{sec12} and \ref{sec13} treat several properties of  morphisms of formal schemes (topological and  finiteness properties). We also introduce notions that are essential later in the development of this work, as pseudo closed immersions, quasi-coverings and the fiber of a morphism. The last section is devoted to the study flat morphisms. There we define completion morphisms, a kind of flat morphisms that topologically are closed immersions.

\section[The  category  of formal schemes]{The  category  of locally noetherian formal schemes} \label{sec11}

We will begin by recalling briefly some basic definitions and results about locally noetherian formal schemes. Of course, for a complete treatment  we refer the reader to \cite[\S 10]{EGA1}. We will give some detailed examples of formal schemes, which we will refer along this exposition, like the affine formal scheme and the formal disc. We will also define the (algebraic) dimension of a locally noetherian formal scheme.

\begin{defn}  \label{defschemeformalafin}
Let $A$ be an $I$-adic ring\footnote{A topological ring $A$ is called \emph{$I$-preadic} if $(I^{n})_{n \in \NN}$ is a  fundamental system of neighborhoods of  $0$, and we say that the topology of  $A$ is  \emph{$I$-adic}. 
Moreover if $A$ is separated and complete for the $I$-adic topology we say that $A$ is \emph{$I$-adic} (\emph{cf.}  \cite[(\textbf{0}.7.1.9)]{EGA1}). 
Trivially, every ring $A$ is separated and complete for the  discrete topology, the topology given by any nilpotent ideal of $A$.}. 
The \emph{formal spectrum of $A$}, written $\spf (A)$, is the topological ringed space (\emph{cf.} \cite[(\textbf{0}.4.1.1)]{EGA1}) whose underlying topological space is the closed subset $V (I) \subset \spec (A)$ and whose sheaf of topological rings is \[\invlim {n \in \NN} \CO_{X_{n}}\] where, for all $n \in \NN$, $\CO_{X_{n}}$ is the structural sheaf of rings of the affine scheme $X_{n} = \spec (A/ I^{n+1})$. 

A topologically ringed space $(\FX, \CO_{\FX})$ is an 
\emph{affine formal scheme} if is isomorphic to the formal spectrum
$\spf (A)$ of an $I$-adic ring $A$. If $A$ is noetherian we will say that $\FX$ is an
\emph{affine noetherian formal scheme}. 

Note that an affine scheme is an affine formal scheme for the discrete topology.
\end{defn}

\begin{rem}
The concepts of formal spectrum and affine formal scheme are established in \cite[\S 10.1]{EGA1} for a general class of topological rings (the \emph{admissible} rings \cite[(\textbf{0}.7.1.2)]{EGA1}). What we call here formal spectrum or affine formal scheme is called in \cite{EGA1} adic formal spectrum or adic affine formal scheme, respectively.
\end{rem}

\begin{parraf}
If $\FX = \spf(A)$ is an affine formal scheme with  $A$ an $I$-adic ring, then:
\begin{enumerate} 
\item
$\ga (\FX, \CO_{\FX}) \cong A$ as topological rings \cite[(10.1.3)]{EGA1}. 
\item
Given  $f \in A$,  if $\fD (f):=\D(f) \cap \FX$, with $\D(f) = \spec(A_{f})$, the collection of open subsets  $\{\fD (f)\}_{f \in A}$ form a base for the topology of $\spf(A)$. Besides, the isomorphism
\[
(\fD (f), \CO_{\FX}|_{\fD (f)}) \cong \spf(A_{\{f\}})  \textrm{ \cite[(10.1.4)]{EGA1}},
\]
where $A_{\{f\}}$ denotes the completion of $A_{f}$ respect to the $I_{f}$-adic topology\footnote{\cite[\S \textbf{0}.7.6]{EGA1} Given a $J$-adic ring $A$ and $S \subset A$ a multiplicative subset,  $A\{S^{-1}\}$ is the completion of $S^{-1}A$ respect to the $S^{-1} J$-adic topology. If $M$ is an $A$-module, $M \{S^{-1}\}$ denotes the completion of $S^{-1} M$ with respect to the $S^{-1}J $-adic topology, that is, \[M \{S^{-1}\}= \invlim {n \in \NN} S^{-1} M/ (S^{-1}J)^{n} M.\] In particular, if $S=\{f^{n},/, n \in \NN\}$ then the completion of $A_{f}$ and of $M_{f}$ for the $J_{f}$-adic topology are denoted $S_{\{f\}}$ and $M_{\{f\}}$, respectively.}, shows that $\fD (f) = (\fD (f), \CO_{\FX}|_{\fD (f)})$ is an affine formal scheme. 
\item
For each $x \in \FX$, let $\fp$ be the corresponding open prime ideal in $A$. The \emph{local ring of $\FX$ at $x$} is\footnote{Given $A$ a $J$-adic ring, $\fp \subset A$  a open prime and $M$ a $A$-module, put \[M_{\{\fp\}} := \dirlim {f \notin \fp} M_{\{f\}}.\] Note that $A_{\{\fp\}}$ is a  $J_{\{\fp\}}$-preadic ring.}:
\[
\CO_{\FX,x} = \dirlim {f \notin \fp} \ga(\fD (f), \CO_{\fD (f)}) = \dirlim {f \notin \fp} A_{\{f\}} =: A_{\{\fp\}} 
\]
It holds that $\CO_{\FX,x}$ is a local noetherian ring  with maximal ideal $\fp\CO_{\FX,x}$ and its residual field is $k(x) = A_{\{\fp\}} / \fp A_{\{\fp\}} = A_{\fp} / \fp A_{\fp}$ (\emph{cf.} \cite[(\textbf{0}.7.6.17) and (\textbf{0}.7.6.18)]{EGA1}).
Observe that  $\CO_{\FX,x}$ need not be a complete ring for the $I\CO_{\FX,x}$-adic topology.
\end{enumerate}
\end{parraf}

\begin{ex} \label{afdis}
Let $A$ be a commutative ring and $A[\mathbf{T}]$ the polynomial ring  with coefficients in $A$ in the indeterminates $\mathbf{T} = T_{1},T_{2},\ldots,T_{r}$. 
Let us assume that $A$ is an $I$-adic noetherian ring and let us denote $I[\mathbf{T}] = I \cdot A[\mathbf{T}]$ the set of polynomial with coefficients in $I$.
\begin{enumerate}
\item \label{af} 
The \emph{ring of restricted formal series} $A\{\mathbf{T}\}$ is the completion of $A[\mathbf{T}]$ respect to the  $I[\mathbf{T}]$-adic topology. If  $I\{\mathbf{T}\}$ is the ideal of the restricted formal series with  coefficients in $I$, it holds that $I\{\mathbf{T}\} = I \cdot A\{\mathbf{T}\}$ and that $A\{\mathbf{T}\}$ is an $I\{\mathbf{T}\}$-adic noetherian ring (\emph{cf.} \cite[(\textbf{0}.7.5.2)]{EGA1}). We call \emph{affine formal $r$-space over $A$} or the \emph{affine formal space of  dimension $r$ over $A$} the formal scheme
\[
\BA_{\spf (A)}^{r} = \spf (A\{\mathbf{T}\})
\]
Observe that the underlying topological space of $\BA_{\spf (A)}^{r}$ is the affine  $r$-space over $A/I$, $\BA_{\spec (A/I)}^{r} = \spec(\frac{A}{I}[\mathbf{T}])$.
\item  \label{dis} 
The  \emph{formal power series ring}  $A[[\mathbf{T}]]$ is the completion of $A[\mathbf{T}]$ with respect to the $(I[\mathbf{T}] + \langle\mathbf{T}\rangle)$-adic topology. If $I[[\mathbf{T}]]$ denotes the ideal of formal  power series with  coefficients in $I$ and $[[\mathbf{T}]]$ is the ideal of formal series without independent term, it is well-known that $I[[\mathbf{T}]] =I \cdot A[[\mathbf{T}]]$,  $[[\mathbf{T}]] = \langle\mathbf{T}\rangle \cdot A[[\mathbf{T}]]$ and that $A[[\mathbf{T}]]$ is an  $(I[[\mathbf{T}]] + [[\mathbf{T}]])$-adic noetherian ring (\emph{cf.} \cite[Theorem 3.3 and Exercise 8.6]{ma2}). We define  the \emph{formal $r$-disc over $A$} or \emph{formal disc of dimension $r$ over $A$} as
\[
\BD_{\spf (A)}^{r} = \spf (A[[\mathbf{T}]])
\]
The underlying topological space of $\BD_{\spf (A)}^{r}$ is $\spec (A/I)$. (Note that if in the ring $A[[\mathbf{T}]]$ we consider the $[[\mathbf{T}]]$-adic topology the underlying topological space of the corresponding affine formal scheme would be $\spec(A)$).
\end{enumerate}
\end{ex}

\begin{defn} \cite[(10.3.1)]{EGA1}	\label{idaf}
Given $A$ an $I$-adic ring let $\FX = \spf(A)$. If
$J$ is an ideal of definition\footnote{Given $A$ an $I$-preadic ring, an \emph{ideal of definition of $A$} is an ideal that determines the $I$-adic topology in  $A$ (\emph{cf.}  \cite[(\textbf{0}.7.1.2)]{EGA1}).} of $A$ then\footnote{In general, if  $A$ is a $J$-adic ring and $M$ an $A$-module,  in \cite[(10.10.1)]{EGA1} it is defined the following sheaf of modules over $\spf(A)$
\[M^{\tr}= \invlim {n \in \NN} \frac{\widetilde{M}}{\widetilde{J}^{n+1}\widetilde{M}}\,.\]
We will recall some properties of this construction from \ref{deftriangulito} on.} 
\[J ^{\tr} = \invlim {n \in \NN} \widetilde{J} / \widetilde{J ^{n}} \subset  \CO_{\FX}\] is called an \emph{Ideal of definition of $\FX$}. In particular,
$I ^{\tr}$ is an  Ideal of definition of $\FX$. 

Given $J ^{\tr}$ an Ideal of definition of $\FX$ it holds that: 
\[(\FX, \CO_{\FX}/ J ^{\tr}) = \spec(A/J)\]
therefore, any Ideal of definition of $\FX$ determines its underlying topological space. 
\end{defn}

\begin{ex}
With the hypothesis and notations of Example \ref{afdis}, we have that $I\{\mathbf{T}\} ^{\tr}$ is an Ideal of definition of $\BA_{\spf (A)}^{r}$ and $(I[[\mathbf{T}]] + [[\mathbf{T}]])^{\tr}$ is an Ideal of definition of $\BD_{\spf (A)}^{r}$.
\end{ex}

\begin{defn}
A topologically ringed space $(\FX, \CO_{\FX})$ is called a
\emph{formal scheme} if for all $x \in \FX$ there exists an open subset
$\FU$ of $\FX$ with $x \in \FU$ and such that $(\FU, \CO_{\FX}|_{\FU})$ is an affine formal scheme \cite[(10.4.2)]{EGA1}. When  the open subsets $\FU$ are affine noetherian formal schemes we say that $(\FX,\CO_{\FX})$ is a  \emph{locally noetherian formal scheme}. We will often write $\FX$ instead of $(\FX,\CO_{\FX})$.
\end{defn}

\begin{ex} \label{excompl}
Given a scheme $X$ and $X' \subset X$ a closed subscheme determined by a coherent Ideal $\CI$ of $\CO_{X}$, we call \emph{completion of $X$ along $X' $} and write 
$(X_{/X' },\CO_{X_{/X' }})$ or, in short, $X_{/X' }$, the topologically ringed space with underlying topological space $X' $ and whose sheaf of topological rings is 
\[\CO_{X_{/X' }} = \invlim {n \in \NN} \CO_{X}/\CI^{n+1}.\] The topologically ringed space $X_{/X' }$ is a formal scheme \cite[(10.8.3)]{EGA1}. A formal scheme isomorphic to one of this type is called  \emph{algebraizable}. In particular, if $X= \spec(A)$ is an affine scheme  and $\CI =  \widetilde{I}$, with $I$ a finitely generated ideal of $A$, it holds that $X_{/X' } = \spf(\widehat{A})$ where $\HA$ is the separated complete ring of $A$ for the $I$-adic topology. Then every affine noetherian formal scheme is algebraizable and, therefore, every locally noetherian formal scheme is locally algebraizable. In \cite[\S 5]{hima} one can find an example of a locally noetherian formal scheme which is not algebraizable. 
\end{ex}

\begin{defn}
Let $(\FX, \CO_{\FX})$ and $(\FY, \CO_{\FY})$ be formal schemes. A morphism  $f:(\FX, \CO_{\FX}) \to (\FY, \CO_{\FY})$ of topologically ringed spaces such that $\forall x \in \FX$, $f^{\sharp}_{x}:  \CO_{\FY,f(x)} \to \CO_{\FX,x}$ is a local homomorphism is a  \emph{morphism of formal schemes} \cite[(10.4.5)]{EGA1}. We denote by $\Hom(\FX,\FY)$ the set of morphisms of formal schemes from $\FX$ to $\FY$.
\end{defn}

\begin{parraf}
Formal schemes together with morphisms of  formal schemes form a category  that we denote by $\sfa$ and, we'll write $\sfn$ for the full subcategory  of $\sfa$ which objects are the locally noetherian formal schemes. Affine formal schemes form a  full subcategory  of $\sfa$, we will denote it by $\sfaa$. In the same way the affine noetherian formal schemes are a full subcategory  of $\sfn$, written as $\sfna$. 
\end{parraf}

The next result, that we will call ``formal Descartes duality'' provides a contravariant equivalence of categories that generalizes the well-known relation between the categories of rings and schemes (``Descartes duality'').

\begin{propo} \cite[(10.2.2)]{EGA1} 
The functors
\begin{equation} \label{equiv}
A \leadsto \spf(A) \qquad \mathrm{and} \qquad
\FX \leadsto \ga(\FX, \CO_{\FX})
\end{equation}
define a duality between the category of adic rings and $\sfaa$. Moreover, the restriction of these functors to the category of adic noetherian rings and to $\sfna$ is also a duality of categories.
\end{propo}

\begin{parraf} \label{morf}
As a consequence of the proposition we get the following morphisms of formal schemes:
\begin{enumerate}
\item \label{morf1}
Given $A$ an $I$-adic ring, the morphism $A \to A_{\{f\}}$ induces the canonical inclusion in $\sfaa$
\[
\fD (f) \inc \spf(A).
\]
	
\item \label{morf2}
If $A$ is an $I$-adic noetherian ring, the canonical inclusions   \[A \inc  A\{T_{1},T_{2},\ldots,T_{r}\} \inc  A[[T_{1},T_{2},\ldots,T_{r}]]\] are continuous morphisms for the adic topologies established in Example \ref{afdis} and, therefore, induce canonical morphisms in $\sfna$
\[\BD_{\spf(A)}^{r} \to \BA_{\spf(A)}^{r} \to \spf(A).\] 

\item \label{morf3}
With the notations of  Example  \ref{excompl},  the morphisms \[\CO_{X} \to \CO_{X}/\CI^{n+1}\] induce a canonical morphism in $\sfa$ 
\[
X_{/X' } \xto{\kappa} X
\]
called \emph{the completion morphism of $X$ along $X'$} that, as a map of topological spaces is the inclusion of $X'$ in $X$.
Locally, for an affine open subset given by an ideal $I$ of the noetherian ring $A$, it is given by $\spf(\widehat{A}) \to \spec(A)$, where $\widehat{A}$ is the $I$-adic completion of $A$.
\end{enumerate}
\end{parraf}

\begin{defn} \cite[(10.5.1)]{EGA1}
Let $(\FX, \CO_{\FX})$ be a  formal scheme. An Ideal $\CJ$ of $ \CO_{\FX}$ is an \emph{Ideal of definition of $\FX$} if for all $x \in \FX$ there exists an open subset 
$\FU \subset \FX$ with $x \in \FU$ such that $(\FU, \CO_{\FX}|_{\FU})$ is an affine
formal scheme  and $\CJ|_{\FU}$ is an Ideal of definition of $\FU$.
\end{defn}

\begin{ex} \label{excom}
Let $X$ be a (usual) scheme.
\begin{enumerate}
\item
Any locally nilpotent Ideal $\CJ \subset \CO_{X}$ is an Ideal of definition of $X$.
\item
If $X'  \subset X$ is a closed  subscheme defined by a coherent Ideal $\CI \subset \CO_X$, the  Ideal $\CI_{/X' } := \invlim {n \in \NN} \CI / \CI^{n+1}$ is  an Ideal of definition of $X_{/X' }$ and it holds that $\CO_{X_{/X' }}/(\CI_{/X' })^{n+1} = \CO_{X}/\CI^{n+1}$. Moreover, if $X= \spec(A)$, $\CI =  \widetilde{I}$ and $\kappa:X_{/X' }=\spf(\HA) \to X =\spec(A)$ is the completion morphism of $X$ along $X' $, we have that
\[
\CI_{/X' } \underset{\textrm{\ref{idaf}}}= I^{\tr} = \kappa^{*} \widetilde{I}
\]
\end{enumerate}
\end{ex}

The next proposition guarantees the existence of  Ideals of  definition for locally noetherian formal schemes.

\begin{propo} \label{idgr} 
Given $(\FX, \CO_{\FX})$ a locally noetherian formal scheme, there exists $\CJ$ the greatest Ideal of  definition of  $\FX$. Besides, $(\FX, \CO_{\FX}/\CJ)$ is a reduced (usual) scheme.
\end{propo}

\begin{proof}
\cite[(10.5.4)]{EGA1}.
\end{proof}

\begin{parraf} \label{existidefmorf}
As a consequence of this Proposition, it is easily seen  that given $f:\FX \to \FY$ a morphism in $\sfn$, if $\CK \subset \CO_{\FY}$ is an Ideal of  definition, then there exists an Ideal of  definition $\CJ \subset \CO_{\FX}$ such that  $f^{*}(\CK) \CO_{\FX} \subset \CJ$, see \cite[(10.6.10)]{EGA1}.
\end{parraf}

An affine formal scheme supports all the possible structures of closed subscheme corresponding to its underlying topological space. Therefore a formal scheme can be considered a special case of ind-scheme, \ie a direct limit of usual schemes. Next, we elaborate on this idea in the case of locally noetherian formal schemes. 

\begin{parraf} \label{notlim}
Given $(\FX, \CO_{\FX})$ in $\sfa$ and $\CJ$ an Ideal of  definition of  $\FX$  such that  $\CJ/\CJ^{2}$ is a $\CO_{\FX}/\CJ$-module of  finite type\footnote{Whenever $\FX$ is in $\sfn$ for every Ideal of  definition $\CJ \subset \CO_{\FX}$ it holds that $\CJ/\CJ^{2}$ a $\CO_{\FX}/\CJ$-module of  finite type} it holds  that $\CJ^{n+1}$ is  an Ideal of  definition and  that $X_{n}=(\FX, \CO_{\FX}/\CJ^{n+1})$ is a (usual) scheme  with the same topological space as $\FX$, $\forall n\in \NN$ (\emph{cf.} \cite[(10.5.1) and (10.5.2)]{EGA1}).  Let $i_{n}:X_{n} \to \FX$
be the canonical morphism determined by $\CO_{\FX} \to \CO_{\FX}/\CJ^{n+1}$, for each $n \in \NN$. If for $m \ge n \ge 0$, $i_{mn}: X_{n} \to X_{m}$ are  the natural closed immersions defined by $\CO_{\FX}/\CJ^{m+1} \to \CO_{\FX}/\CJ^{n+1}$ we have that 
\[i_{m} \circ i_{mn} = i_{n}\] and, therefore, $\{X_{n}, i_{mn}\}$ is a direct system in $\sfa$. There results that $(\FX, \CO_{\FX})$  is the direct limit  of   $\{X_{n}, i_{mn}\} \textrm{ in } \sfa$. We will say that $\FX$ it is expressed as 
\[\FX = \dirlim  {n \in \NN} X_{n}\] with respect to the ideal of definition $\CJ$
and leave implicit that the schemes $\{X_{n}\}$ are defined by the powers of the Ideal of  definition $\CJ$. Notice that if  $\FX$ is in $\sfn$ the  schemes $X_{n}$ are locally noetherian, for all $n \in \NN$.

Suppose  that  $\FX = \spf(A)$ is in $\sfaa$ and that  $ \CJ=J^{\tr}$ is an Ideal of  definition of  $\FX$ being $J \subset A$ an ideal of  definition such that  $J/J^{2}$ is a $(A/J)$-module of  finite type. Applying \cite[(10.3.6)]{EGA1} there results that $\CJ^{n+1} =  (J^{\tr})^{n+1} = (J^{n+1})^{\tr}$ so we have that $X_{n} = \spec(A/J^{n+1})$, for  all $n \in \NN$ (see Definition \ref{idaf}). In this case,  the expression $\FX  = \dirlim {n \in \NN} X_{n} $ corresponds to the equality  $A=\invlim {n \in \NN} A/J^{n+1}$ through the duality (\ref{equiv}).
\end{parraf}

\begin{parraf} \label{lim}
Given $f:\FX \to \FY$ a morphism in $\sfa$, let $\CJ \subset \CO_{\FX}$ and $\CK \subset \CO_{\FY}$ be Ideals of  definition such that $\CJ/\CJ^{2}$  is a $(\CO_{\FX}/\CJ)$-Module of finite type and  $\CK/\CK^{2}$ is  a $\CO_{\FY}/\CK$-Module of  finite  type, and satisfying that $f^{*}(\CK)\CO_{\FX} \subset \CJ$. Let us write  
\[\FX = \dirlim  {n \in \NN} X_{n}\qquad \textrm{ andÊ}\qquadÊ\FY = \dirlim  {n \in \NN} Y_{n}\]
where, for each $n \in \NN$, $X_{n}:=(\FX, \CO_{\FX}/\CJ^{n+1})$ and $Y_{n}:=(\FY, \CO_{\FY}/\CK^{n+1})$. For each $n \in \NN$,, there  exists a unique morphism of  schemes $f_{n}: X_{n} \to Y_{n}$ such that, for $m \ge n \ge 0$ the diagrams:
\begin{diagram}[height=2em,w=2em,labelstyle=\scriptstyle]
\FX		&		\rTto^{f} 		&	\FY\\
\uTto_{ i_{m}} 	&				&	\uTto^{ j_{m}} \\
X_{m}		&	\rTto^{f_{m}}	&	Y_{m}\\
\uTto_{i_{mn}}	&				&\uTto^{j_{mn}}\\
X_{n}		&	\rTto^{f_{n}}	&	Y_{n}, 	&	   \\
\end{diagram} commute.
The morphism $f$ is the direct limit of the system $\{f_{n},i_{mn},j_{mn}\}$ associated  to the Ideals of  definition $\CJ \subset \CO_{\FX}$ and $\CK \subset \CO_{\FY}$ and we will write 
\[f = \dirlim  {n \in \NN} f_{n}\]
\cite[(10.6.7), (10.6.8) and (10.6.9)]{EGA1}.

If $f$ is in $\sfn$ and $\CK \subset \CO_{\FY}$ is an Ideal of  definition, there always exists an Ideal of  definition $\CJ \subset \CO_{\FX}$ such that  $f^{*}(\CK) \CO_{\FX} \subset \CJ$ (see \ref{existidefmorf}). The corresponding morphisms $f_{n}$ are in the  category   of locally noetherian schemes,  for all $n \in \NN$.

Suppose  that  $f:\FX = \spf(A)\to \FY= \spf(B)$  is in $\sfaa$ and that there exist $\CJ=J^{\tr} \subset \CO_{\FX}$ and $ \CK= K^{\tr} \subset \CO_{\FY}$ Ideals of  definition such that $J/J^{2}$ is a $A/J$-module of  finite type and $K/K^{2}$ is a $B/K$-module of  finite type. Let $\phi: B \to A$ be the continuous homomorphism  of adic rings that  corresponds to $f$ through the duality (\ref{equiv}). Then $f^{*}(\CK)\CO_{\FX} \subset \CJ$ if, and only if, $KA \subset J$, and the morphisms $f_{n}$  correspond  with the canonical morphisms $B/K^{n+1} \to A/J^{n+1}$ induced by $\phi$.
\end{parraf}

\begin{parraf} \label{complim}
Let $f:\FX \to \FY$ and  $g:\FY \to \FS$ be two morphisms in $\sfa$ and consider  $\CJ \subset \CO_{\FX}$,  $\CK \subset \CO_{\FY}$ and $\CL \subset \CO_{\FS}$ Ideals of  definition such that $\CJ/\CJ^{2}$ is  a $\CO_{\FX}/\CJ$-Module of finite type, $\CK/\CK^{2}$ is a $\CO_{\FY}/\CK$-Module of  finite type and $\CL/\CL^{2}$ is a $\CO_{\FS}/\CL$-Module of  finite type, and satisfying that $f^{*}(\CK)\CO_{\FX} \subset \CJ$ and $g^{*}(\CL)\CO_{\FX} \subset \CK$. By \ref{lim},  we get
\[f = \dirlim  {n \in \NN} f_{n},\quad g = \dirlim  {n \in \NN} g_{n},\] where $f_{n}:(\FX, \CO_{\FX}/\CJ^{n+1}) \to (\FY, \CO_{\FY}/\CK^{n+1})$ and $g_{n}:(\FY, \CO_{\FY}/\CK^{n+1}) \to (\FS, \CO_{\FS}/ \CL^{n+1})$ are the induced morphisms by $f$ and $g$, respectively, for all $n \in \NN$. Then, it holds that
\[g \circ f = \dirlim  {n \in \NN}g_{n}  \circ  f_{n}\]
\end{parraf}

\begin{notaci}
Henceforth we will use systematically the following notation:
\begin{enumerate}
\item
Given $\FX$ in $\sfa$ and $\CJ \subset \CO_{\FX}$ an Ideal of  definition such that  $\CJ/\CJ^{2}$ is a $\CO_{\FX}/\CJ$-Module of  finite type, for all $n \in \NN$, $X_{n}$ will denote the scheme $(\FX, \CO_{\FX}/\CJ^{n+1})$.
\item
If $f:\FX \to \FY$ is a morphism in $\sfa$, given $\CJ \subset \CO_{\FX}$ and $\CK \subset \CO_{\FY}$ Ideals of  definition such that $\CJ/\CJ^{2}$  is a $(\CO_{\FX}/\CJ)$-Module of  finite type and $\CK/\CK^{2}$ is  a $\CO_{\FY}/\CK$-Module of  finite type, and satisfying that $f^{*}(\CK)\CO_{\FX} \subset \CJ$, $f_{n}:X_{n}:= (\FX, \CO_{\FX}/\CJ^{n+1}) \to Y_{n}:=(\FY, \CO_{\FY}/\CK^{n+1})$ will be the morphism of  schemes induced by $f$, for all $n \in \NN$.
\end{enumerate}
\end{notaci}

\begin{defn}
Let $\FS$ and $\FX$ in $\sfa$. If there exists a  morphism $\FX \to \FS$ we say that $\FX$ is a \emph{$\FS$-formal scheme}. If $\FX=X$ is a (usual) scheme we say that $X$ is a \emph{$\FS$-scheme}. Notice that, every formal scheme $\FX$ in $\sfa$ is a $\spec(\ZZ)$-formal scheme.

Given $\FX \to \FS$ and  $\FY \to \FS$ morphisms in $\sfa$, a morphism $f:\FX \to \FY$ in $\sfa$ is a $\FS$-morphism if the diagram
\begin{diagram}[height=2em,w=2em,labelstyle=\scriptstyle]
\FX 	&	\rTto^f  &	 \FY\\
	&    	\rdTto	& \dTto\\
 &	&	    	\FS\\
\end{diagram} commutes. 
The $\FS$-formal schemes with the  $\FS$-morphisms are a  category  and we will denote by $\Hom_{\FS}(\FX,\FY)$ the set of  $\FS$-morphisms between $\FX$ and $\FY$.
\end{defn}

\begin{propo} \cite[(10.7.3)]{EGA1}
Let $\FX$ and $\FY$ be two  $\FS$-formal schemes. Then there exists $\FX \times_{\FS} \FY$, the fiber product  of  $\FX$ and $\FY$ in the  category  of   $\FS$-formal schemes.
Moreover, if $\FX = \spf(A),\, \FY= \spf(B),\, \FS= \spf(C)$, with $A,\, B$ and $C$  adic rings, the fiber product $\FX \times_{\FS} \FY$ is given by $\spf(A \tc_{C} B)$\footnote{If $J \subset A ,\, K \subset B$ and $L\subset C$ are ideals of  definition such that $LB \subset K$ and $LA \subset J$, $A \tc_{C} B$ is the completion  of  $A \otimes_{C} B$ with respect to the  $(J(A \otimes_{C} B)+ K(A \otimes_{C} B))$-adic topology or, equivalently, \[A \tc_{C} B =\invlim {n \in \NN} A/J^{n+1} \otimes_{C/L^{n+1}} B/K^{n+1}\] \cite[(0.7.7.1) and (0.7.7.2)]{EGA1}. }.
\end{propo}

\begin{rem}
Like the category of  locally noetherian schemes  (\emph{cf.} \cite[(3.2.5)]{EGA1}), $\sfn$ is not stable for fiber products.
\end{rem}

\begin{parraf} \label{cclim}
\cite [(10.7.4)]{EGA1} Given $f:\FX \to \FS$ and $g:\FY \to \FS$ in $\sfn$, let $\CJ \subset \CO_{\FX},\, \CK \subset \CO_{\FY}$ and $\CL \subset \CO_{\FS}$ be Ideals of  definition such that $f^{*}(\CK)\CO_{\FX} \subset \CJ$ and $g^{*}(\CL)\CO_{\FX} \subset \CK$. Express 
\[f = \dirlim  {n \in \NN} (f_{n}:X_{n} \to S_{n}) \text{ and } 
g = \dirlim  {n \in \NN} (g_{n}:Y_{n} \to S_{n}).\] Then the cartesian diagram 
\begin{diagram}[height=2em,w=2em,labelstyle=\scriptstyle] 
\FX \times_{\FS} \FY &	\rTto &	\FY\\
\dTto		&			&	\dTto\\
\FX		&			\rTto	& \FS\\
\end{diagram}
is the direct limit of  the cartesian  diagrams 
\begin{diagram}[height=2em,w=2em,labelstyle=\scriptstyle]
X_{n} \times_{S_{n}} Y_{n} &	\rTto	&	Y_{n}\\
\dTto			&				&	\dTto^{g_{n}}\\
X_{n}			&			\rTto^{f_{n}}	&	S_{n}\\
\end{diagram}
where $n \in \NN$ and, therefore, 
\[\FX \times_{\FS} \FY	=  \dirlim {n \in \NN} X_{n} \times_{S_{n}} Y_{n}.\]
\end{parraf}

\begin{defn} \label{afdis2}
Let $\FX$ be in $\sfn$. 
\begin{enumerate}
\item \label{af2} 
We call \emph{affine formal $r$-space over $\FX$} or \emph{affine formal space of  dimension $r$ over $\FX$} to  $\BA_{\FX}^{r} = \FX \times_{\spec(\ZZ)} \BA_{\spec(\ZZ)}^{r}$
\begin{diagram}[height=2em,w=2em,labelstyle=\scriptstyle]
\BA_{\FX}^{r} &	\rTto^{p}&	\FX\\
\dTto		&			&	\dTto\\
\BA_{\spec(\ZZ)}^{r}	&			\rTto&	\spec(\ZZ)\\
\end{diagram}
Given $\FU=\spf(A) \subset \FX$  an affine open set, then 
\[\BA_{\FX}^{r}|_{\FU} = \BA_{\spf(A)}^{r} \underset{\textrm{\ref{afdis}.(\ref{af})}} = \spf(A\{\mathbf{T}\})\]
where $\mathbf{T}$ denotes the indeterminates $T_{1},T_{2},\ldots,T_{r}$ and it follows that $\BA_{\FX}^{r}$ is in $\sfn$.
\item \label{dis2}
We call \emph{formal $r$-disc over $\FX$} or \emph{formal disc of  dimension $r$ over $\FX$} to $\BD_{\FX}^{r} = \FX \times_{\spec(\ZZ)} \BD_{\spec(\ZZ)}^{r}$
\begin{diagram}[height=2em,w=2em,labelstyle=\scriptstyle]
\BD_{\FX}^{r}	 &	\rTto^{p}&	\FX\\
\dTto		&			&	\dTto\\
\BD_{\spec(\ZZ)}^{r}			&			\rTto&	\spec(\ZZ)
\end{diagram}
If $\FU=\spf(A) \subset \FX$ is an affine open set, we have that \[\BD_{\FX}^{r}|_{\FU} = \BD_{\spf(A)}^{r} \underset{\textrm{\ref{afdis}.(\ref{dis})}} = \spf(A[[\mathbf{T}]])\]
so, $\BD_{\FX}^{r}$ is in $\sfn$.
\end{enumerate}
\end{defn}

\begin{rem}
From now on and, except for exceptions that will be indicated, every formal scheme will be in $\sfn$. We will assume that every ring is noetherian  and, therefore, that every complete ring and every complete module for an adic topology are also separated. 
\end{rem}

\begin{defn}
Let $\FX$ be in $\sfn$, $\CJ \subset \CO_{\FX}$ an Ideal of  definition and $x \in \FX$. We define the  \emph{topological dimension  of   $\FX$ at $x$} as
\[
\dimtop_{x}\FX = \dim_{x}X_{0}
\]
The definition does not depend of the chosen  Ideal of  definition of  $\FX$. As for it, we may assume that  $\FX= \spf(A)$. If $J$ and  $J' $ are ideals of  definition of  $A$, there exist $k,\, l \in \NN$ such that  $J^{k} \subset J'$ and $J'^{l} \subset J$ and, therefore, we have that $\dim A/J = \dim A/J'$. 
We define the \emph{topological dimension  of  $\FX$} as
\[\dimtop\FX = \sup_{x \in \FX} \dimtop_{x} \FX = \sup_{x \in \FX} \dim_{x}X_{0} =  \dim X_{0}.\]
If $\FX = \spf(A)$ with $A$ an $I$-adic noetherian ring, then $ \dimtop \FX = \dim A / I $.  For example,
\[
\begin{array}{ccccc}
 \dimtop \BA_{\spf(A)}^{r} &= &\dim \BA_{\spec(A/I)}^{r} &= &\dim A / I + r \\
\dimtop \BD_{\spf(A)}^{r}  &= &\dim \spec(A/I)           &= &\dim A / I 
\end{array}
\]
\end{defn}

Given $A$ an $I$-adic ring, put $X = \spec(A)$ and $\FX= \spf(A)$. Despite that the only ``visible part'' of  $\FX$ in $X = \spec(A)$ is $V(I)$, it happens that $X \setminus V(I)$ has a deep effect in the behavior of  $\FX$ as we will see along this work. So, apart from the topological dimension of  $\FX$ it is necessary  to consider a dimension that expresses part of the ``hidden" information: the algebraic dimension.

\begin{defn}
Let  $\FX$ be in $\sfn$ and $\CJ$ an Ideal of  definition of  $\FX$. Given $x \in \FX$ we define the \emph{algebraic dimension of  $\FX$ at $x$} as
\[
\dim_{x}\FX = \dim \CO_{\FX,x} 
\]
The  \emph{algebraic dimension of  $\FX$} is \[\dim\FX = \sup_{x \in \FX} \dim_{x} \FX.\]
\end{defn}

\begin{lem} \label{dimha}
Given $A$ a noetherian ring and $I$ an  ideal of  $A$, denote by $\HA$ the completion of  $A$ for the $I$-adic topology. It holds that \[\dim \HA = \underset{I \subset \fp} \sup \dim A_{\fp}\] and,  therefore,  $\dim \HA \le \dim A$. Furthermore, if $A$ is $I$-adic, $\dim A = \dim \HA$.
\end{lem}

\begin{proof}
It is known that there exists a bijective correspondence between the maximal ideals  of  $A$ that  contain $I$ and the maximal ideals of $\HA$ given by $\fm \to \hm$ \cite[III, \S3.4, Proposition 8]{b}. Then
\[
\dim \HA = \underset{ I \subset \fm \in \Max(A)} \sup \dim \HA_{\hm} = \underset{ I \subset \fm \in \Max(A)} \sup \dim A_{\fm} = \underset{ I \subset \fp \in \spec(A)} \sup \dim A_{\fp} 
\]
since $A_{\fm} \inc \HA_{\hm}$  is a flat extension of  local rings with the same residual fields. The other assertions  follow from this  equality. 
 \end{proof}
 \begin{cor} \label{diaf}
If $\FX = \spf(A)$ with $A$ an $I$-adic noetherian ring then
\[\dim \FX = \dim A\]
\end{cor}
\begin{proof}
For each  $x \in \FX$, if $\fp_{x}$ is the corresponding open prime ideal  in $A$ we have that \[\dim_{x} \FX =  \dim A_{\{\fp_{x}\}} = \dim A_{\fp_{x}}\] since $ A_{\fp_{x}} \inc A_{\{\fp_{x}\}} $ is a flat extension of   local rings with the same residual field. Then,  from Lemma \ref{dimha} we deduce the result.
\end{proof}

\begin{cor} \label{dimha2}
Let $A$ be an $I$-adic noetherian ring and denote by $\mathbf{T}$ the set of  indeterminates $T_{1},T_{2},\ldots,T_{r}$. Then:
\[\dim A\{\mathbf{T}\} = \dim A[[\mathbf{T}]] = \dim A + r\]
\end{cor}

\begin{proof}
By Lemma \ref{dimha} it holds that 
\[\dim A[[\mathbf{T}]] \le \dim A\{\mathbf{T}\} \le \dim A[\mathbf{T}] = \dim A + r\]
So, it suffices to show  that $\dim A[[\mathbf{T}]]  \ge \dim A + r$. We will make induction on $r$. Suppose that  $r=1$. Given $\fp_{0} \subset \fp_{1} \subset  \ldots \subset \fp_{l}$  a chain of  prime ideals  in $A$ we have that $\fp_{0} [[T]] \subset \fp_{1} [[T]] \subset \fp_{l}[[T]] \subset \fp_{l}[[T]] + [[T]]$ is a chain of  prime ideals  in $A[[T]]$ and, therefore, $\dim A\{T\} = \dim A[[T]] = \dim A + 1$. Given  $r \in \NN$, assume that the equality holds for  $i<r$. Since $A[[T_{1},T_{2},\ldots,T_{r}]] = A[[T_{1},T_{2},\ldots,T_{r-1}]] [[T_{r}]]$,  by the induction hypothesis we get the result.
\end{proof}

\begin{ex} \label{exdimafdis}
Let $A$ be an   $I$-adic noetherian ring. Then
\[
\begin{array}{cccc}
 \dim \BA_{\spf(A)}^{r} \underset{\textrm{\ref{diaf}}} = &\dim A\{\mathbf{T}\} \underset{\textrm{\ref{dimha2}}} =& \dim A + r  \underset{\textrm{\ref{diaf}}} = &\dim \spf(A) + r\\
 \dim \BD_{\spf(A)}^{r}  \underset{\textrm{\ref{diaf}}}=& \dim A[[\mathbf{T}]]  \underset{\textrm{\ref{dimha2}}} =& \dim A + r  \underset{\textrm{\ref{diaf}}} = &\dim \spf(A) + r
\end{array}
\]
\end{ex}

From this examples, we see that the algebraic dimension of  a formal scheme does not measure  the  dimension of the underlying topological space. 
In general, for $\FX$ in $\sfn$, $\dim_{x} \FX \ge \dimtop_{x} \FX$, for any $x \in \FX$ and, therefore
\[\dim \FX \ge \dimtop \FX.\] 
Moreover, if $\FX= \spf(A)$ with  $A$ an  $I$-adic  ring then $\dim \FX \ge \dimtop \FX + \htt (I)$.

\section{Properties of  morphisms of  formal schemes} \label{sec12}
In this section we will discuss certain properties of  morphisms of  locally noetherian formal schemes. First, we recall  adic morphisms of  formal schemes (\emph{cf.}  \cite[\S 10.12]{EGA1}), whose behavior  is similar to that of morphisms in $\sch$. In fact, an adic morphism can be interpreted  as a family  of  ``relative schemes" with base a formal scheme (see \ref{fibra} below). In our exposition they just play an auxiliary role. Next, we define the  fiber of  a morphism of  formal schemes $f:\FX \to \FY$ in a point $y \in \FY$ and  the relative algebraic dimension of  $f$ in a point $x \in \FX$. The section ends with the definition  and the study of the pseudo closed immersions, that is an important  class  of  non adic morphisms; we also deal with separated morphisms (\emph{cf.}  \cite[\S 10.15]{EGA1}) and radical morphisms. These two last properties behave as the corresponding ones for ordinary schemes, but they will play a decisive  role in the rest of this work.\\


\begin{parraf}  \label{defadic}
A morphism $f:\FX \to \FY$ in $\sfn$ is \emph{adic} (or simply \emph{$\FX$ is a $\FY$-adic} formal scheme) if there exists an Ideal of  definition $\CK$ of  $\FY$ such that $f^{*}(\CK)\CO_{\FX}$ is an Ideal of  definition of  $\FX$ \cite[(10.12.1)]{EGA1}. If $f$ is adic, for each Ideal of  definition $\CK' \subset \CO_{\FY}$ it holds that $f^{*}(\CK')\CO_{\FX}\subset \CO_{\FX}$ is an Ideal of  definition.

Whenever  $f:\FX=\spf(A) \to \FY= \spf(B)$ is in $\sfna$ we have that $f$ is adic if given $K \subset B$ an ideal of  definition, $KA$ is an ideal of  definition of  $A$.
\end{parraf}

If $f:\FX \to \FY$ is an adic morphism, then the topology of $\CO_{\FY}$ determines the topology of $\CO_{\FX}$.

\begin{ex} 
Let  $\FX = \spf(A)$ be in $\sfna$. Then:
\begin{enumerate}
\item
The canonical morphism  $\BA_{\FX}^{r} \to \FX$ is adic. 
\item
For each $f \in A$ the inclusion $\fD (f) \inc \FX$ is adic.
\item
However, the projection $\BD_{\FX}^{r} \to \FX$ is only adic when $r=0$, in which case it is the identity of  $\FX$.
\end{enumerate}
\end{ex}

\begin{ex}
Given a locally noetherian scheme $X$ and $X'  \subset X$ a closed  subscheme, the morphism of  completion of  $X$ along  $X' $, $\kappa: X_{/X' } \to X$
 is adic only if $X$ and $X' $ have the same underlying topological space hence,  $X_{/X'} = X$ and $\kappa=1_{X}$.
\end{ex}

In  Proposition \ref{ccad} we will recall a useful   characterization  of  the adic morphisms. Previously, we need to establish some notations.
 
\begin{parraf} \label{sia}
Let $\FY$ be in $\sfn$ and $\CK\subset \CO_{\FY}$ an  Ideal of  definition such that 
\[\FY =\dirlim {n \in \NN} Y_{n}.\] An \emph{inductive system of  $\{Y_{n}\}_{n \in \NN}$-schemes} is an inductive system of  locally noetherian schemes $\{X_{n}\}_{n \in \NN}$ together with  a collection of  morphisms of  schemes  $\{ f_{n} : X_{n} \to Y_{n} \}_{n \in \NN}$ making commutative    the   diagrams 
\begin{diagram}[height=2em,w=2em,labelstyle=\scriptstyle]
X_{m}		&	\rTto^{f_{m}} &	Y_{m} &\\
\uTto			&	&		\uTto&\\
X_{n}		&	\rTto^{f_{n}} &	Y_{n}& \qquad (m \ge n \ge 0).\\
\end{diagram}
the system $\{f_{n}:X_{n} \to Y_{n}\}_{n \in \NN}$ is \emph{adic} if the above diagrams are cartesian \cite[(10.12.2)]{EGA1}.
\end{parraf}

\begin{propo} \label{ccad}
With the previous notations, the  functor
\[
\{X_{n} \xto{f_{n}}  Y_{n}\}_{n \in \NN}  \leadsto f=\dirlim {n \in \NN} f_{n}
\]
establish a canonical equivalence between the category of adic inductive systems over $\{Y_{n}\}_{n \in \NN}$ and the category of  adic $\FY$-formal schemes, whose inverse is the  functor that takes  an adic locally noetherian $\FY$-formal scheme $f:\FX \to \FY$ to the  inductive system $\{f_{n}:X_{n} \to Y_{n}\}_{n\in \NN}$ determined by the Ideals of  definition $\CK \subset \CO_{\FY}$ and $ f^*(\CK)\CO_{\FY} \subset \CO_{\FY}$.
\end{propo}

\begin{proof}
\cite[(10.12.3)]{EGA1}
\end{proof}

\begin{propo}  \label{adic}
Let $f:\FX \to \FY$ and $g:\FY \to \FS$ be in $\sfn$. Then: 
\begin{enumerate}
\item \label{adic1}
If $f$ and $g$ are adic, so is the  morphism $g \circ f$.
\item \label{adic2}
If $g \circ f$ and $g$ are adic, then $f$ is adic.
\item \label{adic3}
If $f$ is adic, given  $\FY' \to \FY$ in $\sfn$ such that  $\FX_{\FY'} $ is in $\sfn$ it holds that  $f_{\FY'}:\FX_{\FY'} \to \FY'$ is adic.
\end{enumerate}
\end{propo}

\begin{proof}
Let $\CL$ be an Ideal of  definition of  $\FS$.
Assertion (1) follows from  the hypothesis and from the equality $(g \circ f)^{*}(\CL)\CO_{\FX} = f^{*}( g^{*}(\CL)\CO_{\FY})\CO_{\FX}$. Let us prove (2). As  $g$ is adic then $g^{*}(\CL)\CO_{\FY} \subset \CO_{\FY}$ is an Ideal of  definition and, since $g \circ f$ is adic we have that $f^{*}( g^{*}(\CL)\CO_{\FY})\CO_{\FX}$ is an Ideal of  definition of  $\FX$ and, therefore, $f$ is adic. As for   (3) we can assume that the formal schemes are affine. Suppose that $\FX= \spf(A)$, $\FY=\spf(B)$ and $\FY'=\spf(B')$ with $K \subset B$ and $K' \subset B'$ ideals of  definition such that $KB' \subset K'$. If  we put  $A' = A \tc_{B} B'$ let us show that  $K'A'$ is an ideal of  definition of  $A'$. Since $f$ is adic, $KA$ is an ideal of  definition of  $A$ so, $J'=KA(A \tc_{B}B') + K'(A \tc_{B}B')= KB'(A \tc_{B}B') + K'(A \tc_{B}B') =  K'(A \tc_{B}B')$ is an ideal of  definition of  $A'$. Then,  the assertion follows from \ref{defadic}. 
\end{proof}

\begin{ex}
Given $\FX$ in $\sfn$, $\BA_{\FX}^{r} = \FX \times_{\spec(\ZZ)} \BA_{\spec(\ZZ)}^{r} \to \FX$   is an adic morphism, (see Proposition \ref{adic}.(\ref{adic3})). 
\end{ex}


\begin{defn}
Let $f:\FX \to \FY$ be in $\sfa$ and $y \in \FY$. We define the \emph{fiber of  $f$ at the  point $y$} as the formal scheme
\[f^{-1} (y) = \FX \times_{\FY} \spec(k(y)).\]

For example, if $f:\FX = \spf(B) \to \FY= \spf(A)$ is in $\sfaa$ we have that $f^{-1} (y) = \spf(B \tc_{A} k(y))$.
\end{defn}

\begin{ex}
Let $\FX= \spf(A)$ be  in $\sfna$.
If $p: \BA_{\FX}^{r} \to \FX$ is the canonical projection of the  affine formal $r$-space over $\FX$, for all  $x \in \BA_{\FX}^{r}$ and for all $y=p(x)$ we have that 
\[
p^{-1}(y) = \spf(A\{\mathbf{T}\}\tc_{A}k(y)) = \spec(k(y)[\mathbf{T}]) = \BA_{\spec(k(y))}^{r}.
\]

If $q: \BD_{\FX}^{r}\to\FX$ is  the canonical projection of the formal $r$-disc over $\FX$, given $x \in \BD_{\FX}^{r}$ and $y=q(x)$, there results that 
\[
q^{-1}(y) = \spf(A[[\mathbf{T}]]\tc_{A}k(y)) = \spf(k(y)[[\mathbf{T}]]) = \BD_{\spf(k(y))}^{r}.
\]

\end{ex}
\begin{parraf} \label{fibra}
Let $f:\FX \to \FY$ be in $\sfn$ and let us consider $\CJ \subset \CO_{\FX}$ and $\CK\subset \CO_{\FY}$ Ideals of  definition with $f^{*}(\CK)\CO_{\FX} \subset \CJ$. According to \ref{lim}, the  morphism $f$ respect to the  Ideals  of  definition $\CJ$ and $\CK$  can be written as
\[f = \dirlim  {n \in \NN} (X_{n} \xto{f_{n}} Y_{n})\]
Then, by \ref{cclim} it holds that
\[f^{-1} (y) = \dirlim  {n \in \NN}f_{n}^{-1} (y)\] 
where $f_{n}^{-1} (y) = X_{n} \times_{Y_{n}} \spec(k(y))$, for each $n \in \NN$.

If $f$ is also adic, from  Proposition \ref{adic}.(\ref{adic3})  we deduce that $f^{-1} (y)$ is a (ordinary) scheme and that $f^{-1} (y) =f_{n}^{-1} (y)$, for all $n \in \NN$.
\end{parraf}

\begin{notaci}
Let $f:\FX \to \FY$ be in $\sfn$, $x \in \FX$ and $y= f(x)$ and assume that  $\CJ \subset \CO_{\FX}$ and $\CK\subset \CO_{\FY}$ are Ideals of  definition such that $f^{*}(\CK)\CO_{\FX} \subset \CJ$. From now  and, apart from exceptions that will be indicated, whenever we consider the rings $\CO_{\FX,x}$ and $\CO_{\FY,y}$ we will associate them the  $\CJ \CO_{\FX,x}$ and $ \CK \CO_{\FY,y}$-adic topologies, respectively. And we will denote by $\widehat{\CO_{\FX,x}}$  and $\widehat{\CO_{\FY,y}}$ the completion of  $\CO_{\FX,x}$ and  $\CO_{\FY,y}$ with respect to the  $\CJ \CO_{\FX,x}$ and $ \CK \CO_{\FY,y}$-adic topologies, respectively.
\end{notaci}

\begin{defn}
Let $f:\FX \to \FY$ be in $\sfn$. Given $x \in \FX$ and $y= f(x)$, we define the  \emph{relative algebraic dimension of  $f$ at $x$} as 
\[\dim_{x} f = \dim_{x} f^{-1}(y)\]
If $\CJ \subset \CO_{\FX}$ and $\CK\subset \CO_{\FY}$ are Ideals of  definition such that $f^{*}(\CK)\CO_{\FX} \subset \CJ$, then 
\[\dim_{x} f= \dim \CO_{f^{-1}(y),x} = \dim \CO_{\FX,x} \otimes_{\CO_{\FY,y}} k(y) = \dim \widehat{\CO_{\FX,x} }\otimes_{\widehat{\CO_{\FY,y} }} k(y)\]
since the topology in $\widehat{\CO_{\FX,x} }\otimes_{\widehat{\CO_{\FY,y} }} k(y)$ is the $\CJ \widehat{\CO_{\FX,x}}$-adic and, therefore, $\widehat{\CO_{\FX,x} }\tc _{\widehat{\CO_{\FY,y} }} k(y) = \widehat{\CO_{\FX,x} }\otimes_{\widehat{\CO_{\FY,y} }} k(y)$.

On the other hand, we define the \emph{relative topological dimension of  $f$ in $x \in \FX$} as \[\dimtop_{x} f = \dimtop_{x} f^{-1}(y).\]
If $\CJ \subset \CO_{\FX}$ and $\CK\subset \CO_{\FY}$ are Ideals of  definition such that $f^{*}(\CK)\CO_{\FX} \subset \CJ$, and $f_{0}: X_{0} \to Y_{0}$ is the  induced morphism by these ideals, it holds that:
\[\dimtop_{x} f = \dim_{x}  f_{0}^{-1}(y) =\dim_{x}f_{0}\]
\end{defn}

\begin{parraf} \label{exdimalg}
Given $f:\FX \to \FY$ in $\sfn$ and $x \in \FX$, then $\dim_{x} f  \ge \dimtop_{x} f$. Moreover, if $f$ is an adic morphism the equality holds by \ref{fibra}.  For example: 
\begin{enumerate}
\item \label{exdimalg1}
If $p: \BA_{\FX}^{r} \to \FX$ is the canonical projection of the affine formal  $r$-space over $\FX$, for all $x \in \BA_{\FX}^{r}$ we have that
\[
\dim_{x} p = \dimtop_{x} p = \dim k(y)[\mathbf{T}] = r 
\]
where $y=p(x)$.
In contrast, if  $q: \BD_{\FX}^{r}\to\FX$ is the canonical projection of the formal  $r$-disc over $\FX$,  $x \in \BD_{\FX}^{r}$ and $y=q(x)$ there results that
\[
\dim_{x} q = \dim k(y)[[\mathbf{T}]] \underset{\textrm{\ref{exdimafdis}}} = r > \dimtop_{x} q = \dim k(y) = 0
\]
\item \label{exdimalg2}
If $X$ is a usual noetherian scheme  and $X' $ is a closed  subscheme  of  $X$, recall that the  morphism of  completion of  $X$ along  $X' $, $\kappa: X_{/X' } \to X$  is not adic, in general.  Note however that, for all  $x \in X_{/X' }$, it holds   that
\[\dim_{x} \kappa = \dim k(x) = \dimtop_{x} \kappa = 0\]
\end{enumerate}
\end{parraf}

\begin{defn} \label{defnenccerrado}
Let $\FX$ be a locally noetherian  formal scheme and $\CI \subset \CO_{\FX}$ a coherent Ideal . If we put\footnote{\cite[(\textbf{0}.3.1.5)]{EGA1} If $\CF$ is a sheaf on a topological space $X$, the  support of  $\CF$ is \[\supp(\CF) :=\{x \in X \text{ such that  }\CF_{x} \neq 0\}\]}  $\FX' := \supp(\CO_{\FX}/\CI)$ then $\FX'$ is a closed subset and  $(\FX', (\CO_{\FX}/\CI)|_{\FX'})$ is a locally noetherian formal scheme \cite[(10.14.1)]{EGA1}. We will say that $\FX'$ is the   \emph{closed formal subscheme} of  $\FX$ defined by $\CI$ or, briefly, the  \emph{closed subscheme} of  $\FX$ defined by $\CI$ \cite[(10.14.2)]{EGA1}. 
\end{defn}

\begin{parraf}
Let $\FX$ be in $\sfn$, $\CJ \subset \CO_{\FX}$ an Ideal of  definition and $\FX' \subset \FX$ a  closed subscheme defined by the coherent ideal $\CI \subset \CO_{\FX}$. Then \[\FX' = \dirlim {n \in \NN} X'_{n}\] where $X'_{n} \subset X_{n}=(\FX,\CO_{\FX}/\CJ^{n+1})$ is the  closed subscheme defined by $(\CJ^{n+1} +\CI)/\CJ^{n+1}$, 
for all $n \in \NN$.
\end{parraf}

\begin{parraf}
Analogously  to the case of schemes, given $\FX$ in $\sfn$ there exists a bijective correspondence between coherent Ideals $\CI$ of  $\CO_{\FX}$ and closed subschemes $\FX' \inc \FX$ given by $\FX' = \supp(\CO_{\FX}/\CI)$ and $\CO_{\FX'}= (\CO_{\FX}/\CI)|_{\FX'}$. 

In particular, put $\FX= \spf(A)$ where $A$ is a $J$-adic noetherian ring. Given $I \subset A$ an ideal the ring $A/I$ is  $J (A/I)$-adic and $\FX' = \spf(A/I)$ is a closed subscheme of  $\FX$.
\end{parraf}

\begin{defn}
A morphism $j: \FZ \to \FX$ in $\sfn$ is a \emph{closed immersion} if it factors as $\FZ \xto{f} \FX' \inc \FX$
where $f$ is an isomorphism of  $\FZ$ into a closed subscheme $\FX'$ of  $\FX$ \cite[p. 442]{EGA1}.
\end{defn}

\begin{propo} \label{ec} 
Given $f:\FX \to \FY$ an  \emph{adic} morphism in $\sfn$ the following conditions are equivalent:
\begin{enumerate}
\item
The morphism $f$ is a closed immersion.
\item
Given $\CK \subset \CO_{\FY}$ an  Ideal of  definition and  $\CJ =f^{*}(\CK)\CO_{\FX}$, the induced  morphism $f_{0}:X_{0} \to Y_{0}$ is a closed immersion. 
\item
Given $\CK \subset \CO_{\FY}$ an  Ideal of  definition and  $\CJ =f^{*}(\CK)\CO_{\FX}$ the corresponding Ideal of  definition of  $\FX$, the induced  morphisms $f_{n}:X_{n} \to Y_{n}$, are closed immersions, for all $n \in \NN$. 
\end{enumerate}
\end{propo}
\begin{proof}
The equivalence $(1) \Leftrightarrow (2)$ is \cite[(4.8.10)]{EGA31} and  $(2) \Leftrightarrow (3)$ is immediate.
\end{proof}

\begin{defn}
 A morphism $f:\FX \to \FY$ in $\sfn$  is a  \emph{pseudo  closed immersion} if there exists $\CJ \subset \CO_{\FX}$ and $\CK \subset \CO_{\FY}$ Ideals of  definition such that $f^{*}(\CK)\CO_{\FX} \subset \CJ$ and the morphisms induced  by $f$, $f_{n}:X_{n} \to Y_{n}$, are closed immersions, for all $n \in \NN$. As a consequence, if $f:\FX \to \FY$ is a pseudo  closed immersion, $f(\FX)$ is a closed subset  of  $\FY$.
 
Let us show that this definition does not  depend of  the chosen Ideals of  definition. Being a local question, we can assume that $f: \FX= \spf(A) \to \FY = \spf(B)$ is in $\sfna$ and that $ \CJ = J^{\tr},\, \CK= K^{\tr}$ for ideals of definition $J \subset A$ and $K \subset B$. Then given other ideal of definitions $J' \subset A$ and $K' \subset B$ such that $\CJ' = J'^{\tr}\subset \CO_{\FX},\, \CK' = K'^{\tr}\subset \CO_{\FY}$ verifying that $f^{*}(\CK')\CO_{\FX} \subset \CJ'$, there exists $n_{0} > 0$ such that  $J^{n_{0}} \subset J',\, K^{n_{0}} \subset K'$. The  morphism $B \to A$ induces the following commutative diagrams
\begin{diagram}[height=2em,w=2em,labelstyle=\scriptstyle]
B/K^{n_{0}(n+1)}	 &   \rTto    &    A/J^{n_{0}(n+1)}\\
\dTto            &            &    \dTto\\
B/K'^{n+1}       &	\rTto    &    A/J'^{n+1}\\      
\end{diagram}
and it follows that  $B/K'^{n+1}\to A/J'^{n+1}$ is surjective, for all $n \in \NN$.  Thus, $ (\FX, \CO_{\FX}/\CJ'^{n+1}) \to (\FY,\CO_{\FY}/\CK'^{n+1})$ is a  closed immersion, for all $n \in \NN$.
\end{defn}

\begin{ex}
Given $X$ a noetherian scheme and $X'  \subset X$ a closed subscheme defined by an  Ideal $\CI \subset \CO_{X}$, by \ref{excom} and \ref{lim} it holds that 
\[X_{/X' } \xto{\kappa} X = \dirlim {n \in \NN} \left((X' ,\CO_{X}/\CI^{n+1}) \xto{ \kappa_{n}} (X,\CO_{X})\right).\]
Therefore, the  morphism of  completion of  $X$ along  $X'$ is a pseudo  closed immersion.
\end{ex}

Notice that, an adic pseudo  closed immersion  is a  closed immersion. 
Though for the pseudo closed immersions we do not have an analogous characterization to Proposition \ref{ec}:
\begin{ex}
 Given $K$ a field, consider the canonical projection $p: \BD_{\spec(K)}^{1}  \to \spec(K)$. If we take $[[T]]^{\tr}$ as an Ideal of  definition of  $\BD_{\spec(K)}^{1}$ then $p_{0} = 1_{\spec(K)}$ is a  closed immersion. However,  the morphisms $p_{n}: \spec(K[T] / \langle T \rangle  ^{n+1}) \to \spec(K)$ are not closed immersions, for all $ n > 0$ and, thus, $p$ is not a pseudo  closed immersion.
\end{ex}

\begin{propo} \label{ppec}
Let $f:\FX \to \FY$ and $g:\FY \to \FS$ be two morphisms in $\sfn$. It holds that:
\begin{enumerate}
\item \label{ppec1}
If $f$ and $g$ are (pseudo)  closed immersions then $g \circ f$ is a (pseudo)  closed immersion.

\item \label{ppec2}
If $f$ is a (pseudo)  closed immersion, given $h:\FY' \to \FY$ in $\sfn$  we have that $\FX_{\FY'} = \FX \times_{\FY} \FY'$ is in $\sfn$ and  that $f':\FX_{\FY'} \to \FY'$ is a (pseudo)  closed immersion.

\end{enumerate}
\end{propo}

\begin{proof}
By Proposition \ref{adic} it suffices to show the properties for pseudo closed immersions. 
As for (1) let $\CJ \subset \CO_{\FX}$, $\CK \subset \CO_{\FY}$ and $\CL \subset \CO_{\FS}$ be Ideals of  definition with $f^{*}(\CK)\CO_{\FX} \subset \CJ$, $g^{*}(\CL)\CO_{\FY} \subset \CK$ and consider the corresponding expressions for $f$ and $g$ as direct limit  of   scheme morphisms:
\[
f =\dirlim {n \in \NN} (X_{n} \xto{f_{n}} Y_{n}) \qquad g= \dirlim {n \in \NN} (Y_{n} \xto{g_{n}}S_{n})
\]
From \ref{complim} we have that 
\[g \circ f = \dirlim {n \in \NN}  g_{n} \circ f _{n}\] and then the assertion follows from the stability under composition  of the  closed immersions in $\sch$. 
Let us show (2).  Take $\CK' \subset \CO_{\FY'}$ an Ideal of  definition with $h^{*}(\CK)\CO_{\FY'} \subset \CK'$ and such that, by \ref{lim},  
\[h= \dirlim {n \in \NN} (h_{n}: Y'_{n} \to Y_{n}).\] Then by \ref{cclim} we have that 
\[
\begin{array}{ccc}
\begin{diagram}[height=2.5em,w=2.5em,labelstyle=\scriptstyle]
\FX_{\FY'}	& \rTto^{f'}	& \FY' &   \\
\dTto       &         &\dTto^h   & \\
\FX	        &\rTto^f  &	\FY &	\\
\end{diagram}
& =\dirlim {n \in \NN}&
\left(
\begin{diagram}[height=2.5em,w=2.5em,labelstyle=\scriptstyle]
 X_{n} \times_{Y_{n}} Y'_{n} & \rTto^{f'_{n}}&Y'_{n}\\
  \dTto		&   &\dTto^{h_{n}}\\
  X_{n}	&		\rTto^{f_{n}} & Y_{n}\\
\end{diagram}
\right)
\end{array}
\]
By hypothesis $f_{n}$ is a  closed immersion and since  closed immersions in $\sch$ are stable under base-change we have that $f'_{n}$ is a  closed immersion of  noetherian schemes, $\forall n\in \NN$. Finally, $\FX_{\FY'}$ is in $\sfn$, because this is a particular case of proposition \ref{mtf} (\ref{mtf4}). This result is not used in that proof.
\end{proof}

\begin{defn}
 \cite[(10.15.1)]{EGA1}
Given $f:\FX \to \FY$ in $\sfa$, the  \emph{diagonal morphism} $\Delta_{\FX} := (1_{\FX},1_{\FX})_{\FY}: \FX \to \FX \times_{\FY} \FX$ given by the universal property of the fiber product. The morphism $f$ is \emph{separated} if the subset $\Delta_{\FX} (\FX) \subset \FX \times_{\FY} \FX$ is closed. 
\end{defn}

\begin{parraf} \label{psep}
If $f:\FX \to \FY$ is in $\sfn$ and $\CJ \subset \CO_{\FX}$ and $\CK \subset \CO_{\FY}$ are Ideals of  definition with $f^{*}(\CK)\CO_{\FX} \subset \CJ$  then 
\[
f \text{ is separated } \Leftrightarrow X_{0} \xto{f_{0}} Y_{0} \,\text{ is separated \cite[(10.15.2)]{EGA1}.}
\]
By \cite[(5.3.1)]{EGA1}, \ref{complim} and \ref{cclim} we deduce that separated morphisms in $\sfn$ are stable under composition and base-change. Moreover, every (pseudo)  closed immersion is separated.

Also, in \cite[(10.15.5)]{EGA1} it is shown that, if $\FX \times_{\FY} \FX$ is in $\sfn$, then $f:\FX \to \FY$ is separated if, and only if, $\Delta_{\FX}$ is a  closed immersion. 
\end{parraf}


\begin{defn} \label{defencab}
Given $(\FY,\CO_{\FY})$ in $\sfn$ and $\FU \subset \FY$ open, it holds that $(\FU,\CO_{\FX}|_{\FU})$ is a noetherian formal scheme \cite[(10.4.4)]{EGA1} and we say that $\FU$ is a \emph{open subscheme  of  $\FY$}. A morphism $f:\FX \to \FY$  is an \emph{open immersion} if there exists $\FU\subset \FY$ open such that  $f$ factors as
\[
\FX \xto{g} \FU \inc \FY
\]
where $g$ is a isomorphism. 
\end{defn}

\begin{rem}
It is well known that given $X \xto{f} Y \xto{g} S$ in $\sch$ with $f$ a  closed immersion and $g$ an open immersion the  morphism $g \circ f $ factors as $f' \circ g'$ with  $f'$ a  closed immersion and $g'$ an open immersion. However, in $\sfn$ the  analogous result is not true (see a counterexample at the beginning of \cite{AJL2}).
\end{rem}

\begin{defn}\label{rad}
A morphism $f:\FX \to \FY$ in $\sfn$ is \emph{radical} if given $\CJ \subset \CO_{\FX}$ and $\CK \subset \CO_{\FY}$ Ideals of  definition such that  $f^{*}(\CK)\CO_{\FX} \subset \CJ$ the  induced morphism of  schemes $f_{0}:X_{0} \to Y_{0}$ is radical\footnote{\cite[(3.7.2)]{EGA1} A morphism of  schemes $g: X \to Y$ is radical if satisfies the following equivalent  conditions: 1) it is universally injective 2) it is injective and for all $x \in X$, the field extension $k(x)|k(g(x))$ is radical.}. 

Given $x \in \FX$, the residue field of the local rings $\CO_{\FX,x}$ and  $\CO_{X_{0},x}$ are the same and analogously for $\CO_{\FY,f(x)}$ and  $\CO_{Y_{0},f_0(x)}$. Therefore the definition of radical morphisms does not depend on the chosen Ideals of definition of $\FX$ and $\FY$.
\end{defn}

\begin{parraf} \label{prad}
From \cite[(3.7.3) and (3.7.6)]{EGA1} it follows that:
\begin{enumerate}
\item \label{prad1}
The radical morphisms are stable under composition  and noetherian base-change.
\item \label{prad2}
Every monomorphism is a radical morphism. Concretely, open immersions,   closed immersions and pseudo  closed immersions are radical morphisms. 
\end{enumerate}
\end{parraf}

\section{Finiteness conditions of  morphisms  of  formal schemes} \label{sec13}
This section deals with finiteness conditions for morphisms in $\sfn$, that generalize the analogous properties  in $\sch$. In  the class of  adic morphisms we recall the notions of finite type morphisms  and finite morphisms, already defined in \cite[\S 10.13]{EGA1} and \cite[\S 4.8]{EGA31}. In the wide class of  non adic morphisms we will study morphisms of pseudo finite type and pseudo finite morphisms (introduced in \cite[p. 7]{AJL1}\footnote{Morphisms of pseudo finite type have been also introduced independently by Yekutiely in \cite{y} under the name ``formally finite type morphisms''}).
Generalizing  quasi-finite  morphisms in $\sch$, we will define  pseudo quasi-finite  morphisms and  quasi-coverings and establish  its basic  properties.
At the end of  the section  we provide sorites for a property $\CP$ of  morphisms in $\sfn$ (Proposition \ref{sori}).

\begin{defn} \label{defmtf}
A morphism $f:\FX \to \FY$ in $\sfn$ is of \emph{pseudo finite type} (\emph{pseudo finite})   if there exist $\CJ \subset \CO_{\FX}$ and  $\CK\subset \CO_{\FY}$ Ideals of  definition with  $f^{*}(\CK)\CO_{\FX} \subset \CJ$ and such that the  induced morphism of  schemes, $f_{0}: X_{0}  \to Y_{0}$ is of  finite type (finite, respectively). If $f$ is of pseudo finite type (pseudo finite) and adic we say that $f$ is of \emph{finite type} (\emph{finite}, respectively). 

Observe that if $f$ is in $\sch$ the above definitions coincide with  that of    finite type morphism and  finite morphism.
\end{defn}

\begin{propo} \label{ptf}
Let $f:\FX \to \FY$ be in $\sfn$.
\begin{enumerate}
\item \label{ptf1}
The morphism $f$ is of pseudo finite type if, and only if,  for each  $y \in \FY$, there exist  affine open subsets $\FV \subset \FY$ and $\FU \subset \FX$ with $x \in \FU$ and $f(\FU) \subset \FV$ such that $f|_{\FU}$ factors as
\[
\FU \overset{j} \to  \BD_{\BA_{\FV}^{r}}^{s} \xto{p} \FV
\]
where $r,\, s \in \NN$, $j$ is a  closed immersion and $p$ is  the canonical projection. 
\item \label{ptf2}
The morphism $f$ is of finite type if, and only if,  for each $y \in \FY$, there exist affine open subsets $\FV \subset \FY$ and $\FU \subset \FX$ with $x \in \FX$ and $f(\FU) \subset \FV$ such that $f|_{\FU}$ factors as
\[
\FU \overset{j} \to  \BA_{\FV}^{r} \xto{p} \FV
\]
where $r \in \NN$, $j$ is a  closed immersion and $p$ is the canonical projection. 

\item \label{ptf4}
The morphism $f$ is pseudo finite if, and only if,  for each $y \in \FY$, there exists an affine open  $\FV \subset \FY$ with $f^{-1}(\FV) =\FU \subset \FX$ affine  and there exists $s \in \NN$ such that  $f|_{\FU}$ factors as
\[
\FU \xto{f'}  \BD_{\FV}^{s} \xto{p} \FV
\]
where $f'$ is finite and $p$ is the canonical projection (\emph{cf.} \cite[p. 15]{AJL1}).
\item \label{ptf3}
The morphism $f$ is finite if, and only if,  for each $y \in \FY$, there exists an affine open subset  $\FV \subset \FY$ with $f^{-1}(\FV) =\FU$ an affine open subset of $\FX$ such that  $\ga(\FU,\CO_{\FX})$ is a $\ga(\FV,\CO_{\FY})$-module of  finite type.
\end{enumerate}
\end{propo}

\begin{proof}
Since this are local properties  we may assume $f: \FX = \spf(A) \to \FY = \spf(B)$ is in $\sfna$. Given $J \subset A$ and $K \subset B$ ideals of  definition such that $KA \subset J$ let $f_{0}: X_{0}=\spec(A/J) \to Y_{0}= \spec(B/K)$ be the induced morphism induced by $f$.

Let us prove  property (1). As $f$ is  pseudo finite type, there exists a presentation 
\[\frac{B}{K} \inc \frac{B}{K} [T_{1},T_{2}, \ldots ,T_{r}] \overset{\varphi_{0}}\epi \frac{A}{J}.\] 
This  morphism lifts to a ring homomorphism 
\[B \inc  B [T_{1},T_{2}, \ldots ,T_{r}] \to A\]
that extends to a continuous morphism 
\begin{equation} \label{factotipofinit}
 B \inc  B\{\mathbf{T}\}[[\mathbf{Z}]]:=B \{T_{1},T_{2}, \ldots ,T_{r}\}[[Z_{1},Z_{2}, \ldots ,Z_{s}]] \xto{\varphi} A
\end{equation}
 such that the images of   $Z_{i}$ in $A$ generate $J$.  It is easily seen that the  morphism of graduated modules associated to $\varphi$  
 \[
 \bigoplus_{n \in \NN} \frac{(K\{\mathbf{T}\}[[\mathbf{Z}]] +[[\mathbf{Z}]])^{n}}{(K\{\mathbf{T}\}[[\mathbf{Z}]] +[[\mathbf{Z}]])^{n+1}} \xto{\mathrm{gr}(\varphi)}  \bigoplus_{n \in \NN} \frac{J^{n}}{J^{n+1}}
\]
is surjective and, therefore, $\varphi$  is also surjective (\cite[III, \S2.8, Corollary 2]{b}).

As for (3), by hypothesis,  there exists an epimorphism of   $B/K$-modules
\[\underset{t}\bigoplus \frac{B}{K} \overset{\varphi_{0}} \epi \frac{A}{J}\]
for some $t \in \NN$.
This epimorphism  lifts to a continuous morphism of  $B$-modules
\[\underset{t}\bigoplus B \to A\]  that extends to a morphism of topological $B$-modules 
\begin{equation} \label{factmorfinito}
\underset{t}\bigoplus B[[Z_{1},Z_{2}, \ldots, Z_{s}]]  \overset{\varphi}\to A
\end{equation}
 such that the images of  $Z_{i}$ in $A$ generate $J$. From the fact  that $\varphi_{0}$ is surjective we deduce  that $\varphi$ also is (\cite[III, \S2.11, Proposition 14]{b}).
 
In the cases (2) and (4) $KA = J$, so we can choose $s=0$. Then, if $f$ is of finite type the factorization (\ref{factotipofinit}) can be written \[B \to  B \{T_{1},T_{2}, \ldots ,T_{r}\} \epi  A\] and corresponds with the one given in \cite[(10.13.1)]{EGA1}. 
Whenever $f$ is  finite the factorization (\ref{factmorfinito}) is expressed
\[\underset{r}\bigoplus B \epi  A\] and corresponds with the one given in \cite[(4.8.1)]{EGA31}.
\end{proof}

\begin{parraf}\label{defmtf1}
Let $f:\FX \to \FY$ be a pseudo finite type (or pseudo finite) morphism in $\sfn$. As a consequence  of  the last Proposition it holds that \emph{for all couple of  Ideals of  definition} $\CJ \subset \CO_{\FX}$ and $\CK \subset \CO_{\FY}$ such that $f^{*}(\CK) \CO_{\FX} \subset \CJ$, the  induced morphism of  schemes $f_{0}:X_{0} \to Y_{0}$ is of finite type (or finite, respectively).
\end{parraf}

\begin{rem}
In the assertion c) of  the characterization of  finite morphisms  of  formal schemes given in  \cite[(4.8.1)]{EGA31} it should be added   the adic hypothesis for  $f$.
\end{rem}

\begin{propo}  \label{mtf}
Given $f:\FX \to \FY$ and $g:\FY \to \FS$ in $\sfn$ we have the following properties: 
\begin{enumerate}
\item  \label{mtf1}
If $f$ and $g$ are (pseudo) finite type morphisms, then  $g \circ f$ is  a   (pseudo) finite type morphism.
\item  \label{mtf2}
If $f$ and $g$ are  (pseudo)  finite morphisms, the  morphism $g \circ f$  is  (pseudo) finite.
\item   \label{mtf3}
If $f:\FX \to \FY$ is a  (pseudo) finite type morphism, given $h:\FY' \to \FY$ a morphism  in $\sfn$ we have that $\FX \times_{\FY} \FY'$ is in $\sfn$ and that $f':\FX_{\FY'} \to \FY'$ is of  (pseudo) finite type.
\item  \label{mtf4}
If $f:\FX \to \FY$ is a  (pseudo)  finite morphism, given $h:\FY' \to \FY$ a morphism  in $\sfn$ then $f':\FX_{\FY'} \to \FY'$ is (pseudo) finite.
\end{enumerate}
\end{propo}

\begin{proof}
By Proposition \ref{adic} it suffices to prove  the assertions for pseudo finite and of pseudo finite type morphisms. First, (1) and (2) are deduced from  \ref{complim} and the corresponding sorites in $\sch$. 

In order to prove (3) and (4), by \ref{cclim} and by the analogous properties in $\sch$, is suffices to show that if $f$ is a pseudo finite type morphism then the formal scheme $\FX \times_{\FY} \FY'$ is locally noetherian. We may suppose that $\FX= \spf(A),\, \FY = \spf(B)$ and $\FY'= \spf(B')$ with  $J \subset A,\, K \subset B$ and $K' \subset B'$ ideals of  definition such that $KA \subset J$ and $KB' \subset K'$. Let us check that $A'= A \tc_{B} B'$ is noetherian. By hypothesis, we have that $B/K \to A/J$ is of finite type. Therefore, by base-change, there results that $B'/K' \to  A/J \otimes_{B/K} B'/K'$ is of finite type and, since 
\[A' =\invlim {n \in \NN} (A/J^{n+1} \otimes_{B/K^{n+1}} B'/K'^{n+1}),\] and analogously to the argument given in the proof  of  Property \ref{ptf}.(\ref{ptf1}), we have  a surjective morphism $B'\{T_{1},T_{2},\ldots,T_{r}]\}[[Z_{1},Z_{2}, \ldots ,Z_{s}]]  \epi A'$. Then $A'$ is a noetherian ring and.
\end{proof}

\begin{defn} \label{defncuasifin}
Let $f:\FX \to \FY$ a pseudo finite type morphism in $\sfn$. We say that $f$ is \emph{pseudo quasi-finite} if there exist $\CJ \subset \CO_{\FX}$ and $\CK \subset \CO_{\FY}$ Ideals of  definition with $f^{*}(\CK)\CO_{\FX} \subset \CJ$ and such that $f_{0}$ is quasi-finite\footnote{\cite[(6.11.3)]{EGA1} A  finite type morphism $f: X \to Y$ of  schemes is quasi-finite if it satisfies the following equivalent  conditions: (1) Every point $x \in X$ is isolated in $f^{-1}(f(x))$. (2) For all  $x \in X$, the  scheme $f^{-1}(f(x))$ is  $k(x)$-finite. (3) For all $x \in X$, $\CO_{X,x} \otimes_{\CO_{Y,f(x)}} k(f(x))$ is a $k(f(x))$-module of finite type.}. And $f$ is \emph{ pseudo quasi-finite in $x \in \FX$} if there exists an open $\FU \subset \FX$  with $x \in \FU$ such that  $f|_{\FU}$ is pseudo quasi-finite.

Notice that if $f:\FX \to \FY$ is a pseudo quasi-finite morphism (in $\sfn$) then,  \emph{for all couples of  Ideals of  definition} $\CJ \subset \CO_{\FX}$ and $\CK \subset \CO_{\FY}$ such that  $f^{*}(\CK)\CO_{\FX} \subset \CJ$, the  induced morphism of  schemes $f_{0}:X_{0} \to Y_{0}$ is quasi-finite.
\end{defn}

\begin{propo} \label{soritpcf} 
The pseudo quasi-finite morphisms  satisfy the following properties: 
\begin{enumerate}
\item \label{soritpcf1}
Closed immersions, pseudo  closed immersions and open immersions are pseudo quasi-finite.
\item \label{soritpcf2}
If $f:\FX \to \FY$ and $g:\FY \to \FS$ are pseudo quasi-finite morphisms, then $g \circ f$ also is.
\item  \label{soritpcf3}
If $f:\FX \to \FY$ is pseudo quasi-finite, given $h:\FY' \to \FY$ a morphism  in $\sfn$ we have that $f':\FX_{\FY'} \to \FY'$ is pseudo quasi-finite.
\end{enumerate}
\end{propo}

\begin{proof}
The proof  is an immediate consequence  of  the analogous properties in $\sch$.
\end{proof}
In $\sch$ it is the case that a morphism is \'etale if, and only if, it is smooth and quasi-finite. Nevertheless, we will show that in $\sfn$ not  every smooth and pseudo quasi-finite morphism is \'etale (see Example \ref{pcf+plnope} in Chapter \ref{cap3}). That is why we introduce a stronger notion  than   pseudo quasi-finite morphism and  that also generalizes  quasi-finite morphisms  in $\sch$: the quasi-coverings.

\begin{defn} \label{defncuasireves}
Let $f:\FX \to \FY$ be a pseudo finite type morphism   in $\sfn$. The  morphism $f$ is a\emph{ quasi-covering} if  $\CO_{\FX,x} \tc_{\CO_{\FY,f(x)}} k(f(x))$ is a finite type $k(f(x))$-module, for all $x \in \FX$. We say that $f$ is a \emph{quasi-covering in $x \in \FX$} if there exists an open $\FU \subset \FX$ with $x \in \FU$ such that  $f|_{\FU}$ is a  quasi-covering.

Observe that given $\CJ \subset \CO_{\FX}$ and $\CK \subset \CO_{\FY}$ Ideals of  definition such that $f^{*}(\CJ) \CO_{\FX}\subset \CK$, for all $x \in \FX$ there results that 
\[\CO_{\FX,x} \tc_{\CO_{\FY,f(x)}} k(f(x)) = \invlim {n \in \NN} \CO_{X_{n},x} \otimes_{\CO_{Y_{n},f(x)}} k(f(x)).\]

We reserve the word ``covering'' for a dominant (\ie with dense image) quasi-covering. These kind of maps will play almost no role in this work but they are important, for instance, in the study of finite group actions on formal schemes.
\end{defn}

\begin{ex}
If $X$ is a locally noetherian scheme and $X' \subset X $ is a closed subscheme  the  morphism of  completion $\kappa: \FX=X_{/X'} \to X$ is a quasi-covering. In fact, for all $x \in \FX$ we have that 
\[\CO_{\FX,x} \tc_{\CO_{X,\kappa(x)}} k(\kappa(x)) = k(\kappa(x)).\]
\end{ex}

\begin{propo} \label{soritcr}
The quasi-coverings verify the following properties: 
\begin{enumerate}
\item  \label{soritcr1}
Closed immersions,  pseudo  closed immersions and open immersions are  quasi-coverings.
\item  \label{soritcr2}
If $f:\FX \to \FY$ and $g:\FY \to \FS$  are quasi-coverings, the  morphism $g \circ f$ is a quasi-covering.
\item   \label{soritcr3}
If $f:\FX \to \FY$ is a quasi-covering, given $h:\FY' \to \FY$ a morphism  in $\sfn$ it holds that $f':\FX_{\FY'} \to \FY'$ is a quasi-covering.
\end{enumerate}
\end{propo}
\begin{proof}
Immediate.
\end{proof}
\begin{propo} \label{cuairevdim0}
If $f:\FX \to \FY$ is a quasi-covering  in $x \in \FX$ then:
\[
\dim_{x} f =0
\]
\end{propo}
\begin{proof}
It is a consequence of the fact that $\CO_{\FX,x} \tc_{\CO_{\FY,f(x)}} k(f(x))$ is an artinian ring.
\end{proof}

\begin{propo}
Let $f:\FX \to \FY$ be in $\sfn$ a pseudo finite type morphism. If $f$ is a quasi-covering, then  is pseudo quasi-finite. Furthermore, if  $f$ is adic the  reciprocal holds.
\end{propo}

\begin{proof}
Suppose that $f$ is a quasi-covering and let $\CJ \subset \CO_{\FX}$ and $\CK\subset \CO_{\FY}$ be Ideals of  definition such that $f^{*}(\CK)\CO_{\FX} \subset \CJ$. For $x \in \FX$ and for all $y=f(x)$, $\CO_{\FX,x} \tc_{\CO_{\FY,y}} k(y)$ is a finite $k(y)$-module  and, therefore, 
\[\frac{\CO_{X_{0},x}}{\fm_{Y_{0},y}\CO_{X_{0},x}} = \frac{\CO_{\FX,x}}{ \CJ \CO_{\FX,x}} \otimes_{\CO_{Y_{0},y}} k(y)\]
is $k(y)$-finite, so it follows that $f$ is pseudo quasi-finite.
 
If  $f$ is a pseudo quasi-finite  adic  morphism, in \ref{fibra} we have shown that  $f^{-1}(y)=f_{0}^{-1}(y)$ and, then, 
\[\CO_{X_{0},x}/\fm_{Y_{0},y}\CO_{X_{0},x} = \CO_{\FX,x} \tc_{\CO_{\FY,y}} k(y)\]
for all $x \in \FX$ with  $y=f(x)$. It follows that   $f$ is a quasi-covering.
\end{proof}
\begin{cor}
Every finite  morphism $f:\FX \to \FY$  in $\sfn$ is a quasi-covering.
\end{cor}
\begin{proof}

Given $\CJ \subset \CO_{\FX}$ and $\CK \subset \CO_{\FY}$ Ideals of  definition such that $f^{*}(\CK)\CO_{\FX} \subset \CJ$ by hypothesis we have that $f_{0}:X_{0} \to Y_{0}$ is  finite therefore, quasi-finite. Since $f$ is adic the  result  is  consequence of the last proposition.
\end{proof}
Nevertheless and, as it is shown in the next example, not every pseudo  finite morphism is a quasi-covering and, therefore, pseudo quasi-finite does not imply quasi-covering.
\begin{ex} 
The canonical projection $p: \BD_{\FX}^{r} \to \FX$  is not  a quasi-covering since $\dim_{x} p \underset{\textrm{\ref{exdimalg}.(\ref{exdimalg1})}} = r>0$, for all  $x \in \FX$. But the  scheme morphism  $p_{0}= 1_{X_{0}}$ is  finite.
\end{ex}

\begin{parraf}
In short, we have  the  following diagram of  strict implications (with the conditions that implies adic morphism in italics):
\[
\begin{array}{ccccc}
\textrm{\emph{closed immersion}}&\Rightarrow &\textrm{\emph{finite}}&\Rightarrow& \textrm{quasi-covering}\\
\Downarrow		      &			   &\Downarrow    &                     &\Downarrow\\
\textrm{pseudo  closed immersion}&\Rightarrow & \textrm{pseudo finite}&\Rightarrow& \textrm{pseudo quasi-finite}\\

\end{array}
\]
\end{parraf}

\begin{cuesab}
In $\sch$ Zariski's Main Theorem says that every quasi-finite  morphism factors as a open immersion followed a  finite morphism. In $\sfn$, is every quasi-covering the composition of  a open immersion and a  pseudo finite morphism? 
\end{cuesab}

In the following proposition we study the basic sorites of  a property $\CP$ of pseudo finite type morphisms   in $\sfn$. 

\begin{propo}\label{sori}
Let $\CP$ a property of pseudo finite type morphisms  in $\sfn$ and consider the following statements:
\begin{enumerate}
\item
Every  closed immersion satisfies $\CP$.
\item
Every  adic morphism $f:\FX \to \FY$ in $\sfn$ such that, for all  $\CK \subset \CO_{\FY}$,  that is an ideals of definition of $\FY$, and $\CJ = f^{*}(\CK) \CO_{\FX} \subset \CO_{\FX}$, the corresponding Ideal of definition of $\FX$, the induced morphism $f_{0}:X_{0} \to Y_{0}$ is an immersion, satisfies $\CP$.
\item
The property $\CP$ is stable under composition of  morphisms and for base-change in $\sfn$.
\item
The property $\CP$ is stable under   product of  morphisms.
\item
Given $f:\FX \to \FY$ a pseudo finite type morphism in $\sfn$  and $g:\FY \to \FS$ a separated pseudo finite type morphism   in $\sfn$, if $g \circ f$ satisfies $\CP$, then so does $f$.
\item
Let $f:\FX \to \FY$ and  $g:\FY \to \FS$ be  pseudo finite type morphisms. If $g \circ f$ satisfies $\CP$, then so does $f$.
\item
Let $f:\FX \to \FY$ in $\sfn$ and consider $\CJ \subset \CO_{\FX}$ and $\CK \subset \CO_{\FY}$ Ideals of  definition with $f^{*}(\CK)\CO_{\FX} \subset \CJ$ such that  \[f = \dirlim {n \in \NN } (f_{n}: X_{n} \to Y_{n}).\] If $f$ satisfies the property $\CP$, then $f_{n}$ satisfies $\CP$ for all $n \in \NN$.
\end{enumerate}
Then, we have  the following implications:
\begin{align*}
\text{\emph{(3)}}      &\Rightarrow \text{\emph{(4)}}\\
\text{\emph{(1), (3)}} &\Rightarrow \text{\emph{(5), (7)}}\\ 
\text{\emph{(2), (3)}} &\Rightarrow \text{\emph{(6), (7)}}
\end{align*}
\end{propo}

\begin{proof}
The proof of  (3) $\Rightarrow$ (4)  is the same as in the  case of  schemes (\emph{cf.}  \cite[(5.5.12)]{GD}). The implication  (1), (3) $\Rightarrow$ (5) is analogous to that of  schemes using \ref{mtf}.(\ref{mtf1}),  \ref{mtf}.(\ref{mtf3}) and \cite[(10.15.4)]{EGA1}. In order to prove  (1), (3) $\Rightarrow$ (7) consider the following commutative diagrams:
\begin{diagram}[height=2em,w=2em,labelstyle=\scriptstyle]
\FX            &	\rTto^f       & \FY\\
\uTto^{ i_{n}} &	              & \uTto^{ j_{n}}\\
X_{n}          &	\rTto^{f_{n}} &	Y_{n}\\	
\end{diagram}
where $i_{n},\,  j_{n}$ are the  canonical closed immersions, $\forall n \ge 0$. From  (1) and (3) we deduce that   $f \circ i_{n} = j_{n} \circ f_{n} $ satisfies $\CP$ and then, since every  closed immersion is a separated morphism, the  result is consequence of  (5). 

On the other hand, (2) $\Rightarrow$ (1) so, by the aforementioned property we have that (2), (3) $\Rightarrow$ (4), (7). Finally, let us show  (2), (3) $\Rightarrow$  (6). Let $f:\FX \to \FY$ and $g:\FY \to \FS$ be pseudo finite type morphisms. Let $\ga_{f}$ be the graph of $f$
\begin{diagram}[height=2.5em,w=2.5em,labelstyle=\scriptstyle]
\FX 		&			&		&		&\\
		&\rdTto^{\ga_{f}}	\rdTto(4,2)^f  \rdTto(2,4)_{1_{\FX}}	&    &   &   \\		&			&\FX \times_{\FS} \FY	&\rTto^h	&	\FY\\
		&				&	\dTto &	& \dTto^{g}\\
		&		&	\FX	&\rTto^{g \circ f}&	\FS.\\
\end{diagram}
Assume that $\CJ \subset \CO_{\FX},\, \CK \subset \CO_{\FY}$ and $\CL \subset \CO_{\FS}$ are Ideals of  definition such that $f^{*}(\CK)\CO_{\FX} \subset \CJ,\, g^{*}(\CL)\CO_{\FY} \subset \CK$ and, with respect to them,  put   \[f = \dirlim {i \in \NN } f_{n}: X_{n} \to Y_{n} \text{, } g = \dirlim {i \in \NN } g_{n}: Y_{n} \to S_{n}.\] 
If for all $n \in \NN$, $\ga_{f_{n}} $ denotes the graph of $f_n$, it holds that
\[\ga_{f}= \dirlim {n \in \NN} \ga_{f_{n}}\] 
The  morphism $\ga_{f}$ is adic and, since $\ga_{f_{n}}$ is an immersion, for all $n \in \NN$ (\emph{cf. }\cite[(5.1.4)]{EGA1}) using (2) we have that $\ga_{f}$ satisfies $\CP$.
However, since $g \circ f$ possesses property $\CP$ then by (3) $h$ also has $\CP$  and, again by (3) there results that $f = h \circ \ga_{f}$ satisfies $\CP$.
\end{proof}

\begin{cor}\label{cor1sor}
With  the hypothesis of  Proposition \ref{sori},
if $\CP$ is one of the properties of being \emph{(pseudo)  closed immersion or (pseudo) finite} then (1) and (3) hold for so $\CP$, so that (4), (5) and (7) hold for $\CP$ also.
\end{cor}

\begin{proof}
Every (pseudo) closed immersion is (pseudo) finite trivially and, therefore, $\CP$ satisfies (1). Moreover, in  \ref{ppec}.(\ref{ppec1}) we have shown that (pseudo)  closed immersions are stable under composition and in \ref{mtf}.(\ref{mtf2})  we proved that  (pseudo) finite  morphisms are stable under composition, so $\mathcal{P}$ verifies (3).
\end{proof}

\begin{cor} \label{cor2sor}
With the hypothesis of  Proposition \ref{sori},
if $\CP$ is one of the properties of being \emph{separated, radical, (pseudo) finite type, pseudo quasi-finite or quasi-covering}, then (2) and (3) hold for $\CP$, and, therefore, (4), (6) and (7) hold for $\CP$, too.
\end{cor}

\begin{proof}
If $f$ is a morphism with property $\mathcal{P}$ it is easily seen that satisfies (2). Besides  morphisms with property $\mathcal{P}$ are stable under composition (\emph{cf.}  \ref{psep}, \ref{prad}.(\ref{prad1}),  \ref{mtf}.(\ref{mtf1}),  \ref{soritpcf} and \ref{soritcr}).
\end{proof}

\section{Flat morphisms and completion morphisms} \label{sec14}

In  the first part  of  this section we discuss  flat morphisms in $\sfn$. Whenever a  morphism \[f = \dirlim {n \in \NN} f_{n}\] is adic, the  local criterion of flatness for  formal schemes (Proposition \ref{clp}) relates the flat character  of  $f$ and that of  the  morphisms $f_{n}$, for all $n \in \NN$. Though, in absence of  the adic hypothesis this relation does not hold (Example \ref{excompl1}). In the second part we study  the morphisms of  completion in $\sfn$, a  class of flat morphisms that are pseudo  closed immersions (so, they are closed immersions as topological maps) and that play an essential  role  in  the main theorems of Chapter  \ref{cap3} (Theorem \ref{tppalnr}, Theorem \ref{tppalet} and Theorem \ref{tppall}).

\begin{defn}
A morphism $f:\FX \to \FY$ is \emph{flat at $x \in \FX$} if $\CO_{\FX,x}$ is a flat $\CO_{\FY,f(x)}$-module. We say that \emph{$f$ is flat} if it is flat at $x$, for all  $x \in \FX$.
\end{defn}

\begin{parraf}\label{caracterizlocalplanos}
Given $\CJ \subset \CO_{\FX}$ and $\CK\subset \CO_{\FY}$ Ideals of  definition with $f^{*}(\CK)\CO_{\FX} \subset \CJ$, if $\widehat{\CO_{\FX,x}}$  and $\widehat{\CO_{\FY,f(x)}}$ are the completions of  $\CO_{\FX,x}$ and  $\CO_{\FY,f(x)}$ for  the $\CJ \CO_{\FX,x}$ and $ \CK \CO_{\FY,f(x)}$-adic topologies, by \cite[III, \S5.4, Proposition 4]{b} the following are equivalent:
\begin{enumerate}
\item
$f$ is flat at $x \in \FX$
\item
$\CO_{\FX,x}$ is a flat $\CO_{\FY,f(x)}$-module
\item
$\widehat{\CO_{\FX,x}}$ is a flat $\CO_{\FY,f(x)}$-module
\item
$\widehat{\CO_{\FX,x}}$ is a flat $\widehat{\CO_{\FY,f(x)}}$-module
\end{enumerate} 
\end{parraf}

\begin{lem} \cite[7.1.1]{AJL1}\label{pl}
A morphism $f:\spf(A) \to \spf(B)$ is flat if, and only if,  $A$ is a flat $B$-module. 
\end{lem}
\begin{proof}
If $f$ is flat for all $\fm \subset A$ maximal ideal, and $\fn=f^{-1}(\fm)$ holds that $A_{\{\fm\}}$ is a flat $B_{\{\fn\}}$-module and by \cite[III, \S5.4, Proposition 4]{b} we have that $A_{\fm}$ is a flat $B_{\fn}$-module. Then, $A$ is a flat $B$-module. Reciprocally, if $A$ is a flat $B$-module, then $A_{\{\fp\}}$ is a flat $B_{\{\fq\}}$-module for all  $\fp \subset A$ prime ideal and $\fq=f^{-1}(\fp)$ and, in particular for all  $\fp \subset A$ open prime ideal.
\end{proof}

\begin{parraf} \label{propiplanos}
Since flatness is a local property, the last  lemma implies that:
\begin{enumerate}
\item  \label{propiplanos1}
Open immersions are flat.
\item  \label{propiplanos2}
Flatness is stable under composition of  morphisms, base-change  and    product of  morphisms.
\item  \label{propiplanos3}
If $X$ is a  (usual) noetherian scheme and $X' \subset X$ a closed subscheme  the  morphism of  completion $\kappa:X_{/X'} \to X'$ is flat.
 \end{enumerate}
 
The property ``to be flat"  in general does not satisfy  the  statement (7) of   Proposition \ref{sori}.
\end{parraf}

\begin{ex} \label{excompl1}
If $K$ is a field and $ \BA_{K}^{1}= \spec(K[T])$ consider the  closed $X' =V( \langle T \rangle) \subset  \BA_{K}^{1}$. The  canonical morphism  of  completion of   $\BA_{K}^{1}$ along  $X'$
\[\BD_{K}^{1}  \xto{\kappa} \BA_{K}^{1}\] 
is flat but,  the morphisms 
\[\spec(K[T]/ \langle T \rangle ^{n+1}) \xto{\kappa_{n}} \BA_{K}^{1}\]
are not flat, for every $n \in \NN$.
\end{ex}

\begin{propo}
Let $f: \FX \to Y$ be in $\sfn$ with $Y$ a usual scheme, $\CJ \subset \CO_{\FX}$ an Ideal of  definition and denote by $f_{n}: X_{n} = (\FX,\CO_{\FX}/\CJ^{n+1})\to Y$ the  morphisms induced by $f$, for all $n \in \NN$. If $f_{n}$ is a flat morphism, for all  $n \in \NN$, then $f$ is a flat morphism.
\end{propo}

\begin{proof} (\emph{cf.} \cite[Theorem 22.1]{ma2})
We may assume the following: $\FX= \spf(A)$,  $\CJ= J^{\tr}$ with $A$ a  $J$-adic noetherian ring and $Y= \spec(B)$. Let us show that $A$ is a flat $B$-algebra. For that, take $K \subset B$ be an ideal and let us check that the canonical map $u:K':=K \otimes_{B} A\to A$ is injective. Since $K'$ is a finitely generated $A$-module, is complete for the $J$-adic  topology  and, therefore $\cap_{n \in \NN} J^{n+1} K' = 0$. Take $x \in \ker(u)$. By hypothesis, for all  $n \in \NN$, $A/J^{n+1}$ is a flat $B$-algebra  and, therefore  the morphisms
\[\frac{K'}{J^{n+1}K'}=K'\otimes_{A} \frac{A}{J^{n+1}}= K\otimes_{B} \frac{A}{J^{n+1}} \xto{u \otimes \frac{A}{J^{n+1}}} \frac{A}{J^{n+1}}\]
are injective, for all $n \in \NN$. We deduce that $x \in  J^{n+1}K'$, for all $n \in \NN$, so $x=0$.
\end{proof}

Example \ref{excompl1} shows that the reciprocal of the last proposition does not hold.

\begin{propo}[Local flatness criterion for formal schemes] \label{clp}
Let \\ $f:\FX \to \FY$ be an  \emph{adic} morphism, $\CK$ an Ideal of  definition of  $\FY$, consider $\CJ =  f^{*}(\CK)\CO_{\FX}$ and let $f_{n}:X_{n} \to Y_{n}$ be the morphisms  induced  by $f$ and $\CK$, for all $n \in \NN$. The following assertions are equivalent:
\begin{enumerate}
\item
The morphism $f$ is flat.
\item
The morphism $f_{n}$ is flat, for all $n \in \NN$.
\item
The morphism $f_{0}$ is flat and $\tor_{1}^{\CO_{\FY}} (\CO_{\FX}, \CO_{Y_{0}})=0$.
\item
The morphism $f_{0}$ is flat and  $\CJ =  f^{*}(\CK)$.
\end{enumerate}
Moreover, $f$ is adic and flat if and only if, $f_{0}$ is flat and $\CJ =  f^{*}(\CK)$.
\end{propo}

\begin{proof}
We may suppose that $f: \FX=\spf(A) \to \FY=\spf(B)$ is in $\sfn$. Then if $\CK=K^{\tr}$ for an ideal of definition $K \subset B$, we have that $\CJ= (KA)^{\tr}$  and the proposition is a consequence of Lemma \ref{pl} and of the local flatness criterion for rings (\emph{cf.} \cite[Theorem 22.3]{ma2}).
\end{proof}


In Example \ref{excompl} we have shown how to get  from a (usual) locally noetherian scheme $X$  and  a closed subscheme of  $X' \subset X$  a locally noetherian formal scheme $X_{/X'}$, called completion of  $X$ along  $X'$ and,  a canonical morphism $\kappa: X_{/X'} \to X$. Next, we define the completion of  a formal scheme  $\FX$ in $\sfn$ along a closed subscheme $\FX'\subset \FX$. Despite the construction is natural, it has not been systematically developed in the basic references about formal schemes. Since the morphisms of  completion will  be used later we have opted to include them here. 
 
\begin{defn} \label{defcom}
Let $\FX$ be in $\sfn$ and $\FX' \subset \FX$ a closed formal subscheme defined by an Ideal $\CI$ of  $\CO_{\FX}$. Given $\CJ$ an Ideal of  definition of  $\FX$ we define the \emph{completion of  $\FX$ along  $\FX'$}, and is denoted $\FX_{/ \FX'}$, as the  topological ringed space with topological space  $\FX'$ and whose sheaf  of  topological rings is given by
\[
\CO_{\FX_{/\FX'}} = \invlim {n \in \NN} \frac{\CO_{\FX}}{(\CJ+\CI)^{n+1}}.
\] 

The definition does not depend on the chosen   Ideal of  definition $\CJ$ of  $\FX$. Indeed, if $\CJ' \subset \CO_{\FX}$ is another Ideal of  definition for all affine open set  $\FU \subset \FX$ there there exists $r,\, s \in \NN$ such that $\CJ|_{\FU}^{r} \subset \CJ'|_{\FU}$ and $\CJ'|_{\FU}^{s} \subset \CJ|_{\FU}$ and, therefore 
\[
\begin{array}{ccccc}
(\CJ|_{\FU}+\CI|_{\FU})^{r}& \subset &\CJ|_{\FU}^{r} +\CI|_{\FU} &\subset &\CJ'|_{\FU}+ \CI|_{\FU}\\
(\CJ'|_{\FU}+\CI|_{\FU})^{s} &\subset &\CJ'|_{\FU}^{s} +\CI|_{\FU} &\subset &\CJ|_{\FU}+ \CI|_{\FU},
\end{array}
\]
so 
\[\invlim {n \in \NN} \CO_{\FX}/(\CJ+\CI)^{n+1} =\invlim {n \in \NN} \CO_{\FX}/(\CJ'+\CI)^{n+1}.\]
 
It is an easy exercise to check that $\FX_{/ \FX'}$ satisfies the hypothesis of  \cite[(10.6.3) and (10.6.4)]{EGA1},  from it we deduce that:
\begin{enumerate}
\item The formal scheme
 $\FX_{/ \FX'}$ is locally noetherian.
 \item The Ideal
 \[(\CI+\CJ)_{/\FX'}:=\invlim {n \in \sfn} (\CJ+\CI)/(\CJ+\CI)^{n+1}\] is an Ideal of  definition of  $\FX_{/\FX'}$.
 \item It holds that
 $\CO_{\FX_{/\FX'}}/ ((\CI+\CJ)_{/\FX'})^{n+1} = \CO_{\FX}/ (\CJ+\CI)^{n+1}$.
\end{enumerate}
\end{defn}

\begin{parraf} \label{limcomple}
With the above notations, if $Z_{n} = (\FX',  \CO_{\FX}/ (\CJ+\CI)^{n+1})$  for all $n \in \NN$, by \ref{notlim} we have that \[\FX_{/ \FX'} = \dirlim {n \in \NN}Z_{n}\]
On the other hand, if $X_{n} =(\FX, \CO_{\FX}/\CJ^{n+1})$ and $X'_{n} =(\FX', \CO_{\FX}/(\CJ^{n+1}+\CI))$, the canonical  morphisms 
\[\frac{\CO_{\FX}}{\CJ^{n+1}} \epi \frac{\CO_{\FX}}{(\CJ+\CI)^{n+1}} \epi \frac{\CO_{\FX}}{\CJ^{n+1} + \CI}\] 
provide the  closed immersions of  schemes $X'_{n} \xto{j_{n}} Z_{n} \xto{\kappa_{n}} X_{n}$, for all $n \in \NN$ such that the following diagrams are commutative:
\begin{diagram}[height=2em,w=2em,labelstyle=\scriptstyle]
X'_{m}			&	\rTto^{j_{m}} &	Z_{m}	&\rTto^{\kappa_{m}}  &   X_{m}\\
\uTto			   	&			& \uTto	&			         &     \uTto\\
X'_{n}			&	\rTto^{j_{n}} 	& Z_{n}	&\rTto^{\kappa_{n}}  &  X_{n}\\
\end{diagram}for all $m \ge n \ge 0$.
Then by \ref{lim} we have the canonical morphisms in $\sfn$
\[
\FX' \xto{j} \FX_{/ \FX'} \xto{\kappa} \FX 
\]
where $j$ is a  closed immersion (see Proposition \ref{ec}). The  morphism $\kappa$ as topological map is the inclusion and it is called  \emph{morphism of  completion of  $\FX$ along $\FX'$}. 
\end{parraf}

\begin{rem}
Observe that $\kappa$ is adic only if $\CI$ is contained in a Ideal of  definition of  $\FX$, in which case $\FX= \FX_{/ \FX'}$ and $\kappa=1_{\FX}$.
\end{rem}

\begin{parraf} \label{compl}
If $\FX= \spf(A)$ is in $\sfna$ with $A$ a  $J$-adic noetherian ring, and $\FX'= \spf(A/I)$ is a closed formal scheme  of  $\FX$, then 
\[\ga(\FX_{/ \FX'},\CO_{\FX_{/ \FX'}}) =  \invlim {n \in \NN} \frac{A}{(J+I)^{n+1}}=: \HA\]
and from the equivalence of  categories (\ref{equiv}) we have that 
$\FX_{/ \FX'}=\spf(\HA)$ and the morphisms $\FX' \xto{j} \FX_{/ \FX'} \xto{\kappa} \FX$ correspond to the natural continuous  morphisms $A \to \HA \to A /I $. 
\end{parraf}

\begin{propo} \label{caractcom}
Given $\FX$ in $\sfn$ and $\FX'$ a closed formal subscheme of  $\FX$, the  morphism of  completion  $\kappa: \FX_{/ \FX'} \to \FX$ is a flat pseudo  closed immersion.
\end{propo}
\begin{proof}
With  the notations of  \ref{limcomple} we have that \[\kappa= \dirlim {n \in \NN} \kappa_{n}\] and, since $\kappa_{n}$ is a  closed immersion for all $n \in \NN$, it follows that $\kappa$ is a pseudo  closed immersion. In order to prove that $\kappa$ is a flat morphism we may suppose that
$\FX=\spf(A)$ and  $\FX'= \spf(A/I)$, where $A$ is a $J$-adic noetherian ring. Then $\FX_{/ \FX'}=\spf(\HA)$ where $\HA$ is the  completion  of  $A$ for the $(J+I)$-adic topology and, therefore, is flat over $A$. By   Lemma \ref{pl}  it follows that $\kappa$ is a flat morphism. 
\end{proof}

\begin{rem}
Reciprocally, in Theorem  \ref{caracmorfcompl} we will see that every  flat  pseudo  closed immersion is a morphism of  completion.
\end{rem}

\begin{parraf}  \label{complmorf}
Given $f:\FX \to \FY$ in $\sfn$, let $\FX' \subset \FX$ and  $\FY' \subset \FY$ be closed formal subschemes given by Ideals $\CI \subset \CO_{\FX}$ and $\CL \subset \CO_{\FY}$ such that $f^{*}(\CL)\CO_{\FX} \subset \CI$, that is, $f(\FX')\subset \FY'$. Put 
\begin{align*}
X'_{n} &= (\FX', \CO_{\FX}/\CJ^{n+1}+\CI), &\\
Y'_{n} &= (\FY', \CO_{\FX}/\CK^{n+1}+\CL) & \text{ for all } n \in \NN.
\end{align*}
Further, take $\CJ \subset \CO_{\FX}$ and $\CK \subset \CO_{\FY}$ Ideals of  definition of $\FX$ and $\FY$, respectively, such that $ f^{*}(\CK) \CO_{\FX} \subset \CJ$. Write
\begin{align*}
X_{n} &= (\FX, \CO_{\FX}/\CJ^{n+1}), &\quad 
Y_{n} &= (\FY, \CO_{\FY}/\CK^{n+1}), & \\
Z_{n} &= (\FX',  \CO_{\FX}/ (\CJ+\CI)^{n+1}), &\quad 
W_{n} &= (\FY',  \CO_{\FY}/ (\CK+\CL)^{n+1}) & \text{for all } n \in \NN.
\end{align*}
Then the  morphism $f$ induces the following commutative  diagram of  locally noetherian schemes:
\begin{diagram}[height=1.5em,w=2em,labelstyle=\scriptstyle]
X_{n}&&\rTto^{f_{n}}&&Y_{n}&&&\\
&\rdTto&&&\uTto^{\kappa_{n}}&\rdTto&&\\
\uTto^{\kappa_{n}}&&X_{m}&\rTto^{f_{m}}&\HonV&&Y_{m}&\\
&&\uTto^{\kappa_{m}}&&\vLine&&&\\
Z_{n}&\hLine&\VonH&\rTto^{\hf_{n}}&W_{n}&&\uTto_{\kappa_{m}}&\\
&\rdTto&&&\uTto&\rdTto&&\\
\uTto^{j_{n}}&&Z_{m}&\rTto^{\hf_{m}}&\HonV&&W_{m} & \\
&&\uTto^{j_{m}}&&\vLine^{j_{n}}&&&\\
X'_{n}&\hLine&\VonH&\rTto& Y'_{n}&&\uTto_{j_{m}}&\\
&\rdTto&&&&\rdTto&&\\
&&X'_{m}&&\rTto&& Y'_{m}& \\
\end{diagram}
for all $m \ge n \ge 0$. Applying the direct limit over ${n \in \NN}$ we obtain  a morphism $\hf : \FX_{/ \FX'} \to \FY_{/ \FY'}$ in $\sfn$, that we call  \emph{completion of  $f$ along  $\FX'$ and $\FY'$}, such that  the  following diagram is commutative:
\begin{equation} \label{diagrcom}
\begin{diagram}[height=2em,w=2em,labelstyle=\scriptstyle]
\FX	          & \rTto^{f}   & \FY\\
\uTto{\kappa} &             & \uTto{\kappa}\\
\FX_{/ \FX'}	  &	\rTto^{\hf} & \FY_{/ \FY'} \\
\uTto         &             & \uTto\\
\FX'	          & \rTto^{f|_{\FX'}} & \FY'
\end{diagram}
\end{equation}
\end{parraf}

\begin{parraf} \label{complafin}
Suppose that $f: \FX=\spf(A) \to \FY=\spf(B)$ is in $\sfna$ and that $\FX'= \spf(A/I)$ and $\FY'= \spf(B/L)$ with $LA \subset I$. If $J \subset A$ and $K \subset B$ are Ideals of  definition such that $KA \subset J$, the  morphism $\hf : \FX_{/ \FX'} \to \FY_{/ \FY'}$ corresponds through formal Descartes duality (\ref{equiv}) to the   morphism induced by $B \to A$
\[\widehat{B}\to \HA\]
where $\HA$ is the completion of  $A$ for  the $(I+J)$-adic topology and $\widehat{B}$ denotes the  completion of  $B$ for the  $(K+L)$-adic topology.
\end{parraf}

\begin{propo} \label{cccom}
Given $f:\FX \to \FY$ in $\sfn$, let $\FY' \subset \FY$ be a closed formal subscheme and $\FX' = f^{-1}(\FY')$. Then, 
\[\FX_{/ \FX'} = \FY_{/ \FY'} \times_{\FY} \FX.\]
\end{propo}

\begin{proof}
We may suppose that $\FX=\spf(A), \FY=\spf(B)$ and $\FY'= \spf(B/L)$ are affine formal schemes and that $J \subset A$ and $K \subset B$ are ideals of  definition such that $KA \subset J$. By hypothesis,  $\FX'= \spf(A/LA)$, so $\FX_{/ \FX'} = \spf(\HA)$ where $\HA$ is the  completion of  $A$ for the  $(J+LA)$-adic topology. On the other hand, $\FY_{/ \FY'} = \spf(\HB)$ where $\HB$ denotes the completion of  $B$ for the $(K+L)$-adic  topology and  it holds that 
\[\HB \tc_{B} A= B\tc_{B} A=\HA,\]
since $J+(K+L)A = J+KA+LA= J+LA$ so we get the result.
\end{proof}

\begin{propo} \label{complsorit}
Given $f:\FX \to \FY$ in $\sfn$, consider $\FX' \subset \FX$ and $\FY' \subset \FY$ closed formal subschemes such that  $f(\FX')\subset \FY'$.
\begin{enumerate}
\item
Let $\CP$ be one of  the following properties of  morphisms in $\sfn$: 
\begin{center}
\emph{pseudo finite type,  pseudo finite, pseudo  closed immersion, pseudo quasi-finite, quasi-covering, flat, separated, radical.}
\end{center}
If  $f$ verifies $\CP$, then so does $\hf$. 
\item
Moreover, if $\FX' = f^{-1}(\FY')$, let $\CQ$ be one of  the following properties of  morphisms in $\sfn$: 
\begin{center}
\emph{adic, finite type,  finite,   closed immersion.}
\end{center}
Then, if $f$ verifies $\CQ$, then so does $\hf$. 
\end{enumerate}
\end{propo}
\begin{proof}

Suppose that $f$ is flat and let us prove that $\hf$ is flat. For that, we may take   $f: \FX=\spf(A) \to \FY=\spf(B)$  in $\sfna$, $\FX'= \spf(A/I)$ and $\FY'= \spf(B/L)$ with $LA \subset I$. Let $J \subset A$ and $ K \subset B$ be ideals of  definition such that $KA \subset J$ and, $\HA$ and $\widehat{B}$  the completions of $A$ and $B$ for the topologies given by $(I+J)\subset A$ and $(K+L) \subset B$, respectively. By \cite[III, \S5.4, Proposition 4]{b} we have that the  morphism $\widehat{B}\to \HA$ is  flat and, from  \ref{complafin} and  Lemma \ref{pl} there results that $\hf$ is flat.

Suppose that $f$ satisfies any of the other properties $\CP$ and let us prove that $\hf$ inherits it. From  \ref{psep}, \ref{prad}  and   Proposition \ref{caractcom} we deduce that the  morphism of  completion $\kappa: \FX_{/ \FX'} \to \FX$ satisfies $\CP$ and, by Corollary \ref{cor1sor} and  Corollary \ref{cor2sor} we have that $f \circ \kappa$ satisfies $\CP$. Then by the commutativity of diagram (\ref{diagrcom}) we obtain that the composition \[\FX_{/ \FX'} \xto{\hf} \FY_{/ \FY'} \xto{\kappa} \FY\] also satisfies $\CP$. Applying again  Corollary \ref{cor1sor} and  Corollary \ref{cor2sor} there results that $\CP$ holds for $\hf$.

Finally, if $f$ is adic from  Proposition \ref{cccom} and  Proposition \ref{adic}.(\ref{adic3}) we deduce that $\hf$ is adic. Then, if $\CQ$ is any of  the other properties and $f$ has $\CQ$, by (1) we have that $\hf$ satisfies $\CQ$.
\end{proof}

\begin{cor}
Given $\FX$ in $\sfn$ and $\FX' \subset \FX$ a closed formal subscheme,  the  morphism of  completion  $\kappa: \FX_{/ \FX'} \to \FX$  is  pseudo finite type,  pseudo finite, pseudo  closed immersion, pseudo quasi-finite, quasi-covering, flat, separated and radical.

\end{cor}



\chapter[Infinitesimal properties of  morphisms]{Infinitesimal lifting properties of  morphisms of  formal schemes} \label{cap2}
\setcounter{equation}{0}
We consider here infinitesimal lifting properties in the category $\sfn$. Our basic tool is an appropriate analogue of the  module of differentials, which we will use to obtain properties of the infinitesimal conditions. We will adapt the definitions from \cite[\S17.3]{EGA44} to the case of formal schemes and develop its basic properties. As in \emph{loc. cit.\/} these conditions show the most pleasant behavior under a finiteness condition, that in our case it will be to of pseudo finite type. The chapter closes with a version of Zariski's Jacobian criterion in this context.


\section[Differentials of topological algebras]{The module of differentials of topological algebras} \label{sec22}

The usual module of differentials of a homomorphism $\phi:A \to B$ of topological algebras is not complete, but its completion has the good properties of the module of  differentials in the discrete case. We will recall form \cite[Chapter \textbf{0}]{EGA44} some basic properties  of the completed module of  differentials $\om^{1}_{A/B}$ associated to a continuous morphism $A \to B$ of adic rings. We will add some computations for basic examples as power series and restricted power series algebras.

\begin{defn}
(\emph{cf.} \cite[\S \textbf{0}.20.4, p. 219]{EGA44})
Given $B \to A$  a continuous homomorphism  of  preadic rings and $K \subset B,\, J \subset A$  ideals of  definition such that $KA \subset J$, we denote by $I$  the   kernel of the continuous  epimorphism 
\[
\begin{array}{ccc}
A \otimes_{B} A &\to &A \\
a_{1} \otimes a_{2} &\rightsquigarrow  &a_{1} \cdot a_{2}
\end{array}
\]
where in $A \otimes_{B} A$ we consider  the tensor product topology  (\emph{cf.}  \cite[III, Exercise \S2.28]{b}).
We define the  \emph{ module of  $1$-differentials of  $A$ over $B$} as  the  topological $A$-module 
\[\Omega^{1}_{A/B} = \frac{I}{I^{2}}\]
with  the
$J$-adic topology\footnote{Notice that the $J$-adic topology in $\Omega^{1}_{A/B}$ agrees with  the quotient topology  of  the topology induced by $A \otimes_{B} A$ on $I$ (\emph{cf.} \cite[(\textbf{0}.20.4.5)]{EGA44}) }. We denote by $\om^{1}_{A/B}$, the completion of  $\Omega^{1}_{A/B}$ with respect  to  the $J$-adic topology, that is, 
\[\om^{1}_{A/B} = \invlim {n \in \NN} \frac{\Omega^{1}_{A/B}}{J^{n+1}\Omega^{1}_{A/B}}\]
In particular, if $B \to A$ is a homomorphism of discrete rings we have that $\Omega^{1}_{A/B}= \om^{1}_{A/B}$.
\end{defn}

\begin{defn}\cite[(\textbf{0}.20.3.1)]{EGA44} \label{derivcontalg}
Let $B \to A$ be a continuous homomorphism of  preadic rings and $L$ a topological  $A$-module. A \emph{continuous $B$-derivation of  $A$ in $L$ }is a map $d: A \to L$ satisfying the following conditions:
\begin{enumerate}
\item
$d \in \Homcont_{B}(A,L)$
\item
For all $a_{1},a_{2} \in A$, we have that $d(a_{1} \cdot a_{2})= a_{1} \cdot d(a_{2})+a_{2} \cdot d(a_{1})$.
\end{enumerate}
Let us denote by $\Dercont_{B}(A,L)$ the  set of continuous $B$-derivations of  $A$ in $L$. Obviously $\Dercont_{B}(A,L)$ is a subgroup of the group $\Der_{B}(A,L)$ of all $B$-derivations of  $A$ in $L$.
\end{defn}

\begin{parraf} Given $B \to A$ a continuous homomorphism of  preadic rings, let $j_{1}: A \to A \otimes_{B} A$,  $j_{2}: A \to A \otimes_{B} A$ be the canonical continuous $B$-homomorphisms given  by $j_{1}(a)=a \otimes 1,\, j_{2}(a) = 1 \otimes a$. If we denote by $p_{1}$ and $p_{2}$   the canonical continuous $B$-homomorphisms
\[
A \xto{j_{1}} A \otimes_{B} A \xto{p} (A \otimes_{B} A)/I^{2} \,  \, \, \textrm{ y } \, \, \,  A \xto{j_{2}} A \otimes_{B} A \xto{p} (A \otimes_{B} A)/I^{2}\textrm{,}
\]
respectively, 
where $p$ is the canonical projection, there results that 
\[d_{A/B}:=p_{1}-p_{2}: A \to\Omega^{1}_{A/B}\] 
is a continuous $B$-derivation; and it is called \emph{the canonical derivation of  $A$ over $B$}.
The continuous $B$-derivation  $d_{A/B}$ extends  naturally to a continuous $\HB$-derivation  (which, with an abuse of  terminology,  we will call \emph{canonical complete derivation  of  $\HA$ over $\HB$})
\[\HA \xto{\hd_{A/B}} \om^{1}_{A/B}\]
such that the  following diagram is commutative:
\begin{diagram}[height=2em,w=2em,labelstyle=\scriptstyle]
A 			         &     	\rTto			   & \HA     \\
\dTto^{d_{A/B}}          &   					   & \dTto_{\hd_{A/B}}      \\ 
\Omega^{1}_{A/B}    &	\rTto	       		   	   & \om^{1}_{A/B} \\  
\end{diagram}
In particular, if $B \to A$ is a homomorphism of  discrete rings there results that $d_{A/B}=\hd_{A/B}$.
\end{parraf}
\begin{rem}
Given $B \to A$ a morphism of  noetherian adic rings, in general $\Omega^{1}_{A/B}$ is not a complete $A$-module.
\end{rem}

\begin{notacion}
Given $A$ a $J$-preadic ring, we will denote by $A \com$ the category of  complete $A$-modules for the $J$-adic topology. With this notation, if $B \to A$ is a continuous homomorphism of preadic rings, then $\om^{1}_{A/B} \in  A \com$.
\end{notacion}

\begin{parraf}\label{modifn}
Let $B \to A$ be a morphism of   preadic rings. If $J \subset A$ is an ideal of  definition we denote by $\HA$ the  complete ring of  $A$ for the $J$-adic topology. For all $M \in A \com$ the   isomorphism $\Homcont_{A}(\Omega^{1}_{A/B},M)  \cong \Dercont_{B}(A , M)$ (\emph{cf.}   \cite[(\textbf{0}.20.4.8.2)]{EGA44}) induces the  following  canonical isomorphism of  $B$-modules
 \begin{equation} \label{derivmodulalgeb}
 \begin{array}{ccc}
 \Homcont_{A}(\om^{1}_{A/B},M) & \cong &\Dercont_{B}(\HA , M) \\
 u	&\rightsquigarrow& u \circ \hd_{A/B}.\\
 \end{array}
\end{equation}
That is, $(\om^{1}_{A/B}, \hd_{A/B})$ represents  the functor
\[
M  \in A \com \leadsto \Dercont_{B}(\HA,M).
\]
In particular, if $M$ is a  $A/J$-module we have the  isomorphism
 \begin{equation} \label{derivmodulalgebdis}
 \Hom_{A}(\om^{1}_{A/B},M) \cong \Der_{B}(\HA , M).
\end{equation}
\end{parraf}

\begin{parraf} \label{modifcomplto2}
Given a continuous morphism of  preadic rings $B \to A$, let $K \subset B$ and $J \subset A$ be  ideals of  definition such that $KA \subset J$. For all $n \in \NN$, if we put $A_{n} = A/J^{n+1}$ and $B_{n}=B/K^{n+1}$, there is a canonical isomorphism
\begin{equation*} 
 (\om^{1}_{A/B}, \hd_{A/B}) \cong( \invlim {n \in \NN} \Omega^{1}_{A_{n}/B_{n}}, \invlim {n \in \NN} d_{A_{n}/B_{n}}).
\end{equation*}
Indeed, for each $n \in \NN$ we have that $\hd_{A/B} \otimes_{A} A_{n}: A_{n} \to \om^{1}_{A/B} \otimes_{A} A_{n}$ is a continuous $B_{n}$-derivation  and, therefore by (\ref{derivmodulalgeb})  there exists a morphism of  $A_{n}$-modules $u_{n}: \Omega^{1}_{A_{n}/B_{n}} \to \om^{1}_{A/B} \otimes_{A} A_{n}$ such that  the following diagrams are commutative:
\begin{diagram}[height=1.8em,w=1.5em,labelstyle=\scriptstyle,midshaft]
A_{m} &	\rTto^{\hd_{A/B}\otimes_{A} A_{m}} & & &\om^{1}_{A/B} \otimes_{A} A_{m}\\
&\rdTto^{d_{A_{m}/B_{m}}}& &\ruTto^{u_{m}}& \\
&&\Omega^{1}_{A_{m}/B_{m}}&&\\
\dTonto&&\dTonto& & \dTonto\\
A_{n} &	\rLine^{\hd_{A/B}\otimes_{A} A_{n}} & \VonH &\rTto &\om^{1}_{A/B} \otimes_{A} A_{n}\\
&\rdTto^{d_{A_{n}/B_{n}}}&&\ruTto^{u_{n}}& \\
&&\Omega^{1}_{A_{n}/B_{n}}&&\\
\end{diagram}
for all $m \ge n \ge0$. Applying inverse limit we conclude  that there exists a morphism of  $A$-modules 
\[ u\colon \invlim {n \in \NN}\Omega^{1}_{A_{n}/B_{n}} \to \om^{1}_{A/B} \quad \text{ such that } \quad u \circ \invlim {n \in \NN} d_{A_{n}/B_{n}} =\hd_{A/B}.\] 
Since  
\[\invlim {n \in \NN} d_{A_{n}/B_{n}} \in \Dercont_{B}(A, \invlim {n \in \NN} \Omega^{1}_{A_{n}/B_{n}}),\] 
by uniqueness  (up to canonical  isomorphism) of the pair  $(\om^{1}_{A/B}, \hd_{A/B})$, we deduce the result.
\end{parraf}

\begin{rem}
Given $B \to A$ a morphism of  preadic rings, let $K \subset B,\, J \subset A$  ideals of  definition such that $KA \subset J$. As a consequence of  \ref{modifcomplto2} there results that 
\begin{equation} \label{modifcomplto}
( \om^{1}_{A/B},  \hd_{A/B}) \cong  (\om^{1}_{\HA/\HB}, \hd_{\HA/\HB})
\end{equation}
where $\HA$ and $\HB$ denote the completions of  $A$ and $B$, with respect to the  $J$ and $K$-preadic topologies, respectively and $\om^{1}_{A/B}$ and $\om^{1}_{\HA/\HB}$  denote the completions of  $A$ and $B$, with respect to the $K$-preadic and $K\HB$-adic topologies, respectively.
\end{rem}

\begin{parraf} \label{modifeespaefindisc}
Given a discrete ring $A$ and  $\mathbf{T} = T_{1},\, T_{2},\, \ldots T_{r}$ a finite number  of indeterminates, the canonical derivation of  $A[\mathbf{T}]$ over $A$, $d:=d_{A[\mathbf{T}]/A}$, is given by
\[
\begin{array}{ccc}
A[\mathbf{T}] &\overset{d} \lto & \Omega^{1}_{A[\mathbf{T}]/A}\\
P=P(\mathbf{T}) & \leadsto & \sum_{i=1}^{r}  \frac{\partial P}{\partial T_{i}} d T_{i}\\
\end{array}
\]
and $\Omega^{1}_{A[\mathbf{T}]/A}$ is a free $A[\mathbf{T}]$-module of rank $r$ with basis 
\(\{ dT_{1},\, dT_{2},\, \ldots, dT_{r}\} \).

In the  following example this fact is generalized for the formal power series ring  and for the restricted  power series ring  over an  adic noetherian ring $A$.
\end{parraf}

\begin{ex} \label{modifafindis}
Let $A$ be an  adic noetherian ring and  consider the inclusion  $A \inc A\{\mathbf{T}\}[[\mathbf{Z}]] = A\{T_{1},\, T_{2},\,  \ldots,\, T_{r}\}[[Z_{1},\, Z_{2},\, \ldots,\, Z_{s}]]$. The canonical  derivation $d: A[\mathbf{T},\mathbf{Z}] \to \Omega^{1}_{A[\mathbf{T},\mathbf{Z}]/A}$ induces the complete canonical derivation $\hd:= \hd_{A\{\mathbf{T}\}[[\mathbf{Z}]]/A}$ of  $A\{\mathbf{T}\}[[\mathbf{Z}]]$ over $A$
\[
\begin{array}{ccc}
A\{\mathbf{T}\}[[\mathbf{Z}]]& \overset{\hd} \lto & \om^{1}_{A\{\mathbf{T}\}[[\mathbf{Z}]]/A}\\
P=P(\mathbf{T},\mathbf{Z})& \leadsto & \sum_{i=1}^{r}  \frac{\partial P}{\partial T_{i}} \hd T_{i} + \sum_{j=1}^{s}  \frac{\partial P}{\partial Z_{j}} \hd Z_{j}\\
\end{array}
\]
It holds that $\om^{1}_{A\{\mathbf{T}\}[[\mathbf{Z}]]/A}$ is a free $A\{\mathbf{T}\}[[\mathbf{Z}]]$-module of  rank $r+s$, since 
\begin{equation*}
\om^{1}_{A\{\mathbf{T}\}[[\mathbf{Z}]]/A}  \underset{\textrm{(\ref{modifcomplto})}} \cong \om^{1}_{A[\mathbf{T},\mathbf{Z}]/A}  \underset{\textrm{\ref{modifeespaefindisc}}}  \cong  \bigoplus_{i=1}^{r}  A\{\mathbf{T}\}[[\mathbf{Z}]] \hd T_{i} \oplus \bigoplus_{j=1}^{s}  A\{\mathbf{T}\}[[\mathbf{Z}]] \hd Z_{j} .
\end{equation*}
Hence, we deduce that 
$\{ \hd T_{1},\, \hd T_{2},\, \ldots,\, \hd  T_{r},\, \hd Z_{1},\, \ldots,\, \hd Z_{s} \}$
is a basis of the free module $\om^{1}_{A\{\mathbf{T}\}[[\mathbf{Z}]] /A}$.
\end{ex}

\begin{propo} \label{modifinito}
Let $f: \FX= \spf(A) \to \FY =\spf(B)$ be a  pseudo finite type morphism  in $\sfn$.  Then   $\om^{1}_{A/B}$ is a finite type $A$-module. 
\end{propo}

\begin{proof}
Let $J \subset A$ and $K \subset B$ be ideals of  definition such that $KA \subset J$.
By hypothesis we have that $B_{0}=B/K \to A_{0}=A/J$ is a  finite type morphism and therefore,  $\Omega^{1}_{A_{0}/B_{0}}$ is a finite type $A_{0}$-module. From the Second Fundamental Exact Sequence   (\emph{cf.}   \cite[(\textbf{0}.20.5.12.1)]{EGA44}) associated to  the morphisms of  discrete rings $B \to A \to A_{0}$
\[
\frac{J}{J^{2}} \to \Omega^{1}_{A/B} \otimes_{A} A_{0} \to\Omega^{1}_{A_{0}/B_{0}}\to 0
\]
we deduce that $\Omega^{1}_{A/B} \otimes_{A} A_{0}= \Omega^{1}_{A/B}/J\Omega^{1}_{A/B}$ is a $A_{0}$-module of  finite type and from \cite[(\textbf{0}.7.2.7)]{EGA1} it follows that $\om^{1}_{A/B}$ is a finite type $A$-module. 
\end{proof}

Note, however, that even if $f: \FX= \spf(A) \to \FY =\spf(B)$ is a pseudo finite type morphism  in $\sfn$, the following Example (\ref{modifnocom}) shows that, in general,   $\Omega^{1}_{A/B}$ is not  a finitely generated $A$-module.

\begin{ex} \cite[Exercise 16.14.a]{eis} \label{modifnocom}
Let $\mathbb{Q}$ be the  field of  rational numbers and $R=\mathbb{Q} [[T_{1},\, T_{2},\, \ldots,\, T_{r}]]$. Then $\Omega^{1}_{R/\mathbb{Q}}$ is not a finitely generated $R$-module.
\end{ex}

\begin{lem} \label{propimodifcom}
Let $f: \FX= \spf(A) \to \FY =\spf(B)$ be a pseudo finite type morphism in $\sfna$.  Then:
\begin{enumerate}
\item \label{propimodifcom1}
If $B'$ is an adic ring such that  $B'$ is a topological $B$-algebra and $A' = B' \tc_{B} A$ it holds that
\[
\om^{1}_{A'/B'} \cong \om^{1}_{A/B} \otimes_{A} A'.
\]
\item \label{propimodifcom2}
If $S \subset A$ is a multiplicative subset we have that
\[
\om^{1}_{A/B}\{S^{-1}\}\cong \om^{1}_{A\{S^{-1}\}/B} \cong \om^{1}_{A/B} \otimes_{A} A\{S^{-1}\}.
\]
\item \label{propimodifcom3}
For all open prime ideal  $\fp$ of  $A$ there results that
\[
\om^{1}_{A_{\{\fp\}} /B}\cong\widehat{(\om^{1}_{A/B})_{\{\fp\}}} \cong \om^{1}_{\widehat{A_{\fp}}/B}.
\]	
\end{enumerate}
\end{lem}

\begin{proof}
As for (1) it suffices to observe  that the induced topology in $\Omega^{1}_{A/B} \otimes_{A} A' $ is the  one given by the topology of  $A'$ and therefore as a  consequence of the canonical  isomorphism  of  $A'$-modules $\Omega^{1}_{A'/B'}  \cong \Omega^{1}_{A/B} \otimes_{A} A' $ (\emph{cf.}  \cite[(\textbf{0}.20.5.5)]{EGA44}) it holds  that 
\[\om^{1}_{A'/B'} \cong \Omega^{1}_{A/B} \tc_{A} A' \underset{\textrm{\cite[(\textbf{0}.7.7.1)]{EGA1}}}\cong \om^{1}_{A/B} \tc_{A} A'.\] 
From Proposition \ref{modifinito} we get  that $ \om^{1}_{A/B} $ is a $A$-module of  finite type hence, $\om^{1}_{A'/B'} \cong \om^{1}_{A/B} \otimes_{A} A'$. 

Property (2)  is deduced from  the canonical topological isomorphisms
\begin{equation*}
\om^{1}_{A/B}\{S^{-1}\}\cong \Omega^{1}_{A/B}\{S^{-1}\} \underset{\textrm{\cite[(\textbf{0}.20.5.9)]{EGA44}}}\cong \om^{1}_{S^{-1}A/B}
				 \underset{\textrm{(\ref{modifcomplto})}} \cong \om^{1}_{A\{S^{-1}\}/B}.
\end{equation*}	

Property (3) is consequence of  property (2). Indeed,  for all open prime ideals $\fp$ of  $A$ it holds that:
\begin{equation*}
\begin{split}
\widehat{(\om^{1}_{A/B})_{\{\fp\}}} &\cong(\dirlim {f \notin \fp}(\om^{1}_{A/B})_{\{f\}})^{\widehat{}}  
\underset{\textrm{\ref{modifinito}}}\cong (\dirlim {f \notin \fp}\om^{1}_{A/B} \otimes_{A} A_{\{f\}})^{\widehat{}}\\
&\cong \om^{1}_{A/B} \tc_{A} A_{\{\fp\}}\cong \Omega^{1}_{A/B} \tc_{A} A_{\fp} 
\underset{\textrm{\cite[(\textbf{0}.20.5.9)]{EGA44}}} \cong \widehat{\Omega^{1}_{A_{\fp}/B}}\\
& \cong \om^{1}_{A_{\{\fp\}} /B} \cong \om^{1}_{\widehat{A_{\fp}}/B}.
\end{split}
\end{equation*}
The two last  identifications follow from the isomorphism (\ref{modifcomplto}) keeping in mind  that $\widehat{A_{\{\fp\}}}= \widehat{A_{\fp}}$.
\end{proof}

\section[Differentials of formal schemes]{The module of differentials of a pseudo finite type map of formal schemes} \label{sec221}

Given $f: X \to Y$ a finite type morphism of  schemes, is well-known that $\Omega^{1}_{X/Y}$, the module of  $1$-differentials of  $X$ over $Y$, is an essential tool  for   study the smooth, unramified or \'etale character of  $f$. In this section, we introduce the  module of  $1$-differentials for  a morphism $f:\FX \to \FY$ in $\sfn$ and discuss its fundamental properties, which will be used  in the characterizations of  the infinitesimal conditions  in Section \ref{sec23}. We can not use the general definition for ringed spaces, because it does not take into account the topology in the structure sheaves. Our construction will be based on the affine construction developed in the previous section.

\begin{parraf}  \label{modiftfloc}
Let $f: \FX= \spf(A) \to \FY =\spf(B)$ be a pseudo finite type morphism in $\sfna$.  From the two last results in the previous section we obtain  that:
\begin{enumerate}
\item \label{modiftfloc1}
For all $S \subset A$ multiplicative subset,  $\om^{1}_{A/B}\{S^{-1}\}$ is a $A\{S^{-1}\}$-module of  finite type.
\item \label{modiftfloc2}
For all open prime ideal $\fp \subset A$, $(\om^{1}_{A/B})_{\{\fp\}} = \dirlim {f \notin  \fp} (\Omega^{1}_{A/B})_{\{f\}}$ is a $A_{\{\fp\}}$-module of  finite type.
\end{enumerate} 
\end{parraf}

Next, we establish the needed preliminaries for sheafifying our construction and define the differential pair $(\om^{1}_{\FX/\FY}, \hd_{\FX/\FY})$ of  a morphism $\FX \to \FY$ in $\sfn$.

\begin{parraf} \cite[(10.10.1)]{EGA1} \label{deftriangulito}
Put $\FX=\spf(A)$ with $A$ a  $J$-adic noetherian ring, $X=\spec(A)$ and $X'=\spec(A/J)$, so we have that $\FX= X_{/X'}$. Given $M$  a $A$-module, we define  the  topological $\CO_{\FX}$-Module
\[M^{\tr}:=(\widetilde{M})_{/X'} = \invlim {n \in \NN} \frac{\widetilde{M}}{\widetilde{J}^{n+1}\widetilde{M}}\]
Moreover, a morphism $u: M \to N$ in $A\modu$ corresponds to a morphism of  $\CO_{X}$-Modules $\tilde u:\widetilde M \to \widetilde N$ that induces a morphism of  $\CO_{\FX}$-Modules 
\[M^{\tr} \xto{u^{\tr}} N^{\tr} \quad = \quad \invlim {n \in \NN} (\frac{\widetilde{M}}{\widetilde{J}^{n+1}\widetilde{M}} \xto{\tilde{u}_{n}} \frac{\widetilde{N}}{\widetilde{J}^{n+1}\widetilde{N}}).\]
So we get an additive covariant  functor from  the category of  $A$-modules to the category of  $\CO_{\FX}$-Modules
\begin{equation} \label{triangulito}
\begin{array}{ccc}
A\modu& \overset{\tr} \lto & \Modu(\FX)\\
M					&  \rightsquigarrow    &  M^{\tr}.\\
\end{array}
\end{equation}

If $M \in A\modu$ and $\HM$ denotes the complete  module of  $M$ for the $J$-adic topology, we easily deduce that (see the definition of the functor $\tr$): 
\begin{enumerate}
\item $\ga(\FX,M^{\tr}) =\HM$,
\item for all $f \in A$, $\ga(\fD(f),M^{\tr})=M_{\{f\}}$,
\item and $M^{\tr} =\HM^{\tr}$.
\end{enumerate}
\end{parraf}

\begin{lem} \label{propitrian}
Put $\FX=\spf(A)$ with $A$ a $J$-adic noetherian ring. Then:
\begin{enumerate}
\item
The restriction of the functor $(-)^\tr$  to the category $A\com$  is a full  faithful  functor.
\item \label{trianexac}
The functor $(-)^\tr$  is exact in the category of  $A$-modules of  finite type.
\item \label{equivtrian}
\cite[(10.10.2)]{EGA1} The  functor $(-)^\tr$ defines an equivalence of  categories between finite type $A$-modules   and the category $\coh(\FX)$ of  coherent $\CO_{\FX}$-Modules. 
\end{enumerate}
\end{lem}

\begin{proof}
As for (1), it suffices to show that given $M,\, N \in A\com$, then
\[\Hom_{A}(M,N) \cong \Hom_{\CO_{\FX}}(M^{\tr},N^{\tr})\]
Indeed, given $M \to N$ in $A\com$, by \cite[(1.3.8)]{EGA1} it corresponds biunivocally to  a morphism $\widetilde{M} \to \widetilde{N}$ that induces a morphism of  $\CO_{\FX}$-Modules $M^{\tr} \to N^{\tr}$. Reciprocally, every morphism $M^{\tr} \to N^{\tr}$ in $\Modu(\FX)$ gives a morphism in $A\com$ 
\[\ga(\FX,M^{\tr})=M \to \ga(\FX,N^{\tr})=N.\]

In order to prove (2) let  $0 \to M' \to M \to M''\to 0$ be an exact sequence  of   finite type $A$-modules. Then, given $\fp \subset A$ an open prime ideal, for all $f \in A \setminus \fp$ it holds that 
$0 \to M'_{\{f\}} \to M_{\{f\}} \to M''_{\{f\}} \to 0$
is an exact sequence  of  $A_{\{f\}}$-modules and, therefore,  we have  that
\[0 \to M'_{\{\fp\}} \to M_{\{\fp\}} \to M''_{\{\fp\}} \to 0\] is an exact sequence of  $A_{\{\fp\}}$-modules. 

Finally, by (1) the proof of  (3) amounts to check that given $\CF \in \coh(\FX)$ the  $A$-module $M= \ga(\FX, \CF)$ is finitely generated and $M^{\tr} \cong \CF$ (see \cite[IV, \S 4, Theorem 1]{mcl}). The reader can find this verification  in \cite[(10.10.2.9)]{EGA1}.

\end{proof}

In  the following definition we transfer to the geometric context of  locally noetherian  formal schemes  the definition of  continuous derivation.
\begin{defn}\label{derivcontef}
Let $f:\FX \to \FY$ be in $\sfn$ and $\CF$ a topological $\CO_{\FX}$-Module. A continuous  \emph{$\FY$-derivation of  $\CO_{\FX}$ in $\CF$} is a continuous homomorphism of  sheaves of  abelian groups $d: \CO_{\FX} \to \CF$ such that, for all open $\FU \subset \FX$  if we put   $d_{\FU} = \ga(\FU,d)$,  the following conditions hold:
\begin{enumerate}
\item
For all $a_{1},\, a_{2} \in \ga(\FU, \CO_{\FX})$,  \[d_{\FU}(a_{1} \cdot a_{2})= a_{1} \cdot d_{\FU}(a_{2})+a_{2}\cdot d_{\FU}(a_{1}).\]
\item
For all open $\FV \subset \FY$ such that  $f(\FU) \subset \FV$,
\[d_{\FU}(b|_{\FU} \cdot a) = b|_{\FU} \cdot d_{\FU}(a),\]
for all $b \in \ga(\FV, \CO_{\FY}),\, a \in \ga(\FU, \CO_{\FX})$. 
\end{enumerate} 

The set of  continuous $\FY$-derivations  of  $\CO_{\FX}$ with values in $\CF$ is denoted by $\Dercont_{\FY}(\CO_{\FX},\CF)$.
\end{defn}

\begin{rem}
Definition \ref{derivcontef} generalizes to formal schemes the definition of  $Y$-derivation for a morphism $X \to Y$ of  usual schemes (\emph{cf.}  \cite[(16.5.1)]{EGA44}).
\end{rem}

\begin{parraf}
Let $\FX=\spf(A) \to \FY=\spf(B)$ be in $\sfna$ and consider $J \subset A$ and $K \subset B$ Ideals of  definition such that $KA \subset J$. If $X= \spec(A)$, given $M \in A\modu$ and $d \in \Dercont_{B}(A,M)$ we define $d^{\tr} \in \Dercont_{\FY}(\CO_{\FX},M^{\tr})$ as
\[
\CO_{\FX} \xto{d^{\tr}} M^{\tr} = \invlim {n \in \NN} (\frac{\widetilde{A}}{\widetilde{J}^{n+1}\widetilde{A}} \xto{\widetilde {d} \otimes_{\CO_{X}} \CO_{X_{n}}} \frac{\widetilde{M}}{\widetilde{J}^{n+1}\widetilde{M}}).
\]
In particular, the complete canonical derivation of  $A$ over $B$,   $\hd_{A/B}$, induces  the continuous $\FY$-derivation   
\[\CO_{\FX} \xto{(\hd_{A/B})^{\tr}} (\om^{1}_{A/B})^{\tr} = \invlim {n \in \NN} (\frac{\widetilde{A}}{\widetilde{J}^{n+1}\widetilde{A}} \xto{\widetilde{d_{A/B}} \otimes_{\CO_{X}} \CO_{X_{n}}} \frac{\widetilde{\Omega^{1}_{A/B}}}{\widetilde{J}^{n+1}\widetilde{\Omega^{1}_{A/B}}}).\]
\end{parraf}

\begin{defn} 
Given $f: \FX \to \FY$ in $\sfn$ we call \emph{module of  $1$-differen\-tials of  $f$}  or \emph{module of  $1$-differentials of  $\FX$ over $\FY$} and we will denote it  by $\om^{1}_{f}$ or $\om^{1}_{\FX/\FY}$, the sheaf associated to the presheaf of  topological $\CO_{\FX}$-Modules locally  given  by
$(\om^{1}_{A/B})^{\tr} $, for all open sets $\FU=\spf(A) \subset \FX$ and $\FV=\spf(B) \subset \FY$ with $f(\FU) \subset \FV$. 
Note that $\om^{1}_{\FX/\FY}$ has structure of   $\CO_{\FX}$-Module.

Let  $\hd_{\FX/\FY}: \CO_{\FX} \to \om^{1}_{\FX/\FY}$ be the  morphism of  sheaves locally  defined  by
\[(d_{A/B}: A \to \om^{1}_{A/B})^{\tr}\] for all couple of  affine open sets $\FU=\spf(A) \subset \FX$ and $\FV=\spf(B) \subset \FY$ such that $f(\FU) \subset \FV$. The  morphism $\hd_{\FX/\FY}$ is  a continuous $\FY$-derivation and it is called \emph{canonical derivation of  $\FX$ over $\FY$}. We will refer to   $( \om^{1}_{\FX/\FY}, \hd_{\FX/\FY})$ as the  \emph{differential pair of  $\FX$ over $\FY$}.
\end{defn}

\begin{parraf} \label{defmofif}
If $\spf(A) \to \spf(B)$ is in $\sfna$, the  differential pair of  $\FX$ over $\FY$ is \[((\om^{1}_{A/B})^{\tr},(\hd_{A/B})^{\tr}).\]
In particular, if  $X=\spec(A) \to Y=\spec(B)$ is a morphism of  usual schemes, there results that $(\om^{1}_{X/Y},\hd_{X/Y}) = (\Omega^{1}_{X/Y}, d_{X/Y}) = (\widetilde{\Omega^{1}_{A/B}}, \widetilde{d_{A/B}})$ and it agrees with the definition of the  differential pair of  a morphism of  affine schemes  (\emph{cf.}  \cite[(16.5.1)]{EGA44}).
\end{parraf}

\begin{propo} \label{finitmodif}
Let $f: \FX \to \FY $ be a morphism in $\sfn$ of  pseudo finite type. Then \(\om^{1}_{\FX/\FY}\) is a coherent sheaf.
\end{propo}

\begin{proof}
We may suppose that $f: \FX=\spf(A) \to \FY=\spf(B)$ is in $\sfna$. Therefore, since $\om^{1}_{\FX/\FY} = (\om^{1}_{A/B})^{\tr}$ the result is deduced from Proposition \ref{modifinito} and from  the equivalence of categories  \ref{propitrian}.(\ref{equivtrian}).
\end{proof}

\begin{ex} \label{exmodifdisformafingeom}
Let $\FX$ be in $\sfn$. From  Example \ref{modifafindis}, we have that $\om^{1}_{\BD^{s}_{\BA^{r}_{\FX}}/\FX}$ is a locally free $\CO_{\BD^{s}_{\BA^{r}_{\FX}}}$-Module, in fact, free of (constant) rank  $r+s$. 
\end{ex}

\begin{propo} \label{modifmodifn}
Let $f: \FX \to \FY $ be in $\sfn$ and $\CJ \subset \CO_{\FX}$ and $\CK \subset \CO_{\FY}$ be Ideals of  definition  such that $f^{*}(\CK)\CO_{\FX} \subset \CJ$. Then
\[
(\om^{1}_{\FX/\FY}, \hd_{\FX/\FY} ) \cong (\invlim {n \in \NN} \Omega^{1}_{X_{n}/Y_{n}} , \invlim {n \in \NN} d_{X_{n}/Y_{n}}\, \, ).
\]
\end{propo}

\begin{proof}
For all $n \in \NN$, the  morphism $\hd_{\FX/\FY}$ induces a $Y_{n}$-derivation $\hd_{\FX/\FY}\otimes \CO_{X_{n}}: \CO_{X_{n}} \to \om^{1}_{\FX/\FY} \otimes_{\CO_{\FX}} \CO_{X{n}}$. By \cite[(16.5.3)]{EGA44} there exists a morphism $u_{n}: \Omega^{1}_{X_{n}/Y_{n}} \to \om^{1}_{\FX/\FY} \otimes_{\CO_{\FX}} \CO_{X{n}}$ of  $\CO_{X_{n}}$-Modules such that  for all $m \ge n \ge 0$ the  following diagrams are commutative
\begin{diagram}[height=2em,w=2em,labelstyle=\scriptstyle,midshaft]
\CO_{X_{m}} &	\rTto^{\hd_{\FX/\FY}\otimes \CO_{X_{m}}} & & &\om^{1}_{\FX/\FY} \otimes_{\CO_{\FX}} \CO_{X_{m}}\\
&\rdTto^{d_{X_{m}/Y_{m}}}&&\ruTto^{u_{m}}& \\
&&\Omega^{1}_{X_{m}/Y_{m}}&&\\
\dTonto&&\dTonto& & \dTonto\\
\CO_{X_{n}} &	\rLine^{\hd_{\FX/\FY}\otimes \CO_{X_{n}}} &\VonH &\rTto &\om^{1}_{\FX/\FY} \otimes_{\CO_{\FX}} \CO_{X_{n}}\\
&\rdTto^{d_{X_{n}/Y_{n}}}&&\ruTto^{u_{n}}& \\
&&\Omega^{1}_{X_{n}/Y_{n}}&&\\
\end{diagram}
Therefore, applying inverse limit we get a morphism of  $\CO_{\FX}$-Modules \[\invlim {n \in \NN} \Omega^{1}_{X_{n}/Y_{n}} \xto{u} \om^{1}_{\FX/\FY} \, 
\text{ such that } \, u \circ \invlim {n \in \NN} d_{X_{n}/Y_{n}}= \hd_{\FX/\FY}.\] In order to show that $u$ is an isomorphism we may suppose that $f: \FX=\spf(A) \to \FY=\spf(B)$ is in $\sfna$ and therefore the  result follows from  \ref{modifcomplto2}.
\end{proof}

\begin{rem}
In  \cite[2.6]{LNS} the  module of  differentials  of  a morphism $\FX \to \FY$ in $\sfn$ is defined via the characterization given in Proposition \ref{modifmodifn}. Therefore our definition agrees with the one in \emph{loc.~cit}. In our approach it will be a key fact that it represents the functor of continuous derivation.
\end{rem}
Given $X \to Y$ a morphism of  schemes in \cite[(16.5.3)]{EGA44} it is established that $(\Omega^{1}_{X/Y},d_{X/Y})$ is the  universal pair of the  representable  functor $\CF \in  \Modu(X) \rightsquigarrow \Der_{Y}(\CO_{X}, \CF)$. In  Theorem \ref{modrepr} this result  is generalized   for a morphism $\FX \to \FY$ in $\sfn$. In order to do that, we need to define the  category $\Com(\FX)$  of complete $\CO_{\FX}$-Modules. 

\begin{parraf}
Given $\FX $ in $\sfn$ and $\CJ \subset \CO_{\FX}$ an Ideal of  definition of  $\FX$ we will denote by  $\Com (\FX)$ the category of  $\CO_{\FX}$-Modules  $\CF$ such that 
\[\CF = \invlim {n \in \NN} (\CF \otimes_{\CO_{\FX}} \CO_{X_{n}}).\] 
It is easily seen  that the definition does not depend on the election of the Ideal of  definition of  $\FX$. 

For example:
\begin{enumerate}
\item
Given $\FX= \spf(A)$ in $\sfna$ and $J \subset A$ an ideal of  definition for all $A$-modules $M$, it holds that 
\[M^{\tr}= \invlim {n \in \NN} \frac{\widetilde{M}}{\widetilde{J}^{n+1}\widetilde{M}} \in \Com(\FX).\]
\item
Let $\FX$ be  in $\sfn$. For all  $\CF \in \coh (\FX)$, by \cite[(10.11.3)]{EGA1} we have that \[\CF =\invlim {n \in \NN} (\CF \otimes_{\CO_{\FX}} \CO_{X_{n}})\] and therefore, $\coh (\FX)$ is a full subcategory of  $\Com (\FX)$. In particular, given $f:\FX \to \FY$ a pseudo finite type morphism in $\sfn$, from \ref{finitmodif}, we get that  $\om^{1}_{\FX/\FY} \in \Com(\FX)$.
\end{enumerate}
\end{parraf}

Now we are ready to show that given $\FX \to \FY$ a pseudo finite type morphism in $\sfn$, $(\om^{1}_{\FX/\FY}, \hd_{\FX/\FY})$ is the  universal pair of the representable functor
\[\CF \in \Com (\FX)       \leadsto     \Dercont_{\FY}(\CO_{\FX}, \CF).\]

cor \label{modrepr} 
Let $f: \FX \to \FY $ be a morphism in $\sfn$. Then  given $\CF \in \Com (\FX)$ we have a canonical isomorphism
\[
\begin{array}{ccc}
\Hom(\om^{1}_{\FX/\FY},\CF) &\overset{\overset{\varphi}\sim} \lto & \Dercont_{\FY}(\CO_{\FX}, \CF)\\
u				  &\rightsquigarrow    &u \circ \hd_{\FX/\FY}\\
\end{array}
\]

\begin{proof}
Let $\CJ \subset \CO_{\FX}$ and $\CK \subset \CO_{\FY}$ be Ideals of  definition such that $f^{*}(\CK) \CO_{\FX} \subset \CJ$. Use these choices to express \[f = \dirlim {n \in \NN} (f_{n}: X_{n} \to Y_{n}).\] We are going to define the inverse map
 of  $\varphi$. Given $d \in  \Dercont_{\FY}(\CO_{\FX}, \CF)$, for all $n \in \NN$ it holds  that 
\[ d \otimes_{\CO_{\FX}} \CO_{X_{n}}  \in \Der_{Y_{n}}(\CO_{X_{n}}, \CF \otimes_{\CO_{\FX}} \CO_{X_{n}}) \qquad \textrm{and}\qquad d = \invlim {n \in \NN} (d \otimes_{\CO_{\FX}} \CO_{X_{n}}).\] 
There results that 
\[\Hom(\Omega^{1}_{X_{n}/Y_{n}}, \CF \otimes_{\CO_{\FX}}\CO_{X_{n}}) \cong \Der_{Y_{n}}(\CO_{X_{n}}, \CF \otimes_{\CO_{\FX}} \CO_{X_{n}})\]
and therefore, there exists a unique morphism of  $ \CO_{X_{n}}$-Modules 
\[\Omega^{1}_{X_{n}/Y_{n}} \xto{u_{n}} \CF _{n}\]
such  that $u_{n} \circ d_{X_{n}/Y_{n}} = d \otimes_{\CO_{\FX}} \CO_{X_{n}}$, for all $n \in \NN$. Moreover, for all couples of integers $m \ge n \ge 0$ there exists a unique  morphism of  $\CO_{X_{m}}$-Modules  $v_{nm}: \Omega^{1}_{X_{m}/Y_{m}} \to \Omega^{1}_{X_{n}/Y_{n}}$ such   that the  following diagram 
\begin{diagram}[height=2em,w=2em,p=0.3em,labelstyle=\scriptstyle]
\CO_{X_{m} }&\rTto^{d_{X_{m}/Y_{m}}}&\Omega^{1}_{X_{m}/Y_{m}} & \rTto^{u_{m}}	      &   \CF \otimes_{\CO_{\FX}} \CO_{X_{m}}\\
    \dTto^{can.} &&\dTto^{v_{nm}}				  & 			      &   \dTto^{can.}\\ 
\CO_{X_{n}}& \rTto^{d_{X_{n}/Y_{n}}}&\Omega^{1}_{X_{n}/Y_{n}} & \rTto^{u_{n}}	      &   \CF \otimes_{\CO_{\FX}} \CO_{X_{n}}\\
\end{diagram}
is commutative. That is, $(u_{n})_{n \in \NN}$ is an inverse system  of  morphisms and by Proposition \ref{modifmodifn} it induces a morphism of  $\CO_{\FX}$-Modules
\[
\om^{1}_{\FX/\FY} \xto{u} \CF = \invlim {n \in \NN}   \CF \otimes_{\CO_{\FX}} \CO_{X_{n}}
\]
such  that $u \circ \hd_{\FX/\FY} =d$. If we define 
\[
\begin{array}{ccc}
\Dercont_{\FY}(\CO_{\FX}, \CF)& \overset{\psi} \lto & \Hom(\om^{1}_{\FX/\FY},\CF), \\
d& \leadsto & u
\end{array}
\]
it is easily shown by construction  that $\psi$ and $\varphi $ are inverse maps.

The naturality is obvious.
\end{proof}

\begin{lem} \label{imagdirecompl}
Let $f: \FX \to \FY $ be a morphism in $\sfn$. Then given $\CF \in \Com (\FX)$ we have  that 
\[
f_{*} \CF= \invlim {n \in \NN} (f_{*} \CF \otimes_{\CO_{\FY}} \CO_{Y_{n}})
\]
and therefore, $f_{*} \CF$ is in $\Com(\FY)$.
\end{lem}

\begin{proof}
Let $\CJ \subset \CO_{\FX}$ and $\CK \subset \CO_{\FY}$ be  Ideals of  definition such that $f^{*}(\CK)\CO_{\FX} \subset \CJ$. For all $n \in \NN$ we have the canonical morphisms  $f_{*} \CF \to f_{*} \CF \otimes_{\CO_{\FY}} \CO_{Y_{n}}$ that induce the  morphism of  $\CO_{\FY}$-Modules
\[
f_{*} \CF \lto  \invlim {n \in \NN} (f_{*} \CF \otimes_{\CO_{\FY}} \CO_{Y_{n}})
\]
To see wether it is a isomorphism is a local question, therefore we may assume that $f:\FX=\spf(A) \to \FY= \spf(B)$ is in $\sfna$, $\CJ= J^{\tr}$ and $\CK= K^{\tr}$ with $J \subset A$ and $K \subset B$ ideals of  definition such that $KA \subset J$. Then $M=\ga(\FY, f_{*} \CF)$ is a complete $B$-module for  the $f^{-1}(J)$-adic topology and since $K \subset f^{-1}(J)$ we have that
\[M=\invlim {n \in \NN} \frac{M}{K^{n+1}M}\]
and the  result follows.
\end{proof}

\begin{propo} \label{cbasemodif}
Given a  commutative diagram  in $\sfn$ of  pseudo finite type  morphisms:
\begin{diagram}[height=2em,w=2em,labelstyle=\scriptstyle]
\FX	  	 &  \rTto 		& \FY\\
\uTto^{g}    &    			& \uTto\\ 
\FX' 		&    	  \rTto^{h}	&     \FY' \\
\end{diagram}
there exists a morphism of  $\CO_{\FX'}$-Modules 
$g^{*}\om^{1}_{\FX/\FY} \lto \om^{1}_{\FX'/\FY'}$
locally defined by $\hd_{\FX/\FY} a \otimes 1 \rightsquigarrow \hd_{\FX'/\FY'} g(a)$ (equivalently, there exists a morphism of  $\CO_{\FX}$-Modules 
$\om^{1}_{\FX/\FY} \lto g_{*}\om^{1}_{\FX'/\FY'}$
locally defined by $\hd_{\FX/\FY}  \rightsquigarrow \hd_{\FX'/\FY'} g(a)$).
If the  diagram is cartesian, the  above morphism is a isomorphism.
\end{propo}

\begin{proof}
The morphism $\CO_{\FX} \to g_{*}\CO_{\FX'} \xto{g_{*} \hd_{\FX'/\FY'}} g_{*} \om^{1}_{\FX'/\FY'}$ is a continuous $\FY$-derivation. Applying Proposition \ref{finitmodif} and  Lemma \ref{imagdirecompl} we have that $g_{*} \om^{1}_{\FX'/\FY'} \in \Com(\FX)$ and therefore by  Theorem \ref{modrepr} there exists an unique morphism of  $\CO_{\FX}$-Modules $\om^{1}_{\FX/\FY} \to g_{*} \om^{1}_{\FX'/\FY'}$ such that  the  following diagram is commutative
\begin{diagram}[height=2em,w=2em,labelstyle=\scriptstyle]
\CO_{\FX} 		&    	  \rTto^{\hd_{\FX/\FY}}	&  \om^{1}_{\FX/\FY} \\
\dTto   &    			& \dTto\\ 
g_{*} \CO_{\FX'} 	 &  \rTto^{g_{*} \hd_{\FX/\FY}}			&g_{*} \om^{1}_{\FX'/\FY'}\\ 
\end{diagram}
Equivalently, there exists a morphism of  $\CO_{\FX'}$-Modules $g^{*}\om^{1}_{\FX/\FY} \to \om^{1}_{\FX'/\FY'}$ locally  given by $\hd_{\FX/\FY} a \otimes 1 \rightsquigarrow \hd_{\FX'/\FY'} g(a)$. 

Let us assume that the  square is cartesian.We may suppose  that $\FX=\spf(A),\, \FY= \spf(B),\, \FY'=\spf(B')$ and $\FX' = \spf(A \tc_{B} B')$ are in $\sfna$ and therefore the  result is consequence of   \ref{propimodifcom}.(\ref{propimodifcom1}).
\end{proof}
With  the notations of  the previous proposition, if $\FY=\FY'$ the  morphism $g^{*}\om^{1}_{\FX/\FY} \lto \om^{1}_{\FX'/\FY}$ is denoted by $dg$ and is called \emph{the differential of  $g$ over $\FY$}.

\begin{cor}  \label{modifesqu}
Given $f:\FX \to \FY$ a  \emph{finite type} morphism in $\sfn$ consider $\CK\subset \CO_{\FY}$ and $\CJ =f^{*}(\CK)\CO_{\FX} \subset \CO_{\FX}$ Ideals of  definition such that we express  $f = \dirlim {n \in \NN} (f_{n}: X_{n} \to Y_{n})$. Then for all $n \in \NN$ we have that
\[
\Omega^{1}_{X_{n}/Y_{n}} \cong \om^{1}_{\FX/\FY} \otimes_{\CO_{\FX}} \CO_{X_{n}}.
\]
\end{cor}
\begin{proof}
Since $f$ is an adic  morphism, from Proposition \ref{ccad} we have that the diagrams  
\begin{diagram}[height=2em,w=2em,labelstyle=\scriptstyle]
\FX	&	\rTto^{f} &	\FY \\
\uTto			&	&		\uTto\\
X_{n}		&	\rTto^{f_{n}} &	Y_{n}\\
\end{diagram}
are cartesian, $\forall  n \ge 0$. Then the  corollary is deduced  from  the proposition.
\end{proof}

In  the  following example we show that if the  morphism is not adic, the last  corollary does not hold.
\begin{ex} \label{contrmod0}
Let  $K$ be a field  and $p:\BD^{1}_{K} \to \spec(K)$ the  projection  morphism of  the formal disc of  $1$ dimension  over $\spec(K)$. Given the  ideal of  definition $[[T]] \subset K[[T]]]$ such that \[p = \dirlim {n \in \NN} p_{n}\] it holds that $\Omega^{1}_{p_{0}}=0$ but, 
\[\om^{1}_{p} \otimes_{\CO_{\BD^{1}_{K}}} \CO_{\spec(K)} \underset{\ref{exmodifdisformafingeom} }=(K[[T]]\hd T) ^{\tr} \otimes_{K[[T]]^{\tr}} \widetilde{K} \cong \widetilde{K}\neq 0.\]
\end{ex}

\section[Fundamental Exact Sequences]{Fundamental Exact Sequences of modules of differentials of formal schemes} \label{sec222}

We extend the usual First and Second Fundamental Sequences to our construction of differentials of pseudo finite type morphisms between formal schemes. They will provide a basic tool for applying it to the study of the infinitesimal lifting. Also, we will give a local computation based on the Second Fundamental Exact Sequence.

\begin{propo}{(First Fundamental Exact Sequence)}\label{primersef}
Let $f:\FX \to \FY$  and $g: \FY \to \FS$ be two morphisms in $\sfn$ of  pseudo finite type. There exists an exact sequence  of  coherent $\CO_{\FX}$-Modules 
\begin{equation} \label{sef1}
f^{*}(\om^{1}_{\FY/\FS}) \xto{\Phi} \om^{1}_{\FX/\FS} \xto{\Psi} \om^{1}_{\FX/\FY} \to 0\\
\end{equation}
where $\Phi$ and $\Psi$  are locally defined  by
\[
\hd_{\FY/\FS} b \otimes 1 \leadsto \hd_{\FX/\FS} f(b) \qquad 
\hd_{\FX/\FS} a Ê         \leadsto \hd_{\FX/\FY} a\\
\]
\end{propo}

\begin{proof}
The morphism of  $\CO_{\FX}$-Modules $\Phi: f^{*}\om^{1}_{\FY/\FS}  \to \om^{1}_{\FX/\FS}$ is the differential of  $f$ over $\FS$. Since $\hd_{\FX/\FY}: \CO_{\FX} \to \om^{1}_{\FX/\FY}$ is a continuous $\FS$-derivation, from  Theorem \ref{modrepr} there exists a unique  morphism of  $\CO_{\FX}$-Modules $\Psi: \om^{1}_{\FX/\FS} \to \om^{1}_{\FX/\FY}$ such that  $\Psi \circ \hd_{\FX/\FS} = \hd_{\FX/\FY}$.

As for proving  the exactness, we can argue locally and assume that $f: \FX= \spf(A) \to \FY=\spf(B),\, g: \FY= \spf(B) \to \FS=\spf(C)$ are   pseudo finite type morphisms. The sequence of coherent $\CO_{\FX}$-Modules (\ref{sef1})  corresponds via the equivalence between the category of finite type $A$-modules and $\coh(\FX)$ \ref{propitrian}.(\ref{equivtrian}), to  the sequence 
\[
\om^{1}_{B/C} \otimes_{B} A \xto{\Phi} \om^{1}_{A/C}  \xto{\Psi} \om^{1}_{A/B} \to 0
\]
Therefore, it suffices to prove  that the above sequence is exact. Let $J \subset A,\, K \subset B$ and $L \subset C$ be  ideals of  definition such that $LB \subset K,\, KA \subset J$ and for all $n \in \NN$ put $A_{n}=A/J^{n+1},\, B_{n}=B/K^{n+1}$ and $C_{n}=C/L^{n+1}$. For all couples of  integers $m \ge n\ge0$ we have a commutative diagram of  finite type $A_{m}$-modules  
\begin{diagram}[height=2em,w=2em,labelstyle=\scriptstyle]
\Omega^{1}_{B_{m}/C_{m}} \otimes_{B_{m}} A_{m}& & \rTto^{\Phi_{m}} & &\Omega^{1}_{A_{m}/C_{m}} &\rTto^{\Psi_{m}} &  \Omega^{1}_{A_{m}/B_{m}}& \rTto& 0\\
&\rdTonto&&\ruTinc& &&&&\\
&&\Ker \Psi_{m}&&&&&&\\
\dTonto&&\dTonto& & \dTonto& & \dTonto& & \\
\Omega^{1}_{B_{n}/C_{n}} \otimes_{B_{n}} A_{n}&	\hLine^{\Phi_{n}} & \VonH& \rTto&\Omega^{1}_{A_{n}/C_{n}} &\rTto^{\Psi_{n}} &  \Omega^{1}_{A_{n}/B_{n}}& \rTto& 0\\
&\rdTonto&&\ruTinc& &&&&\\
&&\Ker \Psi_{n}&&&&&&\\
\end{diagram}
where the horizontal  sequences correspond with the First Fundamental Exact  Sequences for the morphisms of  (discrete) rings  $C_{j} \xto{g_{j}}B_{j}\xto{f_{j}} A_{j}$ for $j \in \{n, m\}$. Since the vertical maps are surjective and $\Ima \Phi$ is a   finite type $A$-module we deduce that
\begin{equation*}
\begin{split}
\invlim {n \in \NN}  (&\om^{1}_{B_{n}/C_{n}} \otimes_{B_{n}} A_{n} \xto{\Phi_{n}} \om^{1}_{A_{n}/C_{n}}  \xto{\Psi_{n}} \om^{1}_{A_{n}/B_{n}} \to 0) =\\
=\,\, &\om^{1}_{B/C} \otimes_{B} A \xto{\Phi} \om^{1}_{A/C}  \xto{\Psi} \om^{1}_{A/B} \to 0 
\end{split}
\end{equation*}
is an exact sequence.
\end{proof}

\begin{parraf}
Let $f:\FX \to \FY$ be  a pseudo finite type morphism in $\sfn$, and $\CJ \subset \CO_{\FX}$ and $\CK \subset \CO_{\FY}$ be Ideals of  definition with $f^{*}(\CK) \CO_{\FX} \subset \CJ$ and \[f: \FX \to \FY = \dirlim {n \in \NN} (f_{n}: X_{n} \to Y_{n})\]
the relevant expression for $f$. For all $n \in \NN$, from  the First Fundamental Exact Sequence  (\ref{sef1}) associated to $f_{n}: X_{n} \to Y_{n} \inc \FY$, we deduce that $\om^{1}_{X_{n}/\FY} = \Omega^{1}_{X_{n}/\FY}= \Omega^{1}_{X_{n}/Y_{n}}$.
\end{parraf}

\begin{parraf}
Given $\FX' \overset{i} \inc \FX$ a  closed immersion in $\sfn$ (Definition \ref{defnenccerrado}) we have that the  morphism $i^{\sharp}: i^{-1}(\CO_{\FX})  \to \CO_{\FX'}$ is an epimorphism. If $\CK:= \ker(i^{\sharp})$ we call \emph{the conormal sheaf of  $\FX'$ over $\FX$} to $\CC_{\FX'/\FX}:=\CK/ \CK^{2}$. 
 
It is easily  shown that $\CC_{\FX'/\FX}$ verifies the following properties:
\begin{enumerate}
\item
It is a coherent $\CO_{\FX'}$-module.
\item
If $\FX' \subset \FX$ is a closed subscheme given by a coherent  Ideal $\CI \subset \CO_{\FX}$,
then $\CC_{\FX'/\FX}=i^{*}(\CI/\CI^{2})$.
\end{enumerate}
\end{parraf}

\begin{propo}{(Second Fundamental Exact Sequence)} \label{segunsef}
Consider a pseudo finite type morphism in $\sfn$, $f:\FX \to \FY$ and  a  closed immersion $\FX' \overset{i} \inc \FX$. There exists an exact sequence of coherent $\CO_{\FX'}$-Modules
\begin{equation} \label{sef2}
\CC_{\FX'/\FX} \xto{\delta} i^{*}\om^{1}_{\FX/\FY} \xto{\Phi} \om^{1}_{\FX'/\FY} \to 0
\end{equation}
\end{propo}
\begin{proof}
The morphism $\Phi$ is the differential of $i$, given by Proposition \ref{cbasemodif} and  is defined by $\hd_{\FX/\FY} a \otimes 1 \rightsquigarrow \hd_{\FX'/\FY} i(a)$. If $\CI \subset \CO_{\FX}$ is the  Ideal that defines the closed subscheme $i(\FX') \subset \FX$ the  morphism $\delta$ is the one induced by $\hd_{\FX/\FY}|_{\CI}: \CI \to \om^{1}_{\FX/\FY}$. As for showing  the exactness we may suppose that $\FX' = \spf(A/I) \inc  \FX= \spf(A),\, \FY= \spf(B)$ are in $\sfna$ . With these hypothesis the Second Fundamental Exact  Sequence  corresponds through the equivalence \ref{propitrian}.(\ref{equivtrian})  between  the category of  finite type  $A/I$-modules and $\coh(\FX')$, to the sequence 
\[
\frac{I}{I^{2}} \xto{\delta} \om^{1}_{A/B} \otimes_{A}\frac{A}{I} \xto{\Phi} \om^{1}_{(A/I)/B} \to 0
\]
It suffices to show that this sequence is exact. Let $J \subset A$ and $K \subset B$ be ideals of  definition such that $KA \subset B$  and put $A_{n}=A/J^{n+1},\, B_{n} = B/K^{n+1},\, I_{n}=I+J^{n+1}/J^{n+1}$ and $A'_{n} = A_{n}/I_{n}$, for all $n \in \NN$. Then, for $m \ge n \ge 0$, we have the following commutative diagrams of  finite type $A'_{m}$-modules 
\begin{diagram}[height=2em,w=2.2em,labelstyle=\scriptstyle]
\frac{I_{m}}{I^{2}_{m}} &	 & &\rTto^{\delta_{m}\qquad} &\Omega^{1}_{A_{m}/B_{m}} \otimes_{A_{m}} A'_{m}&\rTto^{\Phi_{m}} &  \Omega^{1}_{A'_{m}/B_{m}}& \rTto& 0\\
&\rdTonto&&\ruTinc& &&&&\\
&&\Ker \Phi_{m}&&&&&&\\
\dTonto&&\dTonto& & \dTonto& & \dTonto& & \\
\frac{I_{n}}{I^{2}_{n}} &\hLine &\VonH & \rTto^{\delta_{n}\qquad}&\Omega^{1}_{A_{n}/B_{n}} \otimes_{A_{n}} A'_{n}&\rTto^{\Phi_{n}} &  \Omega^{1}_{A'_{n}/B_{n}}& \rTto& 0\\
&\rdTonto&&\ruTinc& &&&&\\
&&\Ker \Phi_{n}&&&&&&\\
\end{diagram}
where the horizontal sequences are  precisely the Second Fundamental Exact Sequences of  the homomorphisms of discrete rings  $B_{j} \to A_{j} \epi A_{j}/I_{j}$ for $j \in \{n, m\}$ and the vertical maps are surjective. Since $\Ima \delta$ is a finite type $A/I$-module there   results that 
\begin{equation*}
\begin{split}
\invlim {n \in \NN}  (&\frac{I_{n}}{I^{2}_{n}} 	\xto{\delta_{n}} \Omega^{1}_{A_{n}/B_{n}} \otimes_{A_{n}} A'_{n}\xto{\Phi_{n}}   \Omega^{1}_{A'_{n}/B_{n}} \to 0) =\\
=\,\, &\frac{I}{I^{2}} \xto{\delta} \om^{1}_{A/B} \otimes_{A}\frac{A}{I} \xto{\Phi} \om^{1}_{(\frac{A}{I})/B} \to 0
\end{split}
\end{equation*}
is an exact sequence.
\end{proof}

\begin{rem}
Throughout  this work and following the  usual practice in the case of schemes (see  \cite{EGA44}, \cite{at} and \cite{ha1}), we will use the slight abuse of notation of  writing the Second Fundamental  Exact Sequence  associated to a closed immersion $\FX' \overset{i} \to \FX$ given by an Ideal $\CI \subset \CO_{\FX}$ and a morphism $\FX \to \FY$ in $\sfn$, in  any of  the following ways:
\begin{equation*} 
 \CC_{\FX'/\FX} \xto{\delta} i^{*}\om^{1}_{\FX/\FY} \xto{\Phi} \om^{1}_{\FX'/\FY} \to 0
\end{equation*}
\begin{equation*} 
 \CC_{\FX'/\FX} \xto{\delta} \om^{1}_{\FX/\FY} \otimes_{\CO_{\FX}} \CO_{\FX'} \xto{\Phi} \om^{1}_{\FX'/\FY} \to 0
\end{equation*}
\begin{equation*} 
\frac{\CI}{\CI^{2}} \xto{\delta} \om^{1}_{\FX/\FY} \otimes_{\CO_{\FX}} \CO_{\FX'} \xto{\Phi} \om^{1}_{\FX'/\FY} \to 0
\end{equation*}
whichever is convenient for our exposition.
\end{rem}

As it happens in $\sch$ the Second Fundamental Exact Sequence  leads to a local  description of the  module of  differentials of  a pseudo finite type morphism $f:\FX \to \FY$ in $\sfn$. 
\begin{parraf} \label{genermodifptf}
Let  $f: \FX=\spf(A) \to \FY=\spf(B)$ be a morphism in $\sfna$ of  pseudo finite type, then it  factors as  
\[
\FX= \spf(A) \overset{j} \inc \BD^{s}_{\BA^{r}_{\FY}}= \spf(B\{\mathbf{T}\}[[\mathbf{Z}]] ) \xto{p} \FY= \spf(B)
\]
where $r,\, s \in \NN$, $p$ is the canonical projection and $j$ is a  closed immersion given by an Ideal  $\CI= I^{\tr} \subset \CO_{\BD^{s}_{\BA^{r}_{\FY}}}$ with $I= \langle  P_{1},  P_{2}, \ldots,  P_{k}\rangle  \subset B\{\mathbf{T}\}[[\mathbf{Z}]]$. The Second  Fundamental  Exact Sequence (\ref{sef2}) associated to $\FX \overset{j} \inc \BD^{s}_{\BA^{r}_{\FY}} \xto{p} \FY$ 
corresponds through the equivalence  between the category of finite type $A$-modules and $\coh(\FX)$ \ref{propitrian}.(\ref{equivtrian}), to  the sequence 
\begin{equation} \label{segsucfunafdis}
\frac{I}{I^{2}} \xto{\delta} \om^{1}_{B\{\mathbf{T}\}[[\mathbf{Z}]]/B} \otimes_{B\{\mathbf{T}\}[[\mathbf{Z}]]}A\xto{\Phi} \om^{1}_{A/B} \to 0.
\end{equation}
Let us use the following abbreviation $\hd= \hd_{B\{\mathbf{T}\}[[\mathbf{Z}]] /B}$. We have seen in  Example \ref{modifafindis} that $\{ \hd T_{1},\, \hd T_{2}, \ldots,\, \hd  T_{r},\, \hd Z_{1},\, \ldots,\, \hd Z_{s} \}$ is a basis of the free $B\{\mathbf{T}\}[[\mathbf{Z}]] $-module  $\om^{1}_{B\{\mathbf{T}\}[[\mathbf{Z}]] /B}$. Therefore, if $a_{1},\, a_{2}, \ldots,\, a_{r},\, a_{r+1},\, \ldots,\, a_{r+s} $ are the images of  $T_{1},\, T_{2},\, \ldots,\, T_{r},\, Z_{1},\, \ldots,\, Z_{s}$ in $A$, by the definition of  $\Phi$ we have that 
\[\om^{1}_{A/B} = \langle \hd_{A/B} a_{1},\, \hd_{A/B}a_{2}, \ldots,\, \hd_{A/B}a_{r},\, \hd_{A/B}a_{r+1},\, \ldots,\, \hd_{A/B}a_{r+s}   \rangle\]
and from the exactness of  (\ref{segsucfunafdis}) it holds that 
\[
\om^{1}_{A/B} \cong \frac{\om^{1}_{B\{\mathbf{T}\}[[\mathbf{Z}]]/B}\otimes_{B\{\mathbf{T}\}[[\mathbf{Z}]]} A}{\langle \hd P_{1}\otimes1, \hd P_{2}\otimes1, \ldots, \hd P_{k}\otimes1\rangle}
\]
or, equivalently, since the  functor $(-)^\tr$ is exact, 
\[\om^{1}_{\FX/\FY} \cong \frac{\om^{1}_{\BD^{s}_{\BA^{r}_{\FY}}/\FY} \otimes_{\CO_{\BD^{s}_{\BA^{r}_{\FY}}}} \CO_{\FX}}{\langle \hd P_{1}\otimes1, \hd P_{2}\otimes1, \ldots, \hd P_{k}\otimes1\rangle^{\tr}}.\]
\end{parraf}

\section{Definitions of the infinitesimal lifting properties}\label{sec21}

In this section we extend Grothendieck's classical definition of infinitesimal lifting properties in the category  of  schemes (\emph{cf.}   \cite[(17.1.1)]{EGA44}) to the category of locally noetherian formal schemes and we present some of their basic properties. 

\begin{defn}\label{defcinf}
Let $f:\FX \to \FY$ be a morphism in $\sfn$. We say that $f$ is \emph{formally smooth (formally unramified or formally \'etale)} if it satisfies the following lifting condition:

\emph{For all affine $\FY$-scheme $Z$ and for each closed subscheme $T\inc Z$ given by a square zero Ideal $\CI \subset \CO_{Z}$ the induced map
\begin{equation} \label{condlevant1}
\Hom_{\FY}(Z,\FX) \lto \Hom_{\FY}(T,\FX)
\end{equation}
is surjective (injective  or bijective, respectively).}

So we have that $f$ is formally \'etale if, and only if,  is formally smooth and formally unramified.

We stress  that in the lifting conditions (\ref{condlevant1}) the schemes $T$ and $Z$ are not necessarily  locally noetherian. 
These definitions agree with the classical definitions of  the infinitesimal conditions given in \cite[(17.3.1)]{EGA44} for a morphism  $f$ of ordinary schemes.
\end{defn}

\begin{rem}
If the lifting  condition  (\ref{condlevant1}) holds for a square zero Ideal $\CI \subset \CO_{Z}$, then it also holds true for nilpotent Ideals. Indeed, let $Z$ be an affine $\FY$-scheme  and let $T \subset Z$ be a closed subscheme given by an Ideal $\CI \subset \CO_{Z}$  such that  $\CI^{n+1} =0$, for any $n \in \NN$. If we call $T_{k}$ the closed subscheme of  $Z$ given by the  Ideal $\CI^{k+1} \subset \CO_{Z}$, we have that for all $0 \le k \le n$, the canonical morphisms $T_{k} \inc T_{k+1}$ are  closed immersions given by a square zero Ideal. Then, the map $\Hom_{\FY}(T_{k+1},\FX) \lto \Hom_{\FY}(T_{k},\FX)$ is surjective (injective or bijective), for all $0 \le k \le n$ and therefore, $\Hom_{\FY}(Z,\FX) \lto \Hom_{\FY}(T,\FX)$ is  surjective (injective or bijective, respectively).
\end{rem}

Given $f:\FX \to \FY$ in $\sfn$, in the next proposition it is proved that the lifting  condition  (\ref{condlevant1})  extends to  every affine noetherian $\FY$-formal scheme $\FZ$ and every closed formal subscheme $\FT \inc \FZ$ given by a square zero Ideal $\CI \subset \CO_{\FZ}$.

\begin{propo} \label{levantform}
Let $f:\FX \to \FY$ be in $\sfn$. If the  morphism $f$ is formally smooth (formally unramified or formally \'etale), then for all affine noetherian $\FY$-formal scheme $\FZ$ and for all closed formal subschemes $\FT \inc \FZ$ given by a square zero Ideal $\CI \subset \CO_{\FZ}$, the induced map
\begin{equation} \label{condlevant3}
\Hom_{\FY}(\FZ,\FX) \lto \Hom_{\FY}(\FT,\FX)
\end{equation}
is surjective (injective or bijective, respectively). 
\end{propo}

\begin{proof}
Let   $\FT=\spf(C/I) \inc \FZ=\spf(C)$ be a closed formal subscheme given by a square zero ideal $I \subset C$ . Let $L \subset C$ be an ideal of  definition, then the  embedding  $\FT \inc \FZ$ is expressed as 
\[\dirlim {n \in \NN} (T_{n}=\spec(\frac{C}{I+L^{n+1}}) \overset{j_{n}}\inc Z_{n}=\spec(\frac{C}{L^{n+1}}))\] 
where the morphisms $j_{n}$ are  closed immersions of  affine schemes defined by a square zero Ideal. Given $u: \FT \to \FX$ a $\FY$-morphism,  we will denote by $u'_{n}$  the morphisms $T_{n} \inc \FT \xto{u} \FX$ such that  the diagrams
\begin{diagram}[height=2em,w=2em,labelstyle=\scriptstyle]
T_{n}		&	\rTinc^{j_{n}}	&	Z_{n}\\	
\dTto	^{u'_{n}}		&		&	\dTto\\
\FX			&	\rTto^f 	&	\FY\\
\end{diagram} 
commute for all $n \in \NN$. 

Suppose that $f$ is formally smooth, then we will construct a collection of  $\FY$-morphisms $v'_{n}: Z_{n} \to \FX$ such that  $v'_{n}|_{T_{n}}=u'_{n}$ and $v'_{n}|_{Z_{n-1}}=v'_{n-1}$. The $\FY$-morphism 
\[v :=\dirlim {n \in \NN} (v'_{n}: Z_{n} \to \FX)\] 
will satisfy $v|_{\FT}=u$. Indeed, since $f$ is formally smooth, there exists  a $\FY$-morphism $v'_{0}: Z_{0} \to \FX$ such that  $v'_{0}|_{T_{0}}=u'_{0}$. By induction on $n$,  suppose that we have constructed  $\FY$-morphisms $v'_{n}: Z_{n} \to \FX$ such that $v'_{n}|_{T_{n}}=u'_{n}$ and $v'_{n}|_{Z_{n-1}}=v'_{n-1}$. Let us define $v'_{n+1}$. The morphisms $u'_{n+1}$ and $v'_{n}$ induce a $\FY$-morphism  
\[T_{n+1}\coprod_{T_{n}}  Z_{n}=\spec(\frac{C}{I+L^{n+2}}\times_{\frac{C}{I+L^{n+1}}}\frac{C}{L^{n+1}} ) \xto{w'_{n+1}} \FX\]
such that the diagram, where the composition of the middle horizontal maps is $j_{n+1}$,
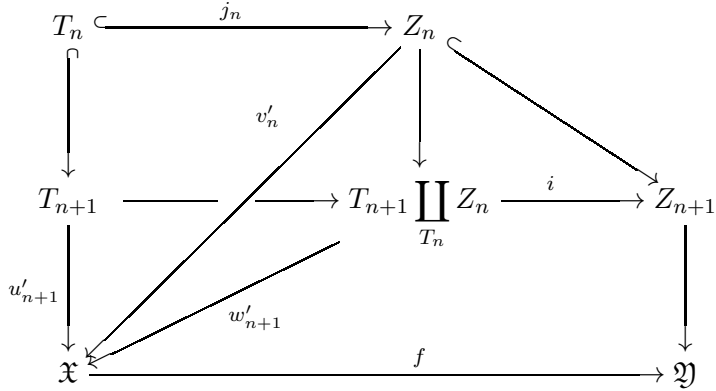
\begin{diagram}[height=3em,w=2.7em,labelstyle=\scriptstyle]
T_{n}   & &\rTinc^{j_{n}}  && Z_{n}   &  &&    \\
\dTinc  & &\qquad\,\,\,\scriptstyle{v'_{n}}&\ldTto(4,4)  & \dTto & \rdTinc(3,2) &&    \\
T_{n+1} &\,\,\,\hLine  & \phantom{\cdot} & \,\,\rTto & T_{n+1} \coprod_{T_{n}} Z_{n} & \rTto^{\phantom{abcd} i} && Z_{n+1} \\
\dTto^{u'_{n+1}}   && & \ldTto(4,2)_{w'_{n+1}} &  &  &&   \dTto \\
\FX     & &&       &  \rTto^{f} &  &&   \FY \\ 
\end{diagram}
is commutative.
Then, since $f$ is formally smooth and $i$ a  closed immersion 
defined by a square zero Ideal, $w'_{n+1}$ lifts to a $\FY$-morphism $v'_{n+1}: Z_{n+1} \to \FX$  that satisfies $v'_{n+1}|_{T_{n+1}}=u'_{n+1}$ and $v'_{n+1}|_{Z_{n}}=w'_{n+1}|_{Z_{n}}=v'_{n}$.

If $f$ is formally unramified, assume there exist $\FY$-morphisms $v: \FZ \to \FX$ and $w: \FZ \to \FX$ such that  $v|_{\FT} = w|_{\FT} = u$ and let us show that $v = w$. With  the notations established at the  beginning of  the proof  consider  
\[v = \dirlim {n \in \NN} (v'_{n}: Z_{n} \to \FX) \,\,\text{ and } \,\, w= \dirlim {n \in \NN} (w'_{n}: Z_{n} \to \FX)\]  such that the   diagram
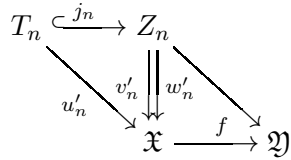
\begin{diagram}[height=2em,w=2em,p=0.3em,labelstyle=\scriptstyle]
T_{n} &   \rTinc^{j_{n}}	    &  Z_{n}				       &  	   	         &      &\\
   & \rdTto_{u'_{n}}&\dTto^{v'_{n}}\dTto_{w'_{n}}& \rdTto 		&     &\\ 
   &			    &\FX	  			       & 	 \rTto^{f} & \FY&\\  
\end{diagram}
commutes.
By hypothesis we have that $v'_{n} = w'_{n}$ for all $n \in \NN$ and we conclude that \[v=\dirlim {n \in \NN}v'_{n}= \dirlim {n \in \NN}  w'_{n}=w.\]
\end{proof}

\begin{parraf} \label{lgfefnr}
In definition \ref{defcinf} the test morphisms for the lifting condition are closed subschemes of \emph{affine} $\FY$-schemes given by square-zero ideals. An easy patching argument gives that the uniqueness of lifting conditions holds for closed subschemes of \emph{arbitrary} $\FY$-schemes given by square-zero ideals. This applies to formally unramified and formally \'etale morphisms.
\end{parraf}

\begin{cor} \label{corlevfor}  
Let $f:\FX \to \FY$ be in $\sfn$. If the  morphism $f$ is formally unramified  (or formally \'etale), then for all noetherian $\FY$-formal schemes $\FZ$ and for each closed formal subscheme $\FT \inc \FZ$  given by a square zero Ideal $\CI \subset \CO_{\FZ}$, the induced map 
\begin{equation} \label{condlevant4}
\Hom_{\FY}(\FZ,\FX) \lto \Hom_{\FY}(\FT,\FX)
\end{equation}
is injective (or bijective, respectively).
\end{cor}

\begin{proof}
Given $\{\FV_{\alpha}\}$ a covering of affine open formal subschemes of  $\FZ$, denote by $\{\FU_{\alpha}\}$ the  covering of affine open formal subschemes of  $\FT$ given by $\FU_{\alpha} = \FV_{\alpha} \cap \FT$, for all $\alpha$.  By \cite[(10.14.4)]{EGA1} $\FU_{\alpha} \inc \FV_{\alpha}$ is a  closed immersion in $\sfn$ determined by a square zero Ideal. Given $u: \FT \to \FX$ a $\FY$-morphism, put $u_{\alpha}= u|_{\FU_{\alpha}}$, for all $\alpha$. If $f$ is formally smooth (formally unramified or formally \'etale) by Proposition \ref{levantform} for all $\alpha$ there exists (there exists at most or there exists a unique, respectively) $\FY$-morphism $v_{\alpha}: \FV_{\alpha} \to \FX$ such that  $v_{\alpha}|_{\FU_{\alpha}} = u_{\alpha}$. In the case  formally unramified, (by uniqueness) glueing the $v_{\alpha}$ we obtain that  there exists at most a $\FY$-morphism $v: \FZ \to \FX$ such that  $v|_{\FT}=u$. And in the  case formally  \'etale (by existence and uniqueness) glueing  the $v_{\alpha}$ we get that there exists an unique  $\FY$-morphism $v: \FZ \to \FX$ such that  $v|_{\FT}=u$.
\end{proof}


\begin{parraf}  \label{condinfalg}
Let $f: \spf(A) \to \spf(B)$  be in $\sfna$. Applying formal Descartes duality (\ref{equiv}), we have that $f$ is formally smooth (formally unramified or formally \'etale) if, and only if,   the topological $B$-algebra $A$ satisfies the following condition of  lifting:

\emph{For each (discrete) ring $C$ and for all continuous homomorphisms $B \to C$, given a square zero ideal $I \subset C$, the induced map
\begin{equation} \label{condlevant2}
\Homcont_{B\alg}(A,C) \to \Homcont_{B\alg}(A,C/I)
\end{equation}
is surjective (injective or bijective, respectively).}
A topological $B$-algebra  $A$ that satisfies the condition (\ref{condlevant2}) is called\emph{ formally smooth (formally unramified or formally \'etale, respectively)} (\emph{cf.}  \cite[(\textbf{0}.19.3.1) and (\textbf{0}.19.10.2)]{EGA41}).
\end{parraf}

A basic reference for the fundamental  properties of  the infinitesimal lifting conditions for preadic rings is   \cite[\S \textbf{0}.19.3 and \S \textbf{0}.19.10]{EGA41}. Let us recall now some of these properties. We will also look at the first basic examples of  morphisms in $\sfna$  that satisfy the infinitesimal conditions (Example \ref{exconforinf}).

\begin{rem} \label{confinfhafin}
Let $B \to A$ be a continuous  morphism of  preadic rings and take $J \subset A,\, K \subset B$ ideals of  definition with $KA \subset J$. Given $J' \subset A,\, K' \subset B$ ideals such that  $K'A \subset J'$, $J \subset J'$ and $K \subset K'$,  if $A$ is a  formally smooth (formally unramified or formally \'etale) $B$-algebra for the  $J$ and $K$-adic topologies then,  we have  that $A$ is a  formally smooth (formally unramified or formally \'etale, respectively) $B$-algebra for the $J'$ and $K'$-adic topologies, respectively.
\end{rem}

\begin{lem} \label{condinfalgcompl}
Let $B \to A$ be a continuous  morphism of  preadic rings, $J \subset A$ and $K \subset B$ ideals of  definition with $KA \subset J$ and we denote by $\HA$ and $\HB$ the respective completions of  $A$ and $B$. The following conditions are equivalent:
\begin{enumerate}
\item
$A$ is a  formally smooth (formally unramified or formally \'etale) $B$-algebra
\item
$\HA$ is a  formally smooth (formally unramified or formally \'etale, respectively) $B$-algebra
\item
$\HA$ is a  formally smooth (formally unramified or formally \'etale, respectively) $\HB$-algebra
\end{enumerate} 
\end{lem}

\begin{proof}
It suffices to notice that for all discrete rings $C$ and for all continuous  homomorphisms $B \to C$, the following equivalences
\[
\Homcont_{B\alg}(A,C) \equiv \Homcont_{B\alg}(\HA,C) \equiv \Homcont_{\HB\alg}(\HA,C) 
\] hold.
\end{proof}

\begin{ex} \label{exconforinf}
Put $\FX= \spf(A)$ with $A$ a $J$-adic  noetherian ring and let $\mathbf{T} = T_{1},\, T_{2},\,  \ldots,\,  T_{r}$ be a   finite number of indeterminates.
\begin{enumerate}
\item \label{restrformliso}
If we take in $A$ the discrete topology, from  the universal property of the polynomial ring it follows that  $A[\mathbf{T}]$ is a formally smooth $A$-algebra. Applying the previous Remark and  Lemma \ref{condinfalgcompl} we have that the  restricted formal series ring  $A\{\mathbf{T}\}$ is a formally smooth $A$-algebra or, equivalently, the  canonical morphism $\BA_{\FX}^{r} \to \FX$ is  formally smooth.

\item  \label{formformliso}
Analogously to the preceding example, we obtain that $A[[\mathbf{T}]]$ is a  formally smooth $A$-algebra, from which  we deduce that projection $\BD_{\FX}^{r} \to \FX$ is formally smooth.

\item  \label{allavesformetale}
If  we take in $A$ the discrete topology it is known that, given $f \in A$, $A_{f}$ is a  formally \'etale $A$-algebra. So, there results that the canonical inclusion $\fD(f)\inc \FX$ is  formally \'etale.

\item  \label{enccerradoformnoram}
Trivially every surjective  morphism of  rings is formally unramified. Therefore, given an ideal $I \subset A$,  the   closed immersion $\spf(A/I) \inc \FX$ is  formally unramified.

\item  \label{compformetale}
If  $\FX'= \spf(A/I)$ is a closed formal subscheme of  $\FX$, in \ref{compl} we have shown that $\kappa: \FX_{/ \FX'} \to \FX$, the  morphism of completion of  $\FX$ along  $\FX'$, corresponds through formal Descartes duality (\ref{equiv}) with the continuous morphism of  rings $A \to \HA$, where $\HA$ is the  complete ring of  $A$ for the $(I+J)$-adic topology  and therefore, $\kappa$ is formally \'etale.
\end{enumerate}
\end{ex}

In the case of usual schemes the study  of the infinitesimal properties using the  module of  differentials leads one to consider the class of finite type morphisms. This hypothesis permits to give precise characterizations of the infinitesimal lifting conditions. In $\sfn$ there are two conditions on morphisms that generalize the property of being of finite type for morphisms  in $\sch$: morphisms of pseudo finite type and their adic counterpart, morphisms of  finite type (Definition \ref{defmtf}). Henceforward, we will focus on morphisms with these properties.

\begin{defn}
Let $f:\FX \to \FY$ be in $\sfn$. The  morphism $f$ is \emph{smooth (unramified or \'etale)}  if, and only if, it is of  pseudo finite type and formally smooth (formally unramified or formally \'etale, respectively). If moreover $f$ is adic, we say that $f$ is \emph{adic smooth (adic unramified or adic \'etale, respectively)}. So $f$ is  adic smooth (adic unramified or adic \'etale) if it is of  finite type and formally smooth (formally unramified or formally \'etale, respectively).
\end{defn}

If $f: X \to Y$ is in $\sch$, both definitions coincide with the one given in \cite[(17.3.1)]{EGA44} and we say  that $f$ is smooth (unramified or \'etale, respectively).

\begin{ex} \label{confinftipof}
Put $\FX= \spf(A)$ with $A$ a $J$-adic noetherian ring and $r \in \NN$. From  Example \ref{exconforinf} we have that 
\begin{enumerate}
\item \label{restrliso}
The canonical morphism $\BA_{\FX}^{r} \to \FX$ is  adic smooth.

\item  \label{formpseliso}
The projection $\BD_{\FX}^{r} \to \FX$ is smooth.

\item  \label{allavesetale}
The canonical  inclusion $\fD(f) \inc \FX$ is  adic \'etale.

\item  \label{enccerradonoram}
Given an ideal $I \subset A$,  the   closed immersion $\spf(A/I) \inc \FX$ is   adic unramified.

\item  \label{comppseetale}
If  $\FX'= \spf(A/I)$ is a closed formal subscheme of  $\FX$, the morphism of  completion $\kappa: \FX_{/ \FX'} \to \FX$ is \'etale.
\end{enumerate}
Notice that in the category of usual schemes, $\sch$, there does not exist non trivial analogues to examples (\ref{formpseliso}) and (\ref{comppseetale}). So, one should expect that the structure of morphisms satisfying infinitesimal lifting properties in $\sfn$ is not a mere extension of the ordinary scheme case.
\end{ex}
From Proposition \ref{infsorit} to  Corollary \ref{corcondinflocal} we study  the basic results about the infinitesimal lifting properties for  morphisms of locally noetherian  formal schemes. Despite their proofs are similar to the ones of  the  analogous results in $\sch$, they do not follow from those. We will prove them with the necessary detail.

\begin{propo} \label{infsorit}
In the category $\sfn$ of  locally noetherian formal schemes  the following properties hold:
\begin{enumerate}
\item  \label{infsorit1}
Composition of smooth (unramified or \'etale)  morphisms   is a smooth  (unramified or \'etale, respectively) morphism.
\item  \label{infsorit2}
Smooth, unramified and \'etale character  is stable under base-change in $\sfn$.
\item  \label{infsorit3}
Product of  smooth (unramified or \'etale) morphisms is a smooth   (unramified or \'etale, respectively) morphism.
\end{enumerate}
\end{propo}

\begin{proof}
By Proposition \ref{sori} it suffices prove (1) and (2). 
In order to prove (1) let $f:\FX \to \FY$ and $g:\FY \to \FS$ be smooth morphisms. Let $Z$ be  an affine  $\FS$-scheme  and $T\inc Z$ a closed subscheme defined by a square zero Ideal $\CI \subset \CO_{Z}$. Given $u: T \to \FX$ a $\FS$-morphism, since $g$ is smooth  the  morphism $f \circ u$ lifts to a $\FS$-morphism $Z  \to \FY$ such that the  following diagram commutes:
\begin{diagram}[height=2.1em,w=2.1em,labelstyle=\scriptstyle]
T    	&        \rTinc  &        Z	  &  &  \\
\dTto^{u}	&          &  \dTto^{w} &		\rdTto(2,2)&   \\
\FX   &   \rTto^f  &   	\FY&	  \rTto^g        & \FS   \\ 
\end{diagram}
And since $f$ is smooth,  there exists a $\FY$-morphism $v: Z \to \FX$ such that  $v|_{T} = u$, from which we deduce that $g \circ f$ is formally  smooth. Applying Proposition \ref{mtf}.(\ref{mtf1}) we deduce that $g \circ f$ is smooth. The  unramified  case is proved analogously.

Now, let us show that given a smooth morphism $f:\FX \to \FY$ and another morphism $\FY' \to \FY$  in $\sfn$, then the morphism $f': \FX \times_{\FS} \FY \to \FY'$ is smooth. Let $Z$ be an affine  $\FY'$-scheme, $T\inc Z$ a closed subscheme defined by a square zero Ideal $\CI \subset \CO_{Z}$ and $u: T \to \FX \times_{\FY} \FY'$ a $\FY'$-morphism. Since $f$ is smooth, we have that the  $\FY$-morphism  $T \xto{u} \FX \times_{\FY} \FY' \to \FX$ lifts to a $\FY$-morphism $Z \to \FX$ that makes the  following diagram commutative
\begin{diagram}[height=2.4em,w=2.1em,labelstyle=\scriptstyle]
T                     &              & \rTinc & & Z\\
\dTto_{u}             &              &  &  \ldTto(4,4) & \dTto \\
\FX \times_{\FY} \FY' & \hLine\,\, & \phantom{\cdot} & \rTto_{f'} & \FY' \\
\dTto                 &              &  & & \dTto \\
\FX                   &              & \rTto^f & & \FY.\\
\end{diagram}
By the universal property of  the fiber product there exists a $\FY'$-morphism $Z \to \FX \times_{\FY} \FY'$ that agrees with $u$ in $T$, so $f'$ is formally smooth. Besides, applying Proposition \ref{mtf}.(\ref{mtf3}) there results that $f'$ is of  pseudo finite type and therefore, $f'$ is smooth. The proof of the unramified case is similar and we leave it to the patient reader.
\end{proof}

\begin{propo}  \label{infsoritcor}
The assertions of the last Proposition hold if  we change the infinitesimal conditions  by the corresponding infinitesimal adic conditions.
 \end{propo}
 
\begin{proof}
By the definition of  the infinitesimal adic conditions, it suffices to  apply the  last result and the sorites of adic morphisms in Proposition \ref{adic}.
\end{proof}

\begin{ex}
Let $\FX$ be in $\sfn$ and $r \in \NN$. From Proposition \ref{infsorit}.(\ref{infsorit2}),  Example  \ref{confinftipof}.(\ref{restrliso}) and    Example \ref{confinftipof}.(\ref{formpseliso})
we get that:
\begin{enumerate}
\item
The morphism of  projection $\BA_{\FX}^{r}= \FX \times_{\spec(\ZZ)} \BA_{\spec(\ZZ)}^{r} \to \FX$ is an adic  smooth morphism.
\item
The canonical  morphism $\BD_{\FX}^{r}=\FX \times_{\spec(\ZZ)} \BD_{\spec(\ZZ)}^{r} \to \FX$ is smooth.
\end{enumerate}
\end{ex}

\begin{propo}\label{infencajes}
The following holds in  the category of  locally noetherian formal schemes:
\begin{enumerate}
\item \label{infencajes1}
A  closed immersion is adic unramified.
\item \label{infencajes2} 
An open immersion is adic \'etale.
\end{enumerate}
\end{propo}

\begin{proof}
Closed and open immersions are monomorphisms and therefore, unramified. On the other hand, open immersions are smooth morphisms.
\end{proof}

\begin{propo} \label{propibasinf}
Let $f: \FX \to \FY ,\, g: \FY \to \FS$ be two morphisms of  pseudo finite type in $\sfn$. 
\begin{enumerate}
\item
If $g \circ f$ is unramified, then so is $f$. 
\item
Let us suppose that $g$ is unramified. If $g \circ f$ is smooth   (or \'etale) then, $f$ is smooth   (or \'etale, respectively).
\end{enumerate}
\end{propo}

\begin{proof}
Item (1) is immediate. In order to prove (2) we consider $Z$ a $\FY$-scheme, $T \inc Z$ a $\FY$-closed subscheme given by a square zero Ideal and $u: T \to \FX$ a $\FY$-morphism. If $g \circ f$ is smooth  there exists a $\FS$-morphism $v: Z \to \FX$  such that  $v|_{T}= u$ and since $g$ is unramified, we have that $u$ is a $\FY$-morphism. Moreover, by  Corollary \ref{cor2sor} $f$ is of  pseudo finite type and therefore, $f$ is smooth. The  \'etale case is  deduced from (1) and from the last  argument.
\end{proof}

\begin{cor} \label{corbasinf}
Let $f:\FX \to \FY$ be a  pseudo finite type morphism and  $g: \FY \to \FS$  an \'etale morphism. The  morphism $g \circ f$ is smooth  (or \'etale) if, and only if,  $f$ is smooth  (or \'etale, respectively).
\end{cor}

\begin{proof}
It is enough to apply  Proposition \ref{infsorit}.(\ref{infsorit1}) and Proposition \ref{propibasinf}. 
\end{proof}

\begin{parraf} (\emph{cf.}  \cite[(\textbf{0}.20.1.1)]{EGA44}) \label{derivversuslevantalg}
Let $B \to A$ be a continuous homomorphism of  preadic rings, $C$ a topological $B$-algebra and $I\subset C$ a square zero ideal. Given $u: A \to C/I$ a continuous $B$-homomorphism of  rings it holds that:
\begin{enumerate}
\item \label{derivversuslevantalg1}
If $v: A \to C,\, w: A \to C$ are two continuous $B$-homomorphisms of  rings such that agree with $u$ modulo $I$, then  $v-w \in \Dercont_{B}(A,I)$.
\item \label{derivversuslevantalg2}
 Fixed $v: A \to C$ a continuous $B$-homomorphism of  rings that  agree with $u$ modulo $I$ and given $d \in \Dercont_{B}(A,I)$ we have that $v+d: A \to C$ is a continuous $B$-homomorphism of  rings that  agrees with $u$ modulo $I$.
\end{enumerate}
\end{parraf}

\begin{parraf} \label{parrfderivlevantcap2}
Given $f:\FX \to \FY$ a  pseudo finite type morphism in $\sfn$, let $Z$  be an affine $\FY$-scheme, $i\colon T \inc Z$ a closed $\FY$-subscheme given by a square zero Ideal $\CI \subset \CO_{Z}$ and $u: T \to \FX$ a $\FY$-morphism. Observe that $\CI=i^{-1}\CI$ has structure of  $\CO_{T}$-Module. From  \ref{derivversuslevantalg} we deduce the following:
\begin{enumerate} 
\item \label{parrfderivlevantcap21}
If $v: Z \to \FX,\, w: Z \to \FX$ are $\FY$-morphisms such that the  diagram
\begin{diagram}[height=2em,w=2em,p=0.3em,labelstyle=\scriptstyle]
T &   \rTinc^{i}	    &  Z				       &  	   	         &      &\\
   & \rdTto_{u}&\dTto^{v}\dTto_{w}& \rdTto 		&     &\\ 
   &			    &\FX	  			       & 	 \rTto^{f} & \FY&\\  
\end{diagram}
is commutative, then  the  morphism
\[
 \CO_{\FX} \xto{v^{\sharp} - w^{\sharp}} u_{*} \CI
\]
is a continuous $\FY$-derivation. 
Lemma  \ref{imagdirecompl} and  Theorem \ref{modrepr}  imply that there exists an unique   morphism of  $\CO_{\FX}$-Modules $\phi: \om^{1}_{\FX/\FY} \to u_{*} \CI$  that makes the  diagram
\begin{diagram}[height=2em,w=2em,p=0.3em,labelstyle=\scriptstyle]
 \CO_{\FX}	&   \rTto^{\hd_{\FX/\FY}}	    & \om^{1}_{\FX/\FY}	\\
\dTto^{v^{\sharp} - w^{\sharp}} &           \ldTto_{\phi}  	    &\\ 
u_{*} \CI	  	&  	    & \\  
\end{diagram} 
commutative.
\item \label{parrfderivlevantcap22}
Fixed a lifting $v: Z \to \FX$ of  $u$ and given $\phi: \om^{1}_{\FX/\FY} \to u_{*} \CI$ a morphism of  $\CO_{\FX}$-Modules, it holds that $w^{\sharp} := v^{\sharp}+ \phi\circ \hd_{\FX/\FY} $ defines another $\FY$-morphism $w: Z \to \FX$ such that  the  diagram 
\begin{diagram}[height=2em,w=2em,p=0.3em,labelstyle=\scriptstyle]
T &  \rTinc^{i}	    &  Z				       &  	   	         &      &\\
   & \rdTto_{u}&\dTdash_{w}& \rdTto 		&     &\\ 
   &			    &\FX	  			       & 	 \rTto^{f} & \FY&\\  
\end{diagram}
commutes.
\end{enumerate}
\end{parraf}

\begin{rem}
Henceforward  and until the end of  this  chapter we will use the  \v Cech  calculus of  cohomology\footnote{Given $X$ a topological space, $\FU$ an open  covering of  $X$ and 
$\CF$ a sheaf  of abelian groups  over $X$, the  $p$-th group of
 \v Cech cohomology of  $\CF$ with respect to the  covering $\FU_{\bullet}$  is  
\[\check{H}^{p} (\FU_{\bullet},\CF) := H^{p}(\check{\C}^{\cdot}(\FU_{\bullet},\CF))=\frac{\check{\Z}^{p} (\FU_{\bullet}, \CF)}{\check{\B}^{p} (\FU_{\bullet}, \CF)}\]
where
$\check{\C}^{\cdot}(\FU_{\bullet}, \CF)$ is the \v Cech complex of  $\CF$ (see, for instance, \cite[Chapter III, p. 219- 220]{ha1} or \cite[\S 5.4]{te}).}. Therefore and whenever it be necessary we'll impose the separation  hypothesis  (each time that we use  cohomology of  degree grater or equal than $2$).
\end{rem}

\begin{propo} \label{condinflocal}
Let  $f:\FX \to \FY$ a morphism in $\sfn$.
\begin{enumerate}
\item
Given  $\{\FU_{\alpha}\}_{\alpha \in L}$ an open covering of  $\FX$,  $ f$  is smooth  (unramified or \'etale)  if, and only if,  for all $\alpha \in L$, $f |_{\FU_{\alpha}}: \FU_{\alpha} \to \FY$ is smooth  (unramified or \'etale, respectively).
\item
If $\{\FV_{\lambda}\}_{\lambda\in J}$ is an open  covering of  $\FY$, $ f$  is smooth  (unramified or \'etale)  if, and only if,  for all $\lambda \in J$, $f |_{f^{-1}(\FV_{\lambda})}: f^{-1}(\FV_{\lambda} ) \to \FV_{\lambda} $ is smooth  (unramified or \'etale, respectively).
\end{enumerate} 
\end{propo}

\begin{proof}
\item
(1) If $f$ is smooth  (unramified or \'etale)  by Proposition  \ref{infencajes}.(\ref{infencajes2}) and Proposition \ref{infsorit}.(\ref{infsorit1}) we have that $f_{\alpha} := f |_{\FU_{\alpha}}$ is smooth  (unramified or \'etale, respectively), $\forall \alpha \in L$. 

In order to prove  the  converse let us establish the following notations: Let $Z$ bet an affine $\FY$-scheme, $T\inc Z$ a closed $\FY$-subscheme given by a square zero Ideal $\CI \subset \CO_{Z}$ of  zero square and $u: T \to \FX$ a $\FY$-morphism. It holds that the collection $\{u^{-1}(\FU_{\alpha})\}_{\alpha \in L}$ is an open  covering of  $T$ and if we denote by  $\{W_{\alpha}\}_{\alpha \in L}$ the corresponding  open  covering  in $Z$, there results  that $u^{-1}(\FU_{\alpha})  \inc W_{\alpha}$  is a closed subscheme defined by the square zero Ideal $\CI_{\alpha}:=\CI|_{W_{\alpha}} \subset \CO_{W_{\alpha}}$, for all $\alpha \in L$. Observe that $\CI_{\alpha}$ has structure of  $\CO_{u^{-1}(\FU_{\alpha})}$-Module.

Assume that  $f_{\alpha}$ is unramified, $\forall \alpha$ and let $v,w: Z \to \FX$ be $\FY$-morphisms such that  $v|_{T} = w|_{T}$. Then $v|_{u^{-1}(\FU_{\alpha})} = w|_{u^{-1}(\FU_{\alpha})}$ and applying \ref{lgfefnr} we have that $v|_{W_{\alpha}} = w|_{W_{\alpha}}$, $\forall \alpha$, that is,  $v=w$ and therefore, $f$ is unramified. 

Let us treat the case when the $\{f_\alpha\}_{\alpha \in L}$ are smooth. For all $\alpha$ we cover $u^{-1}(\FU_{\alpha})$ by affine open subsets  of  $T$, $\{T_{\alpha_{i}}\}_{\alpha_{i} \in L_{\alpha}}$, and take $\{Z_{\alpha_{i}}\}_{\alpha_{i} \in L_{\alpha}}$ the corresponding open subsets in $Z$. Then, simplifying the notation and allowing repetitions in the  covering of  $\FX$, we may assume that the coverings of  $T$ and of  $Z$ are affine and that are indexed by the same set of indexes $L$.
Since the morphisms $\{f_\alpha\}_{\alpha \in L}$ are smooth, for all $\alpha$, there exists a $\FY$-morphism $v'_{\alpha}: Z_{\alpha} \to \FU_{\alpha}$ such that  the  following diagram is commutative:
\begin{diagram}[height=2.5em,w=3em,labelstyle=\scriptstyle]
T_{\alpha}                        &    	\rTinc	     & Z_{\alpha}  \\
\dTto^{u_{\alpha}:=u|_{T_{\alpha}}}& \ldTto^{v'_{\alpha}} & \dTto \\ 
 \FU_{\alpha}	 			      & \rTto^{f_{\alpha}}              & \FY\\  
\end{diagram}
For all $\alpha$, we denote by $v_{\alpha}$ the  $\FY$-morphism $Z_{\alpha} \xto{v'_{\alpha}} \FU_{\alpha} \inc  \FX$. From  $(v_{\alpha})$ we are going to  construct a collection of  $\FY$-morphisms $(w_{\alpha}: Z_{\alpha} \to \FX)$ that glue and that agree with $u_{\alpha}$ in $T_{\alpha}$,  for all $\alpha \in L$. 

In order to do that, we establish the notation
\[\FU_{\alpha  \beta}:= \FU_{\alpha} \cap \FU_{\beta},\, Z_{\alpha \beta}:=Z_{\alpha} \cap Z_{\beta},\, T_{\alpha \beta}:=T_{\alpha} \cap T_{\beta},\, u_{\alpha \beta}:= u|_{T_{\alpha \beta}}\]
for all  couple of  indexes $\alpha,\, \beta$ such that $Z_{\alpha \beta} \neq \varnothing$. 
Observe that $T_{\alpha \beta} \overset{i_{\alpha \beta}} \inc Z_{\alpha \beta}$ is a  closed immersion of  affine schemes given by the square zero Ideal $\CI_{\alpha \beta} \subset \CO_{Z_{\alpha \beta}}$, therefore, let us write $i_{\alpha \beta}^{*} (\CI_{\alpha \beta}) = \CI_{\alpha \beta}$. Fixed $\alpha, \, \beta$ such that $T_{\alpha \beta} \neq \varnothing$, since $v_{\alpha}$ and $v_{\beta}$ agree in $T_{\alpha \beta}$ by \ref{parrfderivlevantcap2}.(\ref{parrfderivlevantcap21}) it holds that  there exists a unique  morphism of  $\CO_{\FX}$-Modules $\phi_{\alpha \beta}: \om^{1}_{\FX/\FY} \to (u_{\alpha \beta})_{*} (\CI_{\alpha \beta})$ such that  the  following diagram is commutative
\begin{diagram}[height=2em,w=2em,p=0.3em,labelstyle=\scriptstyle]
 \CO_{\FX}	&   \rTto^{\hd_{\FX/\FY}\qquad}	    & \om^{1}_{\FX/\FY}	\\
\dTto^{(v_{\alpha}|_{Z_{\alpha\beta}})^{\sharp} -(v_{\beta}|_{Z_{\alpha\beta}})^{\sharp}} &  \ldTto_{\phi_{\alpha \beta}}   &\\ 
(u_{\alpha \beta})_{*} (\CI_{\alpha \beta})	  	&  	    & \\  
\end{diagram}
Let $u_{\alpha \beta}^{*} (\om^{1}_{\FX/\FY}) \to \CI_{\alpha \beta}$ be the  morphism  of  $\CO_{T_{\alpha \beta}}$-Modules adjoint of  $\phi_{\alpha \beta}$, that we will continue to denote by  $\phi_{\alpha \beta}$ making a little abuse  of  notation. The family of  morphisms $\{\phi_{\alpha \beta}\}_{\alpha \beta}$ satisfy the cocycle condition; that is, for  any $\alpha,\, \beta,\, \gamma$ such that $Z_{\alpha \beta \gamma}:= Z_{\alpha } \cap Z_{ \beta } \cap Z_{ \gamma} \neq \varnothing$, we have that
\begin{equation} \label{datosrecoleccap21}
\phi_{\alpha \beta}|_{T_{\alpha \beta \gamma}} - \phi_{\alpha \gamma}|_{T_{\alpha \beta \gamma}}+\phi_{ \beta \gamma}|_{T_{\alpha \beta \gamma}}=0
\end{equation} 
where $T_{\alpha \beta \gamma}:= T_{\alpha } \cap T_{ \beta } \cap T_{ \gamma}$.
Therefore, it defines an element 
\[
\begin{split}
[(\phi_{\alpha \beta})] &\in \check{\h}^{1} (T, \shom_{\CO_{T}}(u^{*} \om^{1}_{\FX/\FY} ,\CI)) = \qquad \textrm{\cite[(5.4.15)]{te}} \\ 
&= \h^{1} (T, \shom_{\CO_{T}}(u^{*} \om^{1}_{\FX/\FY} ,\CI)).
\end{split}
\]
 But,  since $T$ is an affine scheme $\h^{1}(T, \shom_{\CO_{T}}(u^{*} \om^{1}_{\FX/\FY} ,\CI))=0$ and therefore, there exists $\{\phi_{\alpha}\}_{\alpha \in L} \in \check{\C}^{0}(T,\shom_{\CO_{T}}(u^{*} \om^{1}_{\FX/\FY} ,\CI))$ such that, for all couple of indexes $\alpha, \beta$ with $Z_{\alpha \beta} \neq \varnothing$, 
\begin{equation} \label{datosrecoleccap22}
\phi_{\alpha}|_{T_{\alpha \beta}} - \phi_{\beta}|_{T_{\alpha \beta}}=\phi_{\alpha \beta}
\end{equation} 
We again commit an  abuse of  notation and use indifferently  $\phi_{\alpha}$ for the map and its adjoint. For all $\alpha$, let $w_{\alpha}: Z_{\alpha} \to \FX$ be the  morphism that agrees with $v_{\alpha}$ as a map of  topological spaces given as a morphism of  ringed topological spaces  by  $w^{\sharp}_{\alpha}:= v_{\alpha}^{\sharp} - \phi_{\alpha} \circ \hd_{\FX/\FY}$. Therefore, from \ref{parrfderivlevant}.(\ref{parrfderivlevant2}) we have that $w_{\alpha}$ is a lifting of  $u_{\alpha}$ over $\FY$ for all $\alpha \in L$ and  from data (\ref{datosrecoleccap22}) and (\ref{datosrecoleccap21}) we deduce that the morphisms $w_{\alpha}$ glue into a morphism $w: Z \to \FX$ that agrees with $u$ in $T$ and therefore, $f$ is formally smooth. Moreover, since the pseudo finite type condition is   local,   $f$ is of  pseudo finite type, hence, $f$ is smooth.

\item
(2) For all $\lambda \in J$ put $f_{\lambda} = f |_{f^{-1}(\FV_{\lambda})}$ and consider  the  following diagram
\begin{diagram}[height=2.5em,w=3em,labelstyle=\scriptstyle]
f^{-1}(\FV_{\lambda})  & \rTinc_{\quad i_{\lambda}} & \FX  \\
\dTto^{f_{\lambda}}    &                            & \dTto^f \\ 
\FV_{\lambda}          & \rTinc^{\quad j_{\lambda}} & \FY\\  
\end{diagram}
where the horizontal arrows are  open immersions.
If $f$ is smooth  (unramified or \'etale) by Proposition \ref{infsorit}.(\ref{infsorit1})  and Proposition \ref{infencajes}.(\ref{infencajes2}) it holds that $f \circ i_{\lambda} = j_{\lambda} \circ f_{\lambda}$ is smooth  (unramified or \'etale) and since $ j_{\lambda}$ is adic \'etale,  from  Corollary \ref{corbasinf} we get that $f_{\lambda}$ is smooth  (unramified or \'etale, respectively).  Reciprocally, if the morphisms $f_{\lambda}$ are smooth (unramified or \'etale), as the open immersions  $j_{\lambda}$ are \'etale by Proposition \ref {infsorit}.(\ref {infsorit1}) there results that the morphisms $j_{\lambda} \circ f_{\lambda}$ are smooth (unramified or \'etale). Therefore from (1) we get that $f$ is smooth (unramified or \'etale, respectively)
\end{proof}

\begin{cor} \label{corcondinflocal}
The results \ref{propibasinf},  \ref{corbasinf} and  \ref{condinflocal} are true if we replace the infinitesimal lifting properties by their adic counterparts.
\end{cor}
\begin{proof}
It is straightforward  from Proposition \ref{adic},  Proposition \ref{propibasinf},  Corollary \ref{corbasinf} and   Proposition \ref{condinflocal}. 
\end{proof}

\begin{rem}
A consequence of the two last results for the local study of the infinitesimal lifting properties  (with or without the adic hypothesis) for locally noetherian formal schemes, we can restrict to $\sfna$.
\end{rem}

\begin{cor} \label{complpetal}
Let $\FX$ be in $\sfn$ and $\FX' \subset \FX$ a closed formal subscheme. Then the  morphism of   completion of  $\FX$ along  $\FX'$, $\kappa: \FX_{/\FX'} \to \FX$ is \'etale.
\end{cor}
\begin{proof}
Applying Proposition \ref{condinflocal}  we may assume  that $\FX= \spf(A)$ and $\FX'= \spf(A/I)$ are in $\sfna$ with $A$ a $J$-adic noetherian ring and $I \subset A$ an ideal. Then, by  Example \ref{confinftipof}.(\ref{comppseetale}) if $\HA$ is the  completion of  $A$ for the $(I+J)$-adic topology, we have that $\kappa: \FX_{/\FX'}= \spf(\HA) \to \FX=\spf(A)$ is \'etale.
\end{proof}

\begin{propo} \label{complpe}
Given $f:\FX \to \FY$ in $\sfn$, let $\FX' \subset \FX$ and $\FY' \subset \FY$ be closed formal subschemes such that $f(\FX')\subset \FY'$. 
\begin{enumerate}
\item
If $f$ is smooth  (unramified or \'etale) then $\hf: \FX_{/ \FX'} \to \FY_{/ \FY'}$ is smooth  (unramified or \'etale, respectively).
\item
If moreover $\FX' = f^{-1}(\FY')$ and $f$ is adic smooth  (adic unramified or adic \'etale)  we have that   $\hf: \FX_{/ \FX'} \to \FY_{/ \FY'}$ is adic smooth  (adic unramified or adic \'etale, respectively).

\end{enumerate}
\end{propo}
\begin{proof}
By the  diagram (\ref{diagrcom}) we have the  following commutative diagram in $\sfn$
\begin{diagram}[height=2em,w=2em,labelstyle=\scriptstyle]
\FX            &	 \rTto^{f}  & \FY\\
\uTto{ \kappa} &	            & \uTto{ \kappa}\\
\FX_{/ \FX'}   &	\rTto^{\hf}	& \FY_{/ \FY'} \\
\end{diagram}
where the vertical arrows are morphisms of completion that, being slightly imprecise, we denote both by $\kappa$. Let us prove (1). If $f$ is smooth  (unramified or \'etale) by the last Corollary and Proposition \ref{infsorit}.(\ref{infsorit1}) we have that $f \circ \kappa = \kappa \circ \hf$  is also smooth  (unramified or \'etale). Since $\kappa$ is \'etale from Proposition \ref{propibasinf} we deduce that $\hf$ is smooth  (unramified or \'etale, respectively). Assertion (2) is consequence of  (1)  and of  Proposition \ref{complsorit}.
\end{proof}

\begin{cor}
Let $f: X \to Y$ be a morphism of  locally noetherian schemes and let $X' \subset X$ and $Y' \subset Y$ be closed  subschemes  such that $f(X')\subset Y'$. 
If $f$ is smooth  (unramified, \'etale), then $\hf: X_{/ X'} \to Y_{/ Y'}$ is smooth  (unramified, \'etale, respectively), if $X' = f^{-1}(Y')$  we also have that $\hf$ is   adic.
\end{cor}

\begin{defn}
Let $f:\FX \to \FY$ in $\sfn$. We say  that $f$ is \emph{smooth (unramified or \'etale) in $x \in \FX$ }if there exists an open subset $\FU \subset \FX$ such that  $f|_{\FU}$ is smooth  (unramified or \'etale, respectively).

So, by Proposition \ref{condinflocal} it holds that $f$ is smooth  (unramified or \'etale) if, and only if,  $f$ is smooth  (unramified or \'etale, respectively) in $x \in \FX$, $\forall x \in \FX$.

Observe that the  set of  points $x \in \FX$ such that  $f$ is smooth  (unramified, or \'etale) in $x $ is an open subset of  $\FX$.
\end{defn}

\section[Differentials and lifting properties]{The  module of  differentials and the infinitesimal lifting properties} \label{sec23}

In this section we study some characterizations  for  a smooth, unramified and \'etale morphism $f: \FX \to \FY$. Above all, we will focus on the properties related to the  module of  differentials  $\om^{1}_{\FX/\FY}$. We highlight the importance of the Jacobian Criterion for affine formal schemes (Corollary \ref{critjacob}) which allows us to determine when a closed formal subscheme of a smooth formal scheme is smooth rendering Zariski's Jacobian Criterion for topological rings (\emph{cf.} \cite[(\textbf{0}.22.6.1)]{EGA41} into the present context.

\begin{lem} \label{algmodifcero}
Let $B \to A$ a continuous homomorphism of  preadic rings. Then $A$ is a formally unramified $B$-al\-ge\-bra if, and only if, $\om^{1}_{A/B}=0$.
\end{lem}
\begin{proof} \cite[(\textbf{0}.20.7.4)]{EGA41}
Let $A$ be a formally unramified $B$-alge\-bra and $J \subset A$ an ideal of  definition. It is enough to prove that 
\[\Homcont_{A}(\om^{1}_{A/B}, \Omega^{1}_{A/B}/ J^{n+1}\Omega^{1}_{A/B})=0,\] for all $n \in \NN$.  By hypothesis, the  epimorphism of topological  discrete $B$-algebras\footnote{Let $A$ be a ring and $M$ a $A$-module. In  the additive  group $A \oplus M$ we consider  the multiplication
\[
(a,m) \cdot (a',m') = (a a',am'+a'm)
\]
where $a,a' \in A; m,m' \in M$. This operation is bilinear, associative and has a unit. This  ring  structure is called \emph{the semidirect product
 of  $A$ and $M$} and is denoted  $A \ltimes M$. If $A$ is a preadic ring and in $A \ltimes M$ we take the topology induced by the one of  $A$,  we have that the ring  homomorphisms  $A \mono A \ltimes M$ and $A \ltimes M \epi A$ are continuous.}  
\[
C := \frac{A}{J^{n+1}} \ltimes \frac{\Omega^{1}_{A/B}}{J^{n+1}\Omega^{1}_{A/B}} \to \frac{A}{J^{n+1}}\] 
is continuous for all $n \in \NN$ and its kernel $I:= (0,\Omega^{1}_{A/B}/ J^{n+1}\Omega^{1}_{A/B}) \subset C$ is  of  zero square, so, applying \ref{derivversuslevantalg}.(\ref{derivversuslevantalg2}) we deduce  that 
\[0= \Dercont_{B}(A, \frac{\Omega^{1}_{A/B}}{ J^{n+1}\Omega^{1}_{A/B}}) \underset{\textrm{(\ref{modifn})}} \cong \Homcont_{A}(\om^{1}_{A/B}, \frac{\Omega^{1}_{A/B}}{ J^{n+1}\Omega^{1}_{A/B}}).\]

Reciprocally,  let $C$ be topological discrete $B$-algebra, $I \subset C$ an ideal of  zero square, $u: A \to C/I$ a continuous $B$-homomorphism and  suppose  that there exist two continuous  homomorphisms of  topological algebras $v: A \to C$ and $w: A \to C$ such that  the  following diagram is commutative:
\begin{diagram}[height=2em,w=2em,p=0.3em,labelstyle=\scriptstyle]
B & \rTto  & A  \\ 
  & \rdTto & \dTto^{v}\dTto_{w} & \rdTto^{u} & \\ 
  &        & C                  & \rTonto    & \frac{C}{I}. \\
\end{diagram}
By  \ref{modifn} it holds that $\Dercont_{B}(A,I)  \cong \Homcont_{A}(\om^{1}_{A/B},I)=0$ and from \ref{derivversuslevantalg}.(\ref{derivversuslevantalg1}) there results that $u=v$ and therefore, $A$ is a  formally unramified $B$-algebra.
\end{proof}

\begin{propo} \label{modifcero}
Let $f:\FX \to \FY$ be a morphism in $\sfn$ of  pseudo finite type. The  morphism $f$ is unramified  if, and only if,  $\om^{1}_{\FX/\FY}=0$.
\end{propo}

\begin{proof}
By Proposition \ref{condinflocal} we may suppose that $f:\FX \to \FY$ is in $\sfna$ and therefore the  result follows from Lemma \ref{algmodifcero}.
\end{proof}

\begin{cor}
Let $f:\FX \to \FY$ and $g: \FY \to \FS$ be two  pseudo finite type morphisms. Then $f$ is unramified  if, and only if,  the  morphism of  $\CO_{\FX}$-Modules $f^{*}(\om^{1}_{\FY/\FS}) \to \om^{1}_{\FX/\FS}$ is surjective.
\end{cor}

\begin{proof}
Use the last Proposition and the First Fundamental Exact Sequence  (\ref{sef1}) associated to  the morphisms $\FX \xto{f} \FY \xto{g} \FS$.
\end{proof}

\begin{lem}\label{modifll}
Let $f: B \to A$ be a a continuous morphism of  preadic rings. If $A$ is a formally smooth $B$-algebra, then: 
\begin{enumerate}
\item
$A$ is a flat $B$-algebra.
\item
For all ideals of  definition $J \subset A$, $\om^{1}_{A/B}\otimes_{A} A/J$ is a projective $A/J$-module. If moreover,  $\om^{1}_{A/B}$ is a  finite type $A$-module  (for example, if $\spf(A) \to \spf(B)$ is a  pseudo finite type morphism  in $\sfna$), then  $\om^{1}_{A/B}$ is a projective $A$-module.
\end{enumerate}
\end{lem}

\begin{proof}
\item (1) It suffices to show  that for all prime ideals  $\fp \subset A$, $A_{\fp}$ is a flat $B_{\fq}$-module with  $\fq=f^{-1}(\fp)$. Fix $\fp \subset A$ a prime ideal. By  \cite[(\textbf{0}.19.3.5.(iv))]{EGA41}  it holds that $A_{\fp}$ is a formally smooth $B_{\fq}$-algebra  for  the adic topologies and applying Remark  \ref{confinfhafin} there results that $A_{\fp}$ is a formally smooth $B_{\fq}$-algebra  for the topologies given by the maximal  ideals. Then, by \cite[(\textbf{0}.19.7.1)]{EGA41} we have that $A_{\fp}$ is a flat $B_{\fq}$-module.

\item (2) (\emph{cf.}  \cite[Theorem 28.5]{ma2})
We will show that given  $\varphi: M \epi N $ a surjective homomorphism  of  $A/J$-modules,  then the map 
\[\Hom_{\frac{A}{J}}(\om^{1}_{A/B}\otimes_{A} \frac{A}{J},M) \to \Hom_{\frac{A}{J}}(\om^{1}_{A/B}\otimes_{A} \frac{A}{J},N)\]
 is surjective, that is,   $\Hom_{A}(\om^{1}_{A/B},M) \to \Hom_{A}(\om^{1}_{A/B},N)$ is surjective. By \ref{modifn} it  suffices to  prove that the map  induced by $\varphi$
\begin{equation} \label{morfderivover}
\Der_{B}(A,M) \to \Der_{B}(A,N)
\end{equation}
is surjective. If $C$ is the discrete $B$-algebra $(A/J \ltimes M)$ and the ideal $I = (0,\Ker \varphi )\subset C$, given $d \in \Der_{B}(A,N)$ we define a continuous homomorphism of  topological $B$-algebras
\[
\begin{array}{ccc}
A& \lto & \frac{C}{I}=\frac{A}{J} \ltimes N\\
a & \leadsto & (a+J, d(a))\\
\end{array}
\]
and since $I ^{2}=0$ and $A$ is a formally smooth $B$-algebra, there exists a continuous homomorphism of  topological $B$-algebras $u: A \to C$ such that  the  following diagram is commutative
\begin{diagram}[height=3em,w=2em,labelstyle=\scriptstyle]
&B & \rTto& A&\\
&\dTto&\ldTto^{u}&\dTto&\\
\frac{A}{J} \ltimes M=&C & \rTto^{(1, \varphi)} &\frac{C}{I}&= \frac{A}{J} \ltimes N\\
\end{diagram}
The second component of the  morphism $u$ defines a continuous $B$-derivation $d': A \to M$ such that  $\varphi \circ d' = d$ and therefore the map  (\ref{morfderivover}) is surjective, hence, $\om^{1}_{A/B}\otimes_{A} A/J$ is a projective $A/J$-module.

Then, if $\om^{1}_{A/B}$ is a  finite type $A$-module, by \cite[(\textbf{0}.7.2.10)]{EGA1}, it is a projective $A$-module.
\end{proof}

\begin{propo}\label{flplano}
Let $f:\FX \to \FY$ be a smooth morphism. Then $f$ is flat and $\om^{1}_{\FX/\FY}$ is  a locally free $\CO_{\FX}$-module of finite rank.
\end{propo}

\begin{proof}
Since it is a local question, the  result is consequence of the last Lemma and   \cite[(10.10.8.6)]{EGA1}.
\end{proof}

\begin{parraf} \label{formallyescindida}
Given $A$ a $J$-preadic ring, let $\HA$ be the completion  of  $A$ for  the $J$-adic topology and $A_{n}=A/J^{n+1}$, for all $n \in \NN$. Take $M'',\, M'$ and $M$ $A$-modules, denote by $\widehat{M''},\,  \widehat{M' },\,\widehat{M}$ their completions for the  $J$-adic topology and let  $\widehat{M''} \xto{u} \widehat{M' } \xto{v}\widehat{M}$ be
a sequence of  $\HA$-modules. It holds that
\begin{enumerate}
\item  \label{formallyescindida1} If $0 \to \widehat{M''} \xto{u} \widehat{M' }\xto{v} \widehat{M} \to 0$ is a  split  exact sequence   of  $\HA$-modules then, for all $n \in \NN$ 
\[
0 \to M''\otimes_{A} A_{n} \xto{u_{n}} M'\otimes_{A} A_{n} \xto{v_{n}} M\otimes_{A} A_{n} \to 0
\] is a split exact sequence.

\item  \label{formallyescindida2} Reciprocally, if  $M \otimes_{A}A_{n}$ is a projective $A_{n}$-module and \[0 \to M''\otimes_{A} A_{n} \xto{u_{n}} M'\otimes_{A} A_{n} \xto{v_{n}} M\otimes_{A} A_{n} \to 0\] is a split exact sequence  of  $A_{n}$-modules, for all $n \in \NN$,  then
\begin{equation} \label{secformlinv}
0 \to \widehat{M''} \to \widehat{M' }\to \widehat{M} \to 0
\end{equation}
is a split exact sequence of $\HA$-modules. 
\end{enumerate}
Indeed, assertion (1) is immediate. In order to prove (2), for all $n \in \NN$ we have the following commutative diagrams:
\begin{diagram}[height=2em,w=1.8em,labelstyle=\scriptstyle,midshaft]
 0 &\rTto&  M'' \otimes_{A} A_{n+1}& \rTto^{u_{n+1}} &M' \otimes_{A}  A_{n+1} & \rTto{v_{n+1}} &M \otimes_{A} A_{n+1} &  \to &0 \\
 &	&\dTonto^{f_{n}}& & \dTonto^{g_{n}}&	&\dTonto^{h_{n}}& & \\
 0 &\rTto&  M'' \otimes_{A} A_{n}& \rTto^{u_{n}} &M' \otimes_{A}  A_{n} & \rTto{v_{n}} &M \otimes_{A} A_{n} &  \to &0 \\
\end{diagram}
where the rows are split exact sequences. Applying inverse limit we have that the sequence (\ref{secformlinv}) is exact. Let us show that it splits. By hypothesis, for all $n \in \NN$ there exists $t_{n}: M \otimes_{A} A_{n} \to M'\otimes_{A}  A_{n}$ such that  $v_{n} \circ t_{n}= 1$. From  $(t_{n})$ we are going to  define a family of  morphisms $(t'_{n}: M \otimes_{A} A_{n} \to M' \otimes_{A}  A_{n})$ such that 
\begin{equation} \label{conmutatividad escindidas}
v_{n} \circ t'_{n} =1 \qquad g_{n} \circ t'_{n+1} = t'_{n} \circ h_{n}
\end{equation}
for all $n \in \NN$. For  $k=0$ put $t'_{0}:=t_{0}$. Suppose that we have constructed  $t'_{k}$ verifying (\ref{conmutatividad escindidas}) for all $k \le n$ and let us define $t'_{n+1}$. If  $w_{n}:=g_{n} \circ t_{n+1} - t'_{n}\circ h_{n}$ then $v_{n} \circ w_{n} = 0$ and therefore, $\Img w_{n} \subset  \Ker v_{n}= \Img u_{n}$. On the other hand, since $M \otimes_{A}A_{n+1}$ is a projective $A_{n+1}$-module, there exists $\theta_{n+1}: M\otimes_{A}  A_{n+1} \to u_{n+1}(  M''\otimes_{A} A_{n+1})$ such that  the  following diagram is commutative
\begin{diagram}[height=2em,w=1.8em,labelstyle=\scriptstyle,midshaft]
  &			& u_{n+1}( M''\otimes_{A} A_{n+1}) \\
  & \ruTto^{\theta_{n+1}}	& \dTonto\\
M \otimes_{A}  A_{n+1}& \rTto^{w_{n}} & u_{n}( M''\otimes_{A} A_{n}). \\
\end{diagram}
If we put $t'_{n+1}:=t_{n+1} - \theta_{n+1}$, it holds that $v_{n+1} \circ t'_{n+1}=1$ and $g_{n} \circ t'_{n+1} = t'_{n} \circ h_{n}$. The  morphism \[t':= \invlim {n \in \NN} t'_{n}\] satisfies that $v \circ t'=1$ and the sequence (\ref{secformlinv}) splits. 
\end{parraf}

\begin{propo} \label{propslisoimpllocalrot}
Let $f:\FX \to \FY$ be a smooth morphism  in $\sfn$. For all pseudo finite type morphism $\FY \to \FS$  in $\sfn$ the sequence of  coherent $\CO_{\FX}$-modules 
\[
0 \to f^{*}\om^{1}_{\FY/\FS} \xto{\Phi} \om^{1}_{\FX/\FS} \xto{\Psi} \om^{1}_{\FX/\FY} \to 0
\]
is exact and locally split.
\end{propo}

\begin{proof}
It is a local question, so we may assume that $f: \FX=\spf(A) \to \FY=\spf(B)$ and $g: \FY=\spf(B) \to \FS=\spf(C)$ are in $\sfna$ with $A$ a formally smooth $B$-algebra. Given $J \subset A$ an ideal of  definition if we take $A_{n}= A/J^{n+1}$, by  Lemma \ref{modifll} we have that $\om^{1}_{A/B} \otimes_{A} A_{n}$ is a projective $A_{n}$-module, for all $n \in \NN$ and therefore, by \ref{formallyescindida}.(\ref{formallyescindida2}) it suffices to show that for all $n \in \NN$
\begin{equation*} 
0 \to \Omega^{1}_{B/C} \otimes_{B} A_{n} \xto{\Phi_{n}} \Omega^{1}_{A/C} \otimes_{A} A_{n} \xto{\Psi_{n}} \Omega^{1}_{A/B}  \otimes_{A} A_{n} \to 0
\end{equation*}
is a split exact sequence. By the First Fundamental Exact Sequence  of  differential modules associated to $C \xto{g} B \xto{f} A$ it suffices to   prove that for all $n \in \NN$ the map $\Phi_{n}$ is a section or, equivalently, that for all  $A_{n}$-modules $M$, the induced map
\[
\Hom_{A_{n}}(\Omega^{1}_{A/C} \otimes_{A} A_{n},M) \to \Hom_{A_{n}}(\Omega^{1}_{B/C} \otimes_{B} A_{n},M)
\]
is surjective. By (\ref{derivmodulalgebdis}) this  is equivalent to see that  the map 
\begin{equation} \label{morfderivover2}
\Der_{C}(A,M) \to \Der_{C}(B,M)
\end{equation}
is surjective. Let $n \in \NN$ and $M$  be a $A_{n}$-module. Each $d \in \Der_{C}(B,M)$ defines a continuous homomorphism of  topological algebras
\begin{equation*} 
\begin{array}{ccc}
B& \xto{\lambda} & A_{n} \ltimes M\\
b & \leadsto & (f(b) + J^{n+1}, d(b))\\
\end{array}
\end{equation*}
Since $M ^{2}=0$ as an ideal of $A_{n} \ltimes M$ and $A$ is a formally smooth $B$-algebra, there exists a continuous homomorphism of  topological $B$-algebras $v: A \to A_{n} \ltimes M$ such that  the  following diagram is commutative
\begin{diagram}[height=2em,w=2em,labelstyle=\scriptstyle]
B & \rTto{g}& A\\
\dTto{\lambda}&\ldTto{v}&\dTto\\
A_{n} \ltimes M & \rTto &A_{n}.\\
\end{diagram}
The second component of the morphism $v$   gives a $C$-derivation  $d': A \to M$ such that  $ d' \circ g = d$ and therefore the map  (\ref{morfderivover2}) is surjective.
\end{proof}

\begin{cor} \label{corisomodifetal}
Let $f:\FX \to \FY$ be an \'etale morphism  in $\sfn$. For all  pseudo finite type morphism $\FY \to \FS$  in $\sfn$ it holds that
\[
f^{*}\om^{1}_{\FY/\FS} \cong \om^{1}_{\FX/\FS}
\]
\end{cor}

\begin{proof}
It is a  consequence of the last result   and of  Proposition \ref{modifcero}.
\end{proof}
\begin{propo}
Let $f:\FX \to \FY$ be a  pseudo finite type morphism in $\sfn$  and $g:\FY \to \FS$ a smooth morphism  in $\sfn$. The following conditions are equivalent:
\begin{enumerate}
\item
$f$ is smooth 
\item
$g \circ f $ is smooth  and the sequence 
\[
0 \to f^{*}\om^{1}_{\FY/\FS} \to \om^{1}_{\FX/\FS} \to \om^{1}_{\FX/\FY} \to 0
\]
is exact and locally split.
\end{enumerate}
\end{propo}

\begin{proof}
The  implication (1) $\Rightarrow$ (2) is consequence of  Proposition \ref{infsorit} and of  Proposition \ref{propslisoimpllocalrot}.

As for  proving that (2) $\Rightarrow$ (1) we may suppose  that $f: \FX=\spf(A) \to \FY=\spf(B)$ and $g: \FY=\spf(B) \to \FS=\spf(C)$ are in $\sfna$ being $B$ and $A$ formally smooth $C$-algebras. Let us show  that $A$ is a formally smooth $B$-algebra.
Let $E$ be a discrete ring  , $I \subset E$ a square zero ideal and consider the commutative diagram of  continuous homomorphisms of  topological rings
\begin{diagram}[height=2em,w=1.8em,labelstyle=\scriptstyle,midshaft]
C  &	\rTto{f}		&B 				& \rTto{g}		& A\\
     & 		& \dTto^{\lambda}	& 	       		& \dTto^{u}\\
     & 				&E 				& \rTonto^{j}	& E/I.\\
\end{diagram}
Since $A$ is a formally smooth $C$-algebra, there exists a continuous homomorphism of topological $C$-algebras $v: A \to E$ such that  $v \circ g \circ f = \lambda \circ f$ and $j \circ v = u$. Then by  \ref{derivversuslevantalg}.(\ref{derivversuslevantalg1}) we have that $d:= \lambda - v \circ g \in \Dercont_{C}(B,E)$. 
From the hypothesis and  considering the equivalence of categories \ref{propitrian}.(\ref{equivtrian}) we have that the sequence of  finite type $A$-modules  
\[
0 \to \om^{1}_{B/C} \otimes_{B}A \to \om^{1}_{A/C}  \to \om^{1}_{A/B}   \to 0
\]
is exact and split. Besides, since the  morphism $v$ is continuous and $E$ is discrete there exists $n \in \NN$ such that  $E$ is a $A/J^{n+1}$-module.  Therefore the induced map  $\Hom_{A}(\Omega^{1}_{A/C}, E) \to \Hom_{A}(\Omega^{1}_{B/C}, E)$ is surjective and applying  \ref{modifn} we have that  the map 
\[
\Der_{C}(A,E) \to \Der_{C}(B,E)
\]
is surjective too. It follows  that there exists $d' \in \Der_{C} (A,E)$ such that  $d' \circ g =d$. If we put $v':= v+d'$, we have that $v' \circ g = \lambda$ and $j \circ v'  = u$. As a consequence,  $A$ is a formally smooth $B$-algebra.
\end{proof}

\begin{cor} \label{modifisopet}
Let $f:\FX \to \FY$ and $g:\FY \to \FS$ be two  pseudo finite type morphisms such that  $g \circ f$ and $g$ are smooth. Then, $f$ is \'etale  if, and only if,  $f^{*} \om^{1}_{\FY/\FS} \cong \om^{1}_{\FX/\FS}$.
\end{cor}

\begin{proof}
Follows from the last Proposition and Proposition \ref{modifcero}.
\end{proof}

\begin{propo}{(Zariski Jacobian criterion for  preadic rings)} \label{critjacobanillos}
Let $B \to A$ be a continuous morphism of  preadic rings  and suppose that $A$ is a formally smooth $B$-algebra. Given an ideal $I \subset A$. Let us consider in $A' := A/I$ the topology induced by the topology of $A$. The following conditions are  equivalent:
\begin{enumerate}
\item
$A'$ is a  formally smooth $B$-algebra.
\item
Given $J \subset A$ an ideal of  definition of  $A$, define $A'_{n}:= A/(J^{n+1}+I)$. The sequence of  $A'_{n}$-modules 
\begin{equation*}
0 \to \frac{I}{I^{2}} \otimes_{A'} A'_{n} \xto{\delta_{n}} \Omega^{1}_{A/B} \otimes_{A}   A'_{n} \xto{\Phi_{n}} \Omega^{1}_{A'/B} \otimes_{A'} A'_{n}  \to 0
\end{equation*}
is exact and  split,  for all $n \in  \NN$.
\item
The sequence of  $\widehat{A'}$-modules 
\begin{equation*}
0 \to \widehat{ \frac{I}{I^{2}}} \xto{\delta} \Omega^{1}_{A/B} \tc_{A} A' \xto{\Phi} \om^{1}_{A'/B}  \to 0
\end{equation*}
is exact and split.
\end{enumerate}
\end{propo}

\begin{proof}
The equivalence (1) $\Leftrightarrow$ (2) follows from  \cite[(\textbf{0}.22.6.1), (\textbf{0}.19.1.5), (\textbf{0}.19.1.7)]{EGA41} and from the Second Fundamental Exact Sequence associated to the morphisms $B \to A \to A'$. Let us show that (2) $\Leftrightarrow$ (3). Since $A'$ is a formally smooth $B$-algebra,  from Lemma \ref{modifll} we deduce that $\Omega^{1}_{A'/B} \otimes_{A'}A'_{n}$ is a projective $A'_{n}$-module, for all $n \in \NN$  and, therefore, the result follows from \ref{formallyescindida}.(\ref{formallyescindida2}).

\end{proof}

\begin{cor}{(Zariski Jacobian criterion for formal schemes)} \label{critjacob}
Let $f: \FX=\spf(A) \to \FY=\spf(B)$ be a smooth morphism   in $\sfna$ and $\FX' \inc \FX$ a  closed immersion given by an Ideal $\CI=I^{\tr} \subset \CO_{\FX}$. The followings conditions are equivalent:
\begin{enumerate}
\item
The composed morphism  $\FX'  \inc \FX \xto{f} \FY$ is smooth.
\item
Given $\CJ \subset \CO_{\FX}$ an Ideal of  definition, if  $\CO_{X'_{n}}:= \CO_{\FX}/(\CJ^{n+1}+\CI)$, the sequence of  coherent $\CO_{X'_{n}}$-Modules
\[
\qquad \qquad 0 \to \frac{\CI}{\CI^{2}} \otimes_{\CO_{\FX'}} \CO_{X'_{n}} \xto{\delta_{n}} \om^{1}_{\FX/\FY} \otimes_{\CO_{\FX}}   \CO_{X'_{n}} \xto{\Phi_{n}} \om^{1}_{\FX'/\FY} \otimes_{\CO_{\FX'}} \CO_{X'_{n}}  \to 0
\]
is exact  and  locally split,  for all $n \in  \NN$.
\item
The sequence of coherent $\CO_{\FX'}$-Modules
\begin{equation*}
0 \to \frac{\CI}{\CI^{2}} \xto{\delta} \om^{1}_{\FX/\FY} \otimes_{\CO_{\FX}}  \CO_{\FX'} \xto{\Phi} \om^{1}_{\FX'/\FY}  \to 0
\end{equation*}
is exact and locally split.
\end{enumerate}
\end{cor}

\begin{proof}
By the equivalence of categories \ref{propitrian}.(\ref{equivtrian}) it is a consequence of the last Proposition.
\end{proof}


\chapter[Characterization]{Characterization of  the infinitesimal lifting properties}\label{cap3}
\setcounter{equation}{0}


Given $f:\FX \to \FY$ a morphism in $\sfn$ there exist Ideals of  definition of  $\FX$ and $\FY$ such that  
\[f = \dirlim {n \in \NN} f_{n}.\] So it is sensible to ask about the relation that exists between the infinitesimal lifting properties of  $f$ and the infinitesimal lifting properties of the morphisms of ordinary schemes $\{f_{n}\}_{n \in \NN}$. In this chapter two parts are distinguished. 

The first encompasses three sections (\ref{sec31}, \ref{sec32} and \ref{sec33}) we show that there exists  a close relationship between the infinitesimal lifting properties of  an adic  morphism \[f=\dirlim {n \in \NN} f_{n}\] and the infinitesimal lifting properties of the morphism of ordinary schemes $f_{0}$, something that does not happen without the adic hypothesis and that shows the interest of understanding the non adic case. 

The second part is devoted to fit in this framework the completion morphisms. In Section \ref{sec34} we  characterize open immersions and the completion morphisms  in terms of the \'etale property.  One of  the main results  is  Theorem \ref{teorequivet} in which, given 
\[\FY= \dirlim {n \in \NN} Y_{n}\] 
in $\sfn$, it is established an equivalence of categories between  \'etale adic $\FY$-formal schemes  and \'etale $Y_{0}$-schemes. In absence of  the adic hypothesis the  behavior of  the infinitesimal lifting properties can not be reduced to the study of  the infinitesimal lifting properties in  the category of  schemes  (see   Example \ref{peynopen}). The  structure of maps possessing infinitesimal lifting properties in the  non adic case is elucidated in Section \ref{sec35}. Here we obtain some of the main results of this work, namely, Theorem \ref{tppalet} and Theorem \ref{tppall}. They state that every smooth morphism  and every \'etale  morphism factors locally as a completion morphism followed by a smooth adic morphism and an \'etale adic morphism, respectively. These results explain then the local structure of smooth and \'etale morphisms of formal schemes.


\section{Unramified morphisms} \label{sec31}
We  begin  relating the  unramified character of a morphism $f:\FX \to \FY$ of locally noetherian formal schemes 
\[f:\FX \to \FY= \dirlim {n \in \NN} (f_{n}:X_{n} \to Y_{n})\] 
and the one of the underlying ordinary schemes morphisms  $\{f_{n}\}_{n \in \NN}$, where $\CJ \subset \CO_{\FX}$ and $\CK \subset \CO_{\FY}$ are Ideals of  definition such that $f^{*}(\CK) \CO_{\FX} \subset \CJ$ and that allow us to express $f$ as a limit. We will use this notation throughout the section without mentioning it explicitly.

\begin{propo}\label{nrfn}
With the previous notations, the  morphism $f$ is unramified  if, and only if,  $f_{n}:X_{n} \to Y_{n}$  is unramified, for all $n \in \NN$.
\end{propo}

\begin{proof}
Applying  Proposition \ref{modifcero} we have  to show that $\om^{1}_{\FX/\FY} =0$ equivalently  for all $n \in \NN$, $\Omega^{1}_{X_{n}/Y_{n}}=0$. If $\om^{1}_{\FX/\FY} =0$,  by the Second Fundamental Exact Sequence (\ref{sef2}) for the morphisms \[X_{n} \inc \FX \xto{f} \FY,\] we have that $\Omega^{1}_{X_{n}/ \FY}=0$, for all $n \in \NN$.  From  the First Fundamental Exact Sequence (\ref{sef1}) associated to the morphisms \[X_{n} \xto{f_{n}} Y_{n} \inc \FY,\] there results that $\Omega^{1}_{X_{n}/Y_{n}}=0$. Conversely, if for all $n \in \NN$, $\Omega ^{1}_{X_{n}/Y_{n}}=0$, Proposition \ref{modifmodifn} shows that 
\[\om^{1}_{\FX/\FY} =\invlim {n \in \NN} \Omega ^{1}_{X_{n}/Y_{n}} =0.\qedhere\]
\end{proof}

\begin{cor} \label{pecimplnoram}
Let $f:\FX \to \FY$ be a morphism in $\sfn$ and let $\CJ \subset \CO_{\FX}$ and $\CK \subset \CO_{\FY}$ be Ideals of  definition such that $ f^{*}(\CK) \CO_{\FX} \subset \CJ$. If the  induced morphisms $f_{n}:X_{n} \to Y_{n}$ are  immersions, for all $n \in \NN$, then $f$ is unramified. 
\end{cor}

In the class of  adic  morphisms in $\sfn$ the  following proposition provides a criterion, stronger than the last result,  to  determine when a morphism $f$ is unramified.

\begin{propo} \label{fnrf0nr}
Let $f:\FX \to \FY$ be an \emph{adic} morphism in $\sfn$ and $\CK \subset \CO_{\FY}$ an Ideal of  definition. Write \[f=\dirlim {n \in \NN} f_{n}\] by taking Ideals of  definition $\CK \subset \CO_{\FY}$ and $\CJ=f^{*}(\CK)\CO_{\FX} \subset \CO_{\FX}$. The  morphism $f$ is unramified if, and only if,  the induced morphism $f_{0}:X_{0} \to Y_{0}$ is unramified.
\end{propo}

\begin{proof}
If $f$ is unramified by Proposition  \ref{nrfn} we have that $f_{0}$ is unramified. 
Conversely, suppose that  $f_{0}$ is unramified and let us prove that $\om^{1}_{\FX/\FY}=0$. The question is local so we may assume that $f:\FX=\spf(A) \to \FY= \spf(B)$ is in $\sfna$ and that $\CJ=J^{\tr},\, $ with $J \subset A$ an ideal of  definition. By hypothesis $\Omega^{1}_{X_{0}/Y_{0}}=0$  and thus, since $f$ is adic there results that 
\begin{equation} \label{eqmodifceromodiff0cero}
\om^{1}_{\FX/\FY} \otimes_{\CO_{\FX}} \CO_{X_{0}}\underset{\textrm{\ref{modifesqu}}} = \Omega^{1}_{X_{0}/Y_{0}}=0.
\end{equation}
Then by the equivalence  of categories \ref{propitrian}.(\ref{equivtrian}), the last equality says that $\om^{1}_{A/B}/J \om^{1}_{A/B}=0$. On the other hand, since $A$ is a $J$-adic ring it holds that  $J \subset \rad_{A}$ ($\rad_{A}$ is the Jacobson  radical  of  $A$). Moreover, Proposition \ref{modifinito} provides that $\om^{1}_{A/B}$ is a  finite type $A$-module. From Nakayama's  lemma we deduce that $\om^{1}_{A/B}=0$ and therefore, $\om^{1}_{\FX/\FY}= (\om^{1}_{A/B})^{\tr}= 0$. Applying Proposition \ref{modifcero} it follows that $f$ is unramified.
\end{proof}

The following example illustrates that in the  non adic  case  the  analogous of  the last  proposition does not hold.

\begin{ex} \label{framfonram}
Let  $K$ be a field and $p:\BD^{1}_{K} \to \spec(K)$ be the projection morphism of  the formal disc of  dimension $1$ over $\spec(K)$. In    Example \ref{contrmod0} we have seen that $\om^{1}_{p} =(K[[T]] \hd T)^{\tr}$ and therefore, there results that $\BD^{1}_{K}$ is  ramified over $K$ (Proposition \ref{modifcero}).  However, given the  ideal of  definition $[[T]] \subset K[[T]]$ the induced morphism $p_{0} = 1_{\spec(K)}$ is unramified.
\end{ex}

In view of  this example, our next goal will be to determine when a morphism  in $\sfn$
\[f:\FX \to \FY= \dirlim {n \in \NN} (f_{n}:X_{n} \to Y_{n})\] 
with  $f_{0}$ unramified and not necessarily adic, is unramified (Corollary \ref{corf0imfpnr}). In order to do that, we will need some results that describe the local behavior of unramified morphisms. Next, we provide local characterizations of  the unramified  morphisms in $\sfn$, that generalizing the  analogous properties in the category of  schemes (\emph{cf.} \cite[(17.4.1)]{EGA44}).

\begin{propo} \label{caraclocalpnr}
Let $f:\FX \to \FY$ be a morphism in $\sfn$ of  pseudo finite type. Given $x \in \FX$ and $y=f(x)$ the following conditions are equivalent:
\begin{enumerate}
\item[(1)]
 $f$ is unramified  at $x$.
\item[(2)]
$f^{-1}(y)$ is an unramified  $k(y)$-formal scheme  at $x$.
\item[(3)]
$\fm_{\FX,x}\widehat{\CO_{\FX,x}} = \fm_{\FY,y}\widehat{\CO_{\FX,x}}$ and $k(x)|k(y)$ is a finite  separable extension.
\item[(4)]
$\om^{1}_{\CO_{\FX,x}/ \CO_{\FY,y}}  =0$
\item[$(4')$]
$(\om^{1}_{\FX/\FY})_{x}=0$
\item[(5)]
$\CO_{\FX,x}$ is a formally unramified $\CO_{\FY,y}$-algebra  for  the  adic topologies.
\item[$(5')$]
$\widehat{\CO_{\FX,x}}$ is a formally unramified $\widehat{\CO_{\FY,y}}$-algebra   for  the  adic topologies.
\end{enumerate}
\end{propo}

\begin{proof}
Let $\CJ \subset \CO_{\FX}$ and $\CK \subset \CO_{\FY}$ be Ideals of  definition such that $f^{*}(\CK)\CO_{\FX} \subset \CJ$ which allows us to write
\[f:\FX \to \FY= \dirlim {n \in \NN} (f_{n}:X_{n} \to Y_{n})
\]
(1) $\Leftrightarrow$ (2) By Proposition \ref{nrfn}, $f$ is unramified  at $x$ if, and only if, all the morphisms $f_{n}:X_{n} \to Y_{n}$ are unramified at $x$. Applying  \cite[(17.4.1)]{EGA44}, this is equivalent to $f_{n}^{-1}(y)$ be an unramified $k(y)$-scheme at $x$, for all $n \in \NN$, which  is also equivalent   to  \[f^{-1} (y) \underset{\textrm{\ref{fibra}}}= \dirlim {n \in \NN} f_{n}^{-1}(y)\] being an unramified $k(y)$-formal scheme  at $x$ (see Proposition \ref{nrfn}).

(1) $\Leftrightarrow$ (3) The assertion (1) is equivalent to $f_{n}:X_{n} \to Y_{n}$ being unramified at $x$, for all $n \in \NN$, and from \cite[\emph{loc.~ cit.}]{EGA44} there results that  $k(x)|k(y)$  is a finite separable extension, and that $\fm_{X_{n},x} = \fm_{Y_{n},y} \CO_{X_{n},x}$, for all $n \in \NN$. Hence, 
\[
\fm_{\FX,x} \widehat{\CO_{\FX,x}} = \invlim {n \in \NN} \fm_{X_{n},x} =  \invlim {n \in \NN} \fm_{Y_{n},y} \CO_{X_{n},x} = \fm_{\FY,y}\widehat{\CO_{\FX,x}}.
\]

(3) $\Rightarrow$ (4) Since $k(x)|k(y)$  is a finite separable extension  we have that $\Omega^{1}_{k(x)/k(y)}=0$ and therefore, it holds  that \[\om^{1}_{\CO_{\FX,x}/ \CO_{\FY,y}} \otimes_{\widehat{\CO_{\FX,x}}} k(x) \underset{\textrm{\ref{propimodifcom}.(\ref{propimodifcom1})}}= \om^{1}_{(\CO_{\FX,x} \otimes_{\CO_{\FY,y}} k(y))/k(y)} \underset{\textrm{(3)}}=\Omega^{1}_{k(x)/k(y)}=0.\] 
From  \ref{modiftfloc}.(\ref{modiftfloc2}) we deduce that $ \om^{1}_{\CO_{\FX,x}/ \CO_{\FY,y}}$ is a  finite type $\widehat{\CO_{\FX,x}}$-module  and thus, by   Nakayama's lemma, $\om^{1}_{ \CO_{\FX,x}/ \CO_{\FY,y}} =0$.

(4) $\Leftrightarrow  (4')$ By \ref{modiftfloc}.(\ref{modiftfloc2}) it holds that $(\om^{1}_{\FX/\FY})_{x}$ is a finite type $\CO_{\FX,x}$-module   and therefore,  
\[
\om^{1}_{ \CO_{\FX,x}/ \CO_{\FY,y}} \underset{\textrm{\ref{propimodifcom}.(\ref{propimodifcom3})}}= \widehat{(\om^{1}_{\FX/\FY})_{x}} = (\om^{1}_{\FX/\FY})_{x} \otimes_{\CO_{\FX,x}} \widehat{\CO_{\FX,x}}.
\]
Then, since $\widehat{\CO_{\FX,x}}$ is a faithfully flat $\CO_{\FX,x}$-algebra it holds that $\om^{1}_{ \CO_{\FX,x}/ \CO_{\FY,y}}=0$ if, and only if,  $(\om^{1}_{\FX/\FY})_{x}=0$.

(4) $\Leftrightarrow$ (5) It is straightforward from Lemma \ref{algmodifcero}.

(5) $\Leftrightarrow  (5')$ It is enough to apply  Lemma \ref{condinfalgcompl}.

$(4') \Rightarrow$ (1) Since $\om^{1}_{\FX/\FY} \in \coh(\FX)$ (Proposition \ref{finitmodif}), the assertion $(4')$ implies that there exists  an open subset $\FU \subset \FX$  with $x \in \FU$ such that  $(\om^{1}_{\FX/\FY})|_{\FU}=0$ and therefore,  by Proposition \ref{modifcero} we have that $f$ is unramified  at $x$.
\end{proof}

\begin{cor} \label{corcaraclocalpnr}
Let $f:\FX \to \FY$ be a pseudo finite type morphism in $\sfn$. The following conditions are equivalent:
\begin{enumerate}
\item[(1)]
 $f$ is unramified.
\item[(2)]
For all $x \in \FX,\, y=f(x)$, $f^{-1}(y)$ is an unramified $k(y)$-formal scheme  at $x$.
\item [(3)]
For all $x \in \FX,\, y=f(x)$, $\fm_{\FX,x}\widehat{\CO_{\FX,x}} = \fm_{\FY,y}\widehat{\CO_{\FX,x}}$ and $k(x)|k(y)$ is a finite separable extension. 
\item[(4)]
$\om^{1}_{\CO_{\FX,x}/ \CO_{\FY,y}}  =0$, for all $x \in \FX$ with $y=f(x)$. 
\item[$(4')$]
For all $x \in \FX$, $(\om^{1}_{\FX/\FY})_{x} =0$.
\item[(5)]
For all $x \in \FX,\, y=f(x)$, $\CO_{\FX,x}$ is a formally unramified $\CO_{\FY,y}$-algebra for the  adic topologies.
\item[$(5')$]
For all $x \in \FX,\, y=f(x)$, $\widehat{\CO_{\FX,x}}$ is a formally unramified  $\widehat{\CO_{\FY,y}}$-algebra for the  adic topologies.
\end{enumerate}

\end{cor}

\begin{cor} \label{corpnrimplcr}
Let $f:\FX \to \FY$ be a pseudo finite type morphism in $\sfn$. If $f$ is unramified  at $x \in \FX$, then $f$ is a quasi-covering at $x$.
\end{cor}

\begin{proof}
By assertion (3) of  Proposition \ref{caraclocalpnr} we have that \[\CO_{\FX,x} \tc_{\CO_{\FY,f(x)}} k(f(x))=k(x)\] with $k(x)|k(f(x))$ a finite extension and therefore, $f$ is a quasi-covering at $x$ (see Definition \ref{defncuasireves}).
\end{proof}

\begin{cor} \label{pnrdim0} 
Let $f:\FX \to \FY$ be a pseudo finite type morphism in $\sfn$. If $f$ is unramified  at $x \in \FX$, then $\dim_{x} f =0$.
\end{cor}

\begin{proof}
It is straightforward from  Proposition \ref{cuairevdim0} and from the last Corollary.
\end{proof}

\begin{propo} \label{fonrydislocal}
Let $f:\FX \to \FY$ be a pseudo finite type morphism in $\sfn$. Given $x \in \FX$ and $y=f(x)$ the following conditions are equivalent:
\begin{enumerate}
\item
$f$ is unramified  at $x$
\item
$f_{0}:X_{0} \to Y_{0}$ is unramified at $x$ and $\widehat{\CO_{\FX,x}} \otimes_{\widehat{\CO_{\FY,y}}} k(y) = k(x)$
 \end{enumerate}
\end{propo}

\begin{proof}
If $f$ is unramified at $x$, then $f_{0}$ is unramified at $x$ (Proposition \ref{nrfn}). Moreover,  assertion (3) of  Proposition \ref{caraclocalpnr} implies that $\widehat{\CO_{\FX,x}} \otimes_{\widehat{\CO_{\FY,y}}} k(y) = k(x)$ so (1) $\Rightarrow$ (2) holds. Let us prove  that (2) $\Rightarrow$ (1). Since $f_{0}$  is unramified at $x$  we have that $k(x)|k(y)$ is a finite separable extension (\emph{cf.} \cite[(17.4.1)]{EGA44}). From  the equality $\widehat{\CO_{\FX,x}} \otimes_{\widehat{\CO_{\FY,y}}} k(y) = k(x)$ we deduce that $\fm_{\FX,x}\widehat{\CO_{\FX,x}} = \fm_{\FY,y}\widehat{\CO_{\FY,y}}$. Thus,  the  morphism $f$  and the  point $x$ satisfy assertion (3) of Proposition  \ref{caraclocalpnr} and there results that $f$ is unramified  at $x$.
\end{proof}

Now we are ready to state the non adic version of Proposition \ref{fnrf0nr}:

\begin{cor}  \label{corf0imfpnr}
Given $f:\FX \to \FY$ a morphism in $\sfn$ of  pseudo finite type let $\CJ \subset \CO_{\FX}$ and $\CK \subset \CO_{\FY}$ be Ideals of  definition such that $f^{*}(\CK)\CO_{\FX} \subset \CJ$ and let $f_0 \colon X \to Y $ be the induced morphism. The following conditions are equivalent:
\begin{enumerate}
\item
The morphism $f$ is unramified.
\item
The morphism $f_{0}$ is unramified and, for all $x \in \FX$, $f^{-1}(y)=f_{0}^{-1}(y)$ with $y=f(x)$.
\end{enumerate}
\end{cor}

\begin{proof}
Suppose that $f$ is unramified and fix $x \in \FX$ and $y=f(x)$. By the Proposition  \ref{fonrydislocal} we have that $f_{0}$ is unramified and that  $\widehat{\CO_{\FX,x}} \otimes_{\widehat{\CO_{\FY,y}}} k(y) = k(x)$. Therefore, $\CJ \cdot (\widehat{\CO_{\FX,x}} \otimes_{\widehat{\CO_{\FY,y}}} k(y) )= 0$ and applying  Lemma \ref{fibradiscr} we deduce that $f^{-1}(y)=f_{0}^{-1}(y)$. Reciprocally,  suppose that  (2) holds and let us show that given $x \in \FX$, the  morphism $f$ is unramified  at $x$. If $y=f(x)$, we have that $f_{0}^{-1}(y)$ is an unramified   $k(y)$-scheme at $x$ (\emph{cf.} \cite[(17.4.1)]{EGA44}) and since $f^{-1}(y)=f_{0}^{-1}(y)$, from   Proposition \ref{caraclocalpnr} there results that $f$ is unramified  at $x$.
\end{proof}

\begin{lem} \label{fibradiscr} 
Let $A$ be a $J$-adic noetherian ring such that  for all open prime ideals $\fp \subset A$, $J_{\fp}=0$. Then $J=0$ and therefore, the  $J$-adic topology  in $A$ is the discrete topology.  
\end{lem}

\begin{proof}
Since every maximal ideal $\fm \subset A$ is open for the $J$-adic topology, we have that $J_{\fm}=0$, for all maximal ideal $\fm \subset A$, so $J=0$.
\end{proof}

\begin{parraf}
As  consequence of Corollary \ref{corf0imfpnr} there results that: 
\begin{itemize}
\item
If $f:\FX \to \FY$ is an unramified morphism in $\sfn$ then  $f^{-1}(y)$ is a usual scheme for all $x \in \FX$ being  $y=f(x)$.
\item
In  Corollary \ref{corcaraclocalpnr} assertion (2) may be written:
\item[$(2')$]
\emph{For all $x \in \FX,\, y=f(x)$, $f^{-1}(y)$ is a unramified $k(y)$-scheme  at $x$.}

\end{itemize}
\end{parraf}

From Proposition \ref{caraclocalpnr} we obtain the  following result, in which we provide a description of  pseudo  closed immersions that will be used  in the characterization of  completion morphisms (Theorem \ref{caracmorfcompl}).

\begin{cor} \label{pecigf0ecnr} 
Given $f:\FX \to \FY$ in $\sfn$, let $\CJ \subset \CO_{\FX}$ and $\CK \subset \CO_{\FY}$ be Ideals of  definition such that $f^{*}(\CK)\CO_{\FX} \subset \CJ$ and that let us express \[f=\dirlim {n\in \NN} f_{n}.\] The  morphism $f$ is a pseudo  closed immersion if, and only if,  $f$ is unramified  and $f_{0}:X_{0} \to Y_{0}$ is a  closed immersion.
\end{cor}

\begin{proof}
If $f$ is a pseudo  closed immersion, by  Corollary  \ref{pecimplnoram} there results that $f$ is unramified. Reciprocally, suppose that $f$ is unramified and that $f_{0}$ is a  closed immersion and let us show that $f_{n}:X_{n} \to Y_{n}$ is a  closed immersion, for each $n \in \NN$. By \cite[(4.2.2.(ii))]{EGA1} it suffices to prove that, for all $x \in \FX$ with $y=f(x)$, the  morphism  $\CO_{Y_{n},y} \to \CO_{X_{n},x}$ is surjective, for all $n \in \NN$. Fix $x \in \FX$, $ y=f(x) \in \FY$ and $n \in \NN$. Since $f_{0}$ is a  closed immersion, by \cite[\emph{loc. cit.}]{EGA1},   we have that $\CO_{Y_{0},y} \to \CO_{X_{0},x}$ is surjective and therefore, $\spf(\widehat{\CO_{\FX, x}}) \to \spf(\widehat{\CO_{\FY, y}})$ is a  pseudo finite morphism, so, the  morphism $\CO_{Y_{n},y} \to \CO_{X_{n},x}$ is  finite (see \ref{defmtf1}). On the other hand, since $f$ is unramified  by Proposition \ref{nrfn} we get that $f_{n}$ is unramified and therefore, applying Proposition \ref{caraclocalpnr}   we have that $\fm_{Y_{n},y} \CO_{X_{n},x}= \fm_{X_{n},x}$. Then by  Nakayama's lemma we conclude that  $\CO_{Y_{n},y} \to \CO_{X_{n},x}$ is a surjective morphism, for all $n \in \NN$.
\end{proof}

\section{Smooth morphisms} \label{sec32}

The contents of  this section can be structured in two parts. In the first part we study the relationship between the smoothness of  a morphism \[f= \dirlim {n \in \NN} f_{n}\] in $\sfn$ and the smoothness of  the morphisms $f_{n}$ and we locally characterize smooth morphisms. In the second part, we provide a local factorization for smooth morphisms (Proposition  \ref{factpl}). We also prove in Corollary \ref{criteriojacobiano} the  Jacobian criterion, that is a useful explicit condition in terms of a matrix rank for determining whether a closed subscheme in the affine formal space or in the affine formal disc is smooth or not.

\begin{propo} \label{lfn}
Given $f:\FX \to \FY$ in $\sfn$ let $\CJ \subset \CO_{\FX}$ and $\CK \subset \CO_{\FY}$ be Ideals of  definition with $f^{*}(\CK) \CO_{\FX} \subset \CJ$ such that \[f=\dirlim {n\in \NN} f_{n}.\] If $f_{n}:X_{n} \to Y_{n}$ is smooth, for all $n \in \NN$, then $f$ is smooth.
\end{propo}

\begin{proof}
By Proposition \ref{condinflocal} we may assume that $f$ is in $\sfna$. 
Let $Z$ be an affine a scheme, $w: Z \to \FY$ a morphism, $T \inc Z$ a closed $\FY$-subscheme given by a square zero Ideal and consider $u: T \to \FX$ a $\FY$-morphism. Since $f$ and $w$ are morphisms of  affine formal schemes we find  an integer $k \ge 0$ such that  $w^{*}(\CK^{k+1}) \CO_{Z} =0$ and $u^{*}(\CJ^{k+1})\CO_{T}=0$  and therefore $u$ and $w$ factors as $T \xto{u_{k}} X_{k} \xto{i_{k}} \FX $ and $Z \xto{w_{k}} Y_{k} \xto{i_{k}} \FY $, respectively. Since $f_{k}$ is formally smooth, there exists a $Y_{k}$-morphism $v_{k}: Z \to X_{k}$ such that  the  following diagram is commutative
\begin{diagram}[height=2.5em,w=2em,labelstyle=\scriptstyle]
T &    			\rTinc	&   &  Z \\
\dTto^{u_{k}}   &    			 &	\ldTto(3,2)^{v_{k}}  & \dTto^{w_{k}}\\ 
X_{k}	  	 &   	 \rTto^{\qquad f_{k}}   	&	              & Y_{k}\\ 
 \dTto^{i_{k}}    		 &    		         	 &			         & \dTto\\ 
 \FX	  & 		\rTto^{\qquad f} 		&				 & \FY.\\
\end{diagram}
Thus the  $\FY$-morphism $v:=i_{k} \circ v_{k}$ satisfies that $v|_{T}=u$ and then, $f$ is formally smooth. Moreover,  since $f_{0}$ is a finite type morphism, it holds that $f$ is of  pseudo finite type and therefore,  $f$ is smooth. 
\end{proof}

\begin{cor} \label{ladfn}
Let $f:\FX \to \FY$ be an \emph{adic}  morphism  in $\sfn$ and consider $\CK \subset \CO_{\FY}$ an  Ideal of  definition.  
The morphism $f$ is smooth  if, and only if, all the scheme morphisms $\{f_{n}:X_{n} \to Y_{n}\}_{n \in \NN}$, determined by the Ideals of  definition $\CK \subset \CO_{\FY}$ and $\CJ=f^{*}(\CK)\CO_{\FX}$, are  smooth.
\end{cor}

\begin{proof}
If $f$ is  adic, by Proposition  \ref{ccad}, we have that for all $n \in \NN$ the following diagrams are cartesian squares:
\begin{diagram}[height=2em,w=2em,labelstyle=\scriptstyle]
\FX   & \rTto{f}       & \FY\\
\uTto &                & \uTto\\ 
X_{n} &  \rTto^{f_{n}} & Y_{n}.\\ 
\end{diagram}
Then by base-change (Proposition \ref{infsoritcor}) we have that $f_{n}$ is smooth,  for all $n \in \NN$. The  reciprocal follows from  last proposition.
\end{proof}
In the    next example we show that the  reciprocal of  Proposition  \ref{lfn} does not hold in general. 

\begin{ex} \label{peynopen}
Let $K$  a field and $ \BA_{K}^{1}= \spec(K[T])$. Given the  closed subset $X =V( \langle T \rangle) \subset  \BA_{K}^{1}$,   Corollary \ref{complpetal} implies that the  canonical completion morphism 
\[\BD_{K}^{1}  \xto{\kappa} \BA_{K}^{1}\] 
of $\BA_{K}^{1}$ along  $X$ is \'etale. However,  picking in $\BA_{K}^{1}$ the  Ideal of  definition $0$,   the morphisms 
\[\spec(K[T]/ \langle T \rangle ^{n+1}) \xto{\kappa_{n}} \BA_{K}^{1}\]
are not flat, whence it follows that $\kappa_{n}$ can not be smooth for all $n \in \NN$ (see Proposition \ref{flplano}).
\end{ex}

Our next goal will be to determine the relation between smoothness of  a morphism 
\[f= \dirlim {n \in \NN} f_{n}\]  and that of  $f_{0}$ (Corollaries \ref{flf0l} and \ref{corf0imfpl}). In order to do that, we need to characterize smoothness locally.

\begin{propo} \label{pligplafibr}
Let $f:\FX \to \FY$ be a pseudo finite type morphism in $\sfn$. Given $x \in \FX$, $y = f(x)$ the following conditions are equivalent:
\begin{enumerate}
\item
The morphism $f$ is smooth  at $x$.
\item
$\CO_{\FX,x}$ is a  formally smooth $\CO_{\FY,y}$-algebra for the adic topologies.
\item
$\widehat{\CO_{\FX,x}}$ is a formally smooth $\widehat{\CO_{\FY,y}}$-algebra  for the adic topologies.
\item
The morphism $f$ is  flat at $x$ and $f^{-1}(y)$ is a $k(y)$-formal scheme smooth  at $x$.
\end{enumerate}
\end{propo}

\begin{proof}
Since is a local  question and $f$ is of  pseudo finite type, we may assume  that $f: \FX= \spf(A) \to \FY= \spf(B)$ is in $\sfna$, with 
$A = B \{T_{1},\, \ldots ,T_{r}\}[[Z_{1},\ldots, Z_{s}]]/I$ and 
$I \subset B':= B \{\mathbf{T}\}[[\mathbf{Z}]]$ an ideal (Proposition \ref{ptf}). Let $\fp \subset A$ be the open prime  ideal corresponding to $x$, $\fq \subset B'$ the open prime such that  $\fp = \fq/I$ and $\fr \subset B$  the  open prime  ideal corresponding to   $y$. 

(1) $\Rightarrow$ (3) Replacing $\FX$ by a sufficiently small  open neighborhood   of  $x$ we may suppose  that $A$ is a  formally smooth $B$-algebra.  Then, by \cite[(\textbf{0}.19.3.5)]{EGA41} we have that $A_{\fp}$ is a formally smooth $B_{\fr}$-algebra  and    Lemma \ref{condinfalgcompl} implies  that $\widehat{\CO_{\FX,x}}=\widehat{A_{\fp}}$ is a  formally smooth $\widehat{\CO_{\FY,y}}=\widehat{B_{\fr}}$-algebra.

(2) $\Leftrightarrow$ (3) It is a consequence of Lemma \ref{condinfalgcompl}.

(3) $\Rightarrow$ (1) By  Lemma \ref{condinfalgcompl}, assertion (3) is equivalent 
to  $A_{\fp}$ be a formally  smooth $B_{\fr}$-algebra. Then the Zariski's Jacobian criterion  (Proposition \ref{critjacobanillos}) implies that the  morphism of  $\widehat{A_{\fp}}$-modules
\[
\widehat{\frac{I_{\fq}}{I_{\fq}^{2}}} \to \Omega^{1}_{B'_{\fq}/B_{\fr}} \tc_{B'_{\fq}} A_{\fp}
\]
is right invertible and since $\widehat{A_{\fp}}$ is a faithfully  flat $A_{\{\fp\}}$-algebra and $(\om^{1}_{B'/B} \otimes_{B'} A)_{\{\fp\}}$ is a projective $A_{\{\fp\}}$-module (see  Proposition \ref{flplano}), by \cite[(\textbf{0}.19.1.14.(ii))]{EGA41} there results that the  morphism 
\[
(\frac{I}{I^{2}})_{\{\fp\}} \to (\om^{1}_{B'/B} \otimes_{B'} A)_{\{\fp\}}
\]
is right  invertible. From  the equivalence of categories \ref{propitrian}.(\ref{equivtrian}) we find an open  subset $\FU \subset \FX$ with $x \in \FU$ such that  the  morphism
\[
(\frac{I}{I^{2}})^{\tr} \to \om^{1}_{\BD^{s}_{\BA^{r}_{\FY}}/\FY} \otimes _{\CO_{\BD^{s}_{\BA^{r}_{\FY}}}} \CO_{\FX}
\]
is right  invertible in $\FU$. Now, by Zariski's Jacobian criterion for formal schemes (Corollary \ref{critjacob}) it follows that $f$ is smooth  in $\FU$.

(3)  $\Rightarrow$ (4) By  Lemma \ref{modifll} we have that $\widehat{\CO_{\FX,x}}$ is  $\widehat{\CO_{\FY,y}}$-flat and by \ref{caracterizlocalplanos}, $f$ is flat at $x$. Moreover from  \cite[(\textbf{0}.19.3.5)]{EGA41} we deduce that $\widehat{\CO_{\FX,x}}\otimes_{\widehat{\CO_{\FY,y}}} k(y)$ is a formally smooth $k(y)$-algebra  for the adic topologies or, equivalently, by (3) $\dimp$ (1) $f^{-1}(y)$ is a $k(y)$-formal scheme smooth  at $x$.

(4)  $\Rightarrow$ (3) 
By \ref{caracterizlocalplanos} we have that $A_{\fp}$ is a flat $B_{\fr}$-module and therefore, there results that
\begin{equation} \label{ecuacioncita}
0 \to \frac{I_{\fq}}{\fr I_{\fq}} \to \frac{B'_{\fq}}{\fr B'_{\fq}} \to \frac{A_{\fp} }{\fr A_{\fp}} \to 0
\end{equation}
is an exact  sequence. 
On the other hand, since  $f^{-1}(y)$ is a $k(y)$-formal scheme smooth  at $x$, from (1) $\Rightarrow$ (2) we deduce that $\widehat{\CO_{\FX,x}}\otimes_{\widehat{\CO_{\FY,y}}} k(y)$ is a formally smooth $k(y)$-algebra  for the adic topologies or, equivalently by  Lemma \ref{condinfalgcompl}, $A_{\fp}/\fr A_{\fp}$ is a formally smooth $k(\fr)$-algebra  for the adic topologies. Applying Zariski's Jacobian criterion (Proposition \ref{critjacobanillos}), we have that the  morphism 
\[
\widehat{\frac{I_{\fq}}{I_{\fq}^{2}}} \otimes_{B_{\fr}} k(\fr) \to (\om^{1}_{B'/B})_{\fq} \tc_{B'_{\fq}} A_{\fp} \otimes_{B_{\fr}} k(\fr)
\]
is right  invertible. Now, since  $(\om^{1}_{B'/B})_{\fq}$ is a projective $B'_{\fq}$-module (see Proposition \ref{flplano})
by \cite[(\textbf{0}.6.7.2)]{EGA1} we obtain that 
\[
\widehat{\frac{I_{\fq}}{I_{\fq}^{2}}}  \to \om^{1}_{B'_{\fq}/B_{\fr}} \tc_{\widehat{B'_{\fq}}} \widehat{A_{\fp}}
\]
is right invertible. Then, by the Zariski's Jacobian criterion, $A_{\fp}$ is a formally smooth  $B_{\fr}$-algebra for the adic topologies or, equivalently by  Lemma \ref{condinfalgcompl}, $\widehat{A_{\fp}}$ is a  formally smooth $\widehat{B_{\fr}}$-algebra.

\end{proof}

\begin{cor} \label{corpligplafibr}
Let $f:\FX \to \FY$ be a pseudo finite type morphism in $\sfn$. The following conditions are equivalent:
\begin{enumerate}
\item
The morphism $f$ is smooth.
\item
For all $x \in \FX,\, y=f(x)$, $\CO_{\FX,x}$ is a formally smooth $\CO_{\FY,y}$-algebra for the adic topologies.
\item
For all $x \in \FX,\, y=f(x)$, $\widehat{\CO_{\FX,x}}$ is a formally smooth  $\widehat{\CO_{\FY,y}}$-algebra for the adic topologies.
\item
The morphism $f$ is  flat and $f^{-1}(y)$ is a $k(y)$-formal scheme smooth  at $x$, for all $x \in \FX,\, y=f(x)$.
\end{enumerate}
\end{cor}

\begin{cor} \label{flf0l}
Let $f:\FX \to \FY$ be  an  \emph{adic}  morphism in $\sfn$ and $\CK \subset \CO_{\FY}$ an Ideal of  definition. Put \[f=\dirlim {n \in \NN} f_{n}\] using  the Ideals of  definition $\CK \subset \CO_{\FY}$ and $\CJ=f^{*}(\CK)\CO_{\FX} \subset \CO_{\FX}$. Then, the  morphism $f$ is smooth if, and only if, it is  flat and  the  morphism $f_{0}:X_{0} \to Y_{0}$ is smooth.
\end{cor}

\begin{proof}
Since $f$ is  adic,  the  diagram
\begin{diagram}[height=2em,w=2em,labelstyle=\scriptstyle]
\FX        &    	\rTto{f}	&   \FY\\
\uTto   &    			&	\uTto\\ 
X_{0}  &  \rTto^{f_{0}}  	& Y_{0}\\ 
\end{diagram}
is a cartesian square (Proposition  \ref{ccad}).
If $f$ is smooth,
by base-change there results  that $f_{0}$ is smooth. Moreover by Proposition \ref{flplano} we have that $f$ is flat. Reciprocally, if $f$ is adic by \ref{fibra} we have that $f^{-1}(y)=f_{0}^{-1}(y)$, for all $x \in \FX,\,y = f(x)$. Therefore, since $f_{0}$ is smooth, by base-change there  results that $f^{-1}(y)$ is a $k(y)$-scheme smooth at $x$, for all $x \in \FX,\,y = f(x)$ and applying    Corollary \ref{corpligplafibr} we conclude that $f$ is smooth. 
\end{proof}

The following example shows that the  last result is not true without assuming the "adic" hypothesis for the morphism $f$.

\begin{ex} \label{exf0lisonofliso}
Given $K$ a field, let $\mathbb{P}^{n}_{K}$ be the $n$-dimensional projective space and $X \subset \mathbb{P}^{n}_{K}$ a closed subscheme that is not smooth  over $K$. If we denote by $(\mathbb{P}^{n}_{K})_{/ X}$ the completion of  $\mathbb{P}^{n}_{K}$ along  $X$, by Proposition  \ref{complpe} we have that the  morphism 
\[
(\mathbb{P}^{n}_{K})_{/ X} \xto{f} \spec(K)
\]
is smooth  but $f_{0} : X \to \spec(K)$ is not smooth.
\end{ex}

\begin{cor} \label{corf0imfpl}
Given $f:\FX \to \FY$ a morphism  in $\sfn$ let $\CJ \subset \CO_{\FX}$ and $\CK \subset \CO_{\FY}$ be Ideals of  definition such that $f^{*}(\CK) \CO_{\FX} \subset \CJ$ and with this choice let us express \[f=\dirlim {n \in \NN} f_{n}.\] If $f$ is flat, $f_{0}: X_{0} \to Y_{0}$ is a smooth morphism  and $f^{-1}(y)=f_{0}^{-1}(y)$, for all $y = f(x)$ with $x \in \FX$, then $f$ is smooth.
\end{cor}

\begin{proof}
Since $f_{0}$ is smooth and $f^{-1}(y)=f_{0}^{-1}(y)$ for all $y = f(x)$ with $x \in \FX$,  we deduce that $f^{-1}(y)$ is a smooth  $k(y)$-scheme. Besides, by hypothesis $f$ is flat and  Corollary \ref{corpligplafibr} implies that $f$ is smooth.
\end{proof}
Example \ref{exf0lisonofliso} illustrates that the  reciprocal of the last  corollary does not hold.

Every morphism $f: X \to Y$ smooth in $\sch$ is locally a composition of  an \'etale morphism  $U \to \BA^{r}_{Y}$ and of  the projection $\BA^{r}_{Y} \to Y$. Proposition \ref{factpl} generalizes this fact for smooth morphisms in $\sfn$.

\begin{propo} \label{factpl}
Let $f:\FX \to \FY$ be a pseudo finite type morphism in $\sfn$. The  morphism $f$ is smooth  in $x \in \FX$ if, and only if,  there exists an open subset $\FU \subset \FX$ with $x \in \FU$ such that  $f|_{\FU}$ factors as
\[
\FU \xto{g} \BA^{n}_{\FY} \xto{p} \FY
\]
where $g$ is \'etale, $p$ is the canonical projection and $n = \rg (\om^{1}_{\CO_{\FX,x}/\CO_{\FY,f(x)}})$.
\end{propo}

\begin{proof}
As this is a local  question, we may assume that $f: \FX= \spf(A) \to \FY= \spf(B)$ is a smooth morphism   in $\sfna$. By  Proposition \ref{modifinito} and by  Lemma \ref{modifll} we have that $\om^{1}_{A/B}$ is a projective $A$-module of  finite type and therefore, if $\fp \subset A$ is the open prime ideal corresponding to $x$, there exists $h\in A \setminus  \fp$  such that  $\ga(\fD(h), \om^{1}_{\FX/\FY}) = \om^{1}_{A_{\{h\}}/B}$ is a free $A_{\{h\}}$-module of  finite type. Put $\FU = \spf(A_{\{h\}})$. Given $\{\hd a_{1},\hd a_{2},\ldots, \hd a_{n}\}$ a basis of  $\om^{1}_{A_{\{h\}}/B}$ consider the  morphism  of  $\FY$-formal schemes 
\[
\FU \xto{g} \BA^{n}_{\FY}=\spf(B\{T_{1},T_{2},\ldots,T_{n}\} )
\] 
defined  through   the category equivalence (\ref{equiv}) by the continuous morphism of topological $B$-algebras
\[
\begin{array}{ccc}
B\{T_{1},T_{2},\ldots,T_{n}\} & \to & A_{\{h\}}\\
T_{i} &\leadsto  &a_{i}.
\end{array}
\]
The morphism $g$ satisfies  that $f|_{\FU} = p \circ g$. Moreover, we deduce that $g^{*} \om^{1}_{\BA^{n}_{\FY}/\FY}\cong \om^{1}_{\FX/\FY}$ (see the definition of  $g$) and by  Corollary  \ref{modifisopet} we have that $g$ is \'etale. 
\end{proof}

\begin{rem}
This result has appeared in local form in \cite[Proposition 1.11]{y}.
\end{rem}

\begin{cor} \label{dimrango}
Let $f:\FX \to \FY$ be a smooth morphism   at $x \in \FX$ and $y = f(x)$. Then
\[
\dim_{x} f = \rg (\om^{1}_{\CO_{\FX,x}/\CO_{\FY,y}}).
\]
\end{cor}

\begin{proof}
Put  $n= \rg (\om^{1}_{\CO_{\FX,x}/\CO_{\FY,y}})$. By Proposition  \ref{factpl} there exists $\FU \subset \FX$ with $x \in \FU$ such that  $f|_{\FU}$ factors as $\FU \xto{g} \BA^{n}_{\FY} \xto{p} \FY$ where $g$ is an  \'etale  morphism and $p$ is the canonical projection. Applying Proposition \ref{modifll} we have that $f|_{\FU}$ and  $g$ are flat  morphisms and therefore, 
\[
\begin{array}{ccccc}
\dim_{x} f & = &\dim \widehat{\CO_{\FX,x} }\otimes_{\widehat{\CO_{\FY,y} }} k(y) &= &\dim \widehat{\CO_{\FX,x} } - \dim\widehat{\CO_{\FY,y} }\\
\dim_{x} g &=  &\dim \widehat{\CO_{\FX,x} }\otimes_{\widehat{\CO_{\BA^{n}_{\FY},g(x)}}} k(g(x))&= &\dim \widehat{\CO_{\FX,x} } - \dim\widehat{\CO_{\BA^{n}_{\FY},g(x)}}.
\end{array}
\]
Now, since $g$ is unramified  by   Corollary \ref{pnrdim0} we have that $\dim_{x} g =0$ and therefore $\dim_{x} f  = \dim\widehat{\CO_{\BA^{n}_{\FY},g(x)}}- \dim\widehat{\CO_{\FY,y} }=n$.
\end{proof}

\begin{propo} \label{ecppl}
Let $f:\FX \to \FY$ be a morphism of  pseudo finite type and $\FX' \inc \FX$ a  closed immersion given by the  Ideal $\CI \subset \CO_{\FX}$ and put $f'=f|_{\FX'}$. If $f$ is smooth  at $x \in \FX'$, $n = \dim_{x} f$ and $y=f(x)$ the following conditions are equivalent:
\begin{enumerate}
\item
The morphism $f'$ is smooth  at $x$ and $\dim_{x} f'^{-1}(y)=n-m$
\item
The sequence of $\CO_{\FX}$-modules
\[
0 \to \frac{\CI}{\CI^{2}} \to \om^{1}_{\FX/\FY} \otimes_{\CO_{\FX}}  \CO_{\FX'} \to \om^{1}_{\FX'/\FY}  \to 0
\]
is exact\footnote{Let $(X,\CO_{X})$ be a ringed space. We say that the sequence of  $\CO_{X}$-Modules \[0 \to \CF \to \CG \to \CH \to 0\] is exact at $x \in X$ if $0 \to \CF_{x} \to \CG_{x} \to \CH_{x} \to 0$ is exact.} at $x$ and, on a neighborhood of $x$, the Modules are locally free of  ranks $m,\, n$ and $n-m$, respectively.
\end{enumerate}
\end{propo}

\begin{proof}
Since $f:\FX \to \FY$ is a smooth morphism   at $x$, replacing $\FX$, if  necessary, by a smaller neighborhood  of  $x$, we may assume that $f: \FX=\spf(A) \to \FY=\spf(B)$ is a morphism in $\sfna$  smooth at $x$ and that $\FX'= \spf(A/I)$. Therefore, applying Proposition \ref{flplano} and  Corollary \ref{dimrango} we have that $\om^{1}_{\FX/\FY}$ is a locally free $\CO_{\FX}$-Module of  rank $n$. 

Let us prove that (1) $\Rightarrow $ (2). Replacing, again, if it is necessary $\FX'$ by a smaller neighborhood  of  $x$, we may assume that $f':\FX' \to\FY$ is a smooth morphism   and then, applying Proposition \ref{flplano} and  Corollary \ref{dimrango}, there results that $\om^{1}_{\FX'/\FY}$ is a locally free $\CO_{\FX'}$-Module of  rank $n-m$. On the other hand, Zariski's Jacobian criterion for formal schemes (Corollary \ref{critjacob}) implies that the sequence 
\[
0 \to \frac{\CI}{\CI^{2}} \to \om^{1}_{\FX/\FY} \otimes_{\CO_{\FX}}  \CO_{\FX'} \to \om^{1}_{\FX'/\FY}  \to 0
\]
is exact and split, from  where we deduce that $\CI/\CI^{2}$ is a locally free $\CO_{\FX'}$-Module of  rank $m$.

Reciprocally, applying \cite[(\textbf{0}.5.5.4)]{EGA1} to  the Second Fundamental Exact  Sequence (\ref{sef2}) associated to the morphisms $\FX' \inc \FX \xto{f} \FX$, we deduce  that there exists an open formal subscheme $\FU \subset \FX'$ with $x \in \FU$ such that  
\[
0 \to (\frac{\CI}{\CI^{2}})|_{\FU} \to (\om^{1}_{\FX/\FY} \otimes_{\CO_{\FX}}  \CO_{\FX'})|_{\FU} \to (\om^{1}_{\FX'/\FY})|_{\FU} \to 0
\]
is exact  and split. From  Zariski's Jacobian criterion (Corollary \ref{critjacob}) it follows that $f'|_{\FU}$ is smooth  and therefore, $f'$ is smooth  at $x$.
\end{proof}


Locally, a pseudo finite type morphism  $f: \FX \to \FY$,  factors as  $\FU \overset{j} \inc \BD^{r}_{\BA^{s}_{\FY}} \xto{p} \FY$ where $j$ is a  closed immersion (see Proposition \ref{ptf}.(\ref{ptf1})). In  Corollary \ref{criteriojacobiano} we provide a criterion in terms of a rank of a matrix that determines whether  $\FU$ is smooth over $\FY$ or not. Next, we establish  the necessary preliminaries.

\begin{parraf}
Given $\FY= \spf(A)$ in $\sfna$ consider $\FX \subset \BD^{s}_{\BA^{r}_{\FY}}$ a closed formal subscheme given by an Ideal $\CI=I^{\tr}$ with $I=\langle g_{1}\,, g_{2},\, \ldots,\, g_{k} \rangle \subset A\{T_{1},\, \ldots ,T_{r}\}[[Z_{1},\ldots, Z_{s}]]$. In  Example \ref{modifafindis} we have seen that  
\[\{\hd T_{1},\, \hd T_{2},\,\ldots ,\hd T_{r},\, \hd Z_{1},\,\hd Z_{2},\,\ldots, \hd Z_{s}\}\]
is a basis of  $\om^{1}_{A\{\mathbf{T}\}[[\mathbf{Z}]]/A}$ and also,  that given $g \in A\{\mathbf{T}\}[[\mathbf{Z}]]$ it holds that:
\[
\hd g = \sum_{i=1}^{r}  \frac{\partial g}{\partial T_{i}} \hd T_{i} + \sum_{j=1}^{s}  \frac{\partial g}{\partial Z_{j}} \hd Z_{j}
\]
where $\hd$ is the complete canonical derivation of  $A\{\mathbf{T}\}[[\mathbf{Z}]]$ over $A$.
For all $g \in A\{\mathbf{T}\}[[\mathbf{Z}]]$, $s \in \{\hd T_{1},\, \hd T_{2},\,\ldots ,\hd T_{r},\, \hd Z_{1},\,\hd Z_{2},\,\ldots, \hd Z_{s} \}$ and $x \in \FX$,  denote by $\frac{\partial g}{\partial s}(x) $ the image  of  $\frac{\partial g}{\partial s} \in A\{\mathbf{T}\}[[\mathbf{Z}]]$ in $k(x)$. We will call
\begin{equation*}
\Jac_{\FX/\FY}(x)=\begin{pmatrix}
	\frac{\partial g_{1}}{\partial T_{1}}(x) & \frac{\partial g_{1}}{\partial T_{2}}(x) & \ldots & 
	\frac{\partial g_{1}}{\partial Z_{1}}(x) & \frac{\partial g_{1}}{\partial Z_{2}}(x) & \ldots &	\frac{\partial g_{1}}{\partial Z_{s}}(x) \\
	\frac{\partial g_{2}}{\partial T_{1}}(x) & \frac{\partial g_{2}}{\partial T_{2}}(x) & \ldots & 
	\frac{\partial g_{2}}{\partial Z_{1}}(x) & \frac{\partial g_{2}}{\partial Z_{2}}(x) & \ldots &
	\frac{\partial g_{2}}{\partial Z_{s}}(x) \\
	\vdots  & \vdots  & \ddots & \vdots & \vdots  & \ddots  & \vdots\\
	\frac{\partial g_{k}}{\partial T_{1}}(x) & \frac{\partial g_{k}}{\partial T_{2}}(x) & \ldots & 
	\frac{\partial g_{k}}{\partial Z_{1}}(x) & \frac{\partial g_{k}}{\partial Z_{2}}(x) & \ldots &	\frac{\partial g_{k}}{\partial Z_{s}}(x) \\
	\end{pmatrix}.
\end{equation*}
the \emph{Jacobian matrix of  $\FX$ over $ \FY$ at $x$}. This matrix depends on  the chosen generators of  $I$ and therefore, the notation  $\Jac_{\FX/\FY}(x)$ is not completely accurate. 
\end{parraf}

\begin{cor}{(Jacobian criterion for the affine formal space and the affine formal disc).} \label{criteriojacobiano}
Given $\FY= \spf(A)$ in $\sfna$ and an ideal $I=\langle g_{1}\,, g_{2},\, \ldots,\, g_{k} \rangle \subset A\{T_{1},\, \ldots ,T_{r}\}[[Z_{1},\ldots, Z_{s}]]$, let $\FX \subset  \BD^{s}_{\BA^{r}_{\FY}}$ be the closed formal subscheme defined by $\CI=I^{\tr}$. The following assertions  are equivalent:
\begin{enumerate}
\item
The morphism $f: \FX  \to\FY$ is smooth  at $x$ and $\dim_{x} f= r+s-l$.
\item
There exists a set $\{g_{1},\,, g_{2},\, \ldots,\, g_{l}\} \subset \{g_{1}\,, g_{2},\, \ldots,\, g_{k}\}$ such that  $\CI_{x}=\langle g_{1,}\,, g_{2,},\, \ldots,\, g_{l,} \rangle\CO_{\FX,x}$ and $\rg(\Jac_{\FX/\FY}(x)) = l$.
\end{enumerate}
\end{cor}

\begin{proof}
Suppose (1). By Proposition \ref{ecppl} we have that the sequence
\[
0 \to \frac{\CI}{\CI^{2}} \to \om^{1}_{ \BD^{s}_{\BA^{r}_{\FY}}/\FY} \otimes_{\CO_{ \BD^{s}_{\BA^{r}_{\FY}}}} \!\!\! \CO_{\FX} \to \om^{1}_{\FX/\FY}  \to 0
\]
is exact at $x$ and the corresponding $\CO_{\FX}$-Modules are locally free, in a neighborhood of $x$, of ranks $l$, $r+s$ and $r+s-l$, respectively. Therefore, there results that 
\begin{equation} \label{sexcritjacobkev}
0 \to \frac{\CI}{\CI^{2}}\otimes_{\CO_{\FX}} k(x) \to   
\om^{1}_{ \BD^{s}_{\BA^{r}_{\FY}}/\FY} \otimes_{ \CO_{\BD^{s}_{\BA^{r}_{\FY}}}} \!\! k(x) \to  \om^{1}_{\FX/\FY}  \otimes_{\CO_{\FX}} k(x)  \to 0
\end{equation}
is an exact  sequence of  $k(x)$-vector spaces  of  dimension $ l,\, r+s,\, r+s-l$, respectively. Thus, there exists a set $\{g_{1},\, g_{2},\, \ldots,\, g_{l}\} \subset \{g_{1},\, g_{2},\, \ldots,\, g_{k}\}$ such that  $\{g_{1}(x),\, g_{2}(x),\, \ldots,\, g_{l}(x)\}$ provides a basis of  $\CI/\CI^{2}\otimes_{\CO_{\FX}} k(x)$  at $x$ and by   Nakayama's lemma there results that $\CI_{x}=\langle g_{1,}\,, g_{2,},\, \ldots,\, g_{l,} \rangle\CO_{\FX,x}$. Besides, from  the exactness of  the sequence (\ref{sexcritjacobkev}) and from  the equivalence of categories \ref{propitrian}.(\ref{equivtrian}) we deduce that the set 
\[\{\hd g_{1}(x),\, \hd g_{2}(x),\, \ldots,\, \hd g_{l}(x)\} \subset \om^{1}_{ A\{\mathbf{T}\}[[\mathbf{Z}]]/A} \otimes_{A\{\mathbf{T}\} [[\mathbf{Z}]]} k(x)\] 
is linearly independent and, therefore, we have that $\rg(\Jac_{\FX/\FY}(x)) = l$. 

Conversely, from  the Second Fundamental Exact Sequence (\ref{sef2}) associated to the morphisms $\FX \inc \BD^{s}_{\BA^{r}_{\FY}} \to \FY$ we get the exact sequence
\[
 \frac{\CI}{\CI^{2}}\otimes_{\CO_{\FX}} k(x) \to   
\om^{1}_{ \BD^{s}_{\BA^{r}_{\FY}}/\FY} \otimes_{ \CO_{\BD^{s}_{\BA^{r}_{\FY}}}} \!\! k(x) \to  \om^{1}_{\FX/\FY}  \otimes_{\CO_{\FX}} k(x)  \to 0.
\]
On the other hand, since  $\rg(\Jac_{\FX/\FY}(x)) = l$,   we have  that 
\[\{\hd g_{1}(x),\,, \hd g_{2}(x),\, \ldots,\, \hd g_{l}(x)\} \subset \om^{1}_{ A\{\mathbf{T}\}[[\mathbf{Z}]]/A} \otimes_{A\{\mathbf{T}\} [[\mathbf{Z}]]} k(x)\] is a linearly independent set. Extending this set to a basis of the vector space $\om^{1}_{ A\{\mathbf{T}\}[[\mathbf{Z}]]/A} \otimes_{A\{\mathbf{T}\} [[\mathbf{Z}]]} k(x)$, by Nakayama's lemma we find a basis 
\[\mathcal{B} \subset \om^{1}_{ \BD^{s}_{\BA^{r}_{\FY}}/\FY} \,\text{ such that }\, \{\hd g_{1},\, \hd g_{2},\, \ldots, \, \hd g_{l}\} \subset \mathcal{B}\]
and therefore,  
\[\{\hd g_{1},\,, \hd g_{2},\, \ldots,\, \hd g_{l}\} \subset\om^{1}_{ A\{\mathbf{T}\}[[\mathbf{Z}]]/A} \otimes_{A\{\mathbf{T}\} [[\mathbf{Z}]]} A\{\mathbf{T}\}[[\mathbf{Z}]]/I\]
is a linearly independent set  at $x$. Thus the  set $\{g_{1,}\,, g_{2,},\, \ldots,\, g_{l,}\}$ gives a base  of  $\CI/\CI^{2}$ at $x$ and by the equivalence of categories \ref{propitrian}.(\ref{equivtrian}) we have that the sequence of $\CO_\FX$-Modules
\[
0 \to \frac{\CI}{\CI^{2}} \to \om^{1}_{ \BD^{s}_{\BA^{r}_{\FY}}/\FY} \otimes_{\CO_{ \BD^{s}_{\BA^{r}_{\FY}}}} \!\!\! \CO_{\FX} \to \om^{1}_{\FX/\FY}  \to 0
\]
is split exact at $x$ of locally free Modules  of  ranks $l$, $r+s$ and $r+s-l$, respectively. Applying Proposition \ref{ecppl} there results that $f$ is smooth  at $x$ and $\dim_{x} f= r+s-l$.
\end{proof}

Notice that the Jacobian criterion for the  affine formal space  and the affine formal disc (Corollary \ref{criteriojacobiano}) generalize the  Jacobian criterion for the  affine space in $\sch$ (\cite[Ch. VII, Theorem (5.14)]{at}).

The next  corollary will be used later in the  development of  the deformation theory  of  smooth formal schemes (Chapter \ref{cap4}).

\begin{cor} \label{deformloclis}
Let us consider in $\sfn$ a  closed immersion $\FY' \inc \FY$ and a smooth  morphism $f': \FX'  \to \FY'$. For all points $x \in \FX'$ there exists an open subset $\FU' \subset \FX'$  with $x \in \FU'$ and a   locally noetherian formal scheme $\FU$ smooth over $\FY$  such that  $\FU'= \FU\times_{\FY}\FY'$.
\end{cor}
\begin{proof}
Since it is a local question we may assume that the morphisms $\FY'=\spf(B') \inc \FY=\spf(B)$ and $f': \FX'=\spf(A') \to \FY'=\spf(B')$ are in $\sfna$ and that there exist $r,\, s \in \NN$ such that, if we put $B' \{\mathbf{T}\}[[\mathbf{Z}]] :=B'\{T_{1},\, \ldots ,T_{r}\}[[Z_{1},\ldots, Z_{s}]]$, the  morphism $f'$ factors as 
\[\FX' =\spf(A') \inc \BD^{s}_{\BA^{r}_{\FY'}}\!\! = \spf(B' \{\mathbf{T}\}[[\mathbf{Z}]]) \xto{p'} \FY'=\spf(B'),\]
where $\FX'  \inc \BD^{s}_{\BA^{r}_{\FY'}}$ is a closed subscheme given by an Ideal $\CI'=(I')^{\tr} \subset \CO_{\BD^{s}_{\BA^{r}_{\FY'}}}$\!\!, with $I' \subset B' \{\mathbf{T}\}[[\mathbf{Z}]]$ an ideal, and $p'$ is the canonical projection (see Proposition \ref{ptf}.(\ref{ptf1})).  Fix $x \in \FX'$.
As $f'$ is smooth, by  Jacobian criterion for the affine formal space and the affine formal disc (Corollary \ref{criteriojacobiano}),  we have that there exists $\{g'_{1},\, g'_{2},\, \ldots,\, g'_{l}\} \subset I'$ such  that: 
\begin{equation} \label{rangojacobiano}
\langle g_{1,}\,, g_{2,},\, \ldots,\, g_{l,} \rangle\CO_{\FX,x} = I'_{x} \qquad \textrm{and} \qquad \rg (\Jac_{\FX'/\FY'}(x)) = l
\end{equation}
Replacing, if it is necessary, $\FX'$ by a smaller affine open neighborhood  of  $x$ we may assume that $I' =\langle g'_{1},\, g'_{2},\, \ldots,\, g'_{l}\rangle$. Let $\{g_{1},\, g_{2},\, \ldots,\, g_{l}\}\subset B\{\mathbf{T}\}[[\mathbf{Z}]]$ such that  
$
g_{i} \in B\{\mathbf{T}\}[[\mathbf{Z}]] \leadsto g'_{i} \in   B'\{\mathbf{T}\}[[\mathbf{Z}]]
$
through the continuous homomorphism of  rings $B\{\mathbf{T}\}[[\mathbf{Z}]] \epi B'\{\mathbf{T}\}[[\mathbf{Z}]]  $ induced by $B \epi B'$. 
Put $I:= \langle g_{1},\, g_{2},\, \ldots,\, g_{l}\rangle \subset B \{\mathbf{T}\}[[\mathbf{Z}]]$ and  $\FX:= \spf(B \{\mathbf{T}\}[[\mathbf{Z}]]/I)$. It holds  that   $\FX' \subset  \FX$ is a closed subscheme and that in  the   diagram 
\begin{diagram}[height=2.3em,w=2.3em,p=0.3em,labelstyle=\scriptstyle]
\FX	&   \rTinc	 &  \BD^{s}_{\BA^{r}_{\FY}}  &\rTto{p} & \FY	\\
\uTinc &          	&  \uTinc  &	&\uTinc\\ 
\FX'	&   \rTinc	 &  \BD^{s}_{\BA^{r}_{\FY'}}  &\rTto{p'} & \FY'	\\
\end{diagram}
the squares are cartesian. From  (\ref{rangojacobiano}) we deduce  that $\rg (\Jac_{\FX/\FY}(x)) = l$ and, applying the  Jacobian criterion for the affine formal space and the affine formal disc,  there results that $\FX \to \FY$ is smooth at $x \in \FX$.  In order to finish the proof it suffices to take $\FU \subset \FX$ an open neighborhood   of  $x \in \FX$ such that  the  morphism $\FU \to \FY$ is smooth and $\FU' $ the  corresponding open set  in $\FX'$.

\end{proof}


\section{\'Etale morphisms} \label{sec33}
The  main part of  the results of  this section are consequence of  that obtained in Sections \ref{sec31} and \ref{sec32}. These results allow us to characterize in Section \ref{sec34} two important classes  of \'etale morphisms: open immersions and completion morphisms.

\begin{propo} \label{efn}
Given $f:\FX \to \FY$ in $\sfn$  let $\CJ \subset \CO_{\FX}$ and $\CK \subset \CO_{\FY}$ be Ideals of  definition with $f^{*}(\CK) \CO_{\FX} \subset \CJ$. Let us write 
\[f=\dirlim {n \in \NN} f_{n}.\] If $f_{n}:X_{n} \to Y_{n}$ is \'etale, $\forall n \in \NN$, then $f$ is \'etale.
\end{propo}

\begin{proof}
It is a consequence of  the corresponding results for unramified and smooth morphisms (Proposition \ref{nrfn} and Proposition \ref{lfn}).
\end{proof}

\begin{cor} \label{eadfn}
Let $f:\FX \to \FY$ be an  \emph{adic} morphism  in $\sfn$  and  $\CK \subset \CO_{\FY}$ an  Ideal of  definition. Consider $\{f_{n}\}_{n \in \NN}$ the  direct system  of  morphisms of  schemes associated to the Ideals of  definition $\CK \subset \CO_{\FY}$ and $\CJ=f^{*}(\CK)\CO_{\FX} \subset \CO_{\FX}$.
The morphism $f$ is \'etale if, and only if, the morphisms $f_{n}:X_{n} \to Y_{n}$ are \'etale $\forall n \in \NN$.
\end{cor}
\begin{proof}
It follows again from  the results for the unramified and smooth morphisms  (Proposition  \ref{nrfn} and  Corollary \ref{ladfn}).
\end{proof}

\begin{propo} \label{fetf0et}
Let $f:\FX \to \FY$ be an  \emph{adic} morphism  in $\sfn$ and let $f_{0}:X_{0} \to Y_{0}$ be the morphism of schemes associated to the Ideals of definition $\CK \subset \CO_{\FY}$ and $\CJ =f^{*}(\CK)\CO_{\FX} \subset \CO_{\FX}$. Then, $f$ is \'etale if, and only if, $f$ is flat and $f_{0}$ is \'etale.
\end{propo}

\begin{proof}
It is  deduced from  the analogous results for the unramified and smooth morphisms (Proposition  \ref{fnrf0nr} and Corollary   \ref{flf0l}).
\end{proof}

Note that example \ref{peynopen} on page \pageref{peynopen} shows that in the non adic case the last two results do not hold  and also, that in general, the  reciprocal of  Proposition \ref{efn} is not true.


\begin{propo}
Let $f$ be a pseudo finite type morphism in $\sfn$  and $\CJ\subset \CO_{\FX}$ and $\CK \subset \CO_{\FY}$ Ideals of  definition such that $f^{*}(\CK)\CO_{\FX} \subset \CJ$  and, with this choice, let us write \[f=\dirlim {n \in \NN} f_{n}.\] If $f_{0}:X_{0} \to Y_{0}$ is \'etale, $f$ is flat and  $f^{-1}(y)= f_{0}^{-1}(y)$, for all $y=f(x)$ with $x \in \FX$, then $f$ is \'etale.
\end{propo}

\begin{proof}
It follows from Corollary  \ref{corf0imfpnr} and Corollary \ref{corf0imfpl}.
\end{proof}

Example \ref{peynopen} shows that the  reciprocal of the last result  is not true. 
In the next Proposition a local characterization of  \'etale morphisms is given.

\begin{propo} \label{caractlocalpe}
Given $f:\FX \to \FY$ a morphism in $\sfn$ of  pseudo finite type, let  $x \in \FX$ and $y = f(x)$, the following conditions are equivalent:
\begin{enumerate}
\item [(1)]
 $f$ is \'etale  at $x$.
\item [(2)]
$\CO_{\FX,x}$ is a  formally \'etale $\CO_{\FY,y}$-algebra for the adic topologies.
\item [$(2')$]
$\widehat{\CO_{\FX,x}}$ is a formally \'etale $\widehat{\CO_{\FY,y}}$-algebra  for the adic topologies.
\item  [(3)]
$f$ is flat at $x$ and $f^{-1}(y)$ is a $k(y)$-formal scheme \'etale  at $x$.
\item [(4)]
$f$ is flat and unramified  at $x$.
\item [$(4')$]
$f$ is flat at $x$ and $(\om^{1}_{\FX/\FY})_{x}=0$. 
\item [(5)]
$f$ is smooth  at $x$ and a quasi-covering at $x$.

\end{enumerate}
\end{propo}

\begin{proof}
Applying Proposition \ref{caraclocalpnr} and Proposition \ref{pligplafibr} we have  that
\[(5)  \Leftarrow (1) \dimp (2) \dimp (2') \dimp (3) \imp (4) \dimp (4').\]
To show  (4) $\Rightarrow$ (5),  by  Corollary \ref{corpnrimplcr}. it suffices to prove that $f$ is smooth  at $x$. By hypothesis, we have that $f$ is unramified  at $x$ and by Proposition \ref{caraclocalpnr},  there results that $\widehat{\CO_{\FX,x}} \otimes_{\widehat{\CO_{\FY,y}}} k(y) = k(x)$ and $k(x)|k(y)$ is a  finite separable extension, therefore, formally \'etale. Since $f$ is flat at $x$, by Proposition \ref{pligplafibr} we conclude that $f$ is smooth  at $x$. 

Finally, we prove that (5) $\Rightarrow$ (1). It suffices to check that $f$ is unramified  at $x$ or, equivalently by Proposition \ref{caraclocalpnr}, that  $\widehat{\CO_{\FX,x}} \otimes_{\widehat{\CO_{\FY,y}}} k(y) = k(x)$ and that $k(x)|k(y)$ is a finite separable extension. As $f$ is smooth at $x$, applying Proposition \ref{pligplafibr}, we have that $\widehat{\CO_{\FX,x}}$ is a formally smooth $\widehat{\CO_{\FX,x}}$-algebra  for the adic topologies. Then by  base-change there results that $\widehat{\CO_{\FX,x}} \otimes_{\widehat{\CO_{\FY,y}}} k(y)$ is a  formally smooth $k(y)$-algebra. By Remark \ref{confinfhafin} we have that $\widehat{\CO_{\FX,x}} \otimes_{\widehat{\CO_{\FY,y}}}k(y)$ is a formally smooth $k(y)$-algebra  for the topologies given by the maximal ideals and from \cite[Lemma1, p. 216]{ma2} it holds that $\widehat{\CO_{\FX,x}} \otimes_{\widehat{\CO_{\FY,y}}}k(y)$ is a regular local ring. Besides, by hypothesis we have  that  $\widehat{\CO_{\FX,x}} \otimes_{\widehat{\CO_{\FY,y}}}k(y)$ is a finite $k(y)$-module, therefore, an artinian ring, so $\widehat{\CO_{\FX,x}} \otimes_{\widehat{\CO_{\FY,y}}}k(y)=k(x)$. Since   $k(x) = \widehat{\CO_{\FX,x}} \otimes_{\widehat{\CO_{\FY,y}}}k(y)$ is a formally smooth $k(y)$-algebra  we have that $k(x)|k(y)$ is a separable extension (\emph{cf.} \cite[(\textbf{0}.19.6.3)]{EGA41}).  
\end{proof} 

\begin{cor} \label{corcaractlocalpe}
Let $f:\FX \to \FY$ be a  pseudo finite type morphism in $\sfn$. The following conditions are equivalent:
\begin{enumerate}
\item [(1)]
 $f$ is \'etale.
\item [(2)]
For all $x \in \FX,\, y=f(x)$, $\CO_{\FX,x}$ is a formally \'etale $\CO_{\FY,y}$-algebra for  the adic topologies.
\item [$(2')$]
For all $x \in \FX,\, y=f(x)$, $\widehat{\CO_{\FX,x}}$ is a formally \'etale  $\widehat{\CO_{\FY,y}}$-algebra for  the adic topologies.
\item  [(3)]
For all $x \in \FX,\, y=f(x)$, $f^{-1}(y)$ is a $k(y)$-formal scheme \'etale  at $x$ and $f$ is flat. 
 \item [(4)]
$f$ is flat and unramified.
\item [$(4')$]
$f$ is flat  and $\om^{1}_{\FX/\FY}=0$.
\item [(5)]
$f$ is smooth  and a quasi-covering.

\end{enumerate}

\end{cor}

\begin{ex}\label {pcf+plnope}
Given a field, $K$ the  canonical morphism of  projection $\BD^{1}_{K} \to \spec(K)$ is smooth,  pseudo quasi-finite  but  it is not  \'etale.
\end{ex}

In $\sch$ a morphism is \'etale if, and only if, it is smooth and quasi-finite. The  previous example shows that in $\sfn$ there are smooth  and pseudo quasi-finite morphisms that are not \'etale. It is, thus, justified  that  quasi-coverings in $\sfn$ (see Definition \ref{defncuasireves}) is the right generalization of the notion of quasi-finite morphisms in $\sch$.

\section{Special properties of  \'etale morphisms} \label{sec34}

In this  section  we characterize two types  of morphisms: open immersions and completion morphisms in terms of \'etale morphisms. While the first of these results (Theorem \ref{caractencab}) extends a well-known fact in the context of usual schemes, the other  (Theorem \ref{caracmorfcompl}) is a special feature of formal schemes. The section also contains Theorem \ref{teorequivet} in which  we establish an equivalence of  categories between the \'etale adic $\FY$-formal schemes and the  \'etale $Y_{0}$-schemes for a locally noetherian formal scheme  \[\FY= \dirlim {n \in \NN} Y_{n}.\]
These statements are significant by themselves and also play a basic role for the structure theorems on next section.

We begin with two results that will be used in the proof  of  the remainder results of  this  chapter.  
\begin{propo} \label{levantpeusform}
In $\sfn$ let us consider a formally \'etale morphism  $f:\FX \to \FY$, a morphism $g\colon \FS \to \FY$ and $\CL \subset \CO_{\FS}$ an Ideal of definition of $\FS$. Let us write with respect to $\CL$
\[\FS = \dirlim {n \in \NN} S_{n}.\] 
If  $h_{0}: S_{0} \to \FX$ is a morphism in $\sfn$ that makes commutative the  diagram 
\begin{diagram}[height=2em,w=2em,p=0.3em,labelstyle=\scriptstyle]
S_{0}         & \rTinc    & \FS	\\
\dTto^{h_{0}} &           & \dTto_{g}	\\ 
\FX           & \rTto^{f} & \FY\text{,}\\  
\end{diagram}
where $S_{0} \inc \FS$ is the  canonical  closed immersion, then
there exists a unique $\FY$-morphism $l: \FS \to \FX$  in $\sfn$ such that  $l|_{S_{0}}=h_{0}$.
\end{propo}

\begin{proof}
By induction on $n$ we are going to construct a collection of  morphisms $\{h_{n}: S_{n} \to \FX\}_{n \in \NN}$ such that the diagrams
\begin{diagram}[height=2.2em,w=2.5em,labelstyle=\scriptstyle]
S_{n-1}	&			&		&		&\\
		&\rdTto	\rdTto(4,2) \rdTto(2,4)_{h_{n-1}}	&    &      &   \\		&			&S_{n}	&\rTinc & \FS\\
        &       & \dTto_{h_{n}} &	& \dTto_{g}\\
        &       & \FX	&\rTto^{f}&	\FY\\
\end{diagram}
are commutative.
For $n =1$, by  \ref{lgfefnr}  there exists a unique morphism $h_{1}: S_{1}  \to \FX$ such that   $h_{1} |_{S_{0}}= h_{0}$ and $g|_{S_{1}}= f \circ h_{1}$. Take $n \in \NN$ and suppose constructed for all $0< k <n$  morphisms $h_{k}: S_{k} \to \FX$ such that $h_{k} |_{S_{k-1}}= h_{k-1}$ and $g|_{S_{k}}= f \circ h_{k}$. Again by \ref{lgfefnr} there exists an unique morphism $h_{n}: S_{n} \to \FX$ such that $h_{n} |_{S_{n-1}}= h_{n-1}$ and $g|_{S_{n}}= f \circ h_{n}$. It is straightforward that \[l:= \dirlim {n \in \NN} h_{n}\] is a morphism of  formal schemes and is the  unique such that the  diagram
\begin{diagram}[height=2em,w=2em,p=0.3em,labelstyle=\scriptstyle]
S_{0}	&   \rTinc	    & \FS	\\
\dTto^{h_{0}} &           \ldTto_{h}  	    &	\dTto_{g}	\\ 
\FX  	& \rTto 	    &\FY \\  
\end{diagram}
commutes.
\end{proof}

\begin{rem}
The last  result is also true if $f$ is formally smooth.
\end{rem}

\begin{cor} \label{iso0peis}
Let $f:\FX \to \FY$ be an   \'etale morphism in $\sfn$ and  $\CJ \subset \CO_{\FX}$ and $\CK \subset \CO_{\FY}$ Ideals of  definition with $f^{*}(\CK)\CO_{\FX} \subset \FJ$ such that the relevant morphism $f_{0}:X_{0} \to Y_{0}$ is an isomorphism. Then $f$ is an isomorphism. 
\end{cor}

\begin{proof}
By Proposition \ref{levantpeusform} there exists a unique morphism $g: \FY \to \FX$ such that  the  following diagram is commutative
\begin{diagram}[height=2em,w=2.5em,p=0.3em,labelstyle=\scriptstyle]
Y_{0}	           &   \rTinc	    & \FY	\\
\dTto^{f_{0}^{-1}} &   \ldTto(2,4)^{g}  &	\\
 X_{0}             &                 & \dTto_{1_{\FY}}\\
\dTinc             &                 & \\
\FX  	           & \rTto^{f} 	    &\FY\\  
\end{diagram}
Then,  by  Corollary \ref{corbasinf} it holds that $g$ is an \'etale morphism. Thus, applying  Proposition \ref{levantpeusform} we have that there exists an unique morphism $f': \FX \to \FY$ such that   the  following diagram is commutative
\begin{diagram}[height=2em,w=2.5em,p=0.3em,labelstyle=\scriptstyle]
X_{0}	      &   \rTinc	            & \FX \\
\dTto^{f_{0}} &   \ldTto(2,4)^{f'}   &	\\
Y_{0}         &             & \dTto_{1_{\FX}}	\\
 \dTinc       &             & \\
\FY  	      & \rTto^{g}   &\FX\\  
\end{diagram}
so, necessarily $f=f' $ and $f$  is a isomorphism.
 \end{proof}
In $\sch$ open immersions are characterized as being the \'etale morphisms which are radical (see \cite[(17.9.1)]{EGA44}). 
In the  following theorem we  extend this characterization and relate open immersions in formal schemes (Definition \ref{defencab}) with their counterparts in schemes. 
\begin{thm} \label{caractencab}  
Let $f:\FX \to \FY$ be a morphism in $\sfn$. The following conditions are equivalent:
\begin{enumerate}
\item
$f$ is an open immersion.
\item
$f$ is adic, flat and if $\CK \subset \CO_{\FY}$ an Ideal of  definition and $\CJ= f^{*}(\CK) \CO_{\FX} \subset \CO_{\FX}$,  the  associated  morphism of  schemes $f_{0}:X_{0} \to Y_{0}$ is an open immersion.
\item
$f$ is adic \'etale and radical.
\item
There are $\CJ \subset \CO_{\FX}$ and $\CK \subset \CO_{\FY}$ Ideals of  definition satisfying that $ f^{*}(\CK) \CO_{\FX} \subset \CJ$ such that the morphisms $f_{n}:X_{n} \to Y_{n}$ are open immersions, for all $n \in \NN$.
\end{enumerate}
\end{thm}

\begin{proof}
The  implication  (1) $\Rightarrow$ (2) is immediate. Given $\CK \subset \CO_{\FY}$ an Ideal of  definition, suppose  (2)  and let us show (3). Since $f_{0}$ is an open immersion, is radical, so,  $f$ is radical (see Definition \ref{rad} and its attached paragraph). On the other hand, $f_{0}$ is an \'etale morphism and then, by Proposition \ref{fetf0et} we have that $f$ is \'etale. 
Let us prove that (3) $\Rightarrow$ (4). Given $\CK \subset \CO_{\FY}$ an Ideal of  definition and $\CJ = f^{*}(\CK) \CO_{\FX}$, by  Corollary \ref{eadfn} the morphisms $f_{n}:X_{n} \to Y_{n}$ are  \'etale, for all $n \in \NN$. The morphisms $f_{n}$ are radical for all $n \in \NN$ (see Definition \ref{rad}) and thus by \cite[(17.9.1)]{EGA44} it follows that $f_{n}$ is an open immersion, for all $n \in \NN$. Finally, suppose  that  (4) holds and let us see that $f$ is an open immersion. With the notations of  (4), there exists an open subset $U_{0} \subset Y_{0}$ such that  $f_{0}$ factors as
\[
X_{0} \xto{f'_{0}} U_{0} \overset{i_{0}} \inc Y_{0}
\]
where $f'_{0}$ is an isomorphism and $i_{0}$ is the canonical inclusion. Let $\FU \subset \FY$ be the  open formal subscheme with underlying topological space $U_{0}$. Since the  open immersion $i: \FU \to \FY$ is \'etale, then Proposition \ref{levantpeusform} implies that there exists a morphism $f': \FX \to \FU$ of formal schemes such that  the  diagram:
\begin{diagram}[height=1.5em,w=2em,labelstyle=\scriptstyle]
\FX	&	     &\rTto^{\qquad f}&			&	   &\FY\\
\uTinc	&\rdTto(3,2)^{f'}&	&  		&\ruTinc^{i}&	\uTinc \\
		&	    &            & \FU&  & \\
		&	    &            &\uTto	&	   & 	          \\
X_{0}	&	     &\hLine^{\qquad f_{0}}&\VonH	&\rTto   &Y_{0}\\
		&    \rdTto(3,2)_{f'_{0}}       &	&	  &\ruTinc^{i_{0}} &    \\
		&	    &		&U_{0}  				&	 &   \\
\end{diagram} 
is commutative.
On the other hand, since the morphisms $f_{n}$ are \'etale, for all $n \in \NN$, Proposition \ref{efn} implies that $f$ is \'etale. Then by  Corollary \ref{propibasinf} we have that $f'$ is \'etale  and applying  Corollary \ref{iso0peis},  $f'$ is an isomorphism and therefore, $f$ is an open immersion.
\end{proof}

Every completion morphism is a pseudo closed immersion that is flat (\emph{cf.} Proposition \ref{caractcom}). Next, we prove that this condition is also sufficient. Thus, we obtain a criterion to determine whether a $\FY$-formal scheme $\FX$ is the completion of  $\FY$ along  a closed formal subscheme. The completion morphisms will play an important role   in the structure theorems  of  unramified morphisms,  \'etale morphisms and  smooth morphisms of section \ref{sec35}.
 
\begin{thm} \label{caracmorfcompl}
Let $f:\FX \to \FY$ be a morphism in $\sfn$ and let $\CJ \subset \CO_{\FX}$ and $\CK \subset \CO_{\FY}$ be Ideals of  definition such that $f^{*}(\CK) \CO_{\FX} \subset \CJ$. Let us write $f_{0}:X_{0} \to Y_{0}$ the corresponding morphism of ordinary schemes. The following conditions are equivalent:
\begin{enumerate}
\item
There exists a closed formal subscheme $\FY' \subset \FY$ such that    $f$ is the  morphism of  completion of  $\FY$ along  $\FY'$ and therefore, $\FX = \FY_{/\FY'}$.
\item
The morphism $f$ is a flat pseudo closed immersion.
\item
The morphism $f$ is  \'etale  and $f_{0}:X_{0} \to Y_{0}$ is a  closed immersion.
\item
The morphism $f$ is a smooth pseudo  closed immersion.
\end{enumerate}
\end{thm}

\begin{proof}
The implication (1) $\Rightarrow$ (2) is Proposition   \ref{caractcom}. Let us show that (2) $\Rightarrow$ (3). Since $f$ is a pseudo  closed immersion, by  Corollary \ref{pecigf0ecnr}  we have that $f$ is unramified. Then as  $f$ is flat,  Corollary \ref{corcaractlocalpe} implies that $f$ is \'etale. The equivalence (3) $\Leftrightarrow$ (4) is consequence of Corollary \ref{pecigf0ecnr}.  Finally, we show that (3) $\Rightarrow$ (1). By hypothesis, the   morphism $f_{0}: X_{0} \to \FY$ is a  closed immersion. Consider  $\kappa: \FY_{/ X_{0}} \to \FY$ the  morphism of  completion of  $\FY$ along  $X_{0}$ and let us prove that $\FX $ and $\FY_{/ X_{0}}$ are $\FY$-isomorphic. By Proposition \ref{complpe} the  morphism $\kappa$ is \'etale  so,  applying Proposition \ref{levantpeusform}, we have that there exists a $\FY$-morphism $\varphi: \FX \to \FY_{/ X_{0}}$ such that  the  following diagram is commutative
\begin{diagram}[height=1.5em,w=2em,labelstyle=\scriptstyle]
\FX	    &       &\rTto^{\qquad f}&			&	   &\FY\\
\uTinc  &\rdTto(3,2)^{\varphi}&	&  		&\ruTto^{\kappa}&	\uTinc \\
		&	    &            & \FY_{/ X_{0}}	&  & \\
		&	    &            &\uTto	&	   & 	          \\
X_{0}	&	    &\hLine^{\qquad f_{0}}&\VonH	&\rTto   &Y_{0}\\
		&    \rdTto(3,2)_{\varphi_{0}=1_{X_{0}}}       &	&	  &\ruTinc^{f_{0}} &    \\
		&	    &		&X_{0}  				&	 &   \\
\end{diagram}
From Proposition \ref{propibasinf} there results that $\varphi$ is \'etale  and then by  Corollary \ref{iso0peis} we get that $\varphi$ is an isomorphism.
\end{proof}

\begin{rem}
In  implication (3) $\Rightarrow$ (1) we have proved that given $\FY$ in $\sfn$, a closed formal subscheme $\FY' \subset \FY$ defined by the  Ideal $\CI \subset \CO_{\FY}$,  then for every Ideal of  definition $\CK \subset \CO_{\FY}$  of  $\FY$, it holds that 
\[\FY_{/\FY'} = \FY_{/Y'_{0}}\] 
where $Y'_{0} =  (\FY',\CO_{\FY}/(\CI+\CK))$.
\end{rem}

\begin{parraf}
Given  an  \'etale scheme $Y$ and $Y_{0} \subset Y$ a closed subscheme with the  same topological space, the  functor $X \leadsto X \times_{Y} Y_{0}$ defines an equivalence between the  category of   \'etale $Y$-schemes and the category of   \'etale $Y_{0}$-schemes  (see \cite[(18.1.2)]{EGA44}). In the  next theorem we extend this equivalence to the category of  locally noetherian formal schemes.
A special case of this theorem, namely $\FY$ is smooth over a noetherian ordinary base scheme, appears in \cite[Proposition 2.4]{y}.
\end{parraf}

\begin{thm} \label{teorequivet}
Let $\FY$ be in $\sfn$ and $\CK \subset \CO_{\FY}$ an Ideal of  definition respect to which $\FY = \dirlim {n \in \NN} Y_{n}$. Then the  functor
\[
\begin{array}{ccc}
\textrm { \'etale adic }\FY\textrm{-formal schemes}& \xto{F} &\textrm{\'etale }Y_{0}\textrm{-schemes}\\
\FX& \leadsto& \FX \times_{\FY} Y_{0}\\
\end{array}
\]
is an equivalence of  categories.
\end{thm}

\begin{proof}
By \cite[IV, \S 4, Theorem 1]{mcl} it suffices to prove that: a) $F$ is full and faithful;   b) Given $X_{0}$ an  \'etale $Y_{0}$-scheme  there exists an  \'etale adic $\FY$-formal scheme $\FX$ such that  $F(\FX)=\FX \times_{\FY}Y_{0} \cong X_{0}$.

The assertion a) is an immediate consequence of  Proposition \ref{levantpeusform}.

Let us show  b). Given  $X_{0}$ an \'etale $Y_{0}$-scheme  in $\sch$ by \cite[(18.1.2)]{EGA44} there exists $X_{1}$ a locally noetherian \'etale $Y_{1}$-scheme   such that   $X_{1}\times_{Y_{1}}Y_{0}  \cong X_{0}$. Reasoning by induction on $n \in \NN$ and using \cite[\emph{loc. cit.}]{EGA44}, we get a family $\{X_{n}\}_{n \in \NN}$ such that, for all $n \in \NN$ $X_{n}$ is a locally noetherian \'etale $Y_{n}$-scheme and  $X_{n}\times_{Y_{n}}Y_{n-1}  \cong X_{n-1}$, for all $n \in \NN$. Then $\FX:= \dirlim {n \in \NN} X_{n}$ is a locally noetherian adic $\FY$-formal scheme (by \cite[(10.12.3.1)]{EGA1}), 
\[\FX \times_{\FY}Y_{0} \underset{\ref{cclim}}= \dirlim {n \in \NN} (X_{n} \times_{Y_{n}} Y_{0}) =X_{0}\] and 
$\FX$ is an  \'etale $\FY$-formal scheme (see Proposition \ref{efn}).
\end{proof}

\begin{cuesab}
It seems plausible that there is a theory of an algebraic fundamental group for formal schemes that classifies \'etale \emph{adic} surjective maps to a noetherian formal scheme $\FX$. If this is the case, the previous theorem would imply that it agrees with the fundamental group of $X_0$. We also conjecture the existence of a bigger fundamental group classifying arbitrary \'etale surjective maps to a noetherian formal scheme $\FX$, that would give additional information about $\FX$.
\end{cuesab}

\begin{cor}
Let $f:\FX \to \FY$ be an \'etale morphism  in $\sfn$. Given $\CJ \subset \CO_{\FX},$ and $\CK \subset \CO_{\FY}$ Ideals of  definition such that $f^{*}(\CK) \CO_{\FX} \subset \CJ$, if the induced  morphism $f_{0}:X_{0} \to Y_{0}$ is \'etale, then $f$ is adic \'etale.
\end{cor}

\begin{proof}
By  Theorem \ref{teorequivet} there is an adic \'etale morphism $f':\FX' \to\FY$ in $\sfn$  such that  $\FX' \times_{\FY} Y_{0}= X_{0}$. Therefore by Proposition \ref{levantpeusform} there exists a morphism of formal schemes $g: \FX \to \FX'$ such that  the  diagram
\begin{diagram}[height=1.5em,w=2em,labelstyle=\scriptstyle]
\FX	&	     &\rTto^{f}&			&	   &\FY\\
\uTinc		&\rdTto(3,2)^{g}&	&  		&\ruTto^{f'}&	\uTinc \\
		&	    &            & \FX'	&  & \\
		&	    &            &\uTto	&	   & 	          \\
X_{0}	&	     &\hLine&\VonH	&\rTto^{f_{0}}   &Y_{0}\\
		&    \rdTto(3,2)_{g_{0}=1_{X_{0}}}       &	&	  &\ruTto^{f'_{0}} &    \\
		&	    &		&X_{0}  				&	 &   \\
\end{diagram} 
is commutative.
Applying Proposition \ref{propibasinf} we have that $g$ is \'etale  and from Corollary \ref{iso0peis} we deduce that $g$ is an isomorphism and therefore, $f$ is adic \'etale.
\end{proof}

\begin{cor} \label{fnetfet}
Let $f:\FX \to \FY$ be a morphism in $\sfn$. The  morphism $f$ is adic \'etale if, and only if,  there exist $\CJ \subset \CO_{\FX}$ and $\CK \subset \CO_{\FY}$ Ideals of  definition with $f^{*}(\CK) \CO_{\FX} \subset \CJ$  such that the induced morphisms   $f_{n}:X_{n} \to Y_{n}$ are \'etale, for all $n \in \NN$.
\end{cor}

\begin{proof}
If $f$ is \'etale adic, given $\CK \subset \CO_{\FY}$ an Ideal of  definition, let us consider $\CJ=f^{*}(\CK) \CO_{\FX}$ the corresponding Ideal of definition of  $\FX$. By  base change we have that the morphisms  $f_{n}:X_{n} \to Y_{n}$ are  \'etale, for all $n \in \NN$. The  reciprocal is a consequence of Proposition \ref{efn} and of the previous Corollary.
\end{proof}

Theorem \ref{teorequivet} says that given \[\FY= \dirlim {n \in \NN} Y_{n}\] in $\sfn$ and $X_{0}$ an  \'etale $Y_{0}$-scheme there exists  a unique (up to isomorphism)  \'etale $\FY$-formal scheme $\FX$ such that  $\FX \times_{\FY} Y_{0}=X_{0}$.  But, what does it happen when $X_{0}$ is a smooth $Y_{0}$-scheme? 

\begin{propo} \label{equivlocal}
Let $\FY$ be in $\sfn$ and with respect to an Ideal of  definition $\CK \subset \CO_{\FY}$ let us write \[\FY = \dirlim {n \in \NN} Y_{n}.\] Given $f_{0}:X_{0} \to Y_{0}$ a  morphism in $\sch$  smooth at $x \in X_{0}$, there exists an open subset  $U_{0} \subset X_{0}$, with $x \in U_{0} $ and an smooth adic  $\FY$-formal scheme   $\FU$  such that  $\FU \times_{\FY} Y_{0} \cong U_{0}$.
\end{propo}

\begin{proof}
 Since this is a local question in $\FY$, we may assume that $\FY= \spf(B)$ is in $\sfna$, $\CK = K^{\tr}$ with $K \subset B$ an  ideal of  definition of the adic ring $B$, $B_{0} = B/K$ and $f_{0}:X_{0}= \spec(A_{0}) \to Y_{0} = \spec(B_{0})$ is a morphism in $\scha$ smooth at $x \in X_{0}$.  By 
 Proposition \ref{factpl} there exists an open subset $U_{0} \subset X_{0}$ with $x \in U_{0}$ such that  $f_{0} |_{U_{0}}$ factors us 
 \[
 U_{0} \xto{f'_{0}} \BA_{Y_{0}}^{n}= \spec(B_{0}[\mathbf{T}]) \xto{p_{0}} Y_{0}
 \]
 where $f'_{0}$ is an \'etale morphism and $p_{0}$ is the canonical projection. The  morphism $p_{0}$ lifts to a morphism of  projection $p: \BA_{\FY}^{n}= \spf(B\{\mathbf{T}\}) \to \FY$ such that the  following diagram is cartesian
 \begin{diagram}[height=2em,w=2em,labelstyle=\scriptstyle]
	 &	     		    &	\BA_{\FY}^{n}		&\rTto^{p}	   	  &\FY\\
	 &	     		    &	\uTinc		  	&  		   	  &\uTinc \\
U_{0}& \rTto^{f'_{0}}	    &   \BA_{Y_{0}}^{n}         & \rTto^{p_{0}}	  & Y_{0}\\
\end{diagram}
Applying  Theorem \ref{teorequivet}, there exists $\FU$ a locally noetherian \'etale adic $\BA_{\FY}^{n}$-formal scheme  such that  $U_{0} \cong \FU \times_{\BA_{\FY}^{n}} \BA_{Y_{0}}^{n}$. Then $\FU$ is an smooth adic $\FY$-formal scheme  such that  $U_{0} \cong \FU \times_{\FY} Y_{0}$.
\end{proof}

\section[Structure theorems]{Structure theorems of the infinitesimal lifting properties} \label{sec35}

In the previous sections we have seen  that the infinitesimal conditions  of  an adic  morphism \[\dirlim {n \in \NN} f_{n}\] in $\sfn$  are determined by the  morphism $f_{0}$ (\emph{cf.} Propositions \ref{fnrf0nr} and \ref{fetf0et} and   Corollary \ref{flf0l}). However, examples \ref{framfonram} and \ref{peynopen} illustrate that, in absence of  the adic hypothesis this correspondence does not hold any more. 
In this section we prove the main results  of  this chapter in which we provide  a local description of  the infinitesimal lifting properties according to the completion morphisms of  (characterized in  Theorem \ref{caracmorfcompl}) and to   the  infinitesimal \emph{adic}  lifting properties. 

The next theorem transfers the local description of unramified morphisms known in the case of schemes (\cite[(18.4.7)]{EGA44}) to the framework  of  formal schemes.

\begin{thm} \label{tppalnr}
Let $f:\FX \to \FY$ be a morphism in $\sfn$ unramified  at $x \in \FX$. Then there exists an open subset $\FU \subset \FX$ with $x \in \FU$ such that  $f|_{\FU}$ factors as
\[
\FU \xto{\kappa} \FX' \xto{f'} \FY
\]
where $\kappa$ is a pseudo  closed immersion  and $f'$ is an \emph{adic} \'etale morphism.
\end{thm}

\begin{proof}
Let $\CJ \subset \CO_{\FX}$ and $\CK \subset \CO_{\FY}$ be Ideals of  definition such that $f^{*}(\CK) \CO_{\FX} \subset \CJ$.
The morphism of  schemes $f_{0}$ associated to these Ideals is unramified at $x$ (Proposition \ref{nrfn}) and by \cite[(18.4.7)]{EGA44} there exists an open set  $U_{0}\subset X_{0}$ with $x \in U_{0}$ such that  $f_{0}|_{U_{0}}$ factors as
\[
U_{0}\rTinc^{\kappa_{0}} X'_{0} \rTto^{f'_{0}} Y_{0}
\]
where $\kappa_{0}$ is a  closed immersion and $f'_{0}$ is an \'etale morphism. Theorem \ref{teorequivet} implies that there exists  an  \'etale adic morphism $f':\FX' \to\FY$ in $\sfn$ such that  $\FX' \times_{\FY} Y_{0} = X'_{0}$. Then if $\FU \subset \FX$ is the  open formal scheme  with underlying topological space  $U_{0}$, by  Proposition \ref{levantpeusform} there exists a morphism $\kappa: \FU \to \FX'$ such that  the  following diagram is commutative:
\begin{diagram}[height=1.5em,w=2em,labelstyle=\scriptstyle]
\FU	&	     &\rTto^{f|_{\FU}}&			&	   &\FY\\
\uTinc	&\rdTto(3,2)^{\kappa}&	&  		&\ruTto^{f'}&	\uTinc \\
		&	    &            & \FX'	&  & \\
		&	    &            &\uTto	&	   & 	          \\
U_{0}	&	     &\hLine&\VonH	&\rTto^{f_{0}|_{U_{0}}}   &Y_{0}\\
		&    \rdTinc(3,2)_{\kappa_{0}}       &	&	  &\ruTto^{f'_{0}} &    \\
		&	    &		&X'_{0}  				&	 &   \\
\end{diagram}
Since $f$ is unramified,  by Proposition \ref{propibasinf} it holds that $\kappa$ is unramified. On the other hand, $\kappa_{0}$ is a  closed immersion and applying  Corollary \ref{pecigf0ecnr}  we have that $\kappa$ is a pseudo  closed immersion.
\end{proof}

As a consequence of the last result we obtain the following  local description for \'etale morphisms. 

\begin{thm} \label{tppalet}
Let $f:\FX \to \FY$ be a morphism in $\sfn$ \'etale at $x \in \FX$. Then there exists an open subset $\FU \subset \FX$ with $x \in \FU$ such that  $f|_{\FU}$ factors as
\[
\FU \xto{\kappa} \FX' \xto{f'} \FY
\]
where $\kappa$ is a completion morphism and $f'$ is an \emph{adic} \'etale morphism.
\end{thm}

\begin{proof}
By the  last theorem  we have that there exists an open formal subscheme $\FU \subset \FX$ with $x \in \FU$ such that  $f|_{\FU}$ factors as 
\[
\FU \xto{\kappa} \FX' \xto{f'} \FY
\]
where $\kappa$ is a pseudo  closed immersion and $f'$ is an adic \'etale morphism. Then since $f|_{\FU}$ is \'etale  and $f'$ is an adic \'etale morphism, by Proposition \ref{propibasinf} we have that $\kappa$ is \'etale  and applying  Theorem  \ref{caracmorfcompl} there results that $\kappa$ is a completion morphism.
\end{proof}

\begin{thm} \label{tppall}
Let $f:\FX \to \FY$ be a morphism in $\sfn$ smooth at $x \in \FX$. Then there exists an open subset $\FU \subset \FX$ with $x \in \FU$ such that  $f|_{\FU}$ factors as
\[
\FU \xto{\kappa} \FX' \xto{f'} \FY
\]
where $\kappa$ is a completion morphism and $f'$ is an \emph{adic} smooth  morphism.
\end{thm}

\begin{proof}
By Proposition \ref{factpl} there exists an open formal subscheme $\FV \subset \FX$ with $x \in \FV$ such that  $f|_{\FV}$ factors as
\[
\FV \xto{g} \BA^{n}_{\FY} \xto{p} \FY
\]
where $g$ is \'etale  and $p$ is the canonical projection. Applying  the last Theorem to the morphism $g$ we conclude that there exists an open subset $\FU \subset \FX$ with $x \in \FU$ such that  $f|_{\FU}$  factors as 
\[
\FU \xto{\kappa} \FX' \xto{f''} \BA_{\FY}^{n} \xto{p}\FY
\]
where $\kappa$ is a completion morphism, $f''$ is an adic \'etale morphism and $p$ is the canonical projection, from where it follows that $f'= f'' \circ p$ is adic smooth.
\end{proof}

\chapter[Basic Deformation Theory]{Basic Deformation Theory  of  smooth formal schemes}\label{cap4} \setcounter{equation}{0}

As  in $\sch$ the case of smooth  morphisms is a basic ingredient  of  a theory of  deformation in $\sfn$. The  problem consists on  constructing morphisms that extend a given morphism over a smooth  formal scheme to a base which is an ``infinitesimal neighborhood" of  the original. Then, questions of existence and uniqueness can be analyzed. The answers are expressed via cohomological invariants. They can be explicitly computed using the  \v{C}ech complex. Another group of questions that are treated is the construction of  schemes over an   infinitesimal neighborhood of  the base. The existence of such lifting will be controlled by an element belonging to a $2^{\text{nd}}$\!-order cohomology  group. We will restrict ourselves to the use of elementary methods. This will force us to restrict to the separated case which suffices for a large class of applications.

All the results of  this chapter generalize the well-known analogous  results  in $\sch$. Our exposition follows the outline given in  \cite[p. 111--113]{ill}. 

\section{Lifting of  morphisms} \label{sec41}
\begin{parraf}
Given $f:\FX \to \FY$ a morphism of  pseudo finite type consider  a commutative diagram in $\sfn$
\begin{equation}\label{hipotesis}
\begin{diagram}[height=2em,w=2em,p=0.3em,labelstyle=\scriptstyle]
\FZ_{0}       & \rTinc^{i} & \FZ	\\
\dTto^{u_{0}} & \ldTdash   & \dTto\\ 
\FX           & \rTto^{f}  & \FY\\  
\end{diagram}
\end{equation}
where $\FZ_{0}\inc \FZ$ is a closed formal subscheme given by a square zero Ideal $\CI \subset \CO_{\FZ}$.  A morphism $u: \FZ \to \FX$ is a \emph{lifting of  $u_{0}$  over $\FY$} if it makes this diagram commutative. For instance, if $f$ is \'etale for all morphisms $u_{0}$ under the last  conditions there always exists a unique  lifting (see  Corollary \ref{corlevfor}). 

So the basic question is: When can we  guarantee uniqueness and existence of a lifting for a $\FY$-morphism $u_{0}: \FZ_{0} \to \FX$? Observe that in the diagram above, $i\colon \FZ_0 \to \FZ$ is the identity as topological map and, therefore, we may identify $i_{*}\CO_{\FZ_{0}} \equiv \CO_{\FZ_{0}}$. Through this  identification we have that  the  Ideal $\CI$ is a  $\CO_{\FZ_{0}}$-module and  $\CI = i_{*} \CI$.

In the next paragraph it is shown that if 
$\Hom_{\CO_{\FZ_{0}}}(u_{0}^{*}\om^{1}_{\FX/\FY}, \CI)=0$, then the  lifting is unique. And Proposition \ref{propobstruclevant} establishes that there exists an obstruction in 
$\ext^{1}_{\CO_{\FZ_{0}}}(u_{0}^{*} \om^{1}_{\FX/\FY}, \CI)$ to the existence of such a lifting. 
\end{parraf} 

\begin{parraf} \label{parrfderivlevant}
Let us continue to consider the situation depicted in diagram (\ref{hipotesis}). Note that if $u: \FZ \to \FX$ is a lifting of $\FU_0$ over $\FY$ then $u_{0}^{*}\CI = u^{*}\CI$. In view of this identification, from \ref{derivversuslevantalg} we deduce the  following results about the uniqueness of lifting:
\begin{enumerate} 
\item \label{parrfderivlevant1}
If $u: \FZ \to \FX,\, v: \FZ \to \FX$ are liftings of  $u_{0}$ over $\FY$, the  morphism
\[
 \CO_{\FX} \xto{u^{\sharp} - v^{\sharp}} u_{0*} \CI
\]
is a continuous $\FY$-derivation. Lemma  \ref{imagdirecompl} and  Theorem \ref{modrepr}  imply that there exists an unique morphism of  $\CO_{\FX}$-modules $\phi: \om^{1}_{\FX/\FY} \to u_{0*} \CI$  that makes the  diagram
\begin{diagram}[height=2em,w=2em,p=0.3em,labelstyle=\scriptstyle]
 \CO_{\FX}	&   \rTto^{\hd_{\FX/\FY}}	    & \om^{1}_{\FX/\FY}	\\
\dTto^{u^{\sharp} - v^{\sharp}} &           \ldTto_{\phi}  	    &\\ 
u_{0*} \CI	  	&  	    & \\  
\end{diagram} 
commutative.
\item \label{parrfderivlevant2}
Fixed a lifting $u: \FZ \to \FX$ of  $u_{0}$ and given $\phi: \om^{1}_{\FX/\FY} \to u_{0*} \CI$ a morphism of  $\CO_{\FX}$-modules, the map $v^{\sharp} := u^{\sharp}+ \phi\circ \hd_{\FX/\FY} $ defines another lifting of  $u_{0}$. 
\end{enumerate}
In summary, if there exists some lifting  $\FZ \to \FX$ of  $u_{0}$ over $\FY$,  then the  set of  liftings of  $u_{0}$ over $\FY$ is an affine space  over $\Hom_{\CO_{\FX}}(\om^{1}_{\FX/\FY} , u_{0 *}\CI)$ (or equivalently over   $\Hom_{\CO_{\FZ_{0}}}(u_{0}^{*}\om^{1}_{\FX/\FY} , \CI) $). 
\end{parraf}

\begin{rem}
Using the  language of  torsor theory the results in \ref{parrfderivlevant} say that 
 the  sheaf in $\FZ_{0}$ such that it associates to the open subset $\FU_{0} \subset \FZ_{0}$ the  set of  liftings $\FU \to \FX$  of  $u_{0}|_{\FU_{0}}$ over $\FY$ ---where $\FU \subset \FZ$ is the open subset corresponding to $\FU_{0}$--- is a pseudo torsor over $\shom_{\CO_{\FZ_{0}}}(u_{0}^{*}\om^{1}_{\FX/\FY} , \CI)$.
\end{rem}

When can we  guarantee for a diagram like (\ref{hipotesis}) the existence of a lifting of  $u_{0}$ over $\FY$? In Proposition \ref{levantform} we have shown that if $f$ is smooth and $\FZ_{0} \inc \FZ$ is in $\sfna$, then there exists lifting of  $u_{0}$ over $\FY$. So, the issue amounts to patching local data to obtain global data.

\begin{propo} \label{propobstruclevant}
Consider the commutative  diagram (\ref{hipotesis}) where $f:\FX \to \FY$ is a smooth morphism.  Then there exists an element (usually called ``obstruction'') $c_{u_{0}} \in \ext^{1}_{\CO_{\FZ_{0}}}(u_{0}^{*} \om^{1}_{\FX/\FY}, \CI)$ such that: $c_{u_{0}}=0$ if, and only if,  there exists $u: \FZ \to \FX$ a lifting of  $u_{0}$ over $\FY$.
\end{propo}

\begin{proof}
Let $\{\FU_{\alpha}\}_{\alpha}$ be an affine open  covering  of  $\FZ$ and $\FU_{\bullet} = \{\FU_{\alpha, 0}\}_{\alpha}$ the corresponding  affine open covering of $\FZ_{0}$ such that, for all $\alpha$, $\FU_{\alpha,0} \inc \FU_{\alpha}$ is a  closed immersion in $\sfna$ given by the square zero Ideal $\CI|_{\FU_{\alpha}}$.
Since $f$ is a smooth morphism, Proposition  \ref{levantform}  implies that for all $\alpha$ there exists a lifting $v_{\alpha}: \FU_{\alpha} \to \FX$ of  $u_{0}|_{\FU_{\alpha,0}}$ over $\FY$.  For all couples of  indexes $\alpha,\, \beta$ such that  $\FU_{\alpha \beta}:= \FU_{\alpha} \cap \FU_{\beta} \neq \varnothing$, if  we call $\FU_{\alpha \beta,0}$ the  corresponding open formal subscheme of $\FZ_{0}$ from \ref{parrfderivlevant}.(\ref{parrfderivlevant1}) we have that there exists a unique morphism of  $\CO_{\FX}$-Modules $\phi_{\alpha \beta}: \om^{1}_{\FX/\FY} \to (u_{0}|_{\FU_{\alpha \beta, 0}})_{*} (\CI|_{\FU_{\alpha \beta,o}})$ such that  the  following diagram is commutative
\begin{diagram}[height=2em,w=2em,p=0.3em,labelstyle=\scriptstyle]
 \CO_{\FX}	&   \rTto^{\hd_{\FX/\FY}\qquad}	    & \om^{1}_{\FX/\FY}	\\
\dTto^{(v_{\alpha}|_{\FU_{\alpha\beta}})^{\sharp} -(v_{\beta}|_{\FU_{\alpha\beta}})^{\sharp}} &  \ldTto_{\phi_{\alpha \beta}}   &\\ 
(u_{0}|_{\FU_{\alpha \beta, 0}})_{*} (\CI|_{\FU_{\alpha \beta,0}})	  	&  	    & \\  
\end{diagram}
Let $u_{0}^{*} \om^{1}_{\FX/\FY}|_{\FU_{\alpha \beta,0}} \to \CI|_{\FU_{\alpha \beta,0}}$ be the  morphism  of  $\CO_{\FU_{\alpha \beta,0}}$-Modules adjoint to  $\phi_{\alpha \beta}$ that,  with an abuse of notation we continue denoting by $\phi_{\alpha \beta}$. The family of  morphisms $\phi_{\FU_{\bullet}}:=(\phi_{\alpha \beta})$ satisfies the cocycle condition; that is, for any $\alpha,\, \beta,\, \gamma$ such that $\FU_{\alpha \beta \gamma,0}:= \FU_{\alpha,0} \cap \FU_{\beta,0} \cap \FU_{\gamma,0} \neq \varnothing$, we have that
\begin{equation} \label{datosrecolec1}
\phi_{\alpha \beta}|_{\FU_{\alpha \beta \gamma,0}} - \phi_{\alpha \gamma}|_{\FU_{\alpha \beta \gamma,0}}+\phi_{ \beta \gamma}|_{\FU_{\alpha \beta \gamma,0}}=0
\end{equation} 
so, $\phi_{\FU_{\bullet}} \in \check{\Z}^{1} (\FU_{\bullet}, \shom_{\CO_{\FZ_{0}}}(u_{0}^{*} \om^{1}_{\FX/\FY} ,\CI))$. 
Moreover, there results that its  class 
\[[\phi_{\FU_{\bullet}}] \in \check{\h}^{1} (\FU_{\bullet}, \shom_{\CO_{\FZ_{0}}}(u_{0}^{*} \om^{1}_{\FX/\FY} ,\CI))\]
does not depend  of  the liftings $\{v_{\alpha}\}_{\alpha}$. Indeed, if for all arbitrary $\alpha$,  $w_{\alpha}: \FU_{\alpha} \to \FX$ is another  lifting of  $u_{0}|_{\FU_{\alpha,0}}$ over $\FY$,  by  \ref{parrfderivlevant}.(\ref{parrfderivlevant1})  there exists a unique $\xi_{\alpha} \in \Hom_{\FX}(\om^{1}_{\FX/\FY}, (u_{0}|_{\FU_{\alpha , 0}})_{*} (\CI|_{\FU_{\alpha , 0}}))$ such that  $v^{\sharp}_{\alpha}-w^{\sharp}_{\alpha} = \xi_{\alpha} \circ \hd_{\FX/\FY}$. Then for all couples of  indexes $\alpha,\beta$ such that $\FU_{\alpha \beta} \neq \varnothing$ given 
\[\psi_{\alpha \beta} \in \Hom_{\FX}(\om^{1}_{\FX/\FY}, (u_{0}|_{\FU_{\alpha \beta, 0}})_{*} (\CI|_{\FU_{\alpha \beta , 0}}))\] 
such that  
\[w_{\alpha}^{\sharp}|_{\FU_{\alpha \beta}}-w_{\beta}^{\sharp}|_{\FU_{\alpha \beta}} = \psi_{\alpha \beta} \circ \hd_{\FX/\FY},\] 
we have that 
\[\psi_{\alpha \beta}=\phi_{\alpha \beta}+ \xi_{\beta}|_{\FU_{\alpha \beta}} -\xi_{\alpha}|_{\FU_{\alpha \beta}}\]
In other words, the cocycles $\psi_{\FU_{\bullet}}=(\psi_{\alpha \beta})$ and $\phi_{\FU_{\bullet}}$ differ  in a coboundary from where we conclude that  $[\phi_{\FU_{\bullet}}] = [\psi_{\FU_{\bullet}}] \in \check{\h}^{1} (\FU_{\bullet},\shom_{\CO_{\FZ_{0}}}(u_{0}^{*} \om^{1}_{\FX/\FY} ,\CI))$. 
With an analogous argument it is possible to prove that,  given a refinement of  $\FU_{\bullet}$, $\FV_{\bullet}$, then
$[\phi_{\FU_{\bullet}}]=[\phi_{\FV_{\bullet}}] \in \check{\h}^{1} (\FZ_{0}, \shom_{\CO_{\FZ_{0}}}(u_{0}^{*} \om^{1}_{\FX/\FY} ,\CI))$. We define:
\begin{align*}
c_{u_{0}}:=  [\phi_{\FU_{\bullet}}] &\in \, 
     \check{\h}^{1} (\FZ_{0},\shom_{\CO_{\FZ_{0}}}(u_{0}^{*} \om^{1}_{\FX/\FY} ,\CI)) = \\
     &= \h^{1} (\FZ_{0}, \shom_{\CO_{\FZ_{0}}}(u_{0}^{*} \om^{1}_{\FX/\FY} ,\CI))
     \tag{\cite[(5.4.15)]{te}}
\end{align*}
Since $f$ is smooth, Proposition \ref{flplano} implies  that  $\om^{1}_{\FX/\FY}$ is a locally free $\CO_{\FX}$-Module of finite rank, so, 
\[c_{u_{0}} \in \h^{1} (\FZ_{0}, \shom_{\CO_{\FZ_{0}}}(u_{0}^{*} \om^{1}_{\FX/\FY}, \CI))=\ext^{1}(u_{0}^{*} \om^{1}_{\FX/\FY}, \CI).\]

The element $c_{u_{0}}$ is the obstruction to the existence of a lifting of  $u_{0}$. If $u_{0}$ admits a lifting then it is clear that $c_{u_{0}}=0$. Reciprocally, suppose that $c_{u_{0}}=0$. From the family of morphisms $\{v_{\alpha} \colon \FU_{\alpha} \to \FX \}_{\alpha}$ we are going to construct a collection of liftings $\{u_{\alpha}: \FU_{\alpha} \to \FX\}_{\alpha}$ of $\{u_{0}|_{\FU_{\alpha,0}}\}_{\alpha}$ over $\FY$ that  will patch into a morphism $u: \FZ \to \FX$. By hypothesis, we have that there  exists 
$\{\varphi_{\alpha}\}_{\alpha} \in \check{\C}^{0}(\FU_{\bullet},\shom_{\CO_{\FZ_{0}}}(u_{0}^{*} \om^{1}_{\FX/\FY} ,\CI))$ such that for all couples of indexes $\alpha, \beta$ with $\FU_{\alpha \beta} \neq \varnothing$, 
\begin{equation} \label{datosrecolec2}
\varphi_{\alpha}|_{\FU_{\alpha \beta}} - \varphi_{\beta}|_{\FU_{\alpha \beta}}=\phi_{\alpha \beta}
\end{equation} 
For all $\alpha$, let $u_{\alpha}: \FU_{\alpha} \to \FX$ be the  morphism that agrees with $u_{0}|_{\FU_{\alpha,0}}$ as a topological map and is given by  
\[u^{\sharp}_{\alpha}:= v_{\alpha}^{\sharp} - \varphi_{\alpha} \circ \hd_{\FX/\FY}\] as a map of  topological ringed spaces. By \ref{parrfderivlevant}.(\ref{parrfderivlevant2}) we have that $u_{\alpha}$ is a lifting of  $u_{0}|_{\FU_{\alpha,0}}$ over $\FY$ for all $\alpha$ and,  from data (\ref{datosrecolec2}) and (\ref{datosrecolec1}) (for $\gamma = \beta$) we deduce that the morphisms $\{u_{\alpha}\}_\alpha$ glue into a morphism $u: \FZ \to \FX$.
\end{proof}

\begin{cor} \label{**}
With the hypothesis of  the last proposition, if $\FZ$ is an affine  formal scheme, then there exists a lifting of  $u_{0}$ over $\FY$.  
\end{cor}
\begin{proof}
Since $\FZ_{0}$ is an affine  formal scheme we have that 
\[\h^{1} (\FZ_{0},\shom_{\CO_{\FZ_{0}}}(u_{0}^{*} \om^{1}_{\FX/\FY} ,\CI))=0 \qquad \tag{ \cite[Corollary 3.1.8]{AJL1}}\] and, therefore, from the last Proposition we conclude. 
\end{proof}

Note that the last result also follows from Proposition \ref{levantform}.

\section{Lifting of smooth formal schemes } \label{sec42}

Given a smooth morphism $f_{0}:\FX_{0} \to \FY_{0}$ and also a closed immersion $\FY_{0} \inc \FY$ defined by a square zero Ideal $\CI$, one can pose the following questions: 
\begin{enumerate}
\item
Suppose that there exists a smooth $\FY$-formal scheme $\FX$  such that  $\FX \times_{\FY} \FY_{0} = \FX_{0}$. When  is $\FX$  unique? 
\item
Does there exist a smooth $\FY$-formal scheme $\FX$ such that it holds that $\FX \times_{\FY} \FY_{0} = \FX_{0}$?
\end{enumerate}

Relating to the first question, in Proposition \ref{deform4} we will show that  if $\ext^{1}(\om^{1}_{\FX_{0}/\FY_{0}}, f_{0}^{*}\CI)=0$ then $\FX$ is  unique  up to isomorphism.

Apropos of the second, for all  $x \in \FX_{0}$ there exists an open $\FU_{0} \subset \FX_{0}$ with $x \in \FU_{0}$ and a  locally noetherian smooth formal scheme $\FU$ over $\FY$  such that  $\FU_{0}= \FU \times_{\FY} \FY_{0}$ (see  Corollary \ref{deformloclis}). In general,  Proposition \ref{obstrext2} provides  an element in $\ext^{2}(\om^{1}_{\FX_{0}/\FY_{0}}, f_{0}^{*}\CI)$ whose vanishing is equivalent to the existence of such a $\FX$. In particular, whenever $f_{0}$ is in $\sfna$   Corollary  \ref{*} asserts the existence of  $\FX$.

\begin{parraf} \label{deform}
Assume that $f_{0}:\FX_{0} \to \FY_{0}$ is a smooth morphism and $i \colon \FY_{0} \inc \FY$ a  closed immersion given by a square zero Ideal $\CI \subset \CO_{\FY}$, hence, $\FY_{0}$ and $\FY$ have the same underlying topological space. If there exists   a  smooth morphism $f:\FX \to \FY$ in $\sfn$ such that  the   diagram
\begin{equation}\label{hipotesis2}
\begin{diagram}[height=2em,w=2em,p=0.3em,labelstyle=\scriptstyle]
\FX_{0}    & \rTto^{f_{0}} & \FY_{0}	\\
\dTinc^{j} &               & \dTinc \\ 
\FX        & \rTto^{f}     & \FY\\  
\end{diagram}
\end{equation}
is cartesian we say that $f:\FX \to \FY$ is a 
\emph{smooth lifting of $\FX_{0}$ over $\FY$}. 

Let us consider the notation in the commutative diagram (\ref{hipotesis2}), let  $f_{0} \colon \FX_{0} \to \FY_{0}$ be a smooth morphism, $i \colon \FY_{0} \inc \FY$ a  closed immersion given by a square zero Ideal $\CI \subset \CO_{\FY}$ and $f \colon \FX \to \FY$ a smooth lifting of  $\FX_{0}$ over $\FY$.
Observe that, since $f$ is flat, $j \colon \FX_{0} \to \FX$ is a  closed immersion given (up to isomorphism) by the square zero Ideal $f^{*}\CI$. The sheaf $\CI$ is a $\CO_{\FY_0}$-Module in a natural way, $f^{*}\CI$ is a $\CO_{\FX_0}$-Module and it is clear that $f^{*}\CI$ agree with $f_0^{*}\CI$ as $\CO_{\FX_0}$-Modules.
\end{parraf}

\begin{parraf} \label{deform1}
Denote by $\aut_{\FX_{0}}(\FX)$ the  group of  $\FY$-automorphisms of  $\FX$ that induce the identity  in $\FX_{0}$. In  particular, we have that $1_{\FX} \in \aut_{\FX_{0}}(\FX)$ and, therefore, $\aut_{\FX_{0}}(\FX) \neq \varnothing$.  There  exists a bijection  
\[\aut_{\FX_{0}}(\FX)  \iso \Hom_{\CO_{\FX}}(\om^{1}_{\FX/\FY}, j_{*} f_{0}^{*}\CI)\cong \Hom_{\CO_{\FX_{0}}}(\om^{1}_{\FX_{0}/\FY_{0}}, f_{0}^{*}\CI)\] 
defined using the  isomorphism (Theorem \ref{modrepr})
\[\Dercont_{\FY}(\CO_{\FX}, j_{*} f_{0}^{*}\CI)\cong\Hom_{\CO_{\FX}}(\om^{1}_{\FX/\FY}, j_{*} f_{0}^{*}\CI)\] 
and the map
\[
\begin{array}{ccc}
\aut_{\FX_{0}}(\FX) & \lto & \Dercont_{\FY}(\CO_{\FX}, j_{*} f_{0}^{*}\CI)\\
g  &\leadsto &g^{\sharp} - 1_{\FX}^{\sharp}.
\end{array}
\]
In fact, given a homomorphism of  $\CO_{\FX}$-Modules $\phi: \om^{1}_{\FX/\FY} \to j_{*}f_{0}^{*}\CI$, by  \ref{parrfderivlevant}.(\ref{parrfderivlevant1}) the  morphism $g: \FX \to \FX$ defined by 
\[g^{\sharp}= 1_{\FX}^{\sharp} + \phi \circ \hd_{\FX/\FY}\] 
which is the identity as topological map, is an automorphism (its inverse $g^{-1}$ is given by $(g^{-1})^{\sharp}=1_{\FX}^{\sharp} - \phi \circ \hd_{\FX/\FY}$).
\end{parraf}

\begin{parraf} \label{deform2}
If $\FX_{0}$ is in $\sfna$ and  $\FX_{0} \overset{\,\, j'}\inc  \FX' \overset{\, f'}\to \FY$ is another smooth lifting of $\FX_0$ over $\FY$, then there exists a $\FY$-isomorphism $g: \FX \xto{\sim} \FX'$ such that  $g|_{\FX_{0}} = j'$. 
Indeed, from Proposition  \ref{levantform} there are morphisms $g: \FX \to \FX'$, $g': \FX' \to \FX$ such that  the  following diagram is commutative:
\begin{diagram}[height=2em,w=2em,p=0.3em,labelstyle=\scriptstyle]
\FX_{0}		&   \rTinc^{j} & \FX					   &            &   \\
			&\rdTinc^{j'}  & \dTto^{g}  \uTto_{g'}  & \rdTto^{f} &   \\ 
			&			   & \FX'	  			   & \rTto^{f'} & \FY\\  
\end{diagram}
Then since   $g' \circ g|_{\FX_{0}}= j$, $g \circ g'|_{\FX_{0}}= j'$, by  \ref{parrfderivlevant}.(\ref{parrfderivlevant1}) it holds that there exists $\phi \in \Hom_{\CO_{\FX}}(\om^{1}_{\FX/\FY}, j_{*} f_{0}^{*}\CI)$ and   $\phi' \in \Hom_{\CO_{\FX'}}(\om^{1}_{\FX'/\FY}, j'_{*} f_{0}^{*}\CI)$ such that 
\[
(g' \circ g)^{\sharp}= 1_{\FX}^{\sharp} + \phi \circ \hd_{\FX/\FY} \qquad 
(g \circ g')^{\sharp}= 1_{\FX'}^{\sharp} + \phi' \circ \hd_{\FX'/\FY}
\]
From \ref{deform1} we deduce that  $g' \circ g
\in \aut_{\FX_{0}}(\FX)$, $g \circ g'\in \aut_{\FX_{0}}(\FX')$, therefore $g$ is an isomorphism.
\end{parraf}

\begin{parraf} \label{deform3}
In the previous setting, the  set of  $\FY$-isomorphisms of  $\FX$ in $\FX'$ that induce the identity in $\FX_{0}$ is an affine space   over $\Hom_{\CO_{\FX_{0}}}(\om^{1}_{\FX_{0}/\FY_{0}},f_{0}^{*} \CI)$ (or, equivalently over $\Hom_{\CO_{\FX'}}(\om^{1}_{\FX'/\FY}, j'_{*} f_{0}^{*}\CI)$, by adjunction). 
In fact, assume that  $g: \FX \to \FX'$ and $h: \FX \to \FX'$ are two $\FY$-isomorphisms such that $g|_{\FX_{0}} =h|_{\FX_{0}}=j'$. From  \ref{parrfderivlevant}.(\ref{parrfderivlevant1}) there exists a unique homomorphism of  $\CO_{\FX'}$-Modules  
\[\phi: \om^{1}_{\FX'/\FY} \to j'_{*} f_{0}^{*} \CI\] 
such that  
\[g^{\sharp}-h^{\sharp} =  \phi \circ \hd_{\FX'/\FY}.\] 
Reciprocally, if 
\[\phi \in  \Hom_{\CO_{\FX_{0}}}(\om^{1}_{\FX_{0}/\FY_{0}},f_{0}^{*} \CI) \cong \Hom_{\CO_{\FX'}}(\om^{1}_{\FX'/\FY}, j'_{*} f_{0}^{*}\CI)\] 
and  $g: \FX \to \FX'$ is a $\FY$-isomorphism with $g|_{\FX_{0}}=j'$, the  $\FY$-morphism $h: \FX \to \FX'$ defined by 
\[h^{\sharp} = g^{\sharp}+ \phi \circ \hd_{\FX'/\FY},\] 
which as topological space  map  is the identity, is an isomorphism. Indeed, it holds that \[(h\circ g^{-1})^{\sharp} - 1_{\FX'}^{\sharp} \in \Dercont_{\FY}(\CO_{\FX'}, j'_{*} f_{0}^{*}\CI)\] and that 
\[(g^{-1} \circ h)^{\sharp} - 1_{\FX}^{\sharp} \in\Dercont_{\FY}(\CO_{\FX}, j_{*} f_{0}^{*}\CI).\] 
From \ref{deform1} there results that $h\circ g^{-1} \in \aut_{\FX_{0}}(\FX')$ and $g^{-1} \circ h \in \aut_{\FX_{0}}(\FX)$. It follows that $h$ is an isomorphism.
\end{parraf}

\begin{propo} \label{deform4}
Let $\FY_{0} \inc \FY$ be a  closed immersion in $\sfn$ defined by a square zero Ideal $\CI \subset \CO_{\FY}$ and  $f_{0}:\FX_{0} \to \FY_{0}$ a smooth morphism  in $\sfn$ and suppose that there exists a smooth lifting of  $\FX_{0}$ over $\FY$. Then 
the set of  isomorphic classes of smooth liftings  of  $\FX_{0}$  over $\FY$  is an affine  space over $\ext^{1}(\om^{1}_{\FX_{0}/\FY_{0}}, f_{0}^{*} \CI)$.
\end{propo}

\begin{proof}
Let $\FX_{0} \overset{\, j} \inc \FX \overset{\, f}\to \FY$ and $\FX_{0} \overset{\, j'} \inc \FX' \overset{\, f'}\to \FY$ be two smooth  liftings  over $\FY$. Given an affine open covering $\FU_{\bullet}=\{\FU_{\alpha,o}\}$ of  $\FX_{0}$,  let $\{\FU_{\alpha}\}$ and  $\{\FU'_{\alpha}\}$ be the corresponding affine open coverings  of  $\FX$ and $\FX'$, respectively. From \ref{deform2}, for all $\alpha$ there exists an isomorphism of  $\FY$-formal schemes $u_{\alpha}: \FU_{\alpha} \xto{\sim} \FU'_{\alpha}$ such that  the  following diagram is commutative:
\begin{diagram}[height=2em,w=2em,p=0.3em,labelstyle=\scriptstyle]
\FU_{\alpha,0}& &  \rTinc^{j|_{\FU_{\alpha,0}}}& \FU_{\alpha}&    &      &   \\
 & \rdTinc(3,2)_{j'|_{\FU_{\alpha,0}}}&&\dTto_{u_{\alpha}}^{\wr} &\rdTto(3,2)& &  \\ 
 &	&	&\FU'_{\alpha} &  & \rTto& \FY\\  
\end{diagram}
From \ref{deform3}, for all couples of  indexes $\alpha, \beta$ such that  $\FU_{\alpha \beta,0}:= \FU_{\alpha,0} \cap \FU_{ \beta,0} \neq \varnothing$  if  $\FU_{\alpha \beta}:= \FU_{\alpha} \cap \FU_{ \beta}$, $\FU'_{\alpha \beta}:= \FU'_{\alpha} \cap \FU'_{ \beta}$, and $\FU'_{\alpha \beta,0}:= \FU'_{\alpha,0} \cap \FU'_{ \beta,0}$, there exists a unique homomorphism of  $\CO_{\FU_{\alpha \beta,0}}$-Modules
\[\phi_{\alpha \beta}: \om^{1}_{\FX_{0}/\FY_{0}}|_{\FU_{\alpha \beta,0}}  \to ( f_{0}^{*} \CI) |_{\FU_{\alpha \beta,0}}\] 
such that its adjoint $\om^{1}_{\FU'_{\alpha \beta}/\FY} \to (j'_{*} f_{0}^{*} \CI) |_{\FU'_{\alpha \beta}}$, which we will continue denoting $\phi_{\alpha \beta}$, satisfies 
\[u_{\alpha}|_{\FU_{\alpha \beta}} - u_{\beta}|_{\FU_{\alpha \beta}} = \phi_{\alpha \beta} \circ \hd_{\FU'_{\alpha \beta}}\]
Then $\phi_{\FU_{\bullet}}:=\{\phi_{\alpha \beta}\} \in  \check{\C}^{1}(\FU_{\bullet},\shom_{\CO_{\FX_{0}}}( \om^{1}_{\FX_{0}/\FY_{0}}, f_{0}^{*}\CI))$ has the property that for all $\alpha,\, \beta,\, \gamma$ such that  $\FU_{\alpha \beta \gamma, 0}:= \FU_{\alpha, 0} \cap \FU_{ \beta, 0} \cap \FU_{ \gamma, 0} \neq \varnothing$, the cocycle condition
\begin{equation} \label{aaaayyyyy}
\phi_{\alpha \beta}|_{\FU_{\alpha \beta \gamma,0}} - \phi_{\alpha \gamma}|_{\FU_{\alpha \beta \gamma,0}}+\phi_{ \beta \gamma}|_{\FU_{\alpha \beta \gamma,0}}=0
\end{equation} 
holds and, therefore,  
\[\phi_{\FU_{\bullet}} \in \check{\Z}^{1} (\FU_{\bullet},\shom_{\CO_{\FX_{0}}}( \om^{1}_{\FX_{0}/\FY_{0}}, f_{0}^{*}\CI)).\] 
The  definition  of the element  
\[c_{\FU_{\bullet}}:=[\phi_{\FU}] \in \check{\h}^{1} (\FU_{\bullet},\shom_{\CO_{\FX_{0}}}( \om^{1}_{\FX_{0}/\FY_{0}} ,f_{0}^{*}\CI))\] 
does not depend of  the isomorphisms $\{u_{\alpha}\}$. Indeed, let $\{v_{\alpha}: \FU_{\alpha} \xto{\sim} \FU'_{\alpha}\}_\alpha$ be another collection of  $\FY$-isomorphisms such that, for all $\alpha$, $v_{\alpha} \circ j|_{\FU_{\alpha,0}}= j'|_{\FU_{\alpha,0}}$ and let $\psi_{\FU_{\bullet}}$ be the element 
\[\{\psi_{\alpha \beta}\} \in \check{\C}^{1}(\FU_{\bullet},\shom_{\CO_{\FX_{0}}}( \om^{1}_{\FX_{0}/\FY_{0}} ,f_{0}^{*}\CI))\] 
such that  through adjunction it holds that
$v_{\alpha}|_{\FU_{\alpha \beta}} - v_{\beta}|_{\FU_{\alpha \beta}} = \psi_{\alpha \beta} \circ \hd_{\FU'_{\alpha \beta}/\FY}$, for all couples of  indexes $\alpha, \beta$ such that  $\FU_{\alpha \beta,0} \neq \varnothing$. From \ref{deform3}  there exists $\{\xi_{\alpha}\} \in \check{\C}^{0}(\FU_{\bullet},\shom_{\CO_{\FX_{0}}}( \om^{1}_{\FX_{0}/\FY_{0}} ,f_{0}^{*}\CI))$ such that, for all $\alpha$, their adjoints are such that $u_{\alpha}-v_{\alpha} =\xi_{\alpha} \circ \hd_{\FU'_{\alpha}/\FY}$, therefore, \[[\phi_{\FU_{\bullet}}] = [\psi_{\FU_{\bullet}}] \in  \check{\h}^{1} (\FU_{\bullet},\shom_{\CO_{\FX_{0}}}(\om^{1}_{\FX_{0}/\FY_{0}}  ,f_{0}^{*}\CI)).\]
If $\FV_{\bullet}$ is an affine open  refinement of  $\FU_{\bullet}$, by what we have already seen, we deduce that $c_{\FU_{\bullet}}=c_{\FV_{\bullet}}$. Let us define 
\begin{align*}
c:=[\phi_{\FU_{\bullet}}] &\in  \check{\h}^{1} (\FX_{0},\shom_{\CO_{\FX_{0}}}(\om^{1}_{\FX_{0}/\FY_{0}}  ,f_{0}^{*}\CI)) = \\
&= \h^{1} (\FX_{0}, \shom_{\CO_{\FX_{0}}}(\om^{1}_{\FX_{0}/\FY_{0}}, f_{0}^{*}\CI)) \tag{\cite[(5.4.15)]{te}}\\ 
&= \ext^{1}(\om^{1}_{\FX_{0}/\FY_{0}}, f_{0}^{*}\CI)  
\end{align*}

Reciprocally, let $f:\FX \to \FY$ be a smooth lifting of  $\FX_{0}$  and  consider  $c \in \ext^{1}(\om^{1}_{\FX_{0}/\FY_{0}}, f_{0}^{*} \CI)$. Given $\FU_{\bullet}=(\FU_{\alpha,0})$ an affine open covering  of  $\FX_{0}$, take $(\FU_{\alpha})$ the  corresponding affine open covering  in $\FX$ and \[\phi_{\FU_{\bullet}}=(\phi_{\alpha \beta}) \in \check{\Z}^{1} (\FU_{\bullet},\shom_{\CO_{\FX_{0}}}( \om^{1}_{\FX_{0}/\FY_{0}} ,f_{0}^{*}\CI))\] 
such that  $c=[\phi_{\FU_{\bullet}}]$. For each couples of  indexes $\alpha, \beta$ such that  $\FU_{\alpha \beta}= \FU_{\alpha}\cap \FU_{\beta} \neq \varnothing$, let us consider the  morphism $u_{\alpha \beta}: \FU_{\alpha \beta} \to \FU_{\alpha \beta}$ that is the identity as topological map and that,  as ringed topological spaces map is defined by \[u^{\sharp}_{\alpha \beta} := 1^{\sharp}_{\FU_{\alpha \beta}} + \phi_{\alpha \beta} \circ \hd_{\FU_{\alpha \beta}/\FY},\] 
where again $\phi_{\alpha \beta}$ denotes also its adjoint $\phi_{\alpha \beta} \colon \om^{1}_{\FU_{\alpha \beta}/\FY} \to (j_{*} f_{0}^{*} \CI) |_{\FU_{\alpha \beta}}$, satisfies:
\begin{itemize}
\item
$u_{\alpha \beta} \in \aut_{\FU_{\alpha \beta,0}}(\FU_{\alpha \beta})$ (by \ref{deform1})
\item
 $u_{\alpha \beta}|_{\FU_{\alpha \beta \gamma}} \circ u^{-1}_{\alpha \gamma}|_{\FU_{\alpha \beta \gamma}}\circ u_{ \beta \gamma}|_{\FU_{\alpha \beta \gamma}}=1_{\FU_{\alpha \beta \gamma}}$, for any  $\alpha,\, \beta,\, \gamma$ such that $ \FU_{\alpha \beta \gamma}:= \FU_{\alpha} \cap \FU_{\beta} \cap \FU_{\gamma} \neq \varnothing$ (because  $\{\phi_{\alpha \beta}\}$ satisfies the cocycle condition(\ref{aaaayyyyy}))
\item
 $u_{\alpha \alpha}= 1_{\FU_{\alpha}}$ 
 and
 $u_{\alpha \beta}^{-1}= u_{\beta \alpha}$
\end{itemize}
Then the $\FY$-formal schemes $\FU_{\alpha}$ glue into a  smooth lifting $f':\FX' \to\FY$ of  $\FX_{0}$ through the  morphisms $\{u_{\alpha \beta}\}$, since the  morphism $f:\FX \to \FY$ is compatible with the family of isomorphisms $\{u_{\alpha \beta}\}$.

We leave to the reader the verification that these correspondences are mutually inverse.
\end{proof}

\begin{rem}
Proposition \ref{deform4} can be rephrased in the language of  torsor theory as follows: 
The sheaf in $\FX_{0}$ that associates to each open  $\FU_{0} \subset \FX_{0}$
the set of  isomorphism classes of smooth  liftings of  $\FU_{0}$  over $\FY$  is a pseudo torsor over $\ext^{1}(\om^{1}_{\FX_{0}/\FY_{0}}, f_{0}^{*} \CI)$.
\end{rem}

\begin{rem}
With the hypothesis of  Proposition \ref{deform4}, if $f_{0}$ and $\FY$ are  in $\sfna$, we have that \[\h^{1} (\FX_{0},\shom_{\CO_{\FX_{0}}}( \om^{1}_{\FX_{0}/\FY_{0}} ,f_{0}^{*}\CI))=0\] (\emph{cf.}  \cite[Corollary 3.1.8]{AJL1}) and, therefore, there exists a unique isomorphism  class of  liftings of  $\FX_{0}$ over $\FY$.
\end{rem}

\begin{propo} \label{obstrext2}
Let us consider in $\sfn$ a closed immersion $\FY_{0} \inc \FY$ given by a square zero Ideal $\CI \subset \CO_{\FY}$   and $f_{0}:\FX_{0} \to \FY_{0}$ a smooth morphism with $\FX_{0}$ a separated formal scheme. Then there is an element $c_{f_{0}} \in \ext^{2}(\om^{1}_{\FX_{0}/\FY_{0}}, f_{0}^{*}\CI)$  such that: $c_{f_{0}}$ vanishes if, and only if,  there exists a smooth lifting $\FX$ of  $\FX_{0}$ over $\FY$.
\end{propo}

\begin{proof}
From  Proposition \ref{condinflocal} and  Corollary \ref{deformloclis} there exists $\FU_{\bullet}=(\FU_{\alpha, 0})$ an affine open covering of  $\FX_{0}$, such that  for all $\alpha$ there exists $\FU_{\alpha}$ a  smooth lifting of  $\FU_{\alpha, 0}$ over $\FY$. As $\FX_{0}$ is a separated formal scheme $\FU_{\alpha \beta, 0}:= \FU_{\alpha,  0} \cap \FU_{ \beta, 0}$ is an affine open set for any $\alpha,\beta$ and, if we call $\FU_{\alpha \beta} \subset \FU_{\alpha}$ and $\FU_{\beta \alpha} \subset \FU_{\beta}$ to the corresponding open sets, from \ref{deform2} there exists an isomorphism $u_{\alpha \beta}: \FU_{\alpha \beta} \xto{\sim} \FU_{\beta\alpha}$ such that  the  following diagram is commutative:
\begin{diagram}[height=2em,w=2em,p=0.3em,labelstyle=\scriptstyle]
\FU_{\alpha \beta, 0}& &  \rTinc&\FU_{\alpha \beta} &    &      &   \\
			&\rdTinc(3,2)&&\dTto_{u_{\alpha \beta}}^{\wr }     &\rdTto(3,2)& &  \\ 
		&	&	&\FU_{\beta\alpha}&  & \rTto& \FY\\  
\end{diagram}
For any $\alpha, \beta, \gamma$ such that  $\FU_{\alpha \beta \gamma,0}:= \FU_{\alpha,0} \cap \FU_{ \beta,0} \cap \FU_{ \gamma,0}\neq \varnothing$, let us write $\FU_{\alpha \beta \gamma}:= \FU_{\alpha \beta} \times_\FY \FU_{ \alpha \gamma}$. It holds that
\[u_{\alpha \beta \gamma}:=  
u^{-1}_{ \alpha \gamma}|_{\FU_{\gamma \beta} \cap \FU_{ \gamma\alpha}} \circ 
u_{ \beta \gamma}|_{\FU_{ \beta\alpha} \cap \FU_{\beta\gamma}} \circ 
u_{\alpha \beta}|_{\FU_{\alpha \beta} \cap \FU_{\alpha \gamma}} 
\in \aut_{\FU_{\alpha \beta \gamma,0}}(\FU_{\alpha \beta \gamma}).\] 
From  \ref{deform1} there exists a unique \[\phi_{\alpha \beta \gamma} \in \ga(\FU_{\alpha \beta \gamma,0}, \shom_{\CO_{\FX_{0}}}(\om^{1}_{\FX_{0}/\FY_{0}},  f_{0}^{*}\CI))\] 
whose adjoint satisfies the relation  
\[u_{\alpha \beta \gamma}^{\sharp}-1^{\sharp}_{\FU_{\alpha \beta \gamma}}= \phi_{\alpha \beta \gamma} \circ \hd_{\FU_{\alpha \beta \gamma}/\FY}.\] 
Then the element 
\[\phi_{\FU_{\bullet}}:=(\phi_{\alpha \beta \gamma}) \in  \check{\C}^{2}(\FU_{\bullet},\shom_{\CO_{\FX_{0}}}( \om^{1}_{\FX_{0}/\FY_{0}}, f_{0}^{*}\CI))\] 
satisfies for any $\alpha,\, \beta,\, \gamma,\, \delta$ such that  $\FU_{\alpha \beta \gamma \delta,0}:= \FU_{\alpha,0 } \cap \FU_{ \beta,0 } \cap \FU_{ \gamma,0} \cap \FU_{ \delta,0}\neq \varnothing$, the cocycle condition 
\begin{equation} \label{aaaayyyyy2}
\phi_{\alpha \beta \gamma}|_{\FU_{\alpha \beta \gamma \delta,0}}  -
\phi_{\alpha \gamma \delta}|_{\FU_{\alpha \beta \gamma \delta,0}} +
\phi_{ \beta \gamma \delta}|_{\FU_{\alpha \beta \gamma \delta,0}} -
\phi_{ \beta \delta \alpha}|_{\FU_{\alpha \beta \gamma \delta,0}} = 0
\end{equation} 
and, therefore,  $\phi_{\FU_{\bullet}} \in \check{\Z}^{2} (\FU_{\bullet},\shom_{\CO_{\FX_{0}}}( \om^{1}_{\FX_{0}/\FY_{0}}, f_{0}^{*}\CI))$. Using \ref{deform3} and  reasoning  in an analogous way as in  the proof of  Proposition \ref{deform4}, it is easily seen that the  definition  of 
\[c_{\FU_{\bullet}}:=[\phi_{\FU_{\bullet}}] \in \check{\h}^{2} (\FU_{\bullet},\shom_{\CO_{\FX_{0}}}( \om^{1}_{\FX_{0}/\FY_{0}} ,f_{0}^{*}\CI))\] 
does not depend on the chosen  isomorphisms $\{u_{\alpha \beta}\}$. 
Furthermore, if $\FV_{\bullet}$ is an affine open refinement  of  $\FU_{\bullet}$, then
\[c_{\FU_{\bullet}}=c_{\FV_{\bullet}} \in \check{\h}^{2} (\FX_{0},\shom_{\CO_{\FX_{0}}}( \om^{1}_{\FX_{0}/\FY_{0}} ,f_{0}^{*}\CI))\]
Proposition \ref{flplano} implies that $\om^{1}_{\FX_{0}/\FY_{0}}$ is a locally free $\CO_{\FX_{0}}$-Module. Since  $\FX_{0}$ is separated\footnote{Given $\FX$ a separated formal scheme in $\sfn$ and $\CF$ a sheaf of abelian  groups over $\FX$ using  \cite[Corollary 3.1.8]{AJL1} and \cite[Ch. III, Exercise 4.11]{ha1}  we have that $\check{\h}^{i} (\FX,\CF)=\h^{i} (\FX, \CF)$, for all $i >0$.}, we set
\begin{align*}
c_{f_{0}}:=[\phi_{\FU_{\bullet}}] &\in  \check{\h}^{2} (\FX_{0},\shom_{\CO_{\FX_{0}}}(\om^{1}_{\FX_{0}/\FY_{0}}  ,f_{0}^{*}\CI)) = \\
&= \h^{2} (\FX_{0}, \shom_{\CO_{\FX_{0}}}(\om^{1}_{\FX_{0}/\FY_{0}}  ,f_{0}^{*}\CI))\\
&= \ext^{2}(\om^{1}_{\FX_{0}/\FY_{0}}, f_{0}^{*}\CI)  
\end{align*}

Let us show that $c_{f_{0}}$ is the obstruction to the existence of a smooth lifting of  $\FX_{0}$ over $\FY$. If there exists $\FX$ a smooth lifting of  $\FX_{0}$ over $\FY$, one could take the isomorphisms $\{u_{\alpha \beta}\}$ above as the identities, then $c_{f_{0}}=0$, trivially.
Reciprocally, let $\FU_{\bullet}=\{\FU_{\alpha,0}\}$ be an affine open  covering of  $\FX_{0}$ and, for each $\alpha$, $\FU_{\alpha}$ a smooth lifting of $\FU_{\alpha,0}$ over $\FY$ such that, with  the notations  established at the beginning of  the proof, $c_{f_{0}}=[\phi_{\FU_{\bullet}}]$ with 
\[\phi_{\FU_{\bullet}}=(\phi_{\alpha \beta \gamma}) \in \check{\Z}^{2} (\FU_{\bullet},\shom_{\CO_{\FX_{0}}}( \om^{1}_{\FX_{0}/\FY_{0}}, f_{0}^{*}\CI)).\] 
In view of $c_{f_{0}}=0$,
we are going to glue the $\FY$-formal schemes   $(\FU_{\alpha})$ into a lifting of  $\FX_{0}$ over $\FY$. By hypothesis, we have that 
\[\phi_{\FU_{\bullet}} \in  \check{\B}^{2} (\FU_{\bullet}, \shom_{\CO_{\FX_{0}}}( \om^{1}_{\FX_{0}/\FY_{0}} ,f_{0}^{*}\CI))\] 
and therefore, there exists 
\[(\phi_{\alpha \beta}) \in \check{\C}^{1}(\FU_{\bullet},\shom_{\CO_{\FX_{0}}}(\om^{1}_{\FX_{0}/\FY_{0}} ,f_{0}^{*}\CI))\]
such that, for any  $\alpha, \beta, \gamma$ with $\FU_{\alpha \beta \gamma, 0} \neq \varnothing$, 
\begin{equation} \label{datosrecolec22}
\phi_{\alpha \beta}|_{\FU_{\alpha \beta \gamma, 0}} - \phi_{\alpha \gamma}|_{\FU_{\alpha \beta \gamma, 0}}+ \phi_{ \beta \gamma}|_{\FU_{\alpha \beta \gamma, 0}}=\phi_{\alpha \beta \gamma}
\end{equation} 
For each couple of  indexes $\alpha, \beta$ such that  $\FU_{\alpha \beta,0} \neq \varnothing$,  let $v_{\alpha \beta}: \FU_{\alpha \beta} \to \FU_{\beta \alpha }$ be the morphism which is the identity as topological map, and that,  as ringed topological spaces map  is given by 
\[v^{\sharp}_{\alpha \beta} := u^{\sharp}_{\alpha \beta} - \phi_{\alpha \beta} \circ \hd_{\FX/\FY}|_{\FU_{\alpha \beta}}.\] 
The constructed family $\{v_{\alpha \beta}\}$, satisfies:
\begin{itemize}
\item
Each map $v_{\alpha \beta}$ is an isomorphism of  $\FY$-formal schemes (by \ref{deform3}).
\item
For any  $\alpha$, $\beta$, $\gamma$ such that $\FU_{\alpha \beta \gamma, 0}\neq \varnothing$,  
\[v^{-1}_{ \alpha \gamma}|_{\FU_{\gamma \beta} \cap \FU_{ \gamma\alpha}} \circ v_{ \beta \gamma}|_{\FU_{ \beta\alpha} \cap \FU_{\beta\gamma}} \circ v_{\alpha \beta}|_{\FU_{\alpha \beta} \cap \FU_{\alpha \gamma}}=1_{\FU_{\alpha \beta} \cap \FU_{\alpha \gamma}}
\]
by (\ref{aaaayyyyy2}) and (\ref{datosrecolec22}).
\item
For any  $\alpha$, $\beta$,
 $v_{\alpha \alpha}= 1_{\FU_{\alpha}}$ 
and
 $v_{\alpha \beta}^{-1}= v_{\beta \alpha}$.
\end{itemize}
The $\FY$-formal schemes $\{\FU_{\alpha}\}$ glue  into a smooth lifting $f:\FX \to \FY$ of  $\FX_{0}$ over $\FY$ through the glueing morphisms $\{v_{\alpha \beta}\}$.

\end{proof}

\begin{cor} \label{*}
With the  hypothesis of   Proposition \ref{obstrext2}, if $\FX_{0}$ and $\FY$ are  affine, 
there exists a  lifting of  $\FX_{0}$ over $\FY$.
\end{cor}

\begin{proof}
Applying \cite[Corollary 3.1.8]{AJL1} we have that \[\h^{2} (\FX_{0},\shom_{\CO_{\FX_{0}}}( \om^{1}_{\FX_{0}/\FY_{0}} ,f_{0}^{*}\CI))=0\] and the  result follows from  the last proposition.
\end{proof}

\begin{cuesab}
If $\FX_{0}$ is not separated, the  argument used in the proof of  Proposition \ref{obstrext2}  is not valid. We conjecture that it is possible to construct an adequate  version of the cotangent complex generalizing the construction of the sheaf of 1-differentials in $\sfn$ providing a more general version of this deformation theory.
\end{cuesab}



\end{document}